\pdfoutput=1
\documentclass[twoside]{book}

\usepackage{pdfpages}
\usepackage{subfigure}
\usepackage{wrapfig}
\usepackage{bbm}
\usepackage{amsmath, amsfonts, amscd, amssymb, amsthm}
\usepackage{tikz}
\usetikzlibrary{matrix,arrows,calc}

\usepackage[bottom]{footmisc}
\usepackage{aliascnt}
\usepackage{xspace}
\usepackage{stmaryrd}

\usepackage{float}
\usepackage{algorithm}
\usepackage{graphicx}
\usepackage{url}
\usepackage{color}
\usepackage[english,dutch]{babel}

\ifx \@lemma \@empty

\fi
\ifx \@theorem \@empty

\fi
\ifx \@definition \@empty

\fi
\ifx \@corollary \@empty

\fi

\newcommand{\thedocumentname}{Geometric structures and Lie algebroids}
\newcommand{\theauthor}{Ralph L.\ Klaasse}

\theoremstyle{plain}
\newtheorem{thm}{Theorem}[section]

\newtheorem*{thm2}{Theorem}
\newaliascnt{lem}{thm}
\newtheorem{lem}[lem]{Lemma}
\aliascntresetthe{lem}

\newaliascnt{prop}{thm}
\newtheorem{prop}[prop]{Proposition}
\aliascntresetthe{prop}               

\newaliascnt{cor}{thm}
\newtheorem{cor}[cor]{Corollary}
\aliascntresetthe{cor}

\theoremstyle{definition}
\newaliascnt{defn}{thm}
\newtheorem{defn}[defn]{Definition}
\aliascntresetthe{defn}

\newaliascnt{rem}{thm}
\newtheorem{rem}[rem]{Remark}
\aliascntresetthe{rem}

\newaliascnt{exa}{thm}
\newtheorem{exa}[exa]{Example}
\aliascntresetthe{exa}

\newaliascnt{que}{thm}
\newtheorem{que}[que]{Question}
\aliascntresetthe{que}

\theoremstyle{plain}
\newtheorem{thmalpha}{Theorem}

\newcommand{\N}{\mathbb{N}}
\newcommand{\Z}{\mathbb{Z}}
\newcommand{\R}{\mathbb{R}}
\newcommand{\C}{\mathbb{C}}
\newcommand{\im}{{\rm im}\,}
\newcommand{\bi}{\begin{itemize}}
\newcommand{\ei}{\end{itemize}}
\newcommand{\be}{\begin{equation*}}
\newcommand{\ee}{\end{equation*}}
\newcommand{\bp}[1][]{\begin{proof}[Proof#1.]}				
\newcommand{\ep}{\end{proof}}
\newcommand{\wt}{\widetilde}
\newcommand{\wh}{\widehat}
\newcommand{\mf}{\mathfrak}
\newcommand{\ol}{\overline}
\newcommand{\mc}{\mathcal}
\newcommand{\mb}{\mathbb}
\newcommand{\symp}{symplectic structure\xspace}

\newcommand{\blog}{log-\symp}

\newcommand{\fsymp}{folded-\symp}
\newcommand{\Fsymp}{Folded-\symp}
\newcommand{\nsymp}{near-\symp}

\newcommand{\acs}{almost-complex structure\xspace}

\newcommand{\gcs}{generalized complex structure\xspace}
\newcommand{\Gcs}{Generalized complex structure\xspace}
\newcommand{\sgcs}{stable \gcs}
\newcommand{\Sgcs}{Stable \gcs}
\newcommand{\bacs}{$b$-\acs}

\newcommand{\alf}{achiral Lefschetz fibration\xspace}

\newcommand{\lf}{Lefschetz fibration\xspace}
\newcommand{\blf}{$b$-\lf}

\renewcommand{\b}[1][]{{^b}{#1}}
\newcommand{\btwo}[1][]{{^b}{#1}}
\renewcommand{\btwo}[1][]{#1}

\newcommand{\placeholder}{boundary \lf}

\newcommand{\thesistitle}
{Geometric structures and\\ Lie algebroids}

\newcommand{\nativetitle}
{Meetkundige structuren en Lie algebro\"iden}

\newcommand{\thesisauthor}
{R.~L.~Klaasse}

\newcommand{\fullthesisauthor}
{Ralph Leonard Klaasse}

\newcommand{\rectorname}
{prof.~dr.~G.~J.~van~der~Zwaan}

\newcommand{\formalpromotiondate}
{dinsdag 29 augustus 2017 des middags te 12.45 uur}

\newcommand{\birthinfo}
{1 april 1990 te Nijmegen}

\newcommand{\committee}
{
Prof.~dr.~Erik~P.~van den Ban, Universiteit Utrecht\\
Dr.~Christian~Blohmann, Max-Planck-Institut f\"ur Mathematik, Bonn\\
Prof.~dr.~Rui~Loja~Fernandes, University of Illinois at Urbana--Champaign\\
Prof.~dr.~Marco~Gualtieri, University of Toronto\\
Prof.~dr.~Eckhard~Meinrenken, University of Toronto
}

\newcommand{\ISBNnumber}
{978-90-393-6813-8}

\newcommand{\promotorname}
{Prof.~dr.~M.~N.~Crainic}

\newcommand{\copromotorname}
{Dr.~G.~R.~Cavalcanti}

\def\@thesisdraftmode{}

\usepackage{multirow}
\usepackage{bibentry}
\nobibliography*
\makeatletter
\renewcommand\bibentry[1]{\nocite{#1}{\frenchspacing
		\@nameuse{BR@r@#1\@extra@b@citeb}}}
\makeatother

\usepackage{emptypage}
\usepackage[nocompress, noadjust]{cite}

\usepackage[T1]{fontenc}

\usepackage{times}

\usepackage[paperheight=24cm,paperwidth=17cm,headheight=12pt]{geometry}
\addto\captionsenglish{}

\usepackage{fancyhdr}
\usepackage{ifthen}
\fancypagestyle{mystyle}{%
\fancyhf{}%
\renewcommand\headrulewidth{0.4pt}%
\fancyhead[RO,LE]{\thepage}%
\fancyhead[LO]{\small\textbf{\rightmark}}%
\fancyhead[RE]{\small\textbf{\leftmark}}%
}
\fancypagestyle{plain}{%
\fancyhf{}%
\renewcommand\headrulewidth{0.4pt}%
\fancyhead[RO,LE]{\thepage}%
}
\newcommand{\myleftmark}{}
\newcommand{\myrightmark}{}
\newcommand{\mymarkleft}[1] {\renewcommand{\myleftmark}{#1}\markboth{\myleftmark}{\myrightmark}}
\newcommand{\mymarkright}[1]{\renewcommand{\myrightmark}{#1}\markboth{\myleftmark}{\myrightmark}}

\renewcommand\sectionmark[1]{%
\ifthenelse{\value{chapter}>0}%
 {%
 \ifthenelse{\value{section}>0}%
  {%
    \mymarkright{\textbf{\thesection}\ --- \textbf{#1}}%
  }%
  {\mymarkright{\textbf{#1}}}%
 }%
 {\mymarkright{\textbf{#1}}}
}

\usepackage[pdftitle={\thedocumentname},pdfauthor={\theauthor},bookmarks,urlcolor=black,colorlinks=true,linkcolor=black,citecolor=black]{hyperref}

\begin{document}
\selectlanguage{english}

\includepdf[pages={1}]{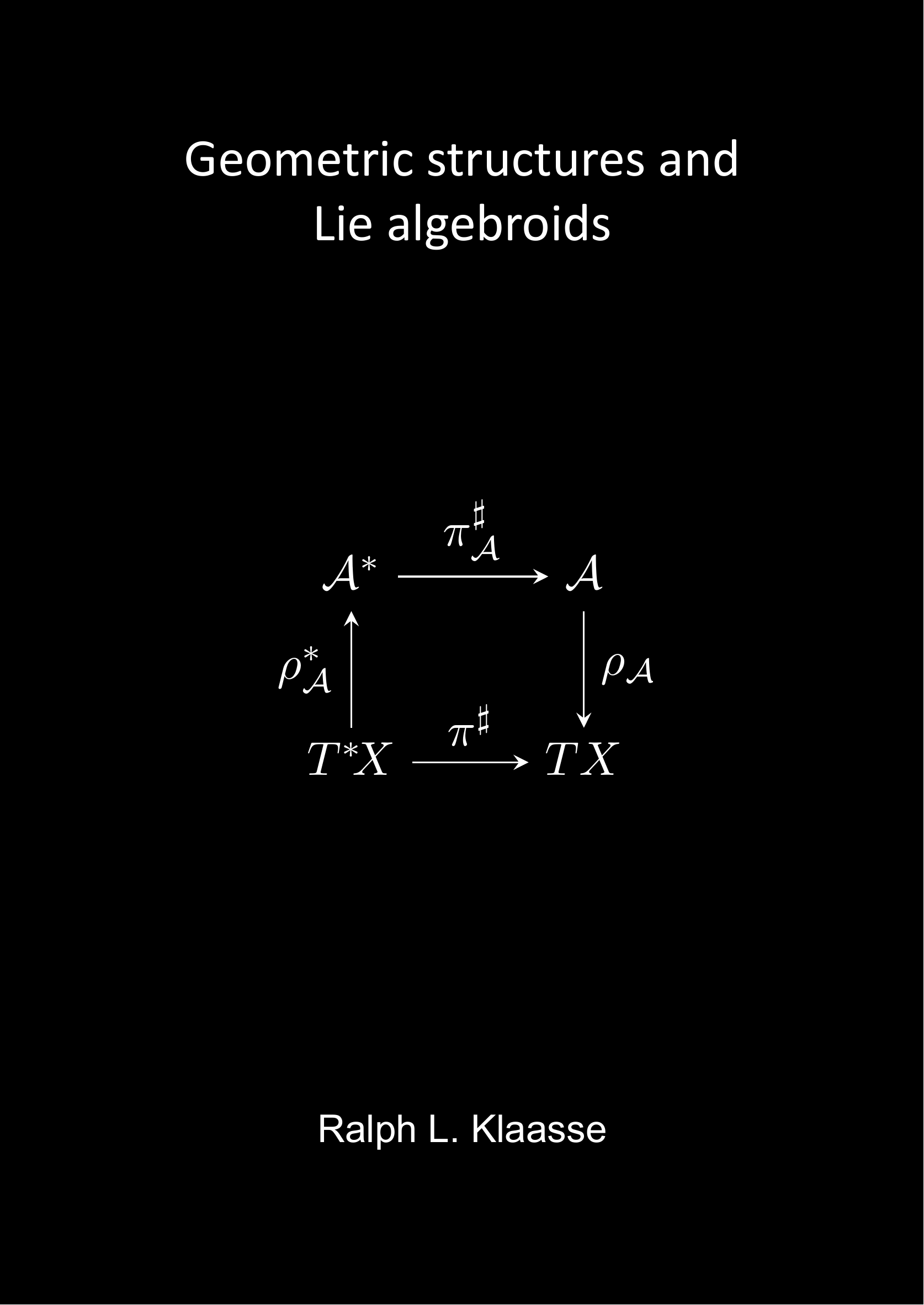}
\newpage

\frontmatter

\thispagestyle{empty}
\hspace{0pt}
\vspace{0.33\textheight}

\begin{center}
{\huge\bf\thesistitle}
\end{center}

\clearpage

\thispagestyle{empty}

\noindent\emph{Thesis committee:}\\
\committee

\vfill

\noindent Geometric structures and Lie algebroids, \thesisauthor\\
\noindent Ph.D.\ thesis Utrecht University, August 2017
\vspace{\baselineskip}

\noindent ISBN: \ISBNnumber\\
\noindent Copyright \copyright\ 2017 by \thesisauthor

\clearpage

\thispagestyle{empty}

\begin{center}
\hspace{0pt}

{\huge\bf\thesistitle}
\vspace{7em}

{\Large\bf\nativetitle}

\vspace{\baselineskip}
(met een samenvatting in het Nederlands)
\vspace{13em}

{\Large Proefschrift}

\vspace{0.5em}
\begin{center}
\begin{minipage}{0.87\textwidth}
\begin{flushleft}
ter verkrijging van de graad van doctor aan de Universiteit Utrecht op gezag van de rector magnificus, \rectorname, ingevolge het besluit van het college van promoties in het openbaar te verdedigen op \formalpromotiondate
\end{flushleft}
\end{minipage}
\end{center}
\vfill

door
\vspace{\baselineskip}

{\large\fullthesisauthor}
\vspace{\baselineskip}

{geboren op \birthinfo}
\end{center}
\vspace{\baselineskip}

\clearpage

\thispagestyle{empty}

\noindent Promotor: \promotorname\\
\noindent Copromotor: \copromotorname

\vspace{\fill}
\noindent{The research in this thesis was supported by VIDI grant number 639.032.221\\ from NWO, the Netherlands Organisation for Scientific Research.}

\clearpage

\thispagestyle{empty}

\hspace{0pt}
\vspace{0.2\textheight}

\begin{center}
{\large \emph{Voor Lide en Leny}}
\end{center}
\cleardoublepage

\newcommand{\resetemptypagestyle}{
\fancypagestyle{empty}{%
\fancyhf{}%
\renewcommand\headrulewidth{0pt}%
}}
\resetemptypagestyle

\chapter*{Preface}
\addcontentsline{toc}{chapter}{Preface}

This thesis is a partial account of research done while being a Ph.D.\ candidate at the Mathematical Institute of Utrecht University. In this thesis we study geometric structures by using Lie algebroids. More precisely, we study the way in which geometric structures (i.e.\ Poisson, Dirac, and generalized complex structures) with singular behavior can have their singularity absorbed in that of a Lie algebroid. Most interesting is the case in which a Lie algebroid can be found for which the lift is now nondegenerate, as this opens the door to application of e.g.\ symplectic techniques.

\vspace{\baselineskip}

This thesis has been written to be self-contained as much as possible to facilitate the reader. A summary giving an overview of its contents for non-specialists both in \hyperref[chap:summary]{English} and in \hyperref[chap:samenvatting]{Dutch} can be found at the end of this thesis.

\vspace{\baselineskip}

\noindent \textbf{Publications.} This thesis is based in part on the research papers	
{\footnotesize\raggedright

\bi
	\item[\cite{CavalcantiKlaasse16}] \bibentry{CavalcantiKlaasse16}.
	\item[\cite{CavalcantiKlaasse17}] \bibentry{CavalcantiKlaasse17}.
\ei
}
	
\noindent as well as the current status of

{\footnotesize\raggedright
\bi
	\item[\cite{KlaasseLanius17two}] \bibentry{KlaasseLanius17two}.
	\item[\cite{KlaasseLanius17}] \bibentry{KlaasseLanius17}.
\ei
}

\noindent Other papers:
	
{\footnotesize\raggedright
\bi
	\item[\cite{BehrensCavalcantiKlaasse17}] \bibentry{BehrensCavalcantiKlaasse17}.
\ei
}

\noindent Some of the other results in this thesis are being prepared for separate publication.

\vspace{\baselineskip}

\noindent I hope you will have a pleasant read,
\vspace{\baselineskip}

Ralph Leonard Klaasse,

August 2017.

\cleardoublepage
\phantomsection
\addcontentsline{toc}{chapter}{Table of Contents}
\tableofcontents
\chapter*{Introduction}
\label{chap:intro}
\mymarkleft{\textbf{Introduction}}
\addcontentsline{toc}{chapter}{Introduction}

This thesis discusses a general procedure by which Lie algebroids can be used to study geometric structures with mild singularities. Such an approach has proven to be successful in several cases, see e.g.\ \cite{Cavalcanti17,CavalcantiGualtieri15, FrejlichMartinezTorresMiranda15,GualtieriLi14,GuilleminMirandaPires14, Lanius16two, Lanius16,MarcutOsornoTorres14,MarcutOsornoTorres14two, Scott16}. We will discuss this process using a concrete example coming from Poisson geometry.
\subsection*{Lifting Poisson structures}
Consider a symplectic structure $\omega$ on a manifold $X$ of dimension $2n$. The symplectic condition means it is both closed, $d\omega = 0$, and nondegenerate, $\omega^n \neq 0$. Due to the Darboux theorem, we know in a neighbourhood of any point $p \in X$ there exists a coordinate system $(x_1,\dots,x_{2n})$ for which
\be
	\omega = dx_1 \wedge dx_2 + \dots + dx_{2n-1} \wedge dx_{2n}.
\ee
A symplectic structure dually can be described by a nondegenerate Poisson bivector $\pi = \omega^{-1}$ which will satisfy $\pi^n \neq 0$ (and the integrability condition $[\pi,\pi] = 0$). In the same coordinate system as above, it would be given by
\be
	\pi = \frac{\partial}{\partial x_1} \wedge \frac{\partial}{\partial x_2} + \dots + \frac{\partial}{\partial x_{2n-1}} \wedge \frac{\partial}{\partial x_{2n}}.
\ee
From the point of view of Poisson geometry, then, it is natural to want to study Poisson structures which are in some sense close to being symplectic. The type of Poisson structures we consider in this example are those we call \emph{log-Poisson structures} (see Section \ref{sec:logpoisson}).\footnote{These also go under the name of $b$-Poisson, $b$-symplectic or log-symplectic structures \cite{Cavalcanti17,FrejlichMartinezTorresMiranda15,GualtieriLi14,GuilleminMirandaPires14,MarcutOsornoTorres14,MarcutOsornoTorres14two, NestTsygan96}, or as topologically stable Poisson structures \cite{Radko02}.} They are defined by demanding that $\pi^n$ need not globally be nonzero, but is allowed to vanish transversally as a section of the line bundle $\wedge^{2n} TX$. Due to this condition, they will fail to be nondegenerate on a smooth hypersurface $Z \subseteq X$. The Weinstein splitting theorem shows we can find coordinates $x_i$ with $Z = \{x_1 = 0\}$ for which
\be
	\pi = x_1 \frac{\partial}{\partial x_1} \wedge \frac{\partial}{\partial x_2} + \dots + \frac{\partial}{\partial x_{2n-1}} \wedge \frac{\partial}{\partial x_{2n}}.
\ee
Its dual $\omega = \pi^{-1}$ then constitutes a ``singular symplectic structure'', viz.
\be
	\omega = \frac1{x_1} dx_1 \wedge dx_2 + \dots + dx_{2n-1} \wedge dx_{2n}.
\ee
Writing $x_1^{-1} dx_1 = d \log x_1$, we see that this form has a logarithmic singularity at $Z$, and is not globally defined on $X$. The goal of the introduction of Lie algebroids to this story is to be able to nevertheless view this as a symplectic, smooth object.

Consider that in the coordinates as above, the sheaf $\mc{V}_X$ of vector fields on $X$ and its dual $\mc{V}_X^*$ locally are given by (where $\langle \, \rangle$ denotes generation as $C^\infty(X)$-module)
\be
	\mc{V}_X = \left\langle \frac{\partial}{\partial x_1}, \dots, \frac{\partial}{\partial x_{2n}} \right\rangle, \qquad \mc{V}_X^* = \langle dx_1, \dots, dx_{2n} \rangle.
\ee
Consider now a new sheaf, which we will denote by $\mc{V}_X(I_Z)$, consisting of all vector fields on $X$ which are tangent to $Z$ (alternatively, all vector fields preserving the vanishing ideal $I_Z$ of $Z$). In the above coordinates, it and its dual are given by
\be
	\mc{V}_X(I_Z) = \left\langle x_1 \frac{\partial}{\partial x_1}, \frac{\partial}{\partial x_2}, \dots, \frac{\partial}{\partial x_{2n}} \right\rangle, \quad \mc{V}_X(I_Z)^* = \left\langle d\log x_1, dx_2, \dots, dx_{2n} \right\rangle.
\ee
This sheaf of $C^\infty(X)$-submodules of $\mc{V}_X$ is locally free, hence is the space of sections of a vector bundle of rank $2n$ by the Serre--Swan theorem. This bundle is called the \emph{log-tangent bundle} (see Section \ref{sec:logtangentbundle}), and will be denoted by $\mc{A}_Z$.\footnote{Another notation we will use is $TX(-\log Z)$. Moreover, it is also called the \emph{$b$-tangent bundle} ${}^b TX$, and was introduced by Melrose \cite{Melrose93} in the case where $Z = \partial X$.} Note while $\mc{V}_X(I_Z)$ is a subsheaf of $\mc{V}_X$, the bundle $\mc{A}_Z$ is not a subbundle of $TX$. The inclusion on sections induces a map called the \emph{anchor}, $\rho_{\mc{A}_Z}\colon \mc{A}_Z \to TX$, which is an isomorphism outside of $Z$. As the sheaf $\mc{V}_X(I_Z)$ is involutive under the Lie bracket on $\mc{V}_X$, the bundle $\mc{A}_Z$ further inherits a bracket from $TX$ which is suitably compatible with the anchor map, turning it into a Lie algebroid (see Section \ref{sec:basicdefinitions}).

The reason for introducing $\mc{A}_Z$ is the following. We can lift the log-Poisson bivector $\pi$ to a bivector in $\mc{A}_Z$. In other words, there exists a section $\pi_{\mc{A}_Z} \in \Gamma(\wedge^2 \mc{A}_Z)$ which satisfies $\rho_{\mc{A}_Z}(\pi_{\mc{A}_Z}) = \pi$. This is a nondegenerate $\mc{A}_Z$-Poisson structure (see Section \ref{sec:apoissonstrs}). The reason for doing this is that now the inverse $\omega_{\mc{A}_Z} := \pi_{\mc{A}_Z}^{-1}$, will be a smooth object. In fact, it is what we call an \emph{$\mc{A}_Z$-symplectic structure}: a section of $\wedge^2 \mc{A}_Z^*$ which satisfies $d_{\mc{A}_Z} \omega_{\mc{A}_Z} = 0$, and $\omega_{{\mc{A}_Z}}^n \neq 0$. Here $d_{\mc{A}_Z}$ is a differential that is defined using the bracket on $\Gamma(\mc{A}_Z)$.

We see that the introduction of the Lie algebroid $\mc{A}_Z$ allows us to view otherwise singular forms as being smooth. In this sense, the goal is to find a Lie algebroid to which a given geometric structure can be \emph{lifted}. Its lift will then be more nondegenerate, with the anchor of the Lie algebroid having absorbed some of the singularities. Once one has found a Lie algebroid to which it lifts nondegenerately, symplectic techniques can be brought to bear to study $\omega_{{\mc{A}_Z}}$, and hence the original log-Poisson structure $\pi$ we started with.
\subsection*{Constructing $\mc{A}$-symplectic structures}
One of the main themes of this thesis is to use the above strategy to construct geometric structures on a given manifold $X$. This is done by relating the existence of such a structure to the existence of a certain type of map the manifold $X$ admits. 

More precisely, we want to extend results from symplectic geometry relating the existence of fibration-like maps to that of $\mc{A}$-symplectic structures. We now give a brief overview of what has been done in this direction.

In \cite{Thurston76}, see also \cite{McDuffSalamon98}, Thurston showed how to equip symplectic fiber bundles with \symp{}s. Gompf \cite{GompfStipsicz99} then showed that Lefschetz fibrations lead to \symp{}s in dimension four, adapting Thurston's methods. Conversely, Donaldson \cite{Donaldson99} proved using approximately holomorphic methods that symplectic forms lead to Lefschetz pencils. Further, Gompf \cite{Gompf04two, Gompf04} introduced hyperpencils in all dimensions and showed they admit \symp{}s, aiming to give a topological characterization of symplectic manifolds.

Indeed, having established this correspondence, the study of Lefschetz pencils can shed light onto symplectic geometry. One can then branch out to other symplectic-like structures, and connect them to suitably generalized Lefschetz-type fibrations. This has been done for so-called near-symplectic and folded-symplectic structures. Indeed, Auroux--Donaldson--Katzarkov \cite{AurouxDonaldsonKatzarkov05} introduced broken Lefschetz pencils to study \nsymp{}s, and established the same correspondence as between Lefschetz pencils and \symp{}s. In another direction, Etnyre--Fuller \cite{EtnyreFuller06} showed that any \alf with a section gives rise to a \nsymp. Moreover, Baykur \cite{Baykur06} showed that \fsymp{}s arise out of \alf{}s. For more background on these results, see also \cite{Honda04,LeBrun97, GayKirby07, Baykur08, Baykur09, Lekili09, AkbulutKarakurt08} and \cite{Baykur06, CannasDaSilva10} respectively.

As we discussed in the previous section, we capture the singular behavior of the geometric structures at hand using Lie algebroids. Consequently, we use morphisms between such Lie algebroids as our maps of choice, and in particular consider what we call Lie algebroid Lefschetz fibrations (see Chapter \ref{chap:constrasymp} for more details). One of the main results of this thesis is the following adaptation of Gompf--Thurston techniques to Lie algebroids. It has appeared before as \cite[Theorem 1.1]{CavalcantiKlaasse17}.
\begin{thmalpha}\label{thm:introasymp} Let $(\varphi,f)\colon \mc{A}^4_{X} \to \mc{A}^2_{\Sigma}$ be a Lie algebroid \lf with connected fibers. Assume that $\mc{A}_\Sigma$ admits a Lie algebroid symplectic structure and there exists a closed $\mc{A}_X$-two-form $\eta$ such that the restriction $\eta|_{\ker \varphi}$ is nondegenerate. Then $X$ admits an $\mc{A}_X$-\symp.
\end{thmalpha}
This result is proven in this thesis as \autoref{thm:thurstonlfibration}. This should be thought of as a tool to construct the geometric structures that are captured by the Lie algebroids one is using. Of course, in order to apply it we must be able to ensure we can satisfy the hypotheses. This we discuss in two cases below.
\subsection*{Constructing log-symplectic structures}
As we discussed before, log-Poisson structures can be viewed as symplectic forms in a Lie algebroid, called the log-tangent bundle. Such \blog{}s can be constructed using \autoref{thm:introasymp}, as we do in Chapter \ref{chap:constructingblogs}. The main result we state in this introduction is the following, whereby we recover a result of Cavalcanti \cite{Cavalcanti17} using our general setup. The following is proven in this thesis as \autoref{thm:alflogsymp}.
\begin{thmalpha}[\cite{Cavalcanti17}] \label{thm:introblogalf} Let $f\colon X^4 \to \Sigma^2$ be an \alf between compact connected manifolds. Assume that the generic fiber $F$ is orientable and $[F] \neq 0 \in H_2(X;\R)$. Then $X$ admits a \blog.
\end{thmalpha}
Note in particular how the assumptions are phrased completely without using the language of Lie algebroids. We mentioned before that one can obtain a \fsymp out of an \alf, but from our point of view this is due to the theorem above and the fact that any \blog induces a \fsymp, see \autoref{thm:blogfolded}. This result is known in the Poisson community, but a proof is included in this thesis for completeness (a similar proof in the orientable case is given in \cite{GuilleminMirandaWeitsman15}).
Further, we show when one can obtain a \blog out of a \fsymp as \autoref{thm:foldedblog}, establishing a converse. In this sense the above result shows that \alf{}s are the correct type of map to be linked to \blog{}s.

\subsection*{Constructing \sgcs{}s}
In a seemingly different vein, other types of geometric structures can also be studied using the language of $\mc{A}$-symplectic geometry. Namely, certain types of \gcs{}s \cite{Hitchin03, Gualtieri11} which are called stable \cite{CavalcantiGualtieri15, Goto16}, are (almost) completely captured by their underlying Poisson bivector. This Poisson structure degenerates on a codimension-two subspace $D \subseteq X$ and is called an elliptic Poisson structure. The subspace allows for the definition of a corresponding Lie algebroid called the elliptic tangent bundle. As in the log-tangent case, an elliptic Poisson structure can be lifted nondegenerately to the elliptic tangent bundle, and we can then apply the machinery of Theorem \ref{thm:introasymp} above. Indeed, we introduce a class of maps called boundary maps, and consider their Lefschetz-analogue. This results in the following theorem, proven in Chapter \ref{chap:constructingsgcs} as \autoref{thm:lfibrationsgcs}, which has appeared before as \cite[Theorem 1.2]{CavalcantiKlaasse17}.
\begin{thmalpha}\label{thm:introsgcs} Let $f\colon (X^4,D) \to (\Sigma,Z)$ be a boundary Lefschetz fibration such that $(\Sigma,Z)$ carries a \blog. Assume that ${[F] \neq 0 \in H_2(X \setminus D;\R)}$ for the generic fiber $F$ and that $D$ is coorientable.. Then $(X,D)$ admits a \sgcs.
\end{thmalpha}
Note again how the statement of the theorem can be phrased without the language of Lie algebroids (we can readily determine topologically which surface pairs $(\Sigma,Z)$ admit a \blog). We classify boundary Lefschetz fibrations over the disk in \cite{BehrensCavalcantiKlaasse17}, and only state the result in this thesis as \autoref{thm:classbdrylefschetz}.
\subsection*{A general framework for lifting}
We hope to have convinced the reader that using Lie algebroids to study geometric structures is a fruitful undertaking. Another goal of this thesis is to provide a general framework in which this can be done. We introduce the concept of a divisor, which serves to capture the desired singularities. We then consider Lie algebroids which adhere to those singularities, and introduce the Lie algebroid analogues of the geometric structures we want to study. We then define a class of geometric structures we call of divisor-type, and show they can be lifted to the aforementioned Lie algebroids.

This type of undertaking invites for a more thorough understanding of the constituent objects themselves. The first six chapters of this thesis provide several results in this direction.
\subsection*{Homotopical obstructions}
Apart from constructing geometric structures using their $\mc{A}$-symplectic reformulations, this point of view can also be used to obstruct their existence. Indeed, in Chapter \ref{chap:homotopicalobstrs} we use the homotopical consequences on the bundle $\mc{A}$ to admit an $\mc{A}$-symplectic structure to determine necessary criteria for existence. In particular, for a given pair $(X,Z)$ where $Z \subseteq X$ is a hypersurface, we define the discrepancy $f(X,Z)$ (measuring the difference in Euler class between $\mc{A}_Z$ and $TX$), and establish the following as \autoref{cor:azsympstr}.
\begin{thmalpha}\label{thm:introhomtopobstr} Let $X$ be a compact oriented four-manifold that admits a log-symplectic structure with singular locus $Z$. Then $b_2^+(X) + b_1(X) + f(X,Z)$ is odd.
\end{thmalpha}
As the techniques directly make use of the Lie algebroids $\mc{A}$, the obstructions we obtain depend not only on the manifold $X$, but also on the prescribed singular loci.
\subsection*{Splitting theorems}
Finally, we study the intermediate objects that are used, namely $\mc{A}$-Lie algebroids, for their own sake. In Chapter \ref{chap:homotopicalobstrs} we apply the techniques of Bursztyn--Lima--Meinrenken \cite{BursztynLimaMeinrenken16} to obtain splitting results for $\mc{A}$-Lie algebroids. Such splitting theorems exist for various kinds of geometric structures, such as Lie algebroids, Poisson structures, and Dirac structures (see the introduction to \cite{BursztynLimaMeinrenken16} for an overview). They provide a useful local model, showing the geometric structure is locally equivalent to its linearized version. The following is proven in this thesis (in more precise form) as \autoref{thm:aanchored} and \autoref{thm:aliealgebroid}.
\begin{thmalpha}\label{thm:introsplitting} Let $E_\mc{A} \to X$ be either an $\mc{A}$-involutive $\mc{A}$-anchored vector bundle, or an $\mc{A}$-Lie algebroid. Moreover, let $i\colon N \hookrightarrow X$ be an $\mc{A}$-transversal for $E_\mc{A}$. Then $E_\mc{A}$ is isomorphic to its $\mc{A}$-linear approximation in a neighbourhood of $N$.
\end{thmalpha}

\section*{Main contributions of this thesis}
\bi
	\item A framework for mildly singular Poisson and Dirac structures (Chapters \ref{chap:divisors} - \ref{chap:diracgeometry});
	\item Adapting the methods of Gompf--Thurston to Lie algebroids (Chapter \ref{chap:constrasymp});
	\item Constructions of log-symplectic structures (Chapter \ref{chap:constructingblogs});
	\item Constructions of \sgcs{}s (Chapter \ref{chap:constructingsgcs});
	\item New obstructions for \blog{}s in dimension four (Chapter \ref{chap:homotopicalobstrs});
	\item A splitting result for $\mc{A}$-Lie algebroids (Chapter \ref{chap:splittingtheorems}).
\ei
\section*{Overview per chapter}
We give a brief overview of the contents of each chapter, emphasizing the main results and concepts introduced.
\subsection*{\hyperref[chap:divisors]{Chapter 1: Divisors}}
In this chapter we introduce the fundamental concept of a divisor. These are line bundles equipped with a section whose zero set is nowhere dense. They specify locally principal ideals we call divisor ideals, which can be thought of as characterizing what type of singularities one is dealing with. We establish several results that are used throughout this thesis, such as local forms for specific examples of divisors, including log and elliptic divisors.

\subsection*{\hyperref[chap:liealgebroids]{Chapter 2: Lie algebroids}}
In this chapter we review the concept of a Lie algebroid $\mc{A} \to X$, and introduce those notions from Lie algebroid theory that are relevant to this thesis, including cohomology, morphisms, and degeneracy loci. In particular we define the isomorphism locus of $\mc{A}$, which is the subset of $X$ over which $\mc{A}$ is isomorphic to $TX$. Moreover, we define the concept of an $\mc{A}$-Lie algebroid, which is a Lie algebroid whose anchor factors through that of $\mc{A}$. We discuss how the process of rescaling can be used to obtain new Lie algebroids, and how in certain cases one can define residue maps for $\mc{A}$-forms.

\subsection*{\hyperref[chap:concreteliealgebroids]{Chapter 3: Concrete Lie algebroids}}
In this chapter we discuss concrete classes of Lie algebroids, focusing on those obtained using rescaling, and on a class we call ideal Lie algebroids. This latter class is defined using their sheaf of sections, as containing those vector fields which interact with a given ideal. These are particularly interesting when using divisor ideals, then adhering to the singularities described by that ideal. Both of these classes of Lie algebroids have dense isomorphism locus. We discuss several important examples of such Lie algebroids, such as the log-tangent bundle and elliptic tangent bundle.

\subsection*{\hyperref[chap:asymplecticstructures]{Chapter 4: $\mc{A}$-symplectic structures}}
In this chapter we introduce the concept of an $\mc{A}$-symplectic structure, which is a symplectic object built using Lie algebroids. When $\mc{A}$ is one of the Lie algebroids from the previous chapter, these can be thought of as singular symplectic structures, with their singularities prescribed by $\mc{A}$ (hence, possibly, by the underlying divisor ideal). We establish Moser-type results for such structures and discuss concrete examples, such as log-symplectic structures.

\subsection*{\hyperref[chap:apoissonstructures]{Chapter 5: $\mc{A}$-Poisson structures}}
In this chapter we consider $\mc{A}$-Poisson structures, which form the Lie algebroid version of a Poisson structure. Moreover, we consider Poisson structures of divisor-type, which have singularities governed by a divisor. We further discuss the proces of lifting ($\mc{A}$-)Poisson structures, whereby these singularities are absorbed into those of the Lie algebroid. We discuss several examples of $\mc{A}$-Poisson structures, including log-Poisson and elliptic Poisson structures. These can be lifted to the log- and elliptic tangent bundles and then become nondegenerate, hence dual to their respective $\mc{A}$-symplectic counterparts.

\subsection*{\hyperref[chap:diracgeometry]{Chapter 6: Dirac structures}}
In this chapter we recall the notion of a Dirac structure, which captures both the study of closed two-forms and Poisson bivectors, as well as that of involutive distributions. This chapter serves as a prerequisite for the next, which uses the language of Dirac geometry we discuss here. In analogy with the previous chapter we moreover introduce the $\mc{A}$-analogue of a Dirac structure, and discuss how they too interact with divisors, and can be lifted to Lie algebroids.

\subsection*{\hyperref[chap:gencomplexstrs]{Chapter 7: Generalized complex structures}}
In this chapter we consider generalized complex structures $\mc{J}$, which are a generalization of both complex and symplectic structures. We introduce the required language and review known results. Underlying any \gcs{} is a Poisson bivector, whose rank determines the type of $\mc{J}$, colloquially measuring the amount of symplectic directions. We are particularly interested in those \gcs{}s which are called stable, being almost-everywhere symplectic. These are defined using the language of divisors. Through their underlying Poisson bivector, these can alternatively be described as certain kinds of symplectic structures in the elliptic tangent bundle.

\subsection*{\hyperref[chap:constrasymp]{Chapter 8: Constructing $\mc{A}$-symplectic structures}}
In this chapter we extend Gompf--Thurston techniques for symplectic Lefschetz fibrations to the world of Lie algebroids. We start by reviewing Lefschetz fibrations and how they correspond to symplectic structures in dimension four. We then introduce the notion of a Lie algebroid Lefschetz fibration and establish general results for when these can be used to construct an $\mc{A}$-symplectic structure on their total space. In particular, we prove \autoref{thm:introasymp}.

\subsection*{\hyperref[chap:constructingblogs]{Chapter 9: Constructing log-symplectic structures}}
In this chapter we use results from the previous chapter to construct \blog{}s, by viewing them as symplectic structures for the log-tangent bundle. The case of the log-tangent bundle is particularly nice as the hypotheses of \autoref{thm:introasymp} can be traced back to not use the language of Lie algebroids. In particular, we prove \autoref{thm:introblogalf} on how \alf{}s provide \blog{}s.

\subsection*{\hyperref[chap:constructingsgcs]{Chapter 10: Constructing \sgcs{}s}}
In this chapter we use results of Chapter \ref{chap:constrasymp} to construct \sgcs{}s. This is done by constructing certain kinds of elliptic symplectic structures. The case of the elliptic tangent bundle is more involved than that of the log-tangent bundle, yet the hypotheses of \autoref{thm:introasymp} are satisfied. In particular, we prove \autoref{thm:introsgcs} relating boundary Lefschetz fibrations to \sgcs{}s. We further discuss how to obtain such maps, and which manifolds admit them.

\subsection*{\hyperref[chap:homotopicalobstrs]{Chapter 11: Homotopical obstructions for $\mc{A}$-\symp{}s}}
In this chapter we tackle a problem that is separate from the previous three chapters, namely obtaining obstructions to the existence of an $\mc{A}$-symplectic structure. This is done using consequences of admitting an $\mc{A}$-\symp, namely the bundle $\mc{A}$ must be both orientable and admit an $\mc{A}$-\acs. We compute the required characteristic classes for certain Lie algebroids including the log-tangent bundle. In particular, we prove \autoref{thm:introhomtopobstr}, which obstructs the existence of log-symplectic structures in dimension four.

\subsection*{\hyperref[chap:splittingtheorems]{Chapter 12: Splitting theorems for $\mc{A}$-Lie algebroids}}
In the final chapter to this thesis we describe work in progress on obtaining splitting results for $\mc{A}$-Lie algebroids following the methods of Bursztyn--Lima--Meinrenken. These splitting results help elucidate the local structure of such objects, and form the stepping stone for similar results for $\mc{A}$-Dirac and $\mc{A}$-Poisson structures. In particular, we prove \autoref{thm:introsplitting} mentioned before.

\mainmatter
\numberwithin{figure}{chapter}
\numberwithin{table}{chapter}
\chapter{Divisors on smooth manifolds}
\label{chap:divisors}
\fancypagestyle{empty}{%
	\fancyhf{}%
	\renewcommand\headrulewidth{0pt}%
	\fancyhead[RO,LE]{\thepage}%
}
In this short chapter we develop the theory of real and complex divisors on smooth manifolds. These are extensions to the smooth setting of the notion of a divisor used in complex geometry. See also \cite{CavalcantiGualtieri15, VanderLeerDuran16}. The theory of divisors permeates much of this thesis. We will see in Section \ref{sec:idealliealgebroids}, Section \ref{sec:poissondivisors} and Section \ref{sec:diracdivisors} that divisors give rise to Lie algebroids, and allow us to define and study interesting classes of Poisson and Dirac structures. Moreover in Section \ref{sec:scgs} we define the notion of a \sgcs. This is done using their spinor description by demanding that its anticanonical bundle is a specific type of divisor. Part of the results in this chapter have appeared before in \cite{CavalcantiKlaasse17} and are joint with Gil Cavalcanti.
\subsection*{Organization of the chapter}
This chapter is built up as follows. In Section \ref{sec:divisors} we define divisors as real or complex line bundles equipped with almost-everywhere nonzero sections. We further discuss the relation between the ideals they give rise to via evaluation of the sections. In Section \ref{sec:morphdivisors} we discuss morphisms of divisors, with Section \ref{sec:divproduct} treating two types of operations one can perform with divisors. We end this chapter with Section \ref{sec:divexamples} in which we discuss important examples of divisors that will be used in later chapters.

\section{Divisors}
\label{sec:divisors}
We start by defining the objects we are interested in. Throughout, let $X$ be a manifold.
\begin{defn} A \emph{divisor} on $X$ is a pair $(U,\sigma)$ where $U \to X$ is a real or complex line bundle and $\sigma \in \Gamma(U)$ is a section whose zero set $Z_\sigma$ is nowhere dense.
\end{defn}
\begin{rem} A nowhere dense subset of a topological space is one whose closure has empty interior. In other words, its intersection with any nonempty open subset is not dense. Often $Z_\sigma$ will be a (union of) submanifold(s) of positive codimension.
\end{rem}
We mostly focus on real divisors. As such, we will often drop the prefix `real', while explicitly stating when divisors are instead complex. Examples of divisors will be discussed in Section \ref{sec:divexamples}. While the zero set $Z_\sigma \subseteq X$ is closed, it is not assumed that $Z_\sigma$ is smooth. Let $(U,\sigma)$ be a divisor. The evaluation $\sigma\colon \Gamma(U^*) \to C^\infty(X)$ determines a specific type of ideal $I_\sigma := \sigma(\Gamma(U^*)) \subseteq C^\infty(X)$.
\begin{defn} An ideal $I \subseteq C^\infty(X)$ on $X$ is called a \emph{divisor ideal} if it is locally principal and locally generated by a function with nowhere dense zero set.
\end{defn}
The name is justified further below (see \autoref{prop:locprincideal} and \autoref{cor:dividealcorresp}). From the definition it is clear that divisors give rise to divisor ideals.
\begin{exa} Let $U$ be the trivial line bundle with $\sigma \in \Gamma(U)$ nonvanishing. Then $Z_\sigma$ is empty, and $(U,\sigma)$ is called the \emph{trivial divisor} on $X$. In this case, $I_\sigma = C^\infty(X)$.
\end{exa}
Depending on whether $U$ is real or complex, the ideal $I_\sigma$ will be either a real or a complex ideal. Letting $\alpha$ be a local trivialization of $U^*$, we have $\alpha(\sigma) = g$ for some local function $g$. Then locally $I_\sigma = \langle \alpha(\sigma) \rangle = \langle g \rangle$. Conversely, out of any such ideal we can construct a divisor, which recovers the ideal via this evaluation process. This extends the correspondence between Cartier divisors and holomorphic line bundles in complex geometry.
\begin{prop}[\cite{VanderLeerDuran16}]\label{prop:locprincideal} Let $I$ be divisor ideal on $X$. Then there exists a divisor $(U_I,\sigma)$ on $X$ such that $I_\sigma = I$.
\end{prop}
In the statement above we are ambiguous about whether we are dealing with real or complex objects. However, the argument given below works in either case. Note that we treat ideals $I$ as specifying ideal sheaves, so that $I(U)$ for an open $U \subseteq X$ denotes all functions defined on $U$ that belong to $I$.
\bp Let $\{U_\alpha\}$ be an open cover of $X$ and $f_\alpha \in I(U_\alpha)$ be generators. Then on $U_\alpha \cap U_\beta$ we have $f_\alpha = g_{\alpha \beta} f_\beta$ with $g_{\alpha \beta} \in C^\infty(U_\alpha \cap U_\beta)$ some nonvanishing function. Consequently $f_\alpha = g_{\alpha \beta} g_{\beta \alpha} f_\alpha$ and similarly $f_\alpha = g_{\alpha \beta} g_{\beta \gamma} g_{\gamma \alpha} f_\alpha$ on $U_\alpha \cap U_\beta \cap U_\gamma$. Using that $f_\alpha$ has nowhere dense zero set, the functions $g_{\alpha \beta}$ are seen to satisfy the cocycle condition. We conclude that $\{(U_{\alpha \beta}, g_{\alpha \beta})\}$ defines a line bundle $U_I$ on $X$. Further, $\sigma|_{U_\alpha} = f_\alpha$ on $U_\alpha$ specifies a section $\sigma$ of $U_I$ with the desired properties.
\ep
The section $\sigma$ constructed in the proof of the above proposition is unique up to multiplication by a smooth nonvanishing function. We now discuss the relation between $I_\sigma$ and $I_{Z_\sigma}$, the vanishing ideal of the zero set $Z_\sigma$ of the divisor.
\begin{prop}\label{prop:divvanishingideal} Let $(U,\sigma)$ be a divisor. Then $I_\sigma \subseteq I_{Z_\sigma}$.
\end{prop}
\bp Let $\alpha^*$ be a local trivialization of $U^*$ and $g \in C^\infty(X)$ such that $\alpha(\sigma) = g$. Then $I_\sigma = \sigma(\Gamma(U^*)) = \langle g \rangle$. As $\sigma^{-1}(0) = Z_\sigma$ we have $\sigma(Z_\sigma) = 0$. We conclude that $g(Z_\sigma) = \sigma(Z_\sigma) = 0$ so that $g \in I_{Z_\sigma}$, hence $I_\sigma \subseteq I_{Z_\sigma}$ as desired. 
\ep
Note that the equality $I_\sigma = I_{Z_\sigma}$ is only possible for $Z_\sigma$ smooth if $Z_\sigma$ is of codimension one. Namely, divisor ideals are locally principal, while vanishing ideals have local generators equal to the associated codimension of the submanifold.
\section{Morphisms of divisors}
\label{sec:morphdivisors}
In this section we discuss morphisms of divisors. We will restrict ourselves to morphisms from real to real, or from complex to complex divisors. Morphisms between divisors are most neatly described in terms of the divisor ideals they give rise to. Denote by $f^*I \subseteq C^\infty(X)$ the ideal generated by the pullback of an ideal $I \subseteq C^\infty(Y)$ along a map $f\colon X \to Y$.
\begin{defn} Let $(U_X,\sigma_X)$ be a divisor on $X$ and $(U_Y,\sigma_Y)$ a divisor on $Y$. A map $f\colon X \to Y$ is a \emph{morphism of divisors} if $f^* I_{\sigma_Y} = I_{\sigma_X}$.
\end{defn}
Alternatively, it suffices to require that $(U_X,  \sigma_X) = (f^* U_Y, g f^* \sigma_Y)$ for some $g \in C^\infty(X;\R^*)$ or $C^\infty(X;\C^*)$, as the sections can be multiplied by nonvanishing functions without changing the ideals they generate. Two divisors are \emph{diffeomorphic} (denoted using $\cong$) if there exists a morphism of divisors between them which is in addition a diffeomorphism. Two divisors on a fixed manifold $X$ are \emph{isomorphic} (denoted using $=$) if the identity map on $X$ is a morphism of divisors. The proper conclusion of \autoref{prop:locprincideal} is the following.
\begin{cor}\label{cor:dividealcorresp} There is a bijective correspondence between divisor ideals and isomorphism classes of divisors.
\end{cor}
There is a relation between morphisms of divisors and their underlying zero sets.
\begin{defn}\label{defn:pairsandmaps} A \emph{pair} $(X,Z)$ is a manifold $X$ with a subset $Z \subseteq X$. A \emph{map of pairs} $f\colon (X,Z) \to (X',Z')$ is a smooth map $f\colon X \to X'$ for which $f(Z) \subseteq Z'$. A \emph{strong map of pairs} is a map of pairs $f\colon (X,Z) \to (X,Z')$ for which $f^{-1}(Z') = Z$.
\end{defn}
\begin{prop}\label{prop:divmorphstrong} Let $f\colon (X,U_X,\sigma_X) \to (Y,U_Y,\sigma_Y)$ be a morphism of divisors. Then $f\colon (X,Z_{\sigma_X}) \to (Y,Z_{\sigma_Y})$ is a strong map of pairs.
\end{prop}
\bp As $f$ is a morphism of divisors, we have that $U_X = f^* U_Y$ and $\sigma_X = g f^* \sigma_Y$ for a nonvanishing function $g$. As $Z_{\sigma_X}$ by definition is the set on which $\sigma_X$ vanishes, we see immediately from this that it is given by $f^{-1}(Z_{\sigma_Y})$.
\ep
\begin{rem} The previous proposition uses something special about the zero sets of divisors. In general for a map $f\colon X \to Y$ and $Z \subseteq Y$ closed, one can only conclude that $f^* I_Z \subseteq I_{f^{-1}(Z)}$. However, in our case we in fact have equality: a rephrasing of \autoref{prop:divmorphstrong} is that a morphism of divisors is also a morphism of ideals between the vanishing ideals of the zero sets, i.e.\ satisfies $f^* I_{Z_{\sigma_Y}} = I_{Z_{\sigma_X}}$.
\end{rem}
Moreover, the Stiefel--Whitney classes (e.g.\ \cite{MilnorStasheff74}) of the divisors are related.
\begin{prop} Let $f\colon (X,U_X,\sigma_X) \to (Y,U_Y,\sigma_Y)$ be a morphism of divisors. Then $w_1(U_X) = f^* w_1(U_Y) \in H^1(X;\Z_2)$.
\end{prop}
When dealing with complex divisors, the same statement is true for their first Chern classes. Note that given a smooth map $f\colon X \to Y$ and a divisor $(U_Y,\sigma_Y)$ on $Y$, one can equip $X$ with a divisor by setting $U_X = f^* U_Y$ and $\sigma_X = f^* \sigma_Y$, as long as $f^{-1}(Z_{\sigma_Y})$ is nowhere dense. It is then automatic that $f$ is a morphism of divisors. We will use this idea later in Section \ref{sec:boundarymaps}.
\section{Products of divisors}
\label{sec:divproduct}
In this section we treat the operation of taking the product of divisors, which is given by the tensor product of bundles and sections. Moreover, we discuss briefly the direct sum of divisors. As with morphisms of divisors, we will only consider products of real divisors, or of complex divisors.
\begin{defn} Let $(U,\sigma)$ and $(U',\sigma')$ be divisors on $X$. Then their \emph{product} is the divisor $(U \otimes U',\sigma \otimes \sigma')$ on $X$, with $Z_{\sigma \otimes \sigma'} = Z_\sigma \cup Z_{\sigma'}$.
\end{defn}
It is immediate that the product is indeed again a divisor, as the union of nowhere dense subsets is nowhere dense. Upon introduction of the shorthand $Z_\sigma = (U,\sigma)$, we will sometimes write the product of divisors additively, i.e.\ we will write $Z_\sigma + Z_{\sigma'}$. This is inspired by the notation for divisors used in algebraic geometry. It is not hard to determine what happens to the associated divisor ideals.
\begin{prop}\label{prop:proddivisor} Let $(U,\sigma)$ and $(U',\sigma')$ be divisors on $X$ with divisor ideals $I_\sigma$ and $I_{\sigma'}$ respectively. Then the product divisor $(U \otimes U',\sigma \otimes \sigma')$ has divisor ideal $I_{\sigma\otimes \sigma'} = I_\sigma \cdot I_{\sigma'}$, the product of the ideals, and $w_1(U\otimes U') = w_1(U) + w_1(U')$.
\end{prop}
\bp Choose local trivializations $\alpha^*, \alpha'^*$ of $U^*$ and $U'^*$ with $\alpha^*(\sigma) = g$ and $\alpha'^*(\sigma') = g'$, so that $I_\sigma = \langle g \rangle$ and $I_{\sigma'} = \langle g' \rangle$ for functions $g, g'$. Then $\alpha^* \otimes \alpha'^*$ trivializes $U \otimes U'$, and it is immediate that $(\alpha^* \otimes \alpha'^*)(\sigma \otimes \sigma') = g \cdot g'$ so that $I_{\sigma \otimes \sigma'} = \langle g \cdot g' \rangle$. The statement regarding $w_1$'s is standard.
\ep
This result shows the product of divisors descends to a product on isomorphism classes of divisors, using \autoref{cor:dividealcorresp}. Products of divisors provide a convenient way to describe divisors with disconnected zero sets, as these can often be decomposed into constituent parts. In particular this is the case when the zero set consists of submanifolds of varying codimension.
\begin{rem} Let $(U,\sigma)$ be a divisor on $X$ and let $(\underline{\R}, \underline{1})$ be the trivial divisor. Then the product $(U \otimes \underline{\R}, \sigma \otimes \underline{1})$ is naturally isomorphic to $(U,\sigma)$ as a divisor. This makes the trivial divisor the unit for the product, as is clear on the level of divisor ideals.
\end{rem}
Divisors admit another type of product. Let $X$ and $X'$ be manifolds and consider the product $X \times X'$ with projections $p_{X}\colon X \times X' \to X$ and $p_{X'}\colon X \times X' \to X'$. Consider two divisors $(U,\sigma) \to X$ and $(U',\sigma') \to X'$.
\begin{defn}\label{defn:divdirectstum} The \emph{external tensor product} of $(U,\sigma)$ and $(U',\sigma')$ is the divisor $(p_X^* U \otimes p_{X'}^* U', p_X^* \sigma \otimes p_X^* \sigma')$ on $X \times X'$, with zero set $Z_\sigma \times X' \cup X \times Z_{\sigma'}$.
\end{defn}
This product will be used in Section \ref{sec:poissondivisors} when discussing Poisson structures.
\section{Examples}
\label{sec:divexamples}
In this section we discuss examples of divisors and establish some of their properties. Of particular interest will be to understand whether we can linearize the divisors, i.e.\ express their divisor ideal as generated by a homogeneous function in normal coordinates to their zero set. This is because of our desire to define Lie algebroids using divisors, as we will explore in Section \ref{sec:idealliealgebroids}. 
\subsection{Log divisors}
\label{sec:logdivisor}
The simplest type of nontrivial divisor is that of a log divisor. These are divisors whose zero set is a smooth codimension-one submanifold, or \emph{hypersurface}. The choice of name is inspired by the name given to the Lie algebroids and Poisson structures it gives rise to (see Sections \ref{sec:logtangentbundle}, \ref{sec:logsympstr} and \ref{sec:logpoisson}).
\begin{defn}\label{defn:logdivisor} A \emph{log divisor} is a divisor $Z = (L,s)$ such that $s$ is transverse to the zero section.
\end{defn}
Given a log divisor $(L,s)$, it follows that its zero set $Z = s^{-1}(0)$ is a smooth hypersurface. Recall that a section $s \in \Gamma(L)$ defines a map $s\colon X \to L$ (with $L$ viewed as a manifold). Hence by applying the tangent functor we obtain a map $Ts\colon TX \to TL$. Over the zero section of $L$, its tangent space canonically splits as $TL|_X \cong TX \oplus L$. The \emph{intrinsic derivative} $ds|_Z\colon NZ \to L|_Z$ is the map obtained from $Ts$ by projecting onto the $L$-component in the splitting of $TL|_X$ above, with $NZ \to Z$ denoting the normal bundle of $Z$. Then, $s$ being transverse to the zero section is equivalent to the intrinsic derivative $ds|_Z$ being surjective. As both $NZ$ and $L|_Z$ are bundles of rank one, this shows the following.
\begin{prop}\label{prop:logdivnormalbundle} Let $Z = (L,s)$ be a log divisor. Then the intrinsic derivative $ds|_Z\colon NZ \to L|_Z$ is an isomorphism.
\end{prop}
We see that $Z$ is coorientable if and only if $L|_Z$ is trivial, and $Z$ having a global defining function is equivalent to $L$ being trivial. In fact, we can immediately determine the first Stiefel--Whitney class of a log divisor using the Poincar\'e dual with $\Z_2$-coefficients, as by definition $s$ is a section of $L$ with transverse zeros.
\begin{prop}\label{prop:logdivstiefelwhitney} Let $Z = (L,s)$ be a log divisor. Then $w_1(L) = {\rm PD}_{\Z_2}([Z])$.
\end{prop}
Let $(L,s)$ be a log divisor. The associated divisor ideal $I_s$ is exactly the vanishing ideal $I_Z$ of the hypersurface $Z$. In other words, we can find local coordinates $(z,x_2,\dots,x_n)$ around points in $Z$ with $\{z = 0\} = Z$ such that $s(z,x_2,\dots,x_n) = z$ and hence $I_Z = \langle z \rangle$. Any hypersurface $Z$ naturally gives rise to a unique log divisor, so that we will often identify a log divisor with the associated submanifold $Z$.
\begin{prop}\label{prop:hypersurfacedivisor} Let $X$ be a manifold and $Z \subset X$ a hypersurface. Then $Z$ carries a unique log divisor structure, i.e.\ there exists a unique isomorphism class of a log divisor $(L,s)$ with $s^{-1}(0) = Z$.
\end{prop}
\bp Apply \autoref{prop:locprincideal} to the vanishing ideal $I_Z$, giving a divisor $(U_{I_Z}, \sigma) =: (L,s)$. The section $s$ vanishes transversely along $Z$ as it is equal to a local defining function for $Z$ in any trivializing open $U_\alpha$ of $L$ containing $Z$. We conclude that $(L,s)$ is a log divisor. Unicity of its isomorphism class is addressed in \autoref{cor:dividealcorresp}.
\ep
Due to the above result we will sometimes talk about \emph{log pairs} $(X,Z)$ as specifying log divisors, and denote the line bundle of a log divisor by $L_Z$.
\begin{rem} In light of the names associated to the Lie algebroid and Poisson structures that are defined using log divisors, it is customary to call log pairs also \emph{$b$-manifolds} \cite{Melrose93}. We will mostly refrain from doing so, except in Chapter \ref{chap:constructingblogs}.
\end{rem}
We can characterise which maps are morphisms between log divisors.
\begin{defn} A map of pairs $f\colon (X,Z_X) \to (Y,Z_Y)$ where $Z_Y \subseteq Y$ is a submanifold is \emph{transverse} if $f$ is transverse to $Z_Y$.
\end{defn}
\begin{prop}\label{prop:bmapideal} Let $(X,Z_X)$ and $(Y,Z_Y)$ be log divisors. Then a map $f\colon X \to Y$ is a morphism of log divisors if and only if $f\colon (X,Z_X) \to (Y,Z_Y)$ is a transverse strong map of pairs.
\end{prop}
\bp We start with the direct implication. If $f$ is a morphism of divisors, it is a strong map of pairs by \autoref{prop:divmorphstrong}. To show that $f$ is transverse to $Z_Y$, let $z_Y$ be a local defining function for $Z_Y$, so that $z_Y \in I_{Z_Y}$ is a local generator. Then $f^* z_Y \in f^* I_{Z_Y}$ is a local generator over $C^\infty(X)$, and by equality of ideals also $f^* z_Y \in I_{Z_X}$ (and generates it). Let $z_X$ be a defining function for $Z_X$. The equality of ideals implies that $z_X = \varphi f^* z_Y$ for $\varphi \in C^\infty(X)$ nonvanishing. Then $d z_X = \varphi f^*(d z_Y) + d \varphi f^* z_Y$. Evaluating at $Z_X$ we see that the second term vanishes as $f^* z_Y = 0$ on $Z_X$. Hence on $Z_X$ we have $0 \neq d z_X = \varphi f^*(dz_Y)$. As $\varphi$ is nonvanishing, this implies $f^* dz_Y = \frac{dz_X}{\varphi} \neq 0$, so $f^* \neq 0$ on $Z_Y$, hence $f$ is transverse to $Z_Y$.

Conversely, assume that $f$ is a transverse strong map of pairs. Let $z_Y$ be a local defining function for $Z_Y$, so that $z_Y \in I_{Z_Y}$ is a local generator. Then $f^* z_Y$ is a local generator for $f^* I_{Z_Y}$. However, transversality gives that $d f^*(z_Y)$ = $d(z_Y \circ f) = dz_Y \circ Tf \neq 0$ in the normal direction, so that $f^* z_Y$ is a local defining function for $f^{-1}(Z_Y)$, hence generates $I_{f^{-1}(Z_Y)}$. We conclude that $f^* I_{Z_Y} = I_{f^{-1}(Z_Y)}$. As $f^{-1}(Z_Y) \subseteq Z_X$, by \autoref{lem:subsetideal} below we have $I_{Z_X} \subseteq I_{f^{-1}(Z_Y)} = f^* I_{Z_Y}$. Now note that also $Z_X \subseteq f^{-1}(Z_Y)$. Again by \autoref{lem:subsetideal}, $I_{f^{-1}(Z_Y)} \subseteq I_{Z_X}$. Using \autoref{lem:pullbackideal} below we have $f^* I_{Z_Y} \subseteq I_{f^{-1}(Z_Y)} \subseteq I_{Z_X}$, hence equality.
\ep
\begin{rem}\label{rem:tsmopbmaps} Transverse strong maps between log pairs are also called \emph{$b$-maps} between $b$-manifolds \cite{Melrose93,GuilleminMirandaPires14}. Using this terminology, \autoref{prop:bmapideal} can be stated by saying that $f$ is a $b$-map if and only if it is a morphism of log divisors.
\end{rem}
\begin{lem}\label{lem:subsetideal} Let $S, S' \subseteq X$ be closed. Then $I_S \subseteq I_{S'}$ if and only if $S' \subseteq S$.
\end{lem}
\bp Let $g \in I_S$ and $x \in S'$. Then $x \in S$, so $g(x) = 0$, so $g \in I_{S'}$, hence $I_S \subset I_{S'}$. Similarly, if $I_S \subset I_{S'}$, take $x \in S'$ and arguing by contradiction, assume that $x \not\in S$. Then as $S$ and $S'$ are closed there exists $h \in I_S$ such that $h(x) \neq 0$. But $I_S \subset I_{S'}$ so $h \in I_{S'}$, hence $0 = h(x) \neq 0$. This is a contradiction, so $x \in S$, so that $S' \subset S$.
\ep
\begin{lem}\label{lem:pullbackideal} Let $f\colon X \to Y$ be a map and $Z \subseteq Y$ closed. Then $f^* I_Z \subseteq I_{f^{-1}(Z)}$.
\end{lem}
\bp Let $g \in f^* I_Z$. Then $g = \sum_i c_i f^*(h_i)$ for some $h_i \in I_Z$ and functions $c_i$. Take $x \in f^{-1}(Z)$. Then $f(x) =: y \in Z$. Hence $g(x) = \sum_i c_i(x) (f^* h_i)(x) = \sum_i c_i(x) h_i(f(x)) = 0$ as $h_i(y) = 0$ for all $i$, so that $g \in I_{f^{-1}(Z)}$.
\ep
\subsection{Complex log divisors}
\label{sec:cplxlogdivisor}
There is also a complex analogue of a (real) log divisor, which will be used in the study of \sgcs{}s (see Section \ref{sec:scgs}). Fix a manifold $X$.
\begin{defn} A \emph{complex log divisor} on $X$ is a complex divisor $D = (U,\sigma)$ whose zero set $D = \sigma^{-1}(0)$ is a smooth codimension-two submanifold with $\sigma$ transverse to the zero section.
\end{defn}
In this case we can immediately determine the first Chern class of $U$.
\begin{prop} Let $D = (U,\sigma)$ be a complex log divisor on $X$. Then $c_1(U) = {\rm PD}_\Z([D]) \in H^2(X;\Z)$, the Poincar\'e dual with $\Z$-coefficients.
\end{prop}
The Poincar\'e dual with $\Z$-coefficients requires the choice of a coorientation for $D$. Now, the intrinsic derivative $d\sigma|_D\colon ND \to U|_D$ (defined as in Section \ref{sec:logdivisor}) is an isomorphism, equipping the normal bundle $ND$ with a complex structure and thus supplying $D$ with a coorientation. The complex ideal $I_\sigma \subseteq C^\infty(X;\C)$ is locally generated by a single complex function $w$ such that $\{w = 0\} = D$. One can view $I_\sigma$ as a complexified version of the vanishing ideal $I_D$ of $D$.

Any complex log divisor $D = (U,\sigma)$ determines a complex conjugate divisor $\overline{D} = (\overline{U}, \overline{\sigma})$ with the same zero locus. This is true for any complex divisor, with the conjugate being obtained (for example) by conjugating its local description in terms of trivializations. In particular, the complex conjugate log divisor has associated ideal $I_{\overline{\sigma}}$ generated locally by the complex function $\overline{w}$ instead.
\subsection{Elliptic divisors}
\label{sec:elldivisor}
In this section we explore a type of divisor with codimension-two zero set which are called elliptic divisors. These will be directly relevant to the study of \sgcs{}s undertaken in Section \ref{sec:scgs}.
\begin{defn} An \emph{elliptic divisor} is a divisor $|D| = (R,q)$ whose zero set $D = q^{-1}(0)$ is a smooth codimension-two critical submanifold of $q$ along which its normal Hessian ${\rm Hess}(q) \in \Gamma(D; {\rm Sym}^2 N^*D \otimes R)$ is definite.
\end{defn}
The normal Hessian of $q$ is the leading term of its Taylor expansion. Note that as $D$ has codimension equal to two, the intrinsic derivative $dq|_D$ is the zero map. We will refer to $|D|$ as the elliptic divisor, yet note there is more information available than just the zero set of $q$. The ideal $I_{|D|} := I_q$ is called an \emph{elliptic ideal}, and both $R$ and $q$ (up to a nonzero smooth function) can be recovered from $I_{|D|}$ by \autoref{prop:locprincideal}.
\begin{rem}\label{rem:ellipticomplexlog} Let $D = (U,\sigma)$ be a complex log divisor. Then its product with its conjugate $(U\otimes \overline{U}, \sigma \otimes \overline{\sigma})$ is invariant under conjugation (i.e.\ naturally attains a real structure, or complex antilinear involution), hence gives rise to a real divisor $(R,q) = ((U\otimes \overline{U})_\R, \sigma \otimes \overline{\sigma})$. This is an elliptic divisor $|D|$. Using \autoref{prop:ellmorsebott} below and the factorization $x^2 + y^2 = (x + iy)(x - iy) = w \overline{w}$, any elliptic divisor with coorientable zero set arises from a complex log divisor in this way, with $(U,\sigma)$ being determined up to diffeomorphism by the choice of coorientation.
\end{rem}
Note that $I_{|D|}$ is not the vanishing ideal $I_D$ of $D$, but instead is locally generated by an even index Morse--Bott function in coordinates normal to $D$, as we now explain.
\begin{defn} Let $g \in C^\infty(X)$ be given. A compact connected submanifold $S \subset X$ is a \emph{nondegenerate critical submanifold} of $g$ if $S \subset {\rm Crit}(g)$ and $\ker {\rm Hess}(g) = T_p S$ for all $p \in S$. If ${\rm Crit}(g)$ consists of nondegenerate critical submanifolds, then $g$ is a \emph{Morse--Bott function}.
\end{defn}
Let $g \in C^\infty(X)$ and $S \subseteq X$ be a nondegenerate critical submanifold of $g$. Consider the exact sequence $0 \to TS \to TX|_S \to NS \to 0$. For $p \in S$ we have that ${\rm Hess}(g)(p) \in {\rm Sym}^2 T_p^* X$, and $T_p S$ is contained in the kernel. But then ${\rm Hess}(g)(p) \in {\rm Sym}^2 N_p^* S$, giving rise to a nondegenerate bilinear form $Q_g \in \Gamma(S; {\rm Sym}^2 N^* S)$. The semi-global version of the Morse--Bott lemma says that $g$ is diffeomorphic to this quadratic approximation in a tubular neighbourhood of $S$.
\begin{lem}[{\cite[Proposition 2.6.2]{Nicolaescu11}}]\label{lem:morsebott} Let $g \in C^\infty(X)$ and $S$ be a nondegenerate critical submanifold of $g$. Then there exists a neighbourhood $U$ of the zero section $S \subset NS$ and an open embedding $\Phi\colon U \to X$ such that $\Phi|_S = {\rm id}_S$ and $\Phi^* g = Q_g$.
\end{lem}
Let $|D| = (R,q)$ be an elliptic divisor. Then $R$ is orientable by $q$ as it is a trivialization away from a codimension-two submanifold. Hence, $R$ is always trivial.
\begin{prop}\label{prop:ellipticdivstiefelwhitney} Let $(R,q)$ be an elliptic divisor on $X$. Then $w_1(R) = 0$.
\end{prop}
If one were to orient $R$ using $q$, the normal Hessian of $q$ along $D$ is positive definite. In other words, let $\alpha$ be a trivialization of $R^*$. Then $g:= \alpha(q) \in C^\infty(X)$ is a function with $g^{-1}(0) = D$ and ${\rm Hess}(g) = \alpha({\rm Hess}(q))$. Moreover, $D$ is a nondegenerate critical submanifold of $g$, and $g$ is locally Morse--Bott around $D$. As $D$ is codimension two, $X \setminus D$ is connected, so that the sign of $g$ on $X \setminus D$ is fixed. Replace $\alpha$ by $-\alpha$ if necessary so that this sign is positive, and then $q$ and $\alpha$ induce compatible orientations. Call such a trivialization $\alpha$ \emph{compatible} with $q$. For compatible trivializations we have $g \geq 0$ so that ${\rm Hess}(g)$ is positive definite. As a consequence of \autoref{lem:morsebott} we obtain the following.
\begin{prop}\label{prop:ellmorsebott} Let $|D| = (R,q)$ be an elliptic divisor. Then there exists a neighbourhood $U$ of the zero section $D \subset ND$ and an open embedding $\Phi\colon U \to X$ such that $\Phi|_D = {\rm id}_D$ and $(\Phi^*\alpha)(\Phi^* q) = Q_g$, where $g = \alpha(q) \in C^\infty(X)$.
\end{prop}
Using a compatible trivialization ensures the bilinear form $Q_g$ is positive definite. Given $p \in D$ we can locally trivialize the bundles $R$ and $ND$, so that using \autoref{prop:ellmorsebott} the section $q$ can be written locally as $q(x_1,\dots, x_n) = \pm(x_1^2 + x_2^2)$ in normal bundle coordinates such that $ND = \langle \partial_{x_1}, \partial_{x_2} \rangle$. Consequently, the elliptic ideal $I_{|D|}$ is locally generated by $r^2$, where $r^2 = x_1^2 + x_2^2$ is the squared radial distance from $D$ inside $ND$.
\begin{rem} While \autoref{prop:hypersurfacedivisor} shows that hypersurfaces carry a unique log divisor structure, the same is not true for codimension-two submanifolds and elliptic divisors. A simple example is provided by $X = \R^2$ with $D = \{(0,0)\}$ and coordinates $(x,y)$. Equip $D$ with the elliptic ideals $I = \langle x^2 + y^2 \rangle$ and $I' = \langle x^2 + 2y^2 \rangle$. As these ideals are distinct, they supply $D$ with two non-isomorphic elliptic divisor structures. However, it is easy to see that these elliptic divisors are diffeomorphic.
\end{rem}
\subsection{Normal-crossing log divisors}
\label{sec:normalcrosslogdivisor}
For our next example we discuss a more general type of log divisor where the zero set is allowed to exhibit normal-crossing behavior.
\begin{defn} A \emph{normal-crossing log divisor} is a divisor $(U,\sigma)$ where $\sigma$ at each point in its zero set $Z_\sigma = \sigma^{-1}(0)$ either vanishes transversally, or has nondegenerate indefinite normal Hessian.
\end{defn}
The zero set $Z_\sigma$ of a normal-crossing log divisor is the union $Z = 
\cup_i Z_i$ of a collection $\{Z_i\}_i$ of smooth hypersurfaces which all intersect transversally. Such a collection is usually called a \emph{normal-crossing divisor} (note the slight difference in names). To typographically distinguish them from log divisors, we sometimes denote normal-crossing log divisors by $\underline{Z} = (U,\sigma)$. In terms of the associated ideal, we have that locally $I_\sigma = \langle \prod_i z_i \rangle$, where each $z_i$ is a local defining function for $Z_i$. This follows from the Morse--Bott lemma, \autoref{lem:morsebott}, applied to $g = \alpha^*(\sigma)$ for $\alpha^*$ a local trivialization of $U$ around $Z_\sigma$. Note the following criterion characterizing normal-crossing divisors due to Saito \cite{Saito80}, phrased in the smooth setting (see \cite{GualtieriLiPelayoRatiu17}).
\begin{prop}\label{prop:saitocriterion} Let $Z = \cup_i Z_i$ be a collection of smooth hypersurfaces. Then $Z$ is a normal-crossing divisor if and only if around each point $x \in Z$ there are $n = |\{i\, | \, x \in Z_i\}|$ commuting vector fields $V_i$ tangent to $Z_i$ such that the determinant $V_1 \wedge \dots \wedge V_n$ vanishes precisely on $Z$ and transversally on its smooth locus.
\end{prop}
As alluded to in the local description given above, normal-crossing log divisors with zero set $Z = \cup_i Z_i$ are related to the log divisors of each $Z_i$ by means of the product operation that was discussed in Section \ref{sec:divproduct}.
\begin{prop}\label{prop:logprodnormalcrossing} Let $X$ be a manifold and $Z = (L,s)$ and $Z' = (L',s')$ be log divisors on $X$ such that $Z$ and $Z'$ intersect transversally. Then their product $(U,\sigma) := (L \otimes L', s \otimes s')$ is a normal-crossing log divisor.
\end{prop}
\bp It is immediate that $(U,\sigma)$ is a divisor. To see that it is a normal-crossing log divisor, note that the divisor ideal of the product is given by $I_s \cdot I_{s'}$, which locally around $Z \cap Z'$ (assuming this set is nonempty) is given by $\langle z z' \rangle$. Here $z$ and $z'$ are local defining functions for $Z$ and $Z'$ respectively. By an application of the Morse--Bott lemma, we realize that the divisor $(U,\sigma)$, which up to isomorphism can be recovered from $I_{s \otimes s'}$ using \autoref{prop:locprincideal}, must be a normal-crossing log divisor.
\ep
The above proposition can be applied repeatedly. More precisely, given a normal-crossing log divisor $(U,\sigma)$ with $Z_\sigma = \cup_i Z_i$ and $Z = (L,s)$ a log divisor which transversally intersects $Z_\sigma$, their product is again a normal-crossing log divisor. In this way we see that a normal-crossing log divisor is built up out of a collection of log divisors using the product operation. As it is standard to denote the zero set additively, we see it makes sense to write $Z_\sigma + Z_\sigma'$ as the zero set of the product of divisors.
\begin{cor} A normal-crossing divisor $Z$ carries a unique isomorphism class of normal-crossing log divisor structures.
\end{cor}
\bp This follows by applying \autoref{prop:hypersurfacedivisor} to each constituent hypersurface $Z_i$, to obtain log divisors $(L_i, s_i)$. Then, using \autoref{prop:logprodnormalcrossing} one takes the product of these divisors to form the normal-crossing log divisor $(U,\sigma) = (\otimes_i L_i, \otimes_i s_i)$ associated to $Z$. The resulting isomorphism class is unique by \autoref{cor:dividealcorresp}.
\ep
We can compute the first Stiefel--Whitney class of a normal-crossing log divisor.
\begin{cor} Let $(U,\sigma)$ be a normal-crossing log divisor with zero set $Z_\sigma = \cup_i Z_i$. Then $w_1(U) = \sum_i {\rm PD}_{\Z_2}[Z_i]$.
\end{cor}
\bp Decompose $(U,\sigma)$ into log divisors as $(U,\sigma) \cong (\otimes_i L_i, \otimes_i s_i)$, and use the fact that $w_1$ is additive with respect to tensor product, together with \autoref{prop:logdivstiefelwhitney}.
\ep
Finally, we mention that given two log divisors $(L,s) \to X$ and $(L',s') \to X'$, their direct sum is a normal-crossing log divisor on $X \times X'$. In fact, the same is true for direct sums of normal-crossing log divisors.
\subsection{Morse--Bott type divisors}
\label{sec:mbtypedivisor}
More generally, we can define any Morse--Bott-type divisor. For this we first define the notion of a nondegenerate critical submanifold of a section of a line bundle.
\begin{defn} Let $(U,\sigma)$ be a line bundle with a section. A compact connected submanifold $S \subseteq X$ is a \emph{nondegenerate critical submanifold} of $\sigma$ if $S \subseteq Z_\sigma$ and $\ker {\rm Hess}(\sigma) = T_p S$ for all $p \in S$.
\end{defn}
We say that $\sigma \in \Gamma(U)$ is a \emph{Morse--Bott section} if $Z_\sigma$ consists of nondegenerate critical submanifolds of $\sigma$. Morse--Bott sections $\sigma \in \Gamma(U)$ turn $(U,\sigma)$ into a divisor.
\begin{defn} A divisor $(U,\sigma)$ is of \emph{Morse--Bott type} if $\sigma$ is Morse--Bott.
\end{defn}
\begin{rem} If $U$ is trivial, note that a Morse--Bott section for $U$ is not the same thing as a Morse--Bott function. Namely, we only demand that the zero set consists of nondegenerate critical submanifolds, instead of demanding this for the set of all critical points. Indeed, the only condition on a divisor section is put on its zero set. In this sense, an alternate way of phrasing the above definition is saying that the section must specify the germ of a Morse--Bott function around its zero set.
\end{rem}
We see immediately that elliptic divisors are of Morse--Bott type. Note however that (normal-crossing) log divisors are not of Morse--Bott type. Any divisor $(U,\sigma)$ of Morse--Bott type will satisfy $w_1(U) = 0$, as its zero set consists of submanifolds of codimension at least two (c.f.\ \autoref{prop:ellipticdivstiefelwhitney}). It should be clear from the discussion above \autoref{prop:ellmorsebott} that its conclusion also holds for any divisor of Morse--Bott type.
\begin{prop} Let $(U,\sigma)$ be a divisor of Morse--Bott type and $\alpha$ a compatible trivialization of $U^*$. Then there exists a neighbourhood $U$ of the zero section $Z_\sigma \subset NZ_\sigma$ and an open embedding $\Phi\colon U \to X$ such that $\Phi|_{Z_\sigma} = {\rm id}_{Z_\sigma}$ and $(\Phi^*\alpha)(\Phi^* q) = Q_g$, where $g = \alpha(q) \in C^\infty(X)$.
\end{prop}
As a consequence of the above result, we see that the divisor ideal $I_\sigma$ will be locally given by $\langle Q_g \rangle$, with $Q_g$ a homogeneous function in normal bundle coordinates.						
\chapter{Lie algebroids}
\label{chap:liealgebroids}
\fancypagestyle{empty}{%
	\fancyhf{}%
	\renewcommand\headrulewidth{0pt}%
	\fancyhead[RO,LE]{\thepage}%
}
In this chapter we introduce and study Lie algebroids. These are vector bundles on a manifold equipped with a Lie bracket on their space of sections, and a compatible bundle map to the tangent bundle of that manifolds. Lie algebroids were first introduced by Pradines \cite{Pradines67} as the infinitesimal analogues of Lie groupoids. We have chosen not to discuss Lie groupoids in this thesis, as we focus on the use of Lie algebroids in desingularizing geometric structures. While Lie algebroids of various kinds are used in the literature, we focus on developing those aspects of Lie algebroid theory that are most relevant to the specific Lie algebroids we are interested in. These are mainly those whose so-called isomorphism locus is dense, in that they are almost-everywhere isomorphic to the tangent bundle of the underlying manifold. Equivalently, they are those whose anchor specifies a divisor. We will not be comprehensive regarding fundamentals, and instead refer the reader to e.g.\ the monograph \cite{Mackenzie05}, or \cite{Marle08}. Most of this chapter reviews the basics of Lie algebroid theory, except Sections \ref{sec:isomorphismloci}, \ref{sec:ladivtype}, \ref{sec:aliealgebroids}, \ref{sec:rescaling} and \ref{sec:residuemaps}, which contain (in part) new material.
\subsection*{Organization of the chapter}
This chapter is built up as follows. In Section \ref{sec:basicdefinitions} we discuss the basic definitions in Lie algebroid theory. In Section \ref{sec:lacohomology} we then note how a Lie algebroid $\mc{A}$ defines an accompanying Lie algebroid de Rham cohomology theory. Section \ref{sec:lamorphisms} discusses morphisms between Lie algebroids, their isomorphism loci, and Lie subalgebroids. Next in Section \ref{sec:aconnareps} we discuss Lie algebroid connections and Lie algebroid representations. Section \ref{sec:aliealgebroids} introduces the notion of an $\mc{A}$-Lie algebroid, which is a Lie algebroid whose anchor factors through that of a fixed Lie algebroid $\mc{A}$. Section \ref{sec:degenloci} discusses degeneraci loci and invariant submanifolds of Lie algebroids. In Section \ref{sec:rescaling} we discuss the process of rescaling, by which a Lie algebroid can be changed along a hypersurface. We discuss algebraic operations on Lie algebroids in Section \ref{sec:algoperations}, and finish this chapter with Section \ref{sec:residuemaps} how residue maps can be used to extract information about $\mc{A}$-forms along certain invariant submanifolds.
\section{Basic definitions}
\label{sec:basicdefinitions}
This section starts with the definition of a Lie algebroid, together with accompanying basic concepts. Before we define the notion of a Lie algebroid, we start by introducing a more fundamental concept, which is that of an anchored vector bundle.
\begin{defn} An \emph{anchored vector bundle} $(E,\rho_E)$ is a vector bundle $E \to X$ and a vector bundle map $\rho_E\colon E \to TX$ called the \emph{anchor} of $E$.
\end{defn}
Anchored vector bundles on their own are perhaps not so interesting. However, many more elaborate types of bundles, such as Lie algebroids, are in particular also anchored vector bundles.
\begin{defn}\label{defn:liealgebroid} A \emph{Lie algebroid} $(\mathcal{A},[\cdot,\cdot]_{\mc{A}},\rho_\mc{A})$ is an anchored vector bundle with a Lie bracket $[\cdot,\cdot]_{\mc{A}}$ on $\Gamma(\mc{A})$ satisfying the Leibniz rule $[v,f w]_{\mc{A}} = f [v,w]_{\mc{A}} + \mc{L}_{\rho_{\mc{A}}(v)} f \cdot w$ for all $v, w \in \Gamma(\mc{A})$, $f \in C^\infty(X)$.
\end{defn}
We will mention here just two examples of Lie algebroids. More examples that are directly relevant for this thesis will be given in Chapter \ref{chap:concreteliealgebroids}.
\begin{exa}\label{exa:tgntbundle} The tangent bundle $\mc{A} = TX$ is a Lie algebroid, with anchor $\rho_\mc{A} = {\rm id}_{TX}$ the identity, and using the standard Lie bracket $[\cdot,\cdot]$ on vector fields.
\end{exa}
\begin{exa} Let $\mf{g}$ be a finite-dimensional Lie algebra. Then $\mf{g}$ defines a Lie algebroid over a point, $X = \{{\rm pt}\}$, and the anchor map is trivial.
\end{exa}
The above two examples are fundamental and are in some sense extreme opposites of each other. In the first example the bundle is trivial (with identity anchor) but the manifold is not, while in the second these roles are reversed.
\begin{defn} Let $\mc{A} \to X$ be a Lie algebroid. An \emph{$\mc{A}$-orientation} on $X$ is an orientation for the bundle $\mc{A}$.
\end{defn}
Note that an $\mc{A}$-orientation exists if and only if $w_1(\mc{A}) = 0$. In this thesis we consider Lie algebroids as generalizations of the tangent bundle $TX$, chosen such that geometric constructions done using $\mc{A}$ are more suitable to the situation at hand. Typically, $\mc{A}$ is in some sense closely related to $TX$ (more precisely, the isomorphism locus of $\mc{A}$ is dense, see Section \ref{sec:isomorphismloci}). For this reason we will not emphasize what happens in the case of Lie algebras at various points.
\begin{defn} Let $\mc{A} \to X$ be a Lie algebroid. The \emph{isotropy} of $\mc{A}$ at $x \in X$ is the subspace $\ker \rho_{\mc{A},x} \subseteq \mc{A}_x$.
\end{defn}
Of course, the above definition makes sense for any anchored vector bundle. However, for Lie algebroids, the isotropy is in fact a Lie algebra: given $v, w \in \Gamma(\mc{A})$ such that $v_x, w_x \in \ker \rho_{\mc{A},x}$ for $x \in X$, we see that $[v,f w]_{\mc{A},x} = f(x) [v,w]_{\mc{A},x}$ for all $f \in C^\infty(X)$ by the Leibniz rule. Hence there exists a well-defined Lie bracket $[\cdot,\cdot]$ on $\ker \rho_{\mc{A},x}$, such that $[v,w]_{\mc{A},x} = [v_x, w_x]$, as any two sections $w, w' \in \Gamma(\mc{A})$ with $w_x = w'_x$ differ locally by a linear combination $\sum_i f_i w_i$ for $f_i \in C^\infty(X)$ and $w_i \in \Gamma(\mc{A})$, where each $f_i$ vanishes at $x$.
\begin{rem} In general, the isotropy Lie algebras $\ker \rho_{\mc{A},x}$ for $x \in X$ need not be of constant dimension. Hence, in general $\ker \rho_\mc{A}$ is not a vector subbundle of $\mc{A}$.
\end{rem} 
The compatibility condition between anchor and bracket has the following consequence (e.g.\ \cite[Proposition 3.1.2]{Marle08}), which used to be included in the axioms.
\begin{prop}\label{prop:anchorliealgmorph} Let $\mc{A} \to X$ be a Lie algebroid. Then $\rho_\mc{A}\colon (\Gamma(\mc{A}),[\cdot,\cdot]_\mc{A}) \to (\Gamma(TX), [\cdot,\cdot])$ is a Lie algebra homomorphism.
\end{prop}
\bp We must show that $\rho_\mc{A}([v,w]_\mc{A}) = [\rho_\mc{A}(v), \rho_\mc{A}(w)]$ for all $v,w \in \Gamma(\mc{A})$. Let $u \in \Gamma(\mc{A})$ be a third section and $f \in C^\infty(X)$. By the Jacobi identity of $[\cdot,\cdot]_\mc{A}$ we have
\be
	[[v,w]_\mc{A},f u]_\mc{A} + [[f u, v]_\mc{A},w]_\mc{A} + [[w, f u]_\mc{A}, v]_\mc{A} = 0.
\ee
Applying the Leibniz rule multiple times we see that
\begin{align*}
	[[v,w]_\mc{A},f u]_\mc{A} &= f [[v,w]_\mc{A},u]_\mc{A} + (\mc{L}_{\rho_\mc{A}([v,w]_\mc{A})} f) \cdot u,\\
	[[f u, v]_\mc{A},w]_\mc{A} &= f[w,[v,u]_\mc{A}]_\mc{A} + (\mc{L}_{\rho_\mc{A}(w)} f) \cdot [v,u]_\mc{A} + (\mc{L}_{\rho_\mc{A}(w)} \mc{L}_{\rho_\mc{A}(v)} f) \cdot u,\\
	[[w, f u]_\mc{A}, v]_\mc{A} &= - f [v,[w,u]_\mc{A}]_\mc{A} - (\mc{L}_{\rho_\mc{A}(v)} f) \cdot [w,u]_\mc{A} - (\mc{L}_{\rho_\mc{A}(w)} f) \cdot [v,u]_\mc{A}\\
	& \qquad - (\mc{L}_{\rho_\mc{A}(v)} \mc{L}_{\rho_\mc{A}(w)} f) \cdot u.
\end{align*}
By combining these four equations and using the Jacobi identity again we obtain
\be
	(\mc{L}_{\rho_\mc{A}([v,w]_\mc{A})} f) \cdot u + (\mc{L}_{\rho_\mc{A}(w)} \mc{L}_{\rho_\mc{A}(v)} f) \cdot u - (\mc{L}_{\rho_\mc{A}(v)} \mc{L}_{\rho_\mc{A}(w)} f) \cdot u = 0.
\ee
This can be rewritten to $(\mc{L}_{\rho_\mc{A}([v,w]_\mc{A}) - [\rho_\mc{A}(v),\rho_\mc{A}(w)]} f) \cdot u = 0$. As both $u$ and $f$ are arbitrary, we conclude that $\rho_\mc{A}([v,w]_\mc{A}) - [\rho_\mc{A}(v),\rho_\mc{A}(w)] = 0$ as desired.
\ep
We remark that the image of the anchor of any Lie algebroid specifies an involutive singular distribution in the sense of Stefan--Sussmann \cite{Sussmann73}. This follows for example from a splitting theorem for Lie algebroids (e.g.\ \cite{Fernandes02, BursztynLimaMeinrenken16}), which we will not discuss here. However, in Chapter \ref{chap:splittingtheorems} we will extend the work in \cite{BursztynLimaMeinrenken16} to obtain in particular a splitting theorem for $\mc{A}$-Lie algebroids (see Section \ref{sec:aliealgebroids}).
\begin{prop} Let $\mc{A} \to X$ be a Lie algebroid. Then $D = \rho_\mc{A}(\mc{A}) \subseteq TX$ is an involutive singular distribution.
\end{prop}
The leaves of the above distribution, i.e.\ the associated maximal immersed submanifolds $\mc{O}$ satisfying $T_x \mc{O} = \im \rho_{\mc{A},x}$ for all $x \in \mc{O}$, are called the \emph{orbits} of $\mc{A}$. Their tangent spaces are spanned locally by the image of the anchor of $\mc{A}$.
\subsection{Complex Lie algebroids}
There is an analogous definition of a complex Lie algebroid, using complex vector bundles. We will not make extensive use of complex Lie algebroids in this thesis, but nevertheless they will play a role in the study of \sgcs{}s in Section \ref{sec:scgs} (see also Section \ref{sec:complexlogtangent}). Given a manifold $X$, denote by $TX_\C = TX \otimes \C$ the complexification of its tangent bundle.
\begin{defn} A \emph{complex anchored vector bundle} $(E_\C, \rho_{E_\C})$ is a complex vector bundle $E_\C \to X$ and a bundle map $\rho_{E_\C}\colon E_\C \to TX_\C$ called the \emph{anchor} of $E_\C$.
\end{defn}
\begin{defn} A \emph{complex Lie algebroid} $(\mc{A}_\C, \rho_{\mc{A}_\C}, [\cdot,\cdot]_{\mc{A}_\C})$ is a complex anchored vector bundle $(\mc{A}_\C,\rho_{\mc{A}_\C})$ with a Lie bracket $[\cdot,\cdot]_{\mc{A}_\C}$ on $\Gamma(\mc{A}_\C)$ satisfying the Leibniz rule $[v,f w]_{\mc{A}_\C} = f [v,w]_{\mc{A}_\C} + \mc{L}_{\rho_{\mc{A}_\C}(v)} f \cdot w$ for all $v, w \in \Gamma(\mc{A}_\C)$, $f \in C^\infty(X;\C)$.
\end{defn}
Essentially all notions we define for Lie algebroids naturally carry over to the complex setting. In particular, a complex Lie algebroid $\mc{A}_\C$ also has an associated Lie algebroid cohomology $H^\bullet(\mc{A}_\C)$ which reduces to complex de Rham cohomology $H^\bullet_{\rm dR}(X;\C)$ if $\mc{A}_\C = TX_\C$ (see the next section). Given a Lie algebroid $\mc{A}$, its complexification $\mc{A} \otimes \C$ (note that we do \emph{not} denote this by $\mc{A}_\C$) is naturally a complex Lie algebroid. There is another operation we can perform, namely that of complex conjugation. This results in another complex Lie algebroid.
\begin{defn}\label{defn:complexconjgla} Let $\mc{A}_\C \to X$ be a complex Lie algebroid. The \emph{complex conjugate} of $\mc{A}_\C$ is the complex Lie algebroid $\overline{\mc{A}}_\C \to X$ given by the complex conjugate complex vector bundle, equipped with the complex conjugate anchor and bracket.
\end{defn}
\section{Lie algebroid cohomology}
\label{sec:lacohomology}
In this section we discuss the differential $\mc{A}$-de Rham complex naturally associated to a Lie algebroid $\mc{A} \to X$. We refer to \cite[Chapter 7]{Mackenzie05} and \cite{Marle08} for more discussion. In analogy with the case when $\mc{A} = TX$, denote by $\Omega^\bullet(\mc{A}) = \Gamma(\wedge^\bullet \mc{A}^*)$ the graded algebra of differential $\mc{A}$-forms. This algebra comes equipped with a differential $d_\mc{A}$, constructed using the bracket $[\cdot,\cdot]_{\mc{A}}$ by means of the usual Koszul formula.
\begin{defn} Let $\mc{A} \to X$ be a Lie algebroid. The \emph{$\mc{A}$-differential} $d_\mc{A}\colon \Omega^\bullet(\mc{A}) \to \Omega^{\bullet+1}(\mc{A})$ is defined by the following expression, given $\eta \in \Omega^k(\mc{A})$ and $v_i \in \Gamma(\mc{A})$,
\begin{align*}
	d_\mc{A} \eta(v_0,\dots,v_k) &= \sum_{i = 0}^k (-1)^{i} \mc{L}_{\rho_\mc{A}(v_i)}(\eta(v_0,\dots,\hat{v}_i,\dots,v_{k}))\\
		&+ \sum_{0 \leq i < j \leq k} (-1)^{i + j} \eta([v_i,v_j]_\mc{A}, v_0, \dots, \hat{v}_i, \dots, \hat{v}_j, \dots, v_{k}),
\end{align*}
where the hat denotes omission of the term.
\end{defn}
\begin{rem} For notational convenience we denote by $\Omega^k_\mc{A}(U)$, for an open subset $U \subset X$, the set of $\mc{A}$-$k$-forms defined on $U$, so that $\Omega^k_\mc{A}(X) = \Omega^k(\mc{A})$.
\end{rem}
The $\mc{A}$-differential satisfies the usual relation with respect to the wedge product, namely that for all $\eta, \xi \in \Omega^\bullet(\mc{A})$ one has
\be
	d_\mc{A}(\eta \wedge \xi) = d_\mc{A} \eta \wedge \xi + (-1)^{|\eta|} \eta \wedge d_\mc{A} \xi,
\ee
where $|\eta|$ is the degree of $\eta$. In fact, the entire Lie algebroid structure on $\mc{A}$ can be recovered from the operator $d_\mc{A}$. Similarly to the usual de Rham differential, the $\mc{A}$-differential squares to zero. The proof is identical to the usual case and is omitted.
\begin{lem} Let $\mc{A} \to X$ be a Lie algebroid. Then $d_\mc{A} \circ d_\mc{A} = 0$.
\end{lem}
The above lemma implies that the pair $(\Omega^\bullet(\mc{A}), d_\mc{A})$ forms a differential complex. Consequently, each Lie algebroid has an associated cohomology theory.
\begin{defn} Let $\mc{A} \to X$ be a Lie algebroid. The \emph{Lie algebroid cohomology} of $\mc{A}$ is given by $H^k_\mc{A}(X) = H^k(\Omega^\bullet(\mc{A}), d_{\mc{A}})$ for $k \in \N \cup \{0\}$.
\end{defn}
We return briefly to our two examples.
\begin{exa} When $\mc{A} = TX$, the $\mc{A}$-differential is the usual de Rham differential $d$, and the Lie algebroid de Rham complex of $TX$ is just its usual de Rham complex.
\end{exa}
\begin{exa}\label{exa:liealgcoh} When $\mc{A} = \mf{g}$, a Lie algebra, the associated Lie algebroid cohomology of $\mf{g}$ is equal to its Chevalley--Eilenberg cohomology.
\end{exa}
Lie algebroid cohomology is generally very hard to compute and need not be finite-dimensional, even on compact manifolds. Nevertheless, the cohomology of most of the explicit Lie algebroids we consider can be computed (see Section \ref{sec:laexamples}).
\begin{defn} Let $v \in \Gamma(\mc{A})$. Then the \emph{$\mc{A}$-Lie derivative} by $v$ is defined as the operator $\mc{L}_v = [\iota_v, d_\mc{A}]$ on $\Omega^\bullet(\mc{A})$, where $[\cdot,\cdot]$ denotes the graded commutator.
\end{defn}
The $\mc{A}$-Lie derivative and $\mc{A}$-differential satisfy the usual Cartan relations as for when $\mc{A} = TX$, and we have chosen not to list them. They can be found in e.g.\ \cite{Marle08}. There is a natural extension of the $\mc{A}$-Lie derivative to $\Gamma(\wedge^\bullet \mc{A})$ as when $\mc{A} = TX$.
\section{Lie algebroid morphisms}
\label{sec:lamorphisms}
In this section we introduce morphisms between Lie algebroids. As Lie algebroids carry a bracket structure on their space of sections, we will see some care must be taken to describe morphisms over different base manifolds. We start with the notion of morphism for anchored vector bundles.
\begin{defn} A \emph{morphism of anchored vector bundles} between $(X,E,\rho_E)$ and $(Y,F,\rho_F)$ is a vector bundle morphism $(\varphi,f)\colon E \to F$ such that $\rho_F \circ \varphi = Tf \circ \rho_E$.
\end{defn}
\begin{exa}\label{exa:anchoredmorphism} Let $E \to X$ be anchored vector bundle. Then its anchor $\rho_E\colon E \to TX$ is a morphism of anchored vector bundles, where $TX$ is given the trivial anchored vector bundle structure.
\end{exa}
Morphisms of anchored vector bundles give rise to maps between the isotropies.
\begin{lem}\label{lem:morphavisotropy} Let $(\varphi,f)\colon (E,X) \to (F,Y)$ be a morphism of anchored vector bundles. Then for each $x \in X$ there is an induced map $\varphi_x\colon \ker \rho_{E,x} \to \ker \rho_{F,f(x)}$.
\end{lem}
\bp Let $v_x \in \ker \rho_{E,x}$ be given. Then as $\varphi$ is a morphism of anchored vector bundles, we have $\rho_{F}(\varphi(v_x)) = T f(\rho_E(v_x)) = 0$, so that $\varphi(v_x) \in \ker \rho_{F,f(x)}$.
\ep
We now turn to Lie algebroid morphisms, where care must be taken in describing bracket compatibility. First we consider Lie algebroid morphisms over the same base.
\begin{defn}\label{def:lamorph} Let $\mc{A}, \mc{A}' \to X$ be Lie algebroids over the same base manifold. A \emph{Lie algebroid morphism} from $\mc{A}$ to $\mc{A}'$ is a morphism of anchored vector bundles $(\varphi, f)\colon \mc{A} \to \mc{A}'$ such that $\varphi[v,w]_{\mc{A}} = [\varphi(v),\varphi(w)]_{\mc{A}'}$ for all $v,w \in \Gamma(\mc{A})$.
\end{defn}
In other words, the extra condition is that $\varphi$ induces a Lie algebra homomorphism on sections. Consequently, we have the following commutative diagrams.
\begin{center}
	\begin{tikzpicture}
	\begin{scope}[xshift=-120pt]
	\matrix (m) [matrix of math nodes, row sep=2.5em, column sep=2.5em,text height=1.5ex, text depth=0.25ex]
	{	\mc{A} & \mc{A}' \\ X & X \\};
	\path[-stealth]
	(m-1-1) edge node [above] {$\varphi$} (m-1-2)
	(m-1-1) edge (m-2-1)
	(m-2-1) edge node [above] {$f$} (m-2-2)
	(m-1-2) edge (m-2-2);
	\end{scope}
	\begin{scope}	
	\matrix (m) [matrix of math nodes, row sep=2.5em, column sep=2.5em,text height=1.5ex, text depth=0.25ex]
	{	\mc{A} & \mc{A}' \\ TX & TX \\};
	\path[-stealth]
	(m-1-1) edge node [above] {$\varphi$} (m-1-2)
	(m-1-1) edge node [left] {$\rho_{\mc{A}}$} (m-2-1)
	(m-2-1) edge node [above] {$Tf$} (m-2-2)
	(m-1-2) edge node [right] {$\rho_{\mc{A}'}$} (m-2-2);
	\end{scope}
	\end{tikzpicture}
\end{center}
As an initial example, anchors of Lie algebroids are Lie algebroid morphisms.
\begin{prop}\label{prop:anchorlamorph} Let $\mc{A} \to X$ be a Lie algebroid. Then $\rho_\mc{A}\colon \mc{A} \to TX$ is a Lie algebroid morphism.
\end{prop}
\bp We saw that $\rho_\mc{A}$ is a morphism of anchored vector bundles (\autoref{exa:anchoredmorphism}). That $\rho_\mc{A}$ is a Lie algebroid morphism then follows from \autoref{prop:anchorliealgmorph}.
\ep
We introduce several adjectives to describe the behavior of the anchor.
\begin{defn} Let $\mc{A} \to X$ be a Lie algebroid. Then $\mc{A}$ is \emph{regular} if $\rho_\mc{A}$ has constant rank, \emph{transitive} if $\rho_\mc{A}$ is surjective, and \emph{totally intransitive} if $\rho_\mc{A} \equiv 0$.
\end{defn}
\begin{defn} Let $\mc{A} \to X$ be a Lie algebroid. Then $\mc{A}$ is \emph{almost-injective} \cite{Debord01} if $\rho_\mc{A}\colon \Gamma(\mc{A}) \to \Gamma(TX)$ is an embedding of sheaves, i.e.\ is injective.
\end{defn}
For regular Lie algebroids (hence for transitive ones), the kernel of $\rho_\mc{A}$ is of constant dimension, so that it is a subbundle of $\mc{A}$. If $\mc{A}$ is transitive, the isotropy Lie algebras together exhibit $\mc{A}$ as an abelian extension of $TX$; there is a short exact sequence as follows
\be
	0 \to \ker \rho_\mc{A} \to \mc{A} \to TX \to 0.
\ee
Returning to our discussion on Lie algebroid morphisms, note that \autoref{def:lamorph} does not immediately extend to varying base, because in general a bundle map does not induce a map on sections. Let $E \to X$ and $F \to X'$ be vector bundles. Recall that a vector bundle morphism $(\varphi, f)\colon E \to F$ factors through the pullback bundle $f^* F$. In other words, there is a unique vector bundle morphism $\varphi^!\colon E \to f^* F$ over $X$ such that $\varphi = {\rm id}_F \circ \varphi^!$, making the following diagram commutative.
\begin{center}
	\begin{tikzpicture}
	\matrix (m) [matrix of math nodes, row sep=2.5em, column sep=2.5em,text height=1.5ex, text depth=0.25ex]
	{	E & f^* F & F \\ X & X & Y \\};
	\path[-stealth]
	(m-1-1) edge node [above] {$\varphi^!$} (m-1-2)
	(m-1-1) edge (m-2-1)
	(m-2-1) edge node [above] {${\rm id}_X$} (m-2-2)
	(m-1-2) edge (m-2-2)
	(m-1-2) edge node [above] {${\rm id}_{F}$} (m-1-3)
	(m-1-3) edge (m-2-3)
	(m-2-2) edge node [above] {$f$} (m-2-3);
	\end{tikzpicture}
\end{center}
As $\varphi^!$ is a vector bundle morphism over the same base, it does induce a morphism on sections $\varphi^!\colon \Gamma(E) \to \Gamma(f^* F)$. Next, note that the $C^\infty(X)$-module $\Gamma(f^* F)$ can be described as the tensor product of modules
\be
	\Gamma(f^* F) = C^\infty(X) \otimes_{C^\infty(Y)} \Gamma(F),
\ee
where $C^\infty(X)$ is a module over $C^\infty(Y)$ using the pullback $f^*$, and the isomorphism is given by the map $g \otimes v \mapsto g f^* v$ for $g \in C^\infty(X)$ and $v \in \Gamma(F)$. In other words, sections of the pullback bundle $f^* F$ can be described nonuniquely as finite $C^\infty(X)$-linear combinations of pullbacks of sections of $F$.
\begin{defn} Given $v \in \Gamma(f^* F)$, we call an expression $v = \sum_i f_i \otimes v_i$ for $f_i \in C^\infty(X)$ and $v_i \in \Gamma(F)$ a \emph{decomposition} of $v$.
\end{defn}
Using such decompositions, we can phrase the bracket compatibility for general Lie algebroid morphisms.
\begin{defn}\label{def:lamorphgen} Let $\mc{A} \to X$ and $\mc{A}' \to X'$ be Lie algebroids. A \emph{Lie algebroid morphism} from $\mc{A}$ to $\mc{A}'$ is a morphism of anchored vector bundles $(\varphi,f)\colon (\mc{A},X) \to (\mc{A}',X')$ such that for all $v, w \in \Gamma(\mc{A})$ and decompositions $\varphi^!(v) = \sum_i f_i \otimes v_i$ and $\varphi^!(w) = \sum_j g_j \otimes w_j$, we have
\be
	\varphi^!([v,w]_\mc{A}) = \sum_{i,j} f_i g_j \otimes [v_i,w_j]_{\mc{A}'} + \sum_j (\rho_{\mc{A}}(v) g_j) \otimes w_j - \sum_i (\rho_{\mc{A}}(w) f_i) \otimes v_i.
\ee
\end{defn}
For this definition to make sense, one must check that the expression on the right hand side is independent of the chosen decompositions. This is indeed so, as can be found in \cite{Mackenzie05} or \cite[Lemma 1.4]{HigginsMackenzie90}.
\begin{prop} Let $(\varphi,f)\colon (\mc{A},X) \to (\mc{A}',X')$ be a morphism of anchored vector bundles between Lie algebroids, and let $v,w \in \Gamma(\mc{A})$. Then the expression for $\varphi^![v,w]_\mc{A}$ in \autoref{def:lamorphgen} is independent of decompositions for $\varphi^!(v)$ and $\varphi^!(w)$.
\end{prop}
\begin{rem} Let $(\varphi,f)\colon (\mc{A},X) \to (\mc{A}',X')$ be a morphism of anchored vector bundles between Lie algebroids. Then the condition that $\varphi$ is a Lie algebroid morphism for sections $v, w \in \Gamma(\mc{A})$ which can be pushed forward to pullback sections, i.e.\ $\varphi^!(v) = f^* v'$ and $\varphi^!(w) = f^* w'$ for $v',w' \in \Gamma(\mc{A}')$, reads $\varphi^!([v,w]_\mc{A}) = f^*([v',w']_{\mc{A}'})$, c.f.\ \autoref{def:lamorph}.
\end{rem}
The reader is invited to verify that the composition of two Lie algebroid morphisms is again a Lie algebroid morphism. This results in a category $\mc{LA}$ of Lie algebroids over smooth manifolds.
\begin{exa}\label{exa:tgntlamorph} Recall from \autoref{exa:tgntbundle} that the tangent bundle is a Lie algebroid. Given a smooth map $f\colon X \to X'$, its differential $T f\colon TX \to TX'$ is a Lie algebroid morphism. See \cite[Proposition 4.3.3]{Mackenzie05}, or use \autoref{prop:duallamorphism} below.
\end{exa}
\begin{rem}
	Another way to phrase the conditions of a Lie algebroid morphism is the following. A bundle map $(\varphi,f)\colon (\mc{A},X) \to (\mc{A}',X')$ is a Lie algebroid morphism if and only if its graph ${\rm Gr}(\varphi) \subseteq \mc{A}' \times \mc{A}$ is a Lie subalgebroid along the graph ${\rm Gr}(f) \subseteq X' \times X$. This allows us to define Lie algebroid morphisms in terms of Lie subalgebroids (to be studied in Section \ref{sec:liesubalgds}). However, to do this we must first specify the Lie algebroid structure on the direct product Lie algebroid $\mc{A}' \times \mc{A} \to X' \times X$. See Section \ref{sec:algoperations} and references therein, and \cite{Meinrenken17}.
\end{rem}
\subsection{Isomorphism loci}
\label{sec:isomorphismloci}
In this section we introduce the isomorphism locus of a Lie algebroid, or more generally, of a Lie algebroid morphism. Namely, for Lie algebroid morphisms covering the identity on a given manifold $X$, we keep track of where it is an isomorphism.
\begin{defn} Let $(\varphi,{\rm id}_X)\colon \mc{A} \to \mc{A}'$ be a Lie algebroid morphism. The \emph{isomorphism locus} of $\varphi$ is the subset $X_\varphi \subseteq X$ where $\varphi$ is an isomorphism.
\end{defn}
This definition similarly works for morphisms of anchored vector bundles covering the identity. Given a Lie algebroid $\mc{A} \to X$, the isomorphism locus of $\rho_\mc{A}$ will be denoted by $X_\mc{A} := X_{\rho_\mc{A}}$. Note that if the isomorphism locus of $\varphi$ is nonempty, we have that ${\rm rank}(\mc{A}) = {\rm rank}(\mc{A}')$. In particular, if $X_\mc{A}$ is nonempty, then ${\rm rank}(\mc{A}) = \dim X$. We denote the complement by $Z_\mc{A} = X \setminus X_\mc{A}$. The following is immediate.
\begin{prop}\label{prop:compisoloci} Let $(\varphi,{\rm id}_X)\colon \mc{A} \to \mc{A}'$ and $(\varphi',{\rm id}_X)\colon \mc{A}' \to \mc{A}''$ be Lie algebroid morphisms. Then $X_{\varphi' \circ \varphi} \supseteq X_\varphi \cap X_{\varphi'}$.
\end{prop}
It can happen that the isomorphism locus of a composition is strictly larger than the intersection of isomorphism loci. This stems from the fact that if a composition of maps $f \circ g$ is bijective, it can only be concluded that $f$ is surjective and $g$ is injective.
\begin{rem}\label{rem:denseisoalmostinj} Let $\mc{A} \to X$ be an almost-injective Lie algebroid. Then $\mc{A}$ has dense isomorphism locus if and only if ${\rm rank}(\mc{A}) = \dim(X)$. Moreover, if $\mc{A}$ has dense isomorphism locus, then $\mc{A}$ is almost-injective.
\end{rem}
Let $\mc{A} \to X$ be a Lie algebroid with dense isomorphism locus. Then the inclusion $i\colon X_{\mc{A}} \hookrightarrow X$ of the isomorphism locus gives a bijection $\rho_{\mc{A}}^*\colon \Omega^k(X_{\mc{A}}) \to \Omega^k_{\mc{A}}(X_{\mc{A}})$ for all $k$. Consequently, we can view $\mc{A}$-forms as smooth forms with certain ``singular'' behavior at $X \setminus X_\mc{A}$. The isomorphism given by $\rho_\mc{A}$ also implies the following, c.f.\ \autoref{rem:denseisoalmostinj}.
\begin{prop}\label{prop:denseisolocusstrongmap} Let $(\varphi,f)\colon \mc{A} \to \mc{A}'$ be a Lie algebroid morphism between Lie algebroids with dense isomorphism loci for which $f\colon (X,Z_\mc{A}) \to (X',Z_{\mc{A}'})$ is a strong map of pairs. Then $\varphi = Tf$ on sections.
\end{prop}
The following is immediate from the previous.
\begin{prop} Let $(\varphi,f)\colon \mc{A} \to \mc{A}'$ be a Lie algebroid morphism between Lie algebroids with dense isomorphism loci for which $f\colon (X,Z_\mc{A}) \to (X',Z_{\mc{A}'})$ is a strong map of pairs. Then $\rho_\mc{A}|_{X_\mc{A}}\colon \ker \varphi \to \ker Tf$ is an isomorphism in $X_\mc{A}$.
\end{prop}
We can sometimes relate the first Stiefel--Whitney class of two Lie algebroids.
\begin{prop}\label{prop:codimtwosw} Let $(\varphi,{\rm id}_X)\colon \mc{A} \to \mc{A}'$ be a Lie algebroid morphism with dense isomorphism locus, such that $X \setminus X_{\varphi}$ is a union of smooth submanifolds of codimension at least two. Then $w_1(\mc{A}) = w_1(\mc{A}')$.
\end{prop}
\bp Note that $w_1(\mc{A}) = w_1(\det(\mc{A}))$ and similarly for $\mc{A}'$. We can consider the determinant $\det \varphi\colon \det(\mc{A}) \to \det(\mc{A}')$ between these line bundles. As $X \setminus X_\varphi$ has codimension at least two, $\det \varphi$ is an isomorphism except on at most a codimension two subset. But then $\det(\mc{A})$ and $\det(\mc{A}')$ are isomorphic.
\ep
\begin{cor}\label{cor:codimtwodegensw} Let $\mc{A} \to X$ be a Lie algebroid such that $Z_\mc{A}$ is a union of smooth submanifolds of codimension at least two. Then $w_1(\mc{A}) = w_1(TX)$.
\end{cor}
\subsection{Lie algebroids of divisor-type}
\label{sec:ladivtype}
We can further relate Lie algebroids with dense isomorphism loci to the divisors of Chapter \ref{chap:divisors} as follows. Given a Lie algebroid morphism $(\varphi,{\rm id}_X)\colon \mc{A} \to \mc{A}'$ with dense isomorphism locus, we can consider its determinant $\det(\varphi)\colon \det(\mc{A}) \to \det(\mc{A}')$. We can view this as specifying a section $\det(\varphi) \in \Gamma(\det(\mc{A}^*) \otimes \det(\mc{A}'))$. The condition that $\varphi$ has dense isomorphism locus is now equivalent to demanding that $(\det(\mc{A}^*) \otimes \det(\mc{A}'), \det(\varphi))$ is a divisor. We call this the \emph{divisor associated to $\varphi$} and denote this by ${\rm div}(\varphi)$. The associated divisor ideal is denoted by $I_\varphi \subseteq C^\infty(X)$.
\begin{defn} Let $\mc{A} \to X$ be a Lie algebroid with dense isomorphism locus. Then ${\rm div}(\rho_\mc{A}) = (\det(\mc{A}^*) \otimes \det(TX), \det(\rho_\mc{A}))$ is the \emph{divisor associated to $\mc{A}$}, and will also be denoted by ${\rm div}(\mc{A})$. The associated divisor ideal is denoted by $I_\mc{A}$.
\end{defn}
We can now instead speak of \emph{Lie algebroids of divisor-type}, as we will do for Poisson and Dirac structures in Sections \ref{sec:poissondivisors} and \ref{sec:diracdivisors} (and for generalized complex structures in Section \ref{sec:agcs}). We have chosen not to constantly use this terminology over explicitly stating the isomorphism locus is dense. Similarly, we can speak of \emph{Lie algebroid morphisms of divisor-type} (only for base-preserving morphisms).
\begin{rem} In light of the definition of the canonical $\mc{A}$-module $Q_\mc{A} = \det(\mc{A}) \otimes \det(T^*X)$ (see \autoref{defn:canamodule}), it is interesting to explore the consequences of the observation that $Q_\mc{A}^* \cong {\rm div}(\mc{A})$ as line bundles, for Lie algebroids of divisor-type.
\end{rem}
It is immediate from \autoref{prop:compisoloci} that divisors of Lie algebroid morphisms behave as follows, after noting that $\det(\varphi' \circ \varphi) = \det(\varphi') \circ \det(\varphi)$.
\begin{prop} Let $(\varphi,{\rm id}_X)\colon \mc{A} \to \mc{A}'$ and $(\varphi',{\rm id}_X)\colon \mc{A}' \to \mc{A}''$ be Lie algebroid morphisms of divisor-type. Then $\varphi' \circ \varphi$ is of divisor-type, and ${\rm div}(\varphi' \circ \varphi) \cong {\rm div}(\varphi') \otimes {\rm div}(\varphi)$, the product of divisors.
\end{prop}
In other words, composition of maps corresponds to taking products of divisors, hence to products of divisor ideals by \autoref{prop:proddivisor}.
\begin{rem}\label{rem:mdivtype} In light of \autoref{defn:adivisortype} and Section \ref{sec:riggedalgebroids}, it makes sense to define a Lie algebroid $\mc{A} \to X$ to be of \emph{$m$-divisor-type} for some $m \geq 0$ if there exists a regular Lie subalgebroid $D \subseteq TX$ of rank $m$ (see Section \ref{sec:liesubalgds}, in particular \autoref{exa:regfoliation}) such that $\rho_\mc{A}$ factors through the inclusion $\rho_D$, i.e.\ $\rho_\mc{A} = \rho_D \circ \varphi_{\mc{A}}$ for some Lie algebroid morphism $\varphi_\mc{A}\colon \mc{A} \to D$, and $\varphi_{\mc{A}}$ is of divisor-type, with divisor given by ${\rm div}(\varphi_\mc{A}) = (\det(\mc{A}^*) \otimes \det(D), \det(\varphi_\mc{A}))$. Similarly for the concept of an $\mc{A}$-Lie algebroid as will be introduced in Section \ref{sec:aliealgebroids}.
\end{rem}
\subsection{Dual description}
There is an alternative description of Lie algebroid morphisms which allows us to circumvent the fact that there is no induced map on sections.

Returning to general vector bundles $E \to X$ and $F \to X'$, recall that there is a one-to-one correspondence between vector bundle morphisms $(\varphi, f)\colon E \to F$, and dually algebra morphisms $\varphi^*\colon \Omega^\bullet(F) \to \Omega^\bullet(E)$, given by $(\varphi^* \xi)(v) = \xi(\varphi(v))$ for $\xi \in \Omega^\bullet(F)$ and $v \in \mf{X}^\bullet(E)$. Recall further that a Lie algebroid $\mc{A}$ equips the space $\Omega^\bullet(\mc{A})$ with a differential $d_\mc{A}$ which encodes the Lie algebroid structure. Using this we can phrase the conditions that $\varphi$ preserves anchors and brackets in terms of $\varphi^*$.
\begin{prop}\label{prop:duallamorphism} Let $(\varphi,f)\colon \mc{A} \to \mc{A}'$ be a bundle morphism between Lie algebroids, with dual $\varphi^*\colon (\mc{A}')^* \to \mc{A}^*$. Then $(\varphi,f)$ is a Lie algebroid morphism if and only if $\varphi^*\colon \Omega^\bullet(\mc{A}') \to \Omega^\bullet(\mc{A})$ is a cochain map, i.e.\ $d_\mc{A} \circ \varphi^* = \varphi^* \circ d_{\mc{A}'}$.
\end{prop}
As is well-known, a chain map between differential complexes descends to a map on cohomology. Consequently, by \autoref{prop:duallamorphism} we obtain the following.
\begin{prop}\label{prop:lalgdmorphcoh} Let $(\varphi,f)\colon \mc{A} \to \mc{A}'$ be a Lie algebroid morphism. Then there is an induced map in Lie algebroid cohomology $\varphi^*\colon H^\bullet(X';\mc{A}') \to H^\bullet(X;\mc{A})$.
\end{prop}
In particular, as the anchor $\rho_\mc{A}\colon \mc{A} \to TX$ is a Lie algebroid morphism by \autoref{prop:anchorlamorph}, we see that there is always a map $\rho_\mc{A}^*\colon H^\bullet_{\rm dR}(X) \to H^\bullet(X;\mc{A})$. In general this map is neither injective nor surjective.

For most of the Lie algebroids we will consider, smooth maps of the underlying manifolds give rise to Lie algebroid morphisms, as long as they intertwine the anchor maps. This is true in general for anchored vector bundle morphisms between Lie algebroids with dense isomorphism loci.
\begin{prop}\label{prop:denselam} Let $\mc{A}_X \to X$ and $\mc{A}_Y \to Y$ be Lie algebroids such that $X_{\mc{A}_X}$ is dense. Suppose that $(\varphi,f)\colon \mc{A}_X \to \mc{A}_Y$ is an anchored bundle morphism and $f^{-1}(Y_{\mc{A}_Y}) = X_{\mc{A}_X}$. Then $(\varphi, f)\colon \mc{A}_X \to \mc{A}_Y$ is a Lie algebroid morphism.
\end{prop}
\bp As $\varphi$ is a vector bundle morphism, $\varphi^*$ is an algebra morphism. In the isomorphism loci, $\varphi$ must equal $Tf$ by \autoref{prop:denseisolocusstrongmap}, and $(Tf,f)$ is a Lie algebroid morphism between $TX$ and $TY$ by \autoref{exa:tgntlamorph}, i.e.\ $f^*$ is a chain map. By density of $X_{\mc{A}_X}$, the map $\varphi^*$ is a chain map everywhere, so that $\varphi$ is a Lie algebroid morphism by \autoref{prop:duallamorphism}.
\ep
Consequently, for such Lie algebroids one can determine whether there is a Lie algebroid morphism $(\varphi, f)$ by checking that $f^*$ extends to a map $\varphi^*$ on forms. This in turn will follow by the universal property of the exterior algebra if it holds on generators, so it suffices to check that $f^*$ extends to $\varphi^*\colon \Omega^1(\mc{A}_Y) \to \Omega^1(\mc{A}_X)$.
\begin{rem}\label{rem:lalgdpoisson} Recall that for a finite dimensional vector space $V$, Lie algebra structures on $V$ are in one-to-one correspondence with linear Poisson structures on $V^*$. This correspondence extends to Lie algebroids \cite{Courant90}, and allows us to in particular define Lie algebroid morphisms in Poisson terms. See e.g.\ \cite{HigginsMackenzie93, Meinrenken17}.
\end{rem}
\begin{rem} There is yet another description of a Lie algebroid morphisms that uses an auxilliary connection. For more information, see one of \cite{CrainicFernandes11,HigginsMackenzie90, Mackenzie05}.
\end{rem}
\subsection{Lie subalgebroids}
\label{sec:liesubalgds}
We next introduce the appropriate subobject in the category of Lie algebroids.
\begin{defn} Let $\mc{A} \to M$ be a Lie algebroid and $N \subseteq X$ a closed embedded submanifold. A \emph{Lie subalgebroid} $\mc{B} \to N$ supported on $N$ is a vector subbundle $\mc{B} \subseteq \mc{A}|_N$ for which the bundle inclusion $(\varphi,f)$ is a Lie algebroid morphism.
\end{defn}
Hence, $(\varphi,f)\colon (\mc{B},N) \to (\mc{A},X)$ is a fiberwise injective Lie algebroid morphism covering an injective immersion.
\begin{rem}\label{rem:liesubalgd} Alternatively but equivalently, one can proceed as follows. The anchor $\rho_\mc{B}$ must be the restriction of $\rho_\mc{A}$ to $\mc{B}$, so that $\rho_\mc{A}(\mc{B}) \subseteq TN$. Given sections $v, w \in \Gamma(\mc{A})$ such that $v|_N, w|_N \in \Gamma(\mc{B})$, we must have $[v,w]_\mc{A}|_N \in \Gamma(\mc{B})$, i.e.\ the space of sections that restrict to $\mc{B}$ is closed under the Lie bracket.
\end{rem}
Sometimes one includes the condition that, given $v,w \in \Gamma(\mc{A})$ with $v|_N = 0$ and $w|_N \in \Gamma(\mc{B})$, we should have $[v,w]_\mc{A}|_N = 0$. However, this condition is automatic, as we now show. This is noted in \cite[Remark 2.1.21]{Li13} and \cite[Proposition 2.15]{Meinrenken17}.
\begin{prop} Let $\mc{A} \to X$ be a Lie algebroid and $Z \subseteq X$ a submanifold. A subbundle $\mc{B}$ of $\mc{A}|_Z$ is a Lie subalgebroid of $\mc{A}$ supported on $Z$ if and only if $\mc{B} \to Z$ is a Lie algebroid.
\end{prop}
\bp Let $v \in \Gamma(\mc{A})$ such that $v|_Z = 0$. Then $v = \sum_i f_i v_i$ for some $v_i \in \Gamma(\mc{A})$ and $f_i \in C^\infty(X)$ which vanish on $Z$. Given $w \in \Gamma(\mc{A})$ such that $w|_Z \in \Gamma(\mc{B})$, we have $\rho_\mc{A}(w)|_Z \in \Gamma(TZ)$. Consequently $\mc{L}_{\rho_\mc{A}(w)} f_i|_Z = 0$ for all $i$, and the Leibniz rule gives $[v,w]_\mc{A}|_Z = 0$, as $[v, w]_\mc{A} = [f v', w]_\mc{A} = \sum_i (f_i [v_i,w]_{\mc{A}} - \mc{L}_{\rho_\mc{A}(w)} f_i \cdot v_i)$.
\ep
\begin{exa}\label{exa:regfoliation} Let $\mc{F}$ be a regular foliation on $X$ and $D = T \mc{F} \subseteq TX$ the associated involutive distribution. Then $\mc{A}_\mc{F} = D$ inherits a Lie algebroid structure from $TX$, and is a Lie subalgebroid. Moreover, if $\mc{A} \to X$ is a regular Lie algebroid with injective anchor, then $\mc{A} = \mc{A}_\mc{F}$ for $\mc{F}$ associated to the distribution $D = \rho_\mc{A}(\mc{A})$.
\end{exa}
\section{Lie algebroid connections and representations}
\label{sec:aconnareps}
In this section we discuss the notion of a Lie algebroid connection on a vector bundle, and moreover discuss Lie algebroid representations, which are vector bundles equipped with a flat Lie algebroid connection. For more information, see e.g.\ \cite{EvensLuWeinstein99, KosmannSchwarzbachLaurentGengouxWeinstein08}.

Let $\mc{A} \to X$ be a Lie algebroid and $E \to X$ a vector bundle.
\begin{defn} An \emph{$\mc{A}$-connection} on $E$ is a bilinear map $\nabla\colon \Gamma(\mc{A}) \times \Gamma(E) \to \Gamma(E)$, $(v,\sigma) \mapsto \nabla_v \sigma$, such that $\nabla_{f v} \sigma = f \nabla_v \sigma$ and $\nabla_v(f\sigma) = f \nabla_v \sigma + (\rho_\mc{A}(v) f) \cdot \sigma$ for all $v \in \Gamma(\mc{A})$, $\sigma \in \Gamma(E)$ and $f \in C^\infty(X)$.
\end{defn}
We start by establishing that $\mc{A}$-connections always exist.
\begin{lem} Let $\nabla'\colon \Gamma(TX) \times \Gamma(E) \to \Gamma(E)$ be a connection on $E$. Then $\nabla_v := \nabla'_{\rho_\mc{A}(v)}$ for $v \in \Gamma(\mc{A})$ defines an $\mc{A}$-connection on $E$.
\end{lem}
\bp It is immediate that $\nabla$ is bilinear and satisfies $\nabla_{f v} \sigma = f \nabla_v \sigma$ because $\rho_\mc{A}$ is a bundle map, which hence is $C^\infty(X)$-linear. For the final property, merely note that $\nabla_v(f \sigma) = \nabla'_{\rho_\mc{A}(v)}(f \sigma) = f \nabla'_{\rho_\mc{A}(v)}(\sigma) + (\rho_\mc{A}(v) f) \cdot \sigma = f \nabla_v \sigma + (\rho_\mc{A}(v) f) \cdot \sigma$, using that $\nabla'$ is a connection.
\ep
Any $\mc{A}$-connection $\nabla$ on $E$ has a \emph{curvature tensor} $F_\nabla \in \Gamma(\wedge^2 \mc{A}^* \otimes {\rm End}(E))$, given by $F_\nabla(v,w) := \nabla_v \circ \nabla_w - \nabla_w \circ \nabla_v - \nabla_{[v,w]_\mc{A}}$ for $v, w \in \Gamma(\mc{A})$. We say $\nabla$ is \emph{flat} if $F_\nabla \equiv 0$. With this we can define the notion of a Lie algebroid representation.
\begin{defn} A pair $(E,\nabla)$ consisting of a vector bundle and an $\mc{A}$-connection $\nabla$ on $E$ is an \emph{$\mc{A}$-representation} or \emph{$\mc{A}$-module} if $\nabla$ is flat, i.e.\ $F_\nabla \equiv 0$.
\end{defn}
Given an $\mc{A}$-module $(E,\nabla)$, we can define a differential $d_{\mc{A},\nabla}$ on $\Omega^\bullet(\mc{A};E)$, the space of $\mc{A}$-forms with values in $E$. Recalling that $\Omega^\bullet(\mc{A};E) = \Gamma(\wedge^\bullet \mc{A}^*) \otimes \Gamma(E)$, we set $d_{\mc{A},\nabla}(\eta \otimes s) := d_{\mc{A}} \eta \otimes s + (-1)^{|\eta|} \eta \otimes \nabla s$, with $|\eta|$ the degree of $\eta$. This differential satisfies $d_{\mc{A},\nabla}(\eta \wedge \xi \otimes s) = d_{\mc{A}} \eta \wedge \xi \otimes s + (-1)^{|\eta|} \eta \wedge d_{\mc{A},\nabla} (\xi \otimes s)$, and $d_{\mc{A},\nabla}$ squares to zero if and only if $\nabla$ is flat.
\begin{defn} Let $\mc{A} \to X$ be a Lie algebroid and $(E,\nabla)$ an $\mc{A}$-module. The \emph{$\mc{A}$-cohomology with values in $E$} is given by $H^k(\mc{A}; E) = H^k(\Omega^\bullet(\mc{A};E), d_{\mc{A},\nabla})$ for $k \in \N \cup \{0\}$.
\end{defn}
\begin{exa} For a Lie algebroid over a point we have $\mc{A} = \mf{g}$, a Lie algebra, and $\mc{A}$-representations are exactly the representations of this Lie algebra.
\end{exa}
\begin{exa} Let $\underline{\R} \to X$ be the trivial line bundle. This carries a trivial $\mc{A}$-representation structure for any Lie algebroid $\mc{A} \to X$, given by $\nabla_v f = \mc{L}_{\rho_{\mc{A}}(v)} f$ for $v \in \Gamma(\mc{A})$ and $f \in \Gamma(\underline{\R}) = C^\infty(X)$. We then have that $H^\bullet(\mc{A};\underline{\R}) = H^\bullet(\mc{A})$.
\end{exa}
Given a Lie algebroid $\mc{A} \to X$, there is always a canonical $\mc{A}$-module (see \cite{EvensLuWeinstein99}).
\begin{defn}\label{defn:canamodule} Let $\mc{A} \to X$ be a Lie algebroid and set $Q_\mc{A} := \det(\mc{A}) \otimes \det(T^*X)$. Then $Q_\mc{A}$ is  the \emph{canonical $\mc{A}$-module}, using the $\mc{A}$-connection on $Q_\mc{A}$ that is defined for $v \in \Gamma(\mc{A})$, $V \in \Gamma(\det(\mc{A}))$ and $\mu \in \Gamma(\det(T^*X))$ by
\be
	\nabla_v(V \otimes \mu) = \mc{L}_{v} V \otimes \mu + V \otimes \mc{L}_{\rho_\mc{A}(v)} \mu.
\ee
\end{defn}
One readily checks this formula indeed defines a flat $\mc{A}$-connection on $Q_\mc{A}$.
\section{Degeneracy loci}
\label{sec:degenloci}
In this section we discuss when one can restrict a Lie algebroid $\mc{A}$ to a submanifold, as a Lie algebroid. Such submanifolds will be called $\mc{A}$-invariant. Moreover, we introduce degeneracy loci of Lie algebroids and show that they are $\mc{A}$-invariant subsets. We follow in part the exposition in \cite{Pym13}.
\begin{defn} Let $\mc{A} \to X$ be a Lie algebroid. Then a submanifold $N \subseteq X$ is \emph{$\mc{A}$-invariant} if $\rho_\mc{A}$ is tangent to $N$, i.e.\ $\rho_\mc{A}|_N$ has image inside $TN$.
\end{defn}
For closed submanifolds $N \subseteq X$, a vector field $V \in \Gamma(TX)$ preserves $I_N$ if and only if it is tangent to $N$. Hence if the submanifold $N$ is $\mc{A}$-invariant, its vanishing ideal $I_N = \{f \in C^\infty(X) \, | \, f|_N = 0\}$ is preserved by the action of $\mc{A}$. In other words, we have $\mc{L}_{\rho_\mc{A}(v)} I_N \subseteq I_N$ for all $v \in \Gamma(\mc{A})$. If $N$ is $\mc{A}$-invariant, the restriction of the anchor $\rho_\mc{A}|_N$ has image inside $TN$, instead of merely $TX|_N$.

We see that $\mc{A}$-invariant submanifolds consist of orbits of $\mc{A}$, i.e.\ are unions of leaves of the singular foliation induced by $\rho_\mc{A}$. The reason for introducing the notion of $\mc{A}$-invariance is the following consequence of the definition.
\begin{prop}\label{prop:ainvariantrestr} Let $\mc{A}$ be a Lie algebroid and $N$ an $\mc{A}$-invariant closed submanifold of $X$. Then $(\mc{A}|_N, \rho_\mc{A}|_N)$ with restriction of $[\cdot,\cdot]_\mc{A}$, is a Lie algebroid over $N$.
\end{prop}
\bp We check that the bracket is well-defined. This follows because $\rho_\mc{A}\colon \Gamma(\mc{A}) \to \Gamma(TX)$ is a Lie algebra homomorphism (\autoref{prop:anchorliealgmorph}), combined with the fact that the Lie bracket of vector fields tangent to $N$ is again tangent to $N$. Moreover, we have that $[v,w]_\mc{A}|_N = 0$ for $v, w \in \Gamma(\mc{A})$ if $w|_N = 0$. Namely, we can write $v = \sum_i f_i w_i$ for $f_i \in I_N$, and compute that $[v,w]_\mc{A}|_N = \sum_i f_i|_N [v,w_i]_\mc{A}|_N + (\mc{L}_{\rho_\mc{A}(v)} f_i)|_N \cdot w_i|_N = 0$, using that $\mc{L}_{\rho_\mc{A}(v)} f_i = 0$ because $\rho_\mc{A}(v)$ is tangent to $N$. We conclude that $\mc{A}|_N$ is a Lie subalgebroid of $\mc{A}$ along $N$.
\ep
\begin{rem}\label{rem:inverseimagesmooth} The proof of \autoref{prop:ainvariantrestr} shows that more generally, given a Lie algebroid $\mc{A} \to X$ and a submanifold $\iota\colon N \hookrightarrow X$ such that $\iota^! \mc{A} := \rho_\mc{A}^{-1}(TN)$ is a smooth subbundle, it is a Lie subalgebroid of $\mc{A}$ along $N$ (\cite[Proposition 2.17]{Meinrenken17}).
\end{rem}
The set of $\mc{A}$-invariant closed submanifolds of $X$ admits unions and intersections. First, we say a subset of $X$ is \emph{$\mc{A}$-invariant} if it is a union of $\mc{A}$-invariant submanifolds.
\begin{lem} Let $N$ and $N'$ be $\mc{A}$-invariant closed submanifolds of $X$. Then their union $N \cup N'$ and intersection $N \cap N'$ are $\mc{A}$-invariant closed subsets.
\end{lem}
\bp The sets $N \cup N'$ and $N \cap N'$ have vanishing ideals given by $I_{N \cup N'} = I_N \cap I_{N'} = I_N \cdot I_{N'}$ and $I_{N \cap N'} = I_N + I_{N'}$. Using the product rule $\mc{L}_{\rho_\mc{A}(v)}(f g) = (\mc{L}_{\rho_\mc{A}(v)}f) \cdot g + f \cdot (\mc{L}_{\rho_{\mc{A}}(v)}g)$, we conclude that $I_{N \cup N'}$ is still $\mc{A}$-invariant. The $\mc{A}$-invariance of $I_{N \cap N'}$ is immediate. Alternatively, this follows from the realization that $\mc{A}$-invariant submanifolds are unions of orbits of $\mc{A}$.
\ep
We next define the degeneracy loci of $\mc{A}$, which will be our main examples of $\mc{A}$-invariant submanifolds. Note that the anchor of $\mc{A}$ is a map $\rho_\mc{A}\colon \mc{A} \to TX$, so that for each integer $k \geq 0$ we obtain a map $\wedge^k \rho_\mc{A}\colon \wedge^k \mc{A} \to \wedge^k TX$. Dually, we can view this as a map $\wedge^k \rho_\mc{A}\colon \Gamma(\wedge^k \mc{A} \otimes \wedge^k T^*X) \to C^\infty(X)$. Its image is an ideal $I_{\mc{A},k}$.
\begin{defn} Let $\mc{A} \to X$ be a Lie algebroid. The $k$th \emph{degeneracy locus} of $\mc{A}$ is the subspace $X_{\mc{A},k} \subseteq X$ defined by the ideal $I_{\mc{A},k+1} \subseteq C^\infty(X)$.
\end{defn}
Hence, $X_{\mc{A},k}$ is the subspace of $X$ where the rank of $\rho_\mc{A}$ is $k$ or less.
\begin{prop}[{\cite[Proposition 2.2.9]{Pym13}}]\label{prop:degenlocus} Let $\mc{A} \to X$ be a Lie algebroid and $k \geq 0$. Then $X_{\mc{A},k}$ is $\mc{A}$-invariant whenever it is smooth.
\end{prop}
\bp Assume that $X_{\mc{A},k}$ is a smooth submanifold and denote the map $\wedge^{k+1} \rho_\mc{A}$ by $\varphi_{k}\colon \wedge^{k+1} \mc{A} \to \wedge^{k+1} TX$. Let $v \in \Gamma(\mc{A})$ and $f \in I_{\mc{A},k}$. Then by definition there exists $\xi \in \Gamma(\wedge^{k+1} \mc{A})$ and $\omega \in \Omega^{k+1}(X)$ such that $f = \langle \varphi_k(\xi), \omega \rangle$. We compute using the compatibility of the bracket and anchor that $\mc{L}_{\rho_\mc{A}(v)} f = \mc{L}_{\rho_\mc{A}(v)} \langle \varphi_k(\xi), \omega \rangle = \langle \mc{L}_{\rho_\mc{A}(v)} \varphi_k(\xi), \omega \rangle + \langle \varphi_k(\xi), \mc{L}_{\rho_\mc{A}(v)} \omega \rangle = \langle \varphi_k(\mc{L}_{v} \xi), \omega\rangle + \langle \varphi_k(\xi), \mc{L}_{\rho_\mc{A}(v)} \omega \rangle \in I_{\mc{A},k}$. As $I_{\mc{A},k}$ is preserved under the action of $\mc{A}$, we conclude that $X_{\mc{A},k}$ is $\mc{A}$-invariant.
\ep
\begin{exa}\label{exa:denseisolocus} Let $\mc{A} \to X$ be a Lie algebroid with dense isomorphism locus, and let $n = \dim X = {\rm rank}(\mc{A})$. Then the complement of the isomorphism locus $X_\mc{A}$ is equal to the $(n-1)$th degeneracy locus of $\mc{A}$, i.e.\ $Z_\mc{A} = X_{\mc{A},n-1}$. In fact, the divisor ideal $I_\mc{A}$ associated to ${\rm div}(\mc{A})$ is equal to $I_{\mc{A},n}$. Note moreover that both $X_\mc{A}$ and $Z_\mc{A}$ are $\mc{A}$-invariant subsets, where $X_\mc{A}$ is open.
\end{exa}
\section{\texorpdfstring{$\mc{A}$}{A}-Lie algebroids}
\label{sec:aliealgebroids}
In this section we introduce the notion of an $\mc{A}$-Lie algebroid, which is a Lie algebroid whose anchor factors through the anchor $\rho_\mc{A}$ of a given Lie algebroid $\mc{A}$. This notion is useful when studying the process of lifting geometric structures, as we will pursue in later chapters. We start by considering the anchored vector bundle version.
\begin{defn} Let $\mc{A} \to X$ be a Lie algebroid. An \emph{$\mc{A}$-anchored vector bundle} is an anchored vector bundle $E_\mc{A} \to X$ with a morphism of anchored vector bundles $\varphi_{E_\mc{A}}\colon E_{\mc{A}} \to \mc{A}$ called the \emph{$\mc{A}$-anchor} satisfying $\rho_{E_\mc{A}} = \rho_\mc{A} \circ \varphi_{E_\mc{A}}$.
\end{defn}
We see that the definition requires that the anchor $\rho_{E_\mc{A}}$ factors through $\rho_\mc{A}$.
\begin{rem} To not have to speak of `Lie algebroid-anchored vector bundles', we will use the name $\mc{A}$-anchored vector bundles even when there are different Lie algebroids $\mc{A}$ to consider, as happens for example in the upcoming \autoref{defn:morphaav}.
\end{rem}
With this notion comes the associated notion of morphism.
\begin{defn}\label{defn:morphaav} A \emph{morphism of $\mc{A}$-anchored vector bundles} between $(X,E_\mc{A},\rho_{E_\mc{A}})$ and $(X',F_{\mc{B}},\rho_{F_\mc{B}})$ is a morphism of anchored vector bundles $(\varphi,f)\colon E_{\mc{A}} \to F_{\mc{B}}$ covering a Lie algebroid morphism $(\varphi',f)\colon \mc{A} \to \mc{B}$, i.e.\ $\varphi' \circ \varphi_{E_{\mc{A}}} = \varphi_{F_{\mc{B}}} \circ \varphi$.
\end{defn}
In other words, we have the following diagram with commutative squares.
\begin{center}
	\begin{tikzpicture}
	\matrix (m) [matrix of math nodes, row sep=2.5em, column sep=2.5em,text height=1.5ex, text depth=0.25ex]
	{	E_\mc{A} & \mc{A} & TX \\ F_\mc{B} & \mc{B} & TX' \\};
	\path[-stealth]
	(m-1-1) edge node [above] {$\varphi_{E_\mc{A}}$} (m-1-2)
	(m-1-1) edge node [left] {$\varphi$} (m-2-1)
	(m-2-1) edge node [above] {$\varphi_{F_\mc{B}}$} (m-2-2)
	(m-1-2) edge node [left] {$\varphi'$} (m-2-2)
	(m-1-2) edge node [above] {$\rho_\mc{A}$} (m-1-3)
	(m-2-2) edge node [above] {$\rho_\mc{B}$} (m-2-3)
	(m-1-3) edge node [left] {$T f$} (m-2-3);
	\draw [->] (m-1-1) [bend left=45] edge node [above] {$\rho_{E_\mc{A}}$} (m-1-3)
	(m-2-1) [bend right=45] edge node [below] {$\rho_{F_\mc{B}}$} (m-2-3)
	;
	\end{tikzpicture}
\end{center}
Alternatively, we can demand for $(\varphi,f)$ to be a vector bundle morphism, as the existence of $\varphi'$ implies that $\varphi$ is also a morphism of anchored vector bundles: we have $\rho_{F_\mc{B}} \circ \varphi = \rho_{\mc{B}} \circ \varphi_{F_\mc{B}} \circ \varphi = \rho_{\mc{B}} \circ \varphi' \circ \varphi_{E_\mc{A}} = Tf \circ \rho_\mc{A} \circ \varphi_{E_\mc{A}} = Tf \circ \rho_{E_\mc{A}}$.
We now define the $\mc{A}$-analogue of a Lie algebroid.
\begin{defn}\label{defn:aliealgebroid} An \emph{$\mc{A}$-Lie algebroid} is an $\mc{A}$-anchored vector bundle $E_\mc{A}$ equipped with a compatible Lie algebroid structure for which $\varphi_{E_\mc{A}}$ is a Lie algebroid morphism.
\end{defn}
Note that $\mc{A}$-Lie algebroids are $\mc{A}$-anchored vector bundles, and are also Lie algebroids. An \emph{$\mc{A}$-Lie algebroid morphism} is a morphism of $\mc{A}$-anchored vector bundles that is also a Lie algebroid morphism. Any Lie algebroid $\mc{A} \to X$ is automatically both an $\mc{A}$- and a $TX$-Lie algebroid. Of course, given an $\mc{A}$-Lie algebroid $\mc{A}'$, any $\mc{A}'$-Lie algebroid is also an $\mc{A}$-Lie algebroid.
\begin{exa}\label{exa:adistribution} Let $\mc{A} \to M$ be a Lie algebroid and consider the notion of an \emph{$\mc{A}$-distribution}, which is a smooth subbundle $D_\mc{A} \subseteq \mc{A}$. Its image under the anchor, $D := \rho_\mc{A}(D_\mc{A}) \subseteq TM$, need not be a distribution. If $D_\mc{A} \subseteq \mc{A}$ is $\mc{A}$-involutive, then it is an $\mc{A}$-Lie algebroid with $\mc{A}$-anchor the natural inclusion.
\end{exa}
We next describe a way to compute Lie algebroid cohomology of $\mc{A}$-Lie algebroids, outlining the general procedure of what is done in \cite{Lanius16,Lanius16two}. Let $\mc{A}'$ be an $\mc{A}$-Lie algebroid for which the $\mc{A}$-anchor $\varphi_{\mc{A}'}$ has dense isomorphism locus, as happens for example when $\mc{A}'$ is a rescaling of $\mc{A}$ (see Section \ref{sec:rescaling}). Considering the dual map $\varphi^*_{\mc{A}'}\colon \mc{A}'^* \to \mc{A}$, we obtain a short exact sequence of complexes
\be
0 \to \Omega^\bullet(\mc{A}) \stackrel{\varphi^*_{\mc{A}'}}{\to} \Omega^\bullet(\mc{A}') \stackrel{p}{\to} \mc{C}^\bullet(\mc{A},\mc{A}') \to 0,
\ee
with $\mc{C}^\bullet(\mc{A},\mc{A}') = \Omega^\bullet(\mc{A}') / \varphi^*_{\mc{A}'} (\Omega^\bullet(\mc{A}))$ and $p$ the projection onto this quotient.
\begin{rem} In our setting, both $\varphi_{\mc{A}'}\colon \Gamma(\mc{A}') \to \Gamma(\mc{A})$ \emph{and} $\varphi^*_{\mc{A}'}\colon \Gamma(\mc{A}^*) \to \Gamma(\mc{A}'^*)$ are injective, because both maps on vector bundle level are almost-everywhere isomorphisms. We remark further that $\Gamma(\mc{A}^*)$ is not the vector space dual of $\Gamma(\mc{A})$.
\end{rem}
There is a differential $d_{\mc{C}(\mc{A},\mc{A}')}$ on $\mc{C}^\bullet(\mc{A},\mc{A}')$ induced by the differential $d_{\mc{A}'}$, given for $\eta \in \mc{C}^\bullet(\mc{A},\mc{A}')$ by $d_{\mc{C}(\mc{A},\mc{A}')}\eta := p(d_{\mc{A}'} \eta')$, for any $\eta' \in \Omega^\bullet(\mc{A}')$ such that $p(\eta') = \eta$. This is well-defined because $\varphi^*_{\mc{A}'}$ is a cochain map due to \autoref{prop:duallamorphism}: given $\eta'$ with $p(\eta'') = 0$, we have by exactness that $\eta' = \varphi_{\mc{A}'}^* \xi$ for some $\xi \in \Gamma(\mc{A}')$. We then compute that $p(d_{\mc{A}'} \eta') = p(d_{\mc{A}'} \varphi_{\mc{A}'}^* \xi)$, but $p \circ d_{\mc{A}'} \circ \varphi_{\mc{A}'}^* = (p \circ \varphi_{\mc{A}'}^*)\circ d_{\mc{A}} = 0$ by exactness. It also follows that $d_{\mc{C}(\mc{A},\mc{A}')}^2 = 0$. Indeed, given $\eta \in \mc{C}^\bullet(\mc{A},\mc{A}')$ and $\eta' \in \Omega^\bullet(\mc{A}')$ with $p(\eta') = \eta$, we compute that $d^2_{\mc{C}(\mc{A},\mc{A}')} \eta = d_{\mc{C}(\mc{A},\mc{A}')} (p(d_{\mc{A}'} \eta'))$. Tautologically $\eta'' := d_{\mc{A}'} \eta'$ is such that $p(\eta'') = p(d_{\mc{A}'} \eta')$, so that by definition of $d_{\mc{C}(\mc{A},\mc{A}')}$ we conclude that $d^2_{\mc{C}(\mc{A},\mc{A}')} \eta = p (d_{\mc{A}'} \eta'') =  p(d_{\mc{A}'}(d_{\mc{A}'} \eta')) = 0$. From this we obtain the following, under the further assumption that the above short exact sequence of cochain complexes is split.
\begin{prop}\label{prop:alacoh} Let $\mc{A}'$ be an $\mc{A}$-Lie algebroid for which $\varphi_{\mc{A}'}$ has dense isomorphism locus, and assume that the projection $p$ admits a splitting. Then $H^\bullet(\mc{A}') \cong H^\bullet(\mc{A}) \oplus H^\bullet(\mc{C}(\mc{A},\mc{A}'))$.
\end{prop}
The above proposition combined with determining $H^\bullet(\mc{C}(\mc{A},\mc{A}'))$ underlies many of the computations of Lie algebroid cohomology presented in Section \ref{sec:laexamples}.
\section{Rescaling}
\label{sec:rescaling}
In this section we recall the process of rescaling \cite{Melrose93,Lanius16}, also referred to as lower elementary modification in \cite{GualtieriLi14}.

We can use Lie subalgebroids supported on hypersurfaces to coherently prescribe extra conditions for sections of a given Lie algebroid to satisfy, such that the sheaf of all such sections is again the sheaf of sections of a Lie algebroid of the same rank. To do this, we first recall the following classical result.
\begin{thm}[Serre--Swan \cite{Serre55, Swan62}]\label{thm:serreswan} Let $X$ be a manifold. There is an equivalence of categories between the category of vector bundles on $X$, and of finitely generated projective $C^\infty(X)$-modules, given by the map $(E \to X) \mapsto \Gamma(E)$.
\end{thm}
The fibers $E_x$ for $x \in X$ can be recovered from $\Gamma(E)$ via the natural isomorphism $E_x \cong \Gamma(E) / I_x \cdot \Gamma(E)$, where $I_x \subseteq C^\infty(X)$ is the ideal of functions vanishing at $x$.
Due to this result, any locally free (or finitely generated projective) $\mc{A}$-involutive $C^\infty(X)$-submodule of the sheaf of sections $\Gamma(\mc{A})$ of a given Lie algebroid, is itself the sheaf of sections of an $\mc{A}$-Lie algebroid. With this in mind, we have the following.
\begin{prop}\label{prop:rescaling} Let $\mc{A} \to X$ be a Lie algebroid, $Z \subseteq X$ a hypersurface and $\mc{B}$ a Lie subalgebroid of $\mc{A}$ supported on $Z$. Then there exists a Lie algebroid ${}^{\mc{B}} \mc{A} \to X$ for which
	\begin{equation*}
		\Gamma({}^{\mc{B}} \mc{A}) = \{v \in \Gamma(\mc{A}) \, | \, v|_Z \in \Gamma(\mc{B})\}.
	\end{equation*}
	The inclusion on sections induces a Lie algebroid morphism $(\varphi,{\rm id}_X)\colon {}^{\mc{B}} \mc{A} \to \mc{A}$.
\end{prop}
\bp Let $\mc{F} = \{v \in \Gamma(\mc{A}) \, | \, v|_Z \in \Gamma(\mc{B})\}$ be the sheaf under consideration. As $\mc{B}$ is a Lie subalgebroid, $\mc{F}$ is closed under the bracket $[\cdot,\cdot]_\mc{A}$ (see \autoref{rem:liesubalgd}). Consequently, by \autoref{thm:serreswan} we need only check that $\mc{F}$ is locally free. This readily follows from $Z$ being of codimension one (see \cite[Proposition 8.1]{Melrose93} or \cite[Theorem 2.2]{Lanius16}): let $n = {\rm rank}(\mc{A})$ and $m = {\rm rank}(\mc{B})$ and trivialize $\mc{B}$ around a point in $Z$ with basis of sections $(w_1, \dots w_m)$. Let $z \in I_Z$ be a local defining function for $Z$. Extend the previous basis to a local basis $(w_1,\dots,w_m,v_{m+1},\dots,v_n)$ of $\Gamma(\mc{A})$. Then $(w_1,\dots,w_m, z v_{m+1}, \dots, z v_n)$ is a local basis of $\Gamma({}^\mc{B} \mc{A})$, so $\mc{F}$ is locally free.
\ep
The above procedure is called \emph{rescaling}; we call ${}^{\mc{B}} \mc{A}$ the \emph{$(\mc{B},Z)$-rescaling} of $\mc{A}$. This is covered in \cite{Melrose93, Lanius16}. Examples of rescaling will be provided in Chapter \ref{chap:concreteliealgebroids}. Note that the induced morphism $\varphi\colon {}^\mc{B} \mc{A} \to \mc{A}$ is an isomorphism outside of $Z$. In \cite{GualtieriLi14}, this operation is called \emph{lower elementary modification} and ${}^\mc{B} \mc{A}$ is denoted by $[\mc{A} : \mc{B}]$. Inspecting the proof of \autoref{prop:rescaling}, we see that ${}^{\mc{B}} \mc{A}$ as a bundle exists because $\mc{B}$ is a vector subbundle of $\mc{A}$ supported on $Z$, i.e.\ because there is an injective vector bundle morphism from $(\mc{B},Z)$ to $(\mc{A},X)$. When this map is further a Lie algebroid morphism, the bundle ${}^\mc{B} \mc{A}$ carries a natural induced Lie algebroid structure.
\begin{rem}\label{rem:larescaling} Let $\mc{A} \to X$ be a Lie algebroid and ${}^\mc{B} \mc{A}$ a $(\mc{B},Z)$-rescaling of $\mc{A}$. Then ${}^\mc{B} \mc{A}$ is an $\mc{A}$-Lie algebroid, as $\Gamma({}^\mc{B} \mc{A})$ is naturally a subsheaf of $\Gamma(\mc{A})$. Its $\mc{A}$-anchor has dense isomorphism locus equal to $X \setminus Z$. Moreover, the divisor ideal of this $\mc{A}$-anchor (c.f.\ \autoref{rem:mdivtype}) is equal to $I_Z^k$, where $k$ is the corank of $\mc{B}$ inside $\mc{A}|_Z$.
\end{rem}
In contrast with the previous remark, note that ${}^\mc{B} \mc{A}$ is \emph{not} a Lie subalgebroid of $\mc{A}$ (unless it makes sense to consider the trivial case $\mc{B} = \mc{A}|_Z$, in which case ${}^\mc{B} \mc{A} = \mc{A}$).
\begin{rem}\label{rem:bbprimerescaling} Let $\mc{A} \to X$ be a Lie algebroid and $\mc{B}, \mc{B}' \subseteq \mc{A}|_Z$ two Lie subalgebroids of $\mc{A}$ supported on $Z$ such that $\mc{B} \subseteq \mc{B}'$. Then ${}^\mc{B} \mc{A}$ is a ${}^{\mc{B}'} \mc{A}$-Lie algebroid.
\end{rem}
The process of rescaling is commutative when rescaling at disjoint hypersurfaces, using the fact that the anchor is an isomorphism outside of these hypersurfaces.
\begin{prop} Let $\mc{A} \to X$ be a Lie algebroid and $(Z,\mc{B})$, $(Z',\mc{C})$ Lie subalgebroids of $\mc{A}$ such that $Z$ and $Z'$ are disjoint hypersurfaces. Then $\mc{B}$ and $\mc{C}$ are naturally Lie subalgebroids of $\, {}^\mc{C} \mc{A}$ and ${}^\mc{B} \mc{A}$ respectively. Further, ${}^\mc{C} ({}^\mc{B} \mc{A}) \cong {}^\mc{B} ({}^\mc{C} \mc{A})$.
\end{prop}
In other words, $(\mc{B},Z)$- and $(\mc{C},Z')$-rescaling commute when $Z \cap Z' = \emptyset$. It is always possible to rescale using $0 \subseteq \mc{A}|_Z$, the trivial Lie subalgebroid. In light of \autoref{rem:bbprimerescaling}, one could say that $0$-rescaling is the harshest, while rescaling using $\mc{B} = \rho_\mc{A}^{-1}(TZ)$ (if it exists, c.f.\ \autoref{rem:inverseimagesmooth}) is the mildest. Note that $0$-rescaling is not idempotent. See \cite{KlaasseLanius17two} for more information.
In fact, the rescaling procedure at disjoint hypersurfaces can be described in terms of the fiber product of Lie algebroids, which will be discussed in Section \ref{sec:algoperations}.
\begin{prop} Let $\mc{A} \to X$ be a Lie algebroid and $(Z, \mc{B})$, $(Z',\mc{C})$ Lie subalgebroids of $\mc{A}$ such that $Z$ and $Z'$ are disjoint hypersurfaces. Then rescaling at $Z$ and $Z'$ is the same as taking fiber product. In other words, ${}^\mc{C} ({}^\mc{B} \mc{A}) \cong {}^\mc{B} \mc{A} \oplus_{\mc{A}} {}^\mc{C} \mc{A} \cong {}^\mc{B} ({}^\mc{C} \mc{A})$.
\end{prop}
In \cite{KlaasseLanius17} we will discuss the relation between the Lie algebroid cohomologies of Lie algebroids ${}^\mc{B} \mc{A}$ and $\mc{A}$, where ${}^\mc{B} \mc{A}$ is a $(\mc{B},Z)$-rescaling of $\mc{A}$.
\section{Algebraic operations}
\label{sec:algoperations}
In this section we briefly discuss two algebraic operations one can perform on Lie algebroids. The first of these is the pullback operation, allowing one (in favorable cases) to transport a Lie algebroid along a smooth map. The second and closely related operation is that of the fiber product, where two Lie algebroids are combined. More information can be found in \cite{HigginsMackenzie90, ChenLiu07, Mackenzie05, KosmannSchwarzbachLaurentGengouxWeinstein08}
\begin{defn} Let $\mc{A} \to X$ be a Lie algebroid. A map $f\colon N \to X$ is \emph{transverse} to $\mc{A}$ if it is transverse to $\rho_\mc{A}$, i.e.\ $Tf(TN) + \rho_\mc{A}(\mc{A}) = TX$ at each point in $f(N)$.
\end{defn}
\begin{defn} Given a map $f\colon N \to X$ transverse to $\mc{A} \to X$, the \emph{pullback} of $\mc{A}$ along $f$ is the Lie algebroid $f^! \mc{A} \to N$ defined as the fiber product
	\be
		f^! \mc{A} = f^* \mc{A} \times_{f^*(TX)} TN = \{(v,X) \in f^* \mc{A} \times TN \, | \, \rho_\mc{A}(v) = Tf(X)\}.
	\ee
	The anchor is the natural projection onto $TN$ and its bracket is generated, for given $g, h \in C^\infty(N)$, $v,w \in \Gamma(\mc{A})$ and $X, Y \in \Gamma(TN)$, by the following expression:
	\be
		[(g v, X), (h w, Y)]_{f^! \mc{A}} = (g h [v,w]_\mc{A} + \mc{L}_X(h) \cdot w - \mc{L}_Y(g) \cdot v, [X,Y]).
	\ee
\end{defn}
\begin{rem}	It is possible for the pullback Lie algebroid to exist along maps $f$ not transverse to $\mc{A}$. All that is required is that $f$ is \emph{admissible}, i.e.\ the fiber product $f^! \mc{A}$ forms a smooth vector bundle, or that $f^* \mc{A} \times_{f^*(TX)} TN$ has constant rank. See \cite{KosmannSchwarzbachLaurentGengouxWeinstein08}.
\end{rem}
The pullback operation constitutes a pullback in the category $\mc{LA}$. First, note that given a bundle morphism $(\varphi,f)\colon (\mc{A},X) \to (\mc{A}',X')$ such that $f$ is admissible, we can define a map $(\varphi^{!!}, {\rm id}_X)\colon \mc{A} \to f^! \mc{A}'$ given by $v \mapsto \varphi(v) \oplus \rho_\mc{A}(v)$ for $v \in \Gamma(\mc{A})$.
\begin{prop}[{\cite[Theorem 4.3.6]{Mackenzie05}}] Let $(\varphi,f)\colon (\mc{A},X) \to (\mc{A}',X')$ be an anchored vector bundle morphism between Lie algebroids such that $f\colon X \to X'$ is admissible, i.e.\ the pullback Lie algebroid $f^! \mc{A}'$ exists. Then $(\varphi,f)$ is a Lie algebroid morphism if and only if $(\varphi^{!!},{\rm id}_X)\colon \mc{A} \to f^! \mc{A}'$ is a Lie algebroid morphism.
\end{prop}
Consequently, if $f$ is admissible we can factor any Lie algebroid morphism as for vector bundles. In other words, there is then the following commutative diagram.
\begin{center}
	\begin{tikzpicture}
	\matrix (m) [matrix of math nodes, row sep=2.5em, column sep=2.5em,text height=1.5ex, text depth=0.25ex]
	{	A & f^! \mc{A}' & \mc{A}' \\ X & X & X' \\};
	\path[-stealth]
	(m-1-1) edge node [above] {$\varphi^{!!}$} (m-1-2)
	(m-1-1) edge (m-2-1)
	(m-2-1) edge node [above] {${\rm id}_X$} (m-2-2)
	(m-1-2) edge (m-2-2)
	(m-1-2) edge node [above] {} (m-1-3)
	(m-1-3) edge (m-2-3)
	(m-2-2) edge node [above] {$f$} (m-2-3);
	\draw[-stealth, bend left=45] (m-1-1) edge node [above] {$\varphi$} (m-1-3);
	\end{tikzpicture}
\end{center}
We next turn to the fiber product of two Lie algebroids over the same base.
\begin{defn}\label{def:fiberproduct} Let $\mc{A}, \mc{A}' \to X$ be Lie algebroids such that $\rho_\mc{A}$ and $\rho_{\mc{A}'}$ are transverse. The \emph{fiber product} of $\mc{A}$ and $\mc{A}'$ is the Lie algebroid defined as the fiber product
\be
	\mc{A} \times_{TX} \mc{A}' = \{(v,w) \in \mc{A} \times \mc{A}' \, | \rho_\mc{A}(v) = \rho_{\mc{A}'}(w) \}.
\ee
The anchor is inherited from either $\mc{A}$ or $\mc{A}'$ and its bracket is given, for $v,w \in \Gamma(\mc{A})$ and $v',w' \in \Gamma(\mc{A})$, by the expression
\be
	[v \oplus v', w \oplus w']_{\mc{A} \times_{TX} \mc{A}'} = [v,w]_{\mc{A}} \oplus [v',w']_{\mc{A}'}.
\ee
\end{defn}
The fiber product is also given by pullback along the diagonal inclusion. There is thus the following commutative diagram over $X$.
\begin{center}
	\begin{tikzpicture}
	\matrix (m) [matrix of math nodes, row sep=2.5em, column sep=2.5em,text height=1.5ex, text depth=0.25ex]
	{	\mc{A} \times_{TX} \mc{A}' & \mc{A}'\\ \mc{A} & TX\\};
	\path[-stealth]
	(m-1-1) edge (m-1-2)
	(m-1-1) edge (m-2-1)
	(m-1-2) edge node [right] {$\rho_{\mc{A}'}$} (m-2-2)
	(m-2-1) edge node [above] {$\rho_\mc{A}$} (m-2-2);
	\end{tikzpicture}
\end{center}
The projections onto $\mc{A}$ and $\mc{A}'$ turn the fiber product $\mc{A} \times_{TX} \mc{A}'$ into both an $\mc{A}$- and an $\mc{A}'$-Lie algebroid. Note that there is a natural isomorphism $\mc{A} \times_{TX} TX \cong \mc{A}$. Consider the case when both $\mc{A}$ and $\mc{A}'$ have dense isomorphism loci. Recall that $Z_\mc{A} = X \setminus X_\mc{A}$  is the degeneracy locus of $\mc{A}$ in this case.
\begin{prop}\label{prop:fiberprodisolocus} Let $\mc{A}, \mc{A}' \to X$ be Lie algebroids with dense isomorphism loci such that $Z_\mc{A}$ and $Z_\mc{A}'$ are disjoint. Then $\rho_\mc{A}$ and $\rho_{\mc{A}'}$ are transverse, and their fiber product has dense isomorphism locus, with complement given by $Z_\mc{A} \cup Z_{\mc{A}'}$.
\end{prop}
\bp As $Z_\mc{A}$ and $Z_{\mc{A}'}$ are disjoint, at each point of $X$, either $\rho_\mc{A}$ or $\rho_{\mc{A}'}$ is an isomorphism. Consequently, it is immediate that $\rho_\mc{A}$ and $\rho_{\mc{A}'}$ are everywhere transverse. The rest of the statement readily follows.
\ep
The above definition immediately generalizes to the more general situation where one has two Lie algebroid morphisms $(\varphi_{\mc{A}'}, {\rm id}_X)\colon \mc{A}' \to \mc{A}$ and $(\varphi_{\mc{A}''}, {\rm id}_X)\colon \mc{A}'' \to \mc{A}$, given three Lie algebroids $\mc{A},\mc{A}',\mc{A}'' \to X$. For example, one can take $\mc{A}'$ and $\mc{A}''$ to both be $\mc{A}$-Lie algebroids. If $\varphi_{\mc{A}'}$ and $\varphi_{\mc{A}''}$ are transverse in $\mc{A}$, i.e.\ $\varphi_{\mc{A}'}(\mc{A}') + \varphi_{\mc{A}''}(\mc{A}'') = \mc{A}$ at each point in $X$, the fiber product $\mc{A}' \times_{\mc{A}} \mc{A}'' \to X$ exists and is both an $\mc{A}'$- and an $\mc{A}''$-Lie algebroid. The following diagram summarizes the situation.
\begin{center}
	\begin{tikzpicture}
	\matrix (m) [matrix of math nodes, row sep=2.5em, column sep=2.5em,text height=1.5ex, text depth=0.25ex]
	{	\mc{A}' \times_{\mc{A}} \mc{A}'' & \mc{A}''\\ \mc{A}' & \mc{A}\\};
	\path[-stealth]
	(m-1-1) edge (m-1-2)
	(m-1-1) edge (m-2-1)
	(m-1-2) edge node [right] {$\varphi_{\mc{A}''}$} (m-2-2)
	(m-2-1) edge node [above] {$\varphi_{\mc{A}'}$} (m-2-2);
	\end{tikzpicture}
\end{center}
There is also an analogue of \autoref{prop:fiberprodisolocus} for this setting, with the same proof.
\begin{rem} There is a more general fiber product of Lie algebroids $\mc{A} \to X$ and $\mc{A}' \to X'$, equipped with transverse Lie algebroid morphisms $(\varphi,f)\colon (\mc{A},X) \to (\mc{B},Y)$ and $(\varphi',f')\colon (\mc{A}',X') \to (\mc{B},Y)$. We will not pursue this here. See \cite{HigginsMackenzie90}.
\end{rem}
\section{Residue maps}
\label{sec:residuemaps}
In this section we discuss residues of Lie algebroid forms. These are a way to extract relevant information of $\mc{A}$-forms over $\mc{A}$-invariant submanifolds where the restriction of $\mc{A}$ is transitive. We further discuss how such residues interact with Lie algebroids morphisms. Residues are used when dealing with geometric structures on $\mc{A}$, in order to extract information along the submanifold (see e.g.\ Chapter \ref{chap:asymplecticstructures} and Section \ref{sec:scgs}).

Let $\mc{A} \to X$ be a Lie algebroid and $D \subseteq X$ an $\mc{A}$-invariant submanifold such that $\mc{A}|_D \to D$ is transitive, i.e.\ the restriction $\rho_\mc{A}|_D$ of the anchor surjects onto $TD$. We will call these \emph{transitive $\mc{A}$-invariant submanifolds}. We obtain a short exact sequence
\be
	0 \to \ker \rho_\mc{A}|_D \to \mc{A}|_D \to TD \to 0.
\ee
Alternatively, we see that $\mc{A}|_D$ is an extension of $TD$ by $\ker \rho_\mc{A}|_D$. The case in which we are mainly interested in is when $D = X \setminus X_\mc{A} = Z_\mc{A}$, the complement of the isomorphism locus of $\mc{A}$ (although in general the restriction need not be transitive). This is an $\mc{A}$-invariant submanifold of $X$ if it is smooth by \autoref{prop:degenlocus}, as $X \setminus X_\mc{A} = X_{\mc{A},n-1}$ where $\dim X = {\rm rank}(\mc{A}) = n$ (see \autoref{exa:denseisolocus}). Dualizing the above sequence, we obtain
\be
	0 \to T^*D \to \mc{A}^*|_D \to (\ker \rho_\mc{A})^*|_D \to 0.	
\ee
We are now in the following more general situation.

Given a short exact sequence $\mc{S}\colon 0 \to E \to W \to V \to 0$ of vector spaces, there is an associated dual sequence $\mc{S}^*\colon 0 \to V^* \to W^* \to E^* \to 0$. For a given $k \in \N$, by taking $k$th exterior powers we obtain a filtration of spaces $\mc{F}^i := \{\rho \in \wedge^k W^* \, | \, \iota_x \rho = 0$ for all $x \in \wedge^i E \}$, for $i = 0,\dots,k+1$. These spaces satisfy $\mc{F}^0 = 0$, $\mc{F}^1 = \wedge^k V^*$, $\mc{F}^i \subset \mc{F}^{i+1}$, and $\mc{F}^{i+1} / \mc{F}^i \cong \wedge^{k-i} V^* \otimes \wedge^i E^*$. Setting $\ell := \dim E$, we have $\mc{F}^{\ell + 1} = \wedge^k W^*$.
\begin{defn} The \emph{residue} of an element $\rho \in \wedge^k W^*$ is defined as its equivalence class ${\rm Res}(\rho) = [\rho] \in \mc{F}^{\ell + 1} / \mc{F}^{\ell} \cong \wedge^{k-\ell} V^* \otimes \wedge^\ell E^*$.
\end{defn}
Upon a choice of trivialization of $\wedge^\ell E^*$, i.e.\ a choice of volume element for $E$, one can view the residue ${\rm Res}(\rho)$ as an element of $\wedge^{k-\ell} V^*$.
\begin{rem}\label{rem:ainvariantresidues} In other words, given a Lie algebroid $\mc{A} \to X$ and a transitive $\mc{A}$-invariant submanifold $D \subseteq X$ with isotropy $E := \ker \rho_\mc{A}|_D$, we obtain a map ${\rm Res}\colon \Omega^\bullet(\mc{A}) \to \Omega^{\bullet-\ell}(D; \wedge^\ell E^*)$, where $\ell = \dim \ker \rho_\mc{A}|_D = {\rm rank}(\mc{A}) - \dim D$.
\end{rem}
There are also lower residues ${\rm Res}_{-m}\colon \wedge^k W^* \to \mc{F}^{\ell+1} / \mc{F}^{\ell-m}$ for $m > 0$. These are always defined but have a better description for forms $\rho \in \wedge^k W^*$ whose higher residues vanish, so that ${\rm Res}_{-m}(\rho) \in \mc{F}^{\ell-m +1} / \mc{F}^{\ell -m} \cong \wedge^{k-\ell+m} V^* \otimes \wedge^{\ell-m} E^*$.
\begin{rem} Under the conditions of \autoref{rem:ainvariantresidues}, these lower residues provide maps ${\rm Res}_{-m}\colon \Omega^\bullet_{1-m}(\mc{A}) \to \Omega^{\bullet - \ell + m}(D; \wedge^{\ell - m} E^*)$, where $\Omega^\bullet_{1-m}(\mc{A})$ is inductively the space of $\mc{A}$-forms all of whose $(1-m)$th or higher residues along $D$ vanish.
\end{rem}
\begin{rem} In case $D = Z_\mc{A}$ for a Lie algebroid $\mc{A}$ with dense isomorphism locus, one should think of these residues as extracting the coefficients in front of the singular parts of an $\mc{A}$-form along $Z_\mc{A}$, as then $(\ker \rho_\mc{A})^*|_D$ consists of ``singular'' generators.
\end{rem}
Given a map of short exact sequences $\Psi\colon \wt{\mc{S}} \to \mc{S}$ with dual map $\Psi^*\colon \mc{S}^* \to \wt{\mc{S}}^*$, there is a corresponding map of filtrations $\Psi^*\colon \mc{F}^i \to \wt{\mc{F}}^i$. Setting $\wt{\ell} := \dim \wt{E}$, we have the following.
\begin{lem}\label{lem:residues} Assume that $\wt{\ell} > \ell$. Then $\wt{\rm Res}(\Psi^* \rho) = 0$ for all $\rho \in \wedge^k W^*$.
\end{lem}
\bp We have $\rho \in \mc{F}^{\ell +1}$ so that $\Psi^* \rho \in \wt{\mc{F}}^{\ell + 1}$. As $\wt{\ell} > \ell$, we have $\wt{\mc{F}}^{\ell+1} \subset \wt{\mc{F}}^{\wt{\ell}} \subset \wt{\mc{F}}^{\wt{\ell}+1}$, so that $\wt{\rm Res}(\Psi^* \rho) = [\Psi^* \rho] \in \wt{\mc{F}}^{\wt{\ell}+1} / \wt{\mc{F}}^{\wt{\ell}}$ vanishes as desired.
\ep
Assuming $\wt{\ell} > \ell$, all lower residues automatically vanish by degree reasons until considering $\Psi^* \rho \in \wt{\mc{F}}^{\ell +1}$. Hence the first possibly nonzero residue is $\wt{\rm Res}_{\ell-\wt{\ell}}(\Psi^* \rho) = [\Psi^* \rho] \in \wt{\mc{F}}^{\ell+1} / \wt{\mc{F}}^{\ell} \cong \wedge^{k-\ell} \wt{V}^* \otimes \wedge^\ell \wt{E}^*$. As in \autoref{lem:residues} we obtain the following.
\begin{lem}\label{lem:residuecommute} In the above setting, assuming $\wt{\ell} \geq \ell$, we have $\Psi^* \circ {\rm Res} = \wt{\rm Res}_{\ell - \wt{\ell}} \circ \Psi^*$.
\end{lem}
Out of the above discussion the following result is immediate.
\begin{prop}\label{prop:ainvsubmfdresidue} Let $(\varphi,f)\colon (\mc{A},X) \to (\mc{A}',X')$ be a Lie algebroid morphism and $D \subseteq X^n$, $D' \subseteq X'$ transitive $\mc{A}$- respectively $\mc{A}'$-invariant submanifolds such that $f\colon (X,D) \to (X',D')$ is a strong map of pairs. Define $\ell = {\rm rank}(\mc{A}) - \dim D$ and $\ell' = {\rm rank}(\mc{A}') - \dim D'$, and assume that $\ell' \geq \ell$. Then ${\rm Res}_{\ell - \ell'} \circ \varphi^* = f^* \circ {\rm Res}$. Moreover, ${\rm Res}_{-m} \circ \varphi^* = 0$ for $m > \ell' - \ell$.
\end{prop}
\bp As $f$ is a strong map of pairs, we have $T f\colon TD \to TD'$. Moreover, by \autoref{lem:morphavisotropy} we see that $\varphi$ restricts to $\varphi\colon \ker \rho_\mc{A}|_D \to \ker \rho_{\mc{A}'}|_{D'}$. Consequently, $\varphi$ induces a map of the relevant short exact sequences defining the residues. The result now follows from \autoref{lem:residuecommute}.
\ep
For Lie algebroids with dense isomorphism loci we obtain the following.
\begin{cor}\label{cor:residuemapsdenseiso} Let $(\varphi,f)\colon (\mc{A},X) \to (\mc{A}',X')$ be a Lie algebroid morphism between Lie algebroids with dense isomorphism loci and transitive degeneracy loci. Assume that $f\colon (X,Z_\mc{A}) \to (X',Z_{\mc{A}'})$ is a strong map of pairs and that ${\rm codim}\, D' \geq {\rm codim}\, D$. Then ${\rm Res}_{\ell-\ell'} \circ \varphi^* = f^* \circ {\rm Res}$ and ${\rm Res}_{-m} \circ \varphi^* = 0$ for $m > \ell' - \ell$, where $\ell = {\rm codim}\, D$ and $\ell' = {\rm codim}\, D'$.
\end{cor}
\bp This follows from \autoref{prop:ainvsubmfdresidue} because degeneracy loci are $\mc{A}$-invariant by \autoref{prop:degenlocus}, and that ${\rm rank}(\mc{A}) = \dim X$ because $X_\mc{A}$ is nonempty.
\ep
\begin{rem} When we speak of transitive degeneracy loci, it is implied that we assume that each degeneracy locus has components being smooth submanifolds.
\end{rem}			
\chapter{Concrete Lie algebroids}
\label{chap:concreteliealgebroids}
\fancypagestyle{empty}{%
	\fancyhf{}%
	\renewcommand\headrulewidth{0pt}%
	\fancyhead[RO,LE]{\thepage}%
}
In this chapter we discuss certain classes of Lie algebroids that are of primary concern to us. The first of these are the Atiyah algebroids, which are transitive Lie algebroids one can associate to a vector bundle. Moreover, there are the Lie algebroids naturally associated to Poisson structures, which we refer to as Poisson algebroids. However, most of this chapter is devoted to studying a class of Lie algebroids we call ideal Lie algebroids. These are defined using an ideal, with the most interesting case being when using a divisor ideal as in Chapter \ref{chap:divisors}. By definition these Lie algebroids then have dense isomorphism locus. Their definition is inspired by specific instances as are used in the literature (see Section \ref{sec:laexamples}).
\subsection*{Organization of the chapter}
This chapter is built up as follows. In Section \ref{sec:atiyahalgebroids} we start with Atiyah algebroids, and relate them to the $\mc{A}$-connections of Section \ref{sec:aconnareps}. Next in Section \ref{sec:poissonalgebroids} we note how Poisson structures interact with Poisson algebroids and their cohomology. We then turn to ideal Lie algebroids in Section \ref{sec:idealliealgebroids}, which come in two flavors. In Section \ref{sec:laexamples} we discuss several explicit examples of Lie algebroids, most of which are ideal Lie algebroids for one of the divisor ideals of Chapter \ref{chap:divisors}, or are obtained using the process of rescaling of Section \ref{sec:rescaling}. We finish with Section \ref{sec:morphilas} by discussing Lie algebroid morphisms between some of the ideal Lie algebroids we considered.
\section{Atiyah algebroids}
\label{sec:atiyahalgebroids}
In this section we discuss a specific kind of Lie algebroid associated to a vector bundle. See e.g.\ \cite{Mackenzie05} for more information. These are relevant in our study of Lie algebroids associated to divisors, as we will pursue in Sections \ref{sec:idealliealgebroids} and \ref{sec:laexamples}.

Let $E \to X$ be a vector bundle. Recall that a \emph{derivation} of $E$ is a first order differential operator $D\colon \Gamma(E) \to \Gamma(E)$ such that there exists a vector field $D_0 \in \mf{X}(X)$ satisfying $D(f\sigma) = f D(\sigma) + D_0(f) \cdot \sigma$ for all $f \in C^\infty(X)$ and $\sigma \in \Gamma(E)$. The vector field $D_0$ associated to $D$ is called the \emph{symbol} of $D$.
\begin{defn} Let $E \to X$ be a vector bundle. The \emph{Atiyah algebroid} ${\rm At}(E) \to X$ of $E$ is the vector bundle of derivations of $E$, with anchor given by the symbol map $\rho_{{\rm At}(E)}(D) = D_0$ for $D \in \Gamma({\rm At}(E))$.
\end{defn}
The Atiyah algebroid is a Lie algebroid with bracket given by the commutator bracket $[D,D'] = D \circ D' - D' \circ D$ for $D, D' \in \Gamma({\rm At}(E))$. One readily checks that this turns ${\rm At}(E)$ into a Lie algebroid. The Atiyah algebroid fits in an exact sequence
\be
	0 \to {\rm End}(E) \to {\rm At}(E) \to TX \to 0,
\ee
so that ${\rm At}(E)$ is a transitive Lie algebroid that is an extension of $TX$ by ${\rm End}(E)$. Here ${\rm End}(E) = E \otimes E^*$ is the space of bundle automorphisms of $E$, which include into ${\rm At}(E)$ as the zeroth order operators. Note that $\dim {\rm End}(E) = {\rm rank}(E)^2$. For later use we record the following explicitly.
\begin{prop}\label{prop:atiyahlarank} Let $E \to X$ be a vector bundle. Then the Atiyah algebroid of $E$ has rank given by ${\rm rank}({\rm At}(E)) = \dim X + {\rm rank}(E)^2$.
\end{prop}
\begin{exa}\label{exa:atiyahnormal} Let $X$ be a manifold and $Y \subseteq X$ a submanifold. Then its normal bundle $NY \to Y$ gives rise to an Atiyah algebroid ${\rm At}(NY) \to Y$. Its sections are all infinitesimal automorphisms of $NY$, which in particular includes the Euler vector field $\mc{E}_{NY}$ of $NY$. We see further that \autoref{prop:atiyahlarank} gives ${\rm rank}({\rm At}(NY)) = \dim Y + {\rm rank}(NY)^2 = \dim Y + {\rm codim}(Y)^2$. If $\dim X = n$ and ${\rm codim}(Y) = k$, this is equal to $n-k + k^2$. Note in particular that if $k = 1$, then ${\rm rank}({\rm At}(NY)) = n$.
\end{exa}
\begin{rem} One can more generally associate an Atiyah algebroid to a principal $G$-bundle $P \to X$ by considering the short exact sequence $0 \to P \times_G \mf{g} \to TP / G \to TX \to 0$, where $P \times_G \mf{g}$ is the associated adjoint bundle (see e.g.\ \cite[Section 3.2]{Mackenzie05}. The specific case presented above is called a \emph{Lie algebroid of derivations} in loc.\ cit.
\end{rem}
Recall that in Section \ref{sec:aconnareps}, given a Lie algebroid $\mc{A} \to X$, we defined the notion of an $\mc{A}$-connection on a vector bundle $E \to X$. An $\mc{A}$-connection has another characterization, in terms of the Atiyah algebroid ${\rm At}(E)$ of $E$, whose proof is immediate.
\begin{prop}\label{prop:aconnbundlemorph} There is a bijection between $\mc{A}$-connections $\nabla$ on $E$ and anchored vector bundle morphisms $(\varphi,f)\colon \mc{A} \to {\rm At}(E)$ given by $\nabla_v = \varphi(v)$ for $v \in \Gamma(\mc{A})$.
\end{prop}
Flatness of $\nabla$ can also be characterized in terms of the bundle morphism description of \autoref{prop:aconnbundlemorph}, as the Lie bracket on ${\rm At}(E)$ is the commutator.
\begin{prop}\label{prop:aconnlamorph} There is a bijection between flat $\mc{A}$-connections $\nabla$ on $E$ and Lie algebroid morphisms $(\varphi,f)\colon \mc{A} \to {\rm At}(E)$ given by $\nabla_v = \varphi(v)$ for $v \in \Gamma(\mc{A})$.
\end{prop}
\bp In light of \autoref{prop:aconnbundlemorph}, we must check that $\varphi\colon \Gamma(\mc{A}) \to \Gamma({\rm At}(E))$ being a Lie algebra homomorphism on sections is equivalent to the condition $F_\nabla \equiv 0$. However, note that $[\varphi(v), \varphi(w)]_{{\rm At}(E)} = \varphi(v) \circ \varphi(w) - \varphi(w) \circ \varphi(v)$ for all $v, w \in \Gamma(\mc{A})$. Consequently, $F_\nabla(v,w) \equiv 0$ is equivalent to $[\varphi(v), \varphi(w)]_{{\rm At}(E)} = \varphi([v,w]_\mc{A})$, from which the claim follows.
\ep

\section{Poisson algebroids}
\label{sec:poissonalgebroids}
In this section we discuss Lie algebroids associated to Poisson bivectors on a manifold $X$, which we refer to as Poisson algebroids. The realization that this can be done is generally attributed to Fuchssteiner \cite{Fuchssteiner82}. The underlying vector bundle is the cotangent bundle $T^*X$, and for this reason these Lie algebroids are also called \emph{cotangent algebroids}. However, this name becomes troublesome when dealing with $\mc{A}$-Poisson structures, as we will do at the end of this section (see \autoref{rem:apoissonalgebroids}).

To not disrupt the flow of this thesis, we assume in this section that the reader is familiar with the basics of Poisson geometry. A short introduction will be provided later, in Section \ref{sec:poissonstructures}. Let $\pi$ be a Poisson structure on $X$, i.e.\ a bivector $\pi \in \mf{X}^2(X)$ for which $[\pi,\pi] = 0$. This determines a Poisson bracket $\{\cdot,\cdot\}_\pi$ by $\{f,g\}_\pi := \pi(df,dg)$ for $f,g\in C^\infty(X)$.

We define an anchored vector bundle structure on $T^*X$ (denoted $\mc{A} = T^*_\pi X$) using the anchor $\rho_\mc{A} = \pi^\sharp\colon v \mapsto \iota_v \pi$ for $v \in T^*X$. Moreover, we define a bracket by
\be
	[v,w]_\mc{A} := \mc{L}_{\pi^\sharp v} w - \mc{L}_{\pi^\sharp w} v - d \pi(v,w), \qquad v, w \in \Gamma(\mc{A}).
\ee
\begin{prop}\label{prop:poissonalgebroid} The bundle $\mc{A} = T^*_\pi X$ is a Lie algebroid if and only if $[\pi,\pi] = 0$.
\end{prop}
We say that $T^*_\pi X$ is the \emph{Poisson algebroid} associated to the Poisson structure $\pi$.
\bp We check that $[\cdot,\cdot]_\mc{A}$ is a compatible Lie bracket. Bilinearity and antisymmetry are clear, and the Leibniz rule holds as for $v,w \in \Gamma(\mc{A})$ and $f \in C^\infty(X)$ we compute
\begin{align*}
	[v,f w]_\mc{A} &= \mc{L}_{\pi^\sharp v}(f w) - \mc{L}_{\pi^\sharp(f w)} v - d \pi(v, f w)\\
	&= f \mc{L}_{\pi^\sharp v} w + (\mc{L}_{\pi^\sharp v} f) \cdot w - f \mc{L}_{\pi^\sharp w} v - f d\pi(v, w)\\
	&= f [v,w]_\mc{A} +  \mc{L}_{\rho_\mc{A}(v)}f \cdot w.
\end{align*}
We omit the verification that the Jacobi identity for $[\cdot,\cdot]_\mc{A}$ is equivalent to $[\pi,\pi] = 0$. Its proof uses the Leibniz rule repeatedly, the fact that $[\pi,\pi] = 0$ if and only if $\{\cdot,\cdot\}_\pi$ satisfies the Jacobi identity, and that as $\mc{A} = T^*_\pi X$ is locally generated by exact one-forms, it suffices to verify the Jacobi identity for elements $u,v,w \in \Gamma(\mc{A})$ given by $v = a df$, $w = b dg$ and $u = c dh$ for $a,b,c,f,g,h \in C^\infty(X)$. See e.g.\ \cite{Vaisman94}.
\ep
\begin{rem} In the proof of \autoref{prop:poissonalgebroid} it is not necessary to use the fact that $T^*X$ is locally generated by exact one-forms. However, it is convenient as it allows one to avoid a more involved computation using Schouten calculus. If one replaces $TX$ by an arbitrary Lie algebroid $\mc{A}$, this shortcut is no longer available in general. See also \autoref{rem:apoissonalgebroids} below.
\end{rem}
\begin{rem} Another characterization of the Lie algebroid $\mc{A} = T^*_\pi X$ is that it is the unique Lie algebroid structure on $T^* X$ such that $\rho_\mc{A}(df) = X_f$ and $[df,dg]_\mc{A} = d \{f,g\}_\pi$ for all $f,g \in C^\infty(X)$ (see e.g.\ \cite{Vaisman94}). Here $X_f = \pi^\sharp(df)$ is the $\pi$-Hamiltonian vector field of $f$, and $\{f,g\}_\pi = \pi(df,dg)$ is the Poisson bracket associated to $\pi$ (see Section \ref{sec:poissonstructures}). The first condition holds by definition for the Lie algebroid structure we consider. Moreover, we have $[df,dg]_\mc{A} = \mc{L}_{X_f} dg - \mc{L}_{X_g} df - d\{f,g\}_\pi = d \iota_{X_f} dg - d \iota_{X_g} df - d \{f,g\}_\pi = d \{f,g\}_\pi$.
\end{rem}
Recall that a Poisson structure defines a cohomology theory, called \emph{Poisson cohomology} $H^*_\pi(X)$ (see Section \ref{sec:poissoncohomology}). Poisson cohomology is isomorphic to the Lie algebroid cohomology of the Poisson algebroid associated to $\pi$. In other words:
\begin{prop}\label{prop:poissoncohlacoh} Let $\pi$ be a Poisson structure on $X$. Then $H^\bullet_\pi(X) \cong H^\bullet(T^*_\pi X)$.
\end{prop}
\begin{rem} In the same way that a Lie algebra $\mf{g}$ carries a canonical linear Poisson structure on its dual, the dual of a Lie algebroid carries a linear Poisson structure (see also \autoref{rem:lalgdpoisson}). Reversing this, a Poisson structure on a Lie algebroid in turn determines a Lie algebroid structure on its dual (see \autoref{rem:apoissonalgebroids} below). There is a notion of when these two are compatible, then forming what is called a \emph{Lie bialgebroid} \cite{MackenzieXu94}. We will not directly use this concept in this thesis.
\end{rem}
\begin{rem}\label{rem:apoissonalgebroids} In Chapter \ref{chap:apoissonstructures} we define $\mc{A}$-Poisson structures on a Lie algebroid $\mc{A}$. Using the same arguments as the ones presented above, it is immediate that $\mc{A}^*$ can be equipped with an analogous $\mc{A}$-Lie algebroid structure $\mc{A}^*_{\pi_\mc{A}} \to \mc{A} \to TX$ using an $\mc{A}$-bivector $\pi_\mc{A} \in \Gamma(\wedge^2 \mc{A})$ if and only if $\pi_\mc{A}$ is $\mc{A}$-Poisson, i.e. satisfies $[\pi_\mc{A},\pi_\mc{A}]_\mc{A} = 0$.
\end{rem}
\section{Ideal Lie algebroids}
\label{sec:idealliealgebroids}
In this section we discuss what we call ideal Lie algebroids. These are Lie algebroids defined using a given ideal of functions. We describe some of their properties before discussing the particular case of ideal Lie algebroids built using divisor ideals. We will be ambiguous about whether we are dealing with real or complex ideals, and similarly with real or complex Lie algebroids.

Let $X$ be a manifold and let $I \subseteq C^\infty(X)$ be an ideal sheaf, and denote by $\mc{V}_X = \Gamma(TX)$ the sheaf of vector fields on $X$. Let $\mc{V}_X(I) := \{v \in \mc{V}_X \, | \, \mc{L}_v I \subseteq I\}$ be the sheaf of derivations preserving the ideal $I$. This is a $C^\infty(X)$-submodule of $\mc{V}_X$.
\begin{lem} Let $I \subseteq C^\infty(X)$ be an ideal. Then $\mc{V}_X(I) \subseteq \mc{V}_X$ is a Lie subalgebra.
\end{lem}
\bp Let $v, w \in \mc{V}_X(I)$ and $f \in I$. Then $\mc{L}_v f = g$ and $\mc{L}_w f = h$ for $g, h \in I$. Hence $\mc{L}_{[v,w]} f = (\mc{L}_v \circ \mc{L}_w - \mc{L}_w \circ \mc{L}_v) f = \mc{L}_v h - \mc{L}_w g$. As $g, h \in I$ and $I$ is closed under addition, we see that $\mc{L}_{[v,w]} f \in I$ for all $f \in I$, so that $[v,w] \in \mc{V}_X(I)$.
\ep
Consequently, $\mc{V}_X(I)$ is a subsheaf of Lie algebras of $\mc{V}_X$. If $\mc{V}_X(I)$ is in addition finitely generated projective (locally free), then $\mc{V}_X(I)$ defines a Lie algebroid $\mc{A}_I \to X$ associated to $I$ with $\Gamma(\mc{A}_I) = \mc{V}_X(I)$ by the Serre--Swan theorem (\autoref{thm:serreswan}).
\begin{defn} Let $I \subseteq C^\infty(X)$ be an ideal. The \emph{ideal Lie algebroid} associated to $I$ is the Lie algebroid $\mc{A}_I$ such that $\Gamma(\mc{A}_I) = \mc{V}_X(I)$ (if it exists, it is unique).
\end{defn}
Of course, it is not guaranteed for a general ideal $I$ that $\mc{V}_X(I)$ is locally free.
\begin{exa} Let $I = C^\infty(X)$ be trivial. Then $\mc{V}_X(I) = \mc{V}_X$ and $\mc{A}_I = TX$.
\end{exa}
\begin{exa}\label{exa:zvx} Let $Z \subseteq X$ be a closed submanifold with vanishing ideal $I_Z$. Then $\mc{V}_X(I_Z)$ consists of all vector fields tangent to $Z$. This is a locally free sheaf if and only if $Z$ has codimension one, in which case $\mc{A}_{I_Z}$ is the $TZ$-rescaling of $TX$. See also Section \ref{sec:logtangentbundle}.
\end{exa}
More nontrivial examples of ideal Lie algebroids will be provided in Section \ref{sec:laexamples}.
\begin{rem} In light of \autoref{defn:secondaryidealla} below, we will sometimes refer to $\mc{A}_I$ as the \emph{primary} ideal Lie algebroid associated to the ideal $I$.
\end{rem}
As $\mc{V}_X(I)$ is a submodule of $\mc{V}_X = \Gamma(TX)$, the anchor of any ideal Lie algebroid is the natural inclusion on sections. However, viewed as vector bundle map, the anchor of $\mc{A}_I$ need not be an isomorphism. More precisely, the isomorphism locus $X_{\mc{A}_I}$ of $\mc{A}_I$ is the complement of the support ${\rm supp}(C^\infty(X) / I)$ of the quotient sheaf, or equivalently, the complement of ${\rm supp}(\mc{V}_X / \mc{V}_X(I))$.

There is another type of Lie algebroid we can consider, by using a different sheaf. Namely, instead consider the sheaf $\mc{W}_X(I) := I \cdot \mc{V}_X \subseteq \mc{V}_X$ consisting of all finite sums of products of elements of $I$ with vector fields on $X$.
\begin{lem}\label{lem:wxisubalgebra} Let $I \subseteq C^\infty(X)$ be an ideal. Then $\mc{W}_X(I) \subseteq \mc{V}_X$ is a Lie subalgebra.
\end{lem}
\bp Given $v, w \in \mc{W}_X(I)$ we have $v = \sum_i f_i v_i$ and $w = \sum_j f'_j w_j$ for some $f_i, f'_j \in I$ and $v_i, w_j \in \mc{V}_X$. Then by repeated application of the Leibniz rule,
\be 
	[v,w] = \sum_{i,j} [f_i v_i, f'_j w_j] = \sum_{i,j} \left( f_i f'_j [v_i, w_j] + f_i (\mc{L}_{v_i} f'_j) \cdot w_j - f'_j (\mc{L}_{w_j} f_i) \cdot v_i\right).
\ee
Each of the three types of terms starts with an element of the ideal $I$, so that we conclude that $[v,w] \in \mc{W}_X(I)$ as desired.
\ep
There is a straightforward relation between $\mc{W}_X(I)$ and $\mc{V}_X(I)$.
\begin{lem}\label{lem:wxvx} Let $I \subseteq C^\infty(X)$ be an ideal. Then $\mc{W}_X(I) \subseteq \mc{V}_X(I)$.
\end{lem}
\bp Let $v \in \mc{W}_X(I)$. Then $v = \sum_i f_i v_i$ for $f_i \in I$ and $v_i \in \mc{V}_X$. For $f \in I$ we then have $\mc{L}_v f = \sum_i f_i (\mc{L}_{v_i} f)$. As $f_i \in I$ for all $i$ and $I$ is an ideal, we conclude that $\mc{L}_v f \in I$ so that $v \in \mc{V}_X(I)$ as desired.
\ep
We thus see that there is a short exact sequence of sheaves given by
\be
	0 \to \mc{W}_X(I) \to \mc{V}_X(I) \to \mc{V}_X(I) / \mc{W}_X(I) \to 0.
\ee
Consider the quotient $C^\infty(X) / I$ and let $Z_I = {\rm supp}(C^\infty(X) / I)$ be its support. It is immediate that $I \subseteq I_{Z_I}$, the vanishing ideal of $Z_I$. Further, as $C^\infty(X)$ and $I$ agree outside of $Z_I$, we conclude that so do $\mc{V}_X$ and $\mc{V}_X(I)$. Moreover, their quotient $\mc{V}_X / \mc{V}_X(I)$ is supported on $Z_I$. Similarly, it also follows that $\mc{W}_X(I)$ and $\mc{V}_X$ agree outside of $Z_I$, and that $\mc{V}_X(I) / \mc{W}_X(I)$ is supported on $Z_I$.

The sheaf $\mc{W}_X(I)$ can again be locally free, in which case it gives rise to a Lie algebroid $\mc{B}_I \to X$ with $\Gamma(\mc{B}_I) = \mc{W}_X(I)$.
\begin{defn}\label{defn:secondaryidealla} Let $I \subseteq C^\infty(X)$ be an ideal. The \emph{secondary ideal Lie algebroid} associated to $I$ is the Lie algebroid $\mc{B}_I$ such that $\Gamma(\mc{B}_I) = \mc{W}_X(I)$.
\end{defn}
As for $\mc{A}_I$, the secondary ideal Lie algebroid $\mc{B}_I$ is unique if it exists.
\begin{exa} Let $I = C^\infty(X)$ be trivial. Then $\mc{W}_X(I) = \mc{V}_X$ and $\mc{B}_I = TX$.
\end{exa}
\begin{exa}\label{exa:zwx} Let $Z \subseteq X$ be a closed submanifold with vanishing ideal $I_Z$. Then $\mc{W}_X(I_Z)$ consists of all vector fields which vanish on $Z$. As in \autoref{exa:zvx} this is a locally free sheaf if and only if $Z$ has codimension one, in which case $\mc{B}_{I_Z}$ is the $0$-rescaling of $TX$. See also Section \ref{sec:zerotangentbundle} for more discussion.
\end{exa}
When $I$ is locally principal and locally generated by a function with nowhere dense zero set, $\mc{W}_X(I)$ is always locally free. Such ideals were called divisor ideals in Chapter \ref{chap:divisors}. This follows from the following noted also by Nistor \cite[Lemma 1.11]{Nistor15}.
\begin{prop}\label{prop:dividealalgebroid} Let $\mc{A} \to X$ be a Lie algebroid and $I \subseteq C^\infty(X)$ a divisor ideal. Then $I \cdot \Gamma(\mc{A})$ is locally free Lie subalgebra of $\Gamma(\mc{A})$, hence there exists a Lie algebroid $\mc{A}'$ such that $\Gamma(\mc{A}') = I \cdot \Gamma(\mc{A})$ by the Serre--Swan theorem (\autoref{thm:serreswan}).
\end{prop}
\bp The argument that $I \cdot \Gamma(\mc{A})$ is a Lie subalgebra is identical to \autoref{lem:wxisubalgebra}. Let $f$ be a local generator for $I$, and $\sigma, \tau \in \Gamma(\mc{A})$. Then by applying the Leibniz rule we obtain $[f \sigma, f \tau]_\mc{A} = f^2 [\sigma, \tau]_\mc{A} + f (\mc{L}_{\rho_\mc{A}(\sigma)} f) \cdot \tau - f (\mc{L}_{\rho_\mc{A}(\tau)} f) \cdot \sigma \in I \cdot \Gamma(\mc{A})$. To see that $I \cdot \Gamma(\mc{A})$ is locally free, note that there is a short exact sequence of sheaves
\be
	0 \to I \cdot \Gamma(\mc{A}) \to \Gamma(\mc{A}) \to \Gamma(\mc{A}) / I \cdot \Gamma(\mc{A}) \to 0.
\ee
As both $\Gamma(\mc{A})$ and $I$ are locally free, so is $\Gamma(\mc{A}) / I \cdot \Gamma(\mc{A}) \cong I^{-1}$ (the inverse of $I$ as a sheaf). The same then holds for $I \cdot \Gamma(\mc{A})$ by \autoref{lem:locfreemodules}.i) below.
\ep
In the previous proof we used the following well-known algebraic lemma, which holds because any short exact sequence ending in a projective module splits.
\begin{lem}\label{lem:locfreemodules} Let $0 \to \mc{E} \to \mc{F} \to \mc{G} \to 0$ be an extension of $C^\infty(X)$-modules.
\bi
	\item[i)]If $\mc{F}$ and $\mc{G}$ are locally free, then so is $\mc{E}$;
	\item[ii)] If $\mc{E}$ and $\mc{G}$ are locally free, then so is $\mc{F}$.
\ei
\end{lem}
\begin{rem} It is not true in general that in the situation of \autoref{lem:locfreemodules}, if we know that $\mc{E}$ and $\mc{F}$ are locally free, the same is true for $\mc{G}$.
\end{rem}
Applying the above proposition to $\mc{A} = TX$ we see that if $I$ is a divisor ideal, we conclude that $\mc{B}_I$ exists. Moreover, in this case you can conclude that $\mc{A}_I$ exists using \autoref{lem:locfreemodules}.ii) and \autoref{lem:wxvx} by considering the quotient sheaf $\mc{V}_X(I) / \mc{W}_X(I)$.
\begin{rem}\label{rem:isoloci} The inclusions $\mc{W}_X(I) \subseteq \mc{V}_X(I) \subseteq \mc{V}_X$ of \autoref{lem:wxvx} translate into a sequence of Lie algebroid morphisms $\mc{B}_I \to \mc{A}_I \to TX$, both with isomorphism locus equal to $X \setminus Z_I$. In other words, the anchor of $\mc{B}_I$ factors through that of $\mc{A}_I$, making $\mc{B}_I$ into an $\mc{A}_I$-Lie algebroid.
\end{rem}
We next study to what extent the sheaf $\mc{V}_X(I)$ depends on the ideal $I$. Recall that the \emph{radical} of an ideal $I \subseteq C^\infty(X)$ is defined as $\sqrt{I} := \{f \in C^\infty(X) \, | \, f^n \in I \text{ for some } n \in \N\}$. Then $I \subseteq \sqrt{I}$ and $Z_I = Z_{\sqrt{I}}$. Moreover, we have the following.
\begin{prop}\label{prop:locprincidealmodule} Let $I$ be a locally principal ideal. Then $\mc{V}_X(I) = \mc{V}_X(\sqrt{I})$.
\end{prop}
\bp This can be verified locally. Let $f \in \sqrt{I}$ be such that $\sqrt{I} = \langle f \rangle$, so that $I = \langle f^n \rangle$ for some $n \in \N$. Take $V \in \mc{V}_X(I)$, so that $\mc{L}_V f^n \in I$, say $\mc{L}_V f^n = g$. Then $g = g' f^n$ for some function $g'$, hence $n f^{n-1} \mc{L}_V f = g' f^n$, from which it follows that $\mc{L}_V f = \frac1n g' f$, so that $\mc{L}_V f \in \sqrt{I}$. As $f$ generates $\sqrt{I}$, we conclude that $V \in \mc{V}_X(\sqrt{I})$, whence $\mc{V}_X(I) \subseteq \mc{V}_X(\sqrt{I})$. For the other inclusion, take $W \in \mc{V}_X(\sqrt{I})$ and $h \in I$. Then $h = h' f^n$ for some function $h'$, hence $\mc{L}_W h = \mc{L}_W h' f^n = n h' f^{n-1} \mc{L}_W f + f^n \mc{L}_W h'$. As $\mc{L}_W f \in \sqrt{I}$, write $\mc{L}_W f = f' f$ for some function $f'$. Then we see that $\mc{L}_W h = (n h' f' + \mc{L}_W h') f^n \in I$, so that $W \in \mc{V}_X(I)$, from which we conclude that $\mc{V}_X(\sqrt{I}) \subseteq \mc{V}_X(I)$ hence equality follows.
\ep
Hence the ideal Lie algebroids $\mc{A}_I$ associated to locally principal ideals $I$ are only sensitive to $I$ up to taking powers. An analogous statement does not hold for the secondary ideal Lie algebroids $\mc{B}_I$. We next discuss the relation between $I$ and the vanishing ideal $I_{Z_I}$ of its support $Z_I$. This is similar to \autoref{prop:divvanishingideal}. We saw in \autoref{rem:isoloci} that $\mc{A}_I$ has isomorphism locus equal to the complement of $Z_I$, i.e.\ $X_{\mc{A}_I} = X \setminus Z_I$. This means that $Z_I$ is the degeneracy locus of $\mc{A}_I$, so that by \autoref{prop:degenlocus} we know that $Z_I$ is $\mc{A}_I$-invariant if it is smooth.
\begin{prop}\label{prop:ziaiinvariant} Let $I \subseteq C^\infty(X)$ be an ideal such that $\mc{A}_I$ exists. Then $Z_I$ is $\mc{A}_I$-invariant if it is smooth, i.e.\ $\Gamma(\mc{A}_I)$ consists of vector fields tangent to $Z_I$.
\end{prop}
We finish by discussing ideal Lie algebroids that are obtained using divisor ideals. All ideal Lie algebroids we will consider in this thesis are of this type. Recall that we discussed with \autoref{prop:dividealalgebroid} another way to obtain a new Lie algebroid out of a given Lie algebroid $\mc{A}$ using a divisor ideal $I$, by considering the sheaf $I \cdot \Gamma(\mc{A})$. Using \autoref{prop:locprincideal} so that $I = I_\sigma$ for a divisor $(U,\sigma)$, evaluation gives a map $\Gamma(\mc{A} \otimes U^*) \to I \cdot \Gamma(\mc{A})$. In case $I = I_Z$ is a log divisor ideal, we obtain the following.
\begin{prop}\label{prop:seclaisotype} Let $Z = (L_Z,s)$ be a log divisor. Then $\mc{B}_{I_Z} \cong TX \otimes L^*_Z$.
\end{prop}
This immediately implies the following.
\begin{cor} We have $\det(\mc{B}_{I_Z}^*) \otimes \det(TX) \cong L_Z^n$, where $n = \dim X$.
\end{cor}
In other words, the divisor of $\mc{B}_{I_Z}$ is equal to the $n$-fold product $(L_Z^n, s^n)$, and $I_{\mc{B}_{I_Z}} = I_Z^n$.
If a divisor is expressed as a product $(U \otimes U', \sigma \otimes \sigma')$, recall from \autoref{prop:proddivisor} that its divisor ideal is given by the product $I_\sigma \cdot I_{\sigma'}$. If $\mc{A}_{I \cdot I'}$ exists, each of the divisors $(U,\sigma)$ and $(U',\sigma')$ then themselves give rise to ideal Lie algebroids $\mc{A}_I$ and $\mc{A}_{I'}$. Moreover, $\mc{A}_{I \cdot I'}$ is obtained from $\mc{A}_I$ and $\mc{A}_{I'}$ using fiber product.
\begin{prop}\label{prop:ilaproduct} Let $I \cdot I'$ be a divisor ideal for which the primary ideal Lie algebroid exists. Then $\mc{A}_{I \cdot I'} = \mc{A}_I \times_{TX} \mc{A}_{I'}$.
\end{prop}
In favorable cases, the divisor of $\mc{A}_{I \cdot I'}$ is the product of divisors of $\mc{A}_I$ and $\mc{A}_{I'}$.
\begin{rem} While $\mc{V}_X(I \cdot I')$ being locally free implies both $\mc{V}_X(I)$ and $\mc{V}_X(I')$ are also locally free, the converse is not necessarily true if $Z_I$ and $Z_{I'}$ have nonempty intersection. For example, take $I = I_Z$ and $I' = I_{Z'}$ for $Z,Z' \subseteq X$ hypersurfaces which do not intersect transversely (see also Section \ref{sec:logtangentbundle}).
\end{rem}
If one expresses a divisor $(U,\sigma)$ as the trivial product $(U \otimes \underline{\R}, \sigma \otimes \underline{1})$, the above proposition is consistent with the natural isomorphism $\mc{A} \times_{TX} TX \cong \mc{A}$. Moreover, we see that $\mc{A}_{I \cdot I'}$ is both an $\mc{A}_I$- and an $\mc{A}_{I'}$-Lie algebroid (again, if it exists).
\section{Examples}
\label{sec:laexamples}
In this section we discuss several examples of Lie algebroids with dense isomorphism loci. They are all obtained as ideal Lie algebroids using divisor ideals, or as rescalings thereof. We discuss their local descriptions, possible residue maps, and Lie algebroid cohomology (if available).
\subsection{Log-tangent bundle}
\label{sec:logtangentbundle}
We start with the log-tangent bundle, which is associated to a log divisor.

Let $Z = (L,s)$ be a log divisor on $X$ with associated ideal $I_Z = s (\Gamma(L^*))$. This ideal is exactly the vanishing ideal of $Z$, and $\mc{V}_X(I_Z)$ consists of all vector fields preserving $I_Z$. But then $\mc{V}_X(I_Z)$ is the sheaf of vector fields tangent to $Z$ \cite{Melrose93}. In local adapted coordinates $(z,x_2,\dots,x_n)$ around $Z = \{z = 0\}$ with $I_z = \langle z \rangle$, one has $\mc{V}_X(I_Z) = \langle z \partial_z, \partial_{x_2}, \dots, \partial_{x_n} \rangle$, which shows this is a locally free sheaf.

\begin{defn} Let $Z = (L,s)$ be a log divisor on $X$. The \emph{log-tangent bundle} $TX(- \log Z)$ is the ideal Lie algebroid on $X$ with $\Gamma(TX(- \log Z)) = \mc{V}_X(I_Z)$.
\end{defn}
We will sometimes denote the log-tangent bundle by $\mc{A}_Z = TX(-\log Z)$ when the underlying manifold $X$ is understood. Hence, $\mc{A}_Z = \mc{A}_{I_Z}$ in the notation of the previous section. It is immediate that the isomorphism locus of $TX(-\log Z)$ is given by $X \setminus Z$, so that $Z$ is an $\mc{A}_Z$-invariant submanifold. In analogy with the holomorphic case, denote $\Omega^k(\log Z) = \Omega^k(TX(-\log Z))$.
\begin{rem} This Lie algebroid is also called the \emph{$b$-tangent bundle} \cite{Melrose93, NestTsygan96, GuilleminMirandaPires14}, especially when $Z = \partial X$, and is then denoted by ${}^b TX$, where the $b$ stands for `boundary'. We use the name log-tangent bundle as it is more consistent with other Lie algebroids we consider, and its notation allows us to keep track of $Z$.
\end{rem}
In terms of a local generator $z$ for $I_Z$ which is extended to a coordinate system $(z,x_2,\dots,x_n)$ of $X$, such generators can be given by
\be
\Gamma(TX(-\log Z)) = \mc{V}_X(I_Z) = \langle z \partial_z, \partial_{x_i} \rangle.
\ee
Its dual is given in these coordinates by
\be
\Gamma(T^*X(\log Z)) = \Gamma(\mc{A}_Z^*) = \langle d\log z, dx_i \rangle,
\ee
with $d\log z = \frac{dz}{z}$ explaining the name of this Lie algebroid.
\begin{rem} It follows immediately from the definition that $TX(-\log Z) = {}^{TZ} TX$, the $TZ$-rescaling of $TX$. Indeed, sections of $TX(-\log Z)$ consist of all vector fields on $X$ which are tangent to $Z$, as they preserve its vanishing ideal.
\end{rem}
Because $Z$ is in fact a transitive $\mc{A}_Z$-invariant submanifold, the anchor gives rise to a line bundle $\mathbb{L}_Z := \ker \rho_{\mc{A}_Z}$ over $Z$. This line bundle is always trivial and has a canonical trivialization \cite[Proposition 4]{GuilleminMirandaPires14}. Locally this is given by the Euler vector field $\mc{E}_{NZ} = z \frac{\partial}{\partial z}$, for $z$ a local defining function for $Z$, i.e.\ a local generator of $I_{Z}$.

As a consequence, the sections of the isotropy of $\mc{A}_Z$ over $Z$ are equal to the endomorphisms of $NZ$. In fact, $\mc{A}_Z$ restricts over $Z$ to the Atiyah algebroid ${\rm At}(NZ)$ of the normal bundle of $Z$ \cite{GualtieriLi14, CannasDaSilvaWeinstein99}.
\begin{rem} In \cite{CannasDaSilvaWeinstein99}, it is asked how much of $\mc{A}_Z$ is determined by ${\rm PD}_{\Z_2}[Z]$ and $TX$. In Section \ref{sec:stablebundleiso} we give an answer to this, namely that in the even-dimensional orientable case, the isomorphism class of $\mc{A}_Z$ is determined by ${\rm PD}_{\Z_2}[Z]$ and $TX$ up to knowing the Euler class of $\mc{A}_Z$.
\end{rem}
As $Z$ is a transitive $\mc{A}_Z$-invariant submanifold, the log-tangent bundle admits a residue map as in \autoref{sec:residuemaps}. The restriction of $TX(-\log Z)$ to $Z$ surjects onto $TZ$ via the anchor map $\rho_{\mc{A}_Z}$, giving the exact sequence
\be
	0 \to \mathbb{L}_Z \to TX(-\log Z)|_Z \stackrel{\rho_{\mc{A}_Z}}{\to} TZ \to 0,
\ee
where $\mathbb{L}_Z \to Z$ is the kernel of $\rho_{\mc{A}_Z}$. Upon dualizing this sequence we get a projection map ${\rm Res}_Z\colon \Omega^k(\log Z) \to \Omega^{k-1}(Z)$, which fits in the residue sequence
\be
	0 \to \Omega^\bullet(X) \stackrel{\rho_{\mc{A}_Z}^*}{\to} \Omega^\bullet(\log Z) \stackrel{{\rm Res}_Z}{\to} \Omega^{\bullet - 1}(Z) \to 0.
\ee
In terms of the local coordinate system above, where we have $\Gamma(TX(-\log Z)) = \langle z \partial_z, \partial_{x_2}, \dots, \partial_{x_n} \rangle$, a given log $k$-form $\alpha \in \Omega^k(\log Z)$ can be expressed as
\be
\alpha = d \log z \wedge \alpha_0 + \alpha_1,
\ee
with $\alpha_i$ smooth forms. The inclusion $j_Z\colon Z \hookrightarrow X$ gives ${\rm Res}_Z(\alpha) = j_Z^* \alpha_0$. The following result referred to as the Mazzeo--Melrose theorem shows the above sequence splits, and identifies the Lie algebroid cohomology of the log-tangent bundle $TX(-\log Z)$ in terms of $X$ and $Z$. Explicitly (see e.g.\ \cite{MarcutOsornoTorres14}), a splitting is given in a tubular neighbourhood of $Z$ by the map $\sigma\colon \Omega^{\bullet-1}(Z) \to \Omega^\bullet(\log Z)$, $\beta \mapsto d\log \lambda \wedge p^*(\beta)$, where $p\colon NZ \to Z$ is the projection, and $\lambda\colon M \setminus Z \to \R_{> 0}$ is given by $\lambda(z) = |z|$ near $Z$ and by $\lambda \equiv 1$ further away.
\begin{thm}[\cite{Melrose93}]\label{thm:loglacohomology} Let $(X,Z)$ be a log pair. Then there is an isomorphism $H^k(TX(-\log Z)) \cong H^k(X) \oplus H^{k-1}(Z)$.
\end{thm}
If $TX(-\log Z)$ carries an orientation, it induces an orientation of $TX$ over the isomorphism locus $X \setminus Z$. However, such an orientation can never come from a global orientation on $TX$. Indeed, either $TX$ itself is nonorientable, or if both $\mc{A}_Z$ and $TX$ are oriented, we know that $Z$ is separating (see \autoref{cor:orlogpairseparating}). We can then decompose $X \setminus Z = X_+ \cup X_-$ according to whether the orientations agree or disagree. This is readily seen from the local description of $TX(-\log Z)$ above, or more intrinsically from the fact that $\det \rho_{\mc{A}_Z}\colon \det(\mc{A}_Z) \to \det(TX)$ vanishes transversally. Viewing $\det \rho_{\mc{A}_Z} \in \Gamma(\det(\mc{A}_Z^*) \otimes \det(TX))$, we have the following.
\begin{prop}\label{prop:logtangentlogdiv} Let $(X,Z)$ be a log pair. Then $(\det(\mc{A}_Z^*) \otimes \det(TX), \det \rho_{\mc{A}_Z})$ is a log divisor. Consequently, we have that $\det(\mc{A}_Z^*) \otimes \det(TX) \cong L_Z$.
\end{prop}
\bp This follows from the local description of $\mc{A}_Z$, together with \autoref{prop:hypersurfacedivisor} on unicity of log divisors with zero set $Z$. In terms of the local coordinate system $(z,x_i)$ used above, we have $\det(\rho_{\mc{A}_Z})\colon z \partial z \wedge \partial_{x_i} \mapsto z (\partial_z \wedge \partial_{x_i})$.
\ep
\begin{rem} In other words, \autoref{prop:logtangentlogdiv} shows that the divisor ideal associated to $\mc{A}_Z$ is $I_Z$, and that $\mc{A}_Z$ is of log divisor-type.
\end{rem}
\subsection{Normal-crossing log-tangent bundle}
There is a variant of the log-tangent bundle which uses a normal-crossing log divisor (see Section \ref{sec:normalcrosslogdivisor}). Given a normal-crossing log divisor $\underline{Z} = (U,\sigma)$, where $Z = \cup_i Z_i$ is a normal-crossing divsor, we can again consider the sheaf $\mc{V}_X(I_Z)$ of all vector fields tangent to $Z$. Due to \autoref{prop:saitocriterion}, this is a locally free sheaf.
\begin{defn}[\cite{GualtieriLiPelayoRatiu17}] Let $\underline{Z} = (U,\sigma)$ be a normal-crossing log divisor on $X$. The \emph{normal-crossing log tangent bundle} $TX(-\log \underline{Z}) \to X$ is the ideal Lie algebroid with $\Gamma(TX(- \log \underline{Z})) = \mc{V}_X(I_Z)$.
\end{defn}
We will also use the notation $\mc{A}_{\underline{Z}}$ for this Lie algebroid. In terms of local generators $z_i$ for $I_{Z_i}$ extended to a coordinate system $(z_i,x_j)$ of $X$, we have that
\be
	\Gamma(TX(-\log \underline{Z})) = \mc{V}_X(I_{\underline{Z}}) = \langle z_i \partial_{z_i}, \partial_{x_j} \rangle.
\ee
It is immediate that $\mc{A}_{\underline{Z}}$ is an $\mc{A}_{Z_i}$-Lie algebroid for all $i$. In fact, using \autoref{prop:logprodnormalcrossing} and \autoref{prop:ilaproduct} repeatedly, we obtain the following.
\begin{prop} Let $\underline{Z} = (U,\sigma)$ be a normal-crossing log divisor on $X$, with $Z = \cup_{i=1}^m Z_i$. Then
\be
	\mc{A}_{\underline{Z}} = \mc{A}_{Z_1} \times_{TX} \dots \times_{TX} \mc{A}_{Z_m}.
\ee
\end{prop}
Letting $Z_\tau = \cap_{i \in \tau} Z_i$ for $\tau \subseteq \{1,\dots,m\}$, we see that $\mc{A}_{\underline{Z}}$ is an $\mc{A}_{\underline{Z}_\tau}$-Lie algebroid for any $\tau$. There are residue maps ${\rm Res}_{Z_\tau}\colon \Omega^\bullet(\log \underline{Z}) \to \Omega^{\bullet - |\tau|}(Z_\tau)$. These allow us to establish the Lie algebroid cohomology of $\mc{A}_{\underline{Z}}$, similarly to \autoref{thm:loglacohomology}.
\begin{thm}[\cite{GualtieriLiPelayoRatiu17}] Let $\underline{Z} = (U,\sigma)$ be a normal-crossing log divisor on $X$, with $Z = \cup_{i=1}^m Z_i$. Then we have the following isomorphism in cohomology:
\be
	H^k(TX(-\log \underline{Z})) \cong H^k(X) \oplus \bigoplus_{|\tau| \leq k} H^{k-|\tau|}(Z_\tau).
\ee
\end{thm}
\subsection{Complex log-tangent bundle}
\label{sec:complexlogtangent}
In this section we discuss the complex log-tangent bundle associated to a complex log divisor as defined in \cite{CavalcantiGualtieri15}. Similar to the real log divisor case, given a complex log divisor $D = (U,\sigma)$, there is an associated complex ideal Lie algebroid.
\begin{defn} Let $D = (U,\sigma)$ be a complex log divisor on $X$. The \emph{complex log-tangent bundle} $TX(-\log D)$ is the complex ideal Lie algebroid on $X$ with $\Gamma(TX(-\log D)) = \mc{V}_X(I_D)$.
\end{defn}
We will also denote the complex log-tangent bundle by $\mc{A}_D = TX(-\log D)$. In terms of a local generator $w \in I_\sigma = I_D$, we can express the sections of the complex log-tangent bundle as
\be
	\Gamma(TX(-\log D)) = \mc{V}_X(I_D) = \langle w \partial_w, \partial_{\overline{w}}, \partial_{x_i} \rangle.
\ee
From this we see that an $\mc{A}_D$-form can be locally expressed as $d \log w \wedge \alpha + \beta$, with unique $\alpha, \beta$ contained in the subalgebra generated by $\langle d\overline{w}, dx_i \rangle$.

We described in Section \ref{sec:cplxlogdivisor} that a complex log divisor has a complex conjugate divisor $\overline{D} = (\overline{U}, \overline{\sigma})$ with $I_{\overline{\sigma}} = \langle \overline{w} \rangle$. Consequently, the complex log-tangent bundle $TX(-\log \overline{D})$ gives rise to the module
\be
	\Gamma(TX(-\log \overline{D})) = \mc{V}_X(I_{\overline{D}}) = \langle \partial_w, \overline{w} \partial_{\overline{w}}, \partial_{x_i} \rangle.
\ee
\begin{prop}\label{prop:complexlogtransverse} Let $D = (U,\sigma)$ be a complex log divisor with complex conjugate $\overline{D} = (\overline{U}, \overline{\sigma})$. Then $TX(-\log D)$ and $TX(-\log \overline{D})$ are transverse Lie algebroids.
\end{prop}
\bp In the local coordinate system used above, sections of the complexified tangent bundle can be written as $\Gamma(TX_\C) = \langle \partial_w, \partial_{\overline{w}}, \partial_{x_i} \rangle$. Transversality follows from the local description of $\mc{A}_D$ and $\mc{A}_{\overline{D}}$.
\ep
Consequently, we can form the fiber product Lie algebroid $\mc{A}_D \times_{TX} \mc{A}_{\overline{D}}$ out of a complex log divisor (see \autoref{def:fiberproduct}). By \autoref{prop:ilaproduct}, this corresponds to the product of the pair of conjugate complex log divisors. This complex Lie algebroid has as local description that
\be
	\Gamma(\mc{A}_D \times_{TX} \mc{A}_{\overline{D}}) = \langle w \partial_w, \overline{w} \partial_{\overline{w}}, \partial_{x_i} \rangle.
\ee
\begin{rem} There is a description of $\mc{A}_D$ as a Lie subalgebroid of the complexified Atiyah algebroid ${\rm At}_\C(U)$ of $U$. See \cite[Proposition 1.16]{CavalcantiGualtieri15}.
\end{rem}
We finish by discussing the Lie algebroid cohomology of the complex log-tangent bundle. As the isomorphism locus of $\mc{A}_D$ is given by $X\setminus D$, the inclusion $i\colon X \setminus D \hookrightarrow X$ induces a map on forms $i^*\colon \Omega^k(\log D) \to \Omega^k(X \setminus D;\C)$.
\begin{thm}[{\cite[Theorem 1.3]{CavalcantiGualtieri15}}] The map $i^*$ provides an isomorphism in cohomology, $H^k(i^*)\colon H^k(\log D) \stackrel{\cong}{\to} H^k(X \setminus D; \C)$.
\end{thm}
Further, there is a residue map which locally is given by ${\rm Res}(d\log w \wedge \alpha + \beta) = j_D^* \alpha$, where $j_D\colon D \hookrightarrow X$ is the inclusion. This is a cochain morphism, so it induces a map in cohomology, ${\rm Res}_*\colon H^k(\log D) \to H^{k-1}(D;\C)$.
\subsection{Elliptic tangent bundle}
\label{sec:elltangentbundle}
In this section we discuss the ideal Lie algebroid constructed out of an elliptic divisor. For more information, see \cite{CavalcantiGualtieri15}. Let $(X,|D|)$ be an elliptic pair with associated ideal $I_{|D|} = q (\Gamma(R^*))$ and elliptic divisor $(R,q)$. Note that $I_{|D|}$ is not the vanishing ideal of $D$. The associated submodule $\mc{V}_X(I_{|D|})$ defines a sheaf of locally constant rank \cite{CavalcantiGualtieri15}, locally generated in appropriate polar coordinates $(r,\theta,x_3,\dots,x_n)$ around $D = \{r = 0\}$ such that $I_{|D|} = \langle r^2 \rangle$ by $\langle r \partial_r, \partial_\theta, \partial_{x_3}, \dots, \partial_{x_n} \rangle$.
\begin{defn} The \emph{elliptic tangent bundle} $TX(-\log|D|) \to X$ is the ideal Lie algebroid on $X$ with $\Gamma(TX(-\log|D|)) = \mc{V}_X(I_{|D|})$.
\end{defn}
The isomorphism locus of $TX(-\log |D|)$ is given by $X \setminus D$. We will sometimes denote the elliptic tangent bundle by $\mc{A}_{|D|}$. As for the log-tangent bundle, we denote $\Omega^k(\log |D|) = \Omega^k(TX(-\log |D|))$.

The elliptic tangent bundle admits several residue maps \cite{CavalcantiGualtieri15}, three of which we will now describe. The \emph{elliptic residue} ${\rm Res}_q$ comes from considering the restriction of $TX(-\log |D|)$ to $D$, which fits in an exact sequence
\begin{equation}
	\label{eqn:ellipticrestriction}
	0 \to \underline{\R} \oplus \mf{k} \to TX(-\log |D|)|_D \to TD \to 0,
\end{equation}
with $\underline{\R}$ generated by the Euler vector field of $ND$, and $\mf{k} \cong \wedge^2 N^*D \otimes R$. Choosing a coorientation for $D$, i.e.\ a trivialization of $ND$, also trivializes $\mf{k}$. Dualizing the above sequence we obtain a projection map ${\rm Res}_q\colon \Omega^k(\log |D|) \to \Omega^{k-2}(D,\mf{k}^*)$, with kernel $\Omega^\bullet_0(\log |D|)$ the subcomplex of $\Omega^\bullet(\log |D|)$ of zero elliptic residue forms. 

Denote by $S^1ND$ the circle bundle of $ND$. The \emph{radial residue} ${\rm Res}_r$ arises from quotienting \eqref{eqn:ellipticrestriction} by the Euler vector field of $ND$, giving the extension
\be
0 \to \mf{k} \to {\rm At}(S^1 ND) \to TD \to 0,
\ee
where ${\rm At}(S^1 ND)$ is the associated Atiyah algebroid of $S^1 ND$. Noting that the restriction $TX(-\log |D|)|_D$ is a trivial extension of ${\rm At}(S^1 ND)$, the elliptic residue factors through the radial residue ${\rm Res}_r\colon \Omega^k(\log |D|) \to \Gamma(D, \wedge^{k-1} {\rm At}(S^1  ND)^*)$. When the elliptic residue vanishes, the radial residue naturally maps to $\Omega^{k-1}(D)$ without needing a coorientation. Finally, there is a \emph{$\theta$-residue} ${\rm Res}_\theta\colon \Omega^k_0(\log|D|) \to \Omega^{k-1}(D)$, which we will only define for forms with zero elliptic residue. We provide a description of these residue maps in local coordinates. In the adapted coordinate system around $D$ as above, where $\Gamma(TX(-\log |D|)) = \langle r \partial_r, \partial_\theta, \partial_{x_3}, \dots, \partial_{x_n} \rangle$, a given elliptic $k$-form $\alpha \in \Omega^k(\log |D|)$ can be locally written as
\be
\alpha = d \log r \wedge d\theta \wedge \alpha_0 + d \log r \wedge \alpha_1 + d\theta \wedge \alpha_2 + \alpha_3,
\ee
where each $\alpha_i$ is a smooth form. The inclusion $j_D\colon D \hookrightarrow X$ gives ${\rm Res}_q(\alpha) = j_D^* \alpha_0$, ${\rm Res}_r(\alpha) = (d\theta \wedge \alpha_0 + \alpha_1)|_D$, and we set ${\rm Res}_\theta(\alpha) = j_D^*\alpha_2$. Moreover, we see that ${\rm Res}_q(\alpha) = \iota_{\partial_\theta} {\rm Res}_r(\alpha)$.

The Lie algebroid cohomology of the elliptic tangent bundle as well as its zero elliptic residue version can be expressed in terms of $X$ and $D$.
\begin{thm}[{\cite[Theorems 1.8, 1.12]{CavalcantiGualtieri15}}]\label{thm:elllacohomology} Let $(X,|D|)$ be an elliptic pair. Then $H^k(\log |D|) \cong H^k(X \setminus D) \oplus H^{k-1}(S^1 ND)$. Further, one has $H^k_0(\log |D|) \cong H^k(X \setminus D) \oplus H^{k-1}(D)$.
\end{thm}
The above isomorphisms are both induced by the maps $(i^*, {\rm Res}_r)$ with $i\colon X\setminus D \hookrightarrow X$ the inclusion of the divisor complement, noting ${\rm Res}_r$ naturally maps to $\Omega^{k-1}(D)$ when the elliptic residue vanishes (in that case, ${\rm Res}_r(\alpha) = \alpha_1|_D$ in the above local coordinates). Similarly to \autoref{prop:logtangentlogdiv} we note the following.
\begin{prop}[{\cite[Lemma 3.4]{CavalcantiGualtieri15}}]\label{prop:elltangentlogdiv} Let $(X,|D|)$ be an elliptic pair. Then we have that $(\det(\mc{A}_{|D|}^*) \otimes \det(TX), \det \rho_{\mc{A}_{|D|}}) \cong (R,q)$.
\end{prop}
In other words, the divisor of $\mc{A}_{|D|}$ is isomorphic to $(R,q)$.
\bp Choose local compatible coordinates $(r,\theta,x_i)$ around a point in $D$ and split $r^2 = x^2 + y^2$. Then we have $r \partial_r = x \partial_x + y \partial_y$ and $\partial_{\theta} = x \partial_y - y \partial_x$, so that $r \partial_r \wedge \partial_\theta = r^2 \partial_x \wedge \partial_y$. This implies that $\det \rho_{\mc{A}_{|D|}}\colon r \partial_r \wedge \partial_{\theta} \wedge \partial_{x_i} \mapsto r^2 (\partial_x \wedge \partial_y \wedge \partial_{x_i})$, from which the statement follows.
\ep
Note this statement is consistent with \autoref{prop:ellipticdivstiefelwhitney} and \autoref{cor:codimtwodegensw}. We next discuss the relation between the elliptic tangent bundle and the complex log tangent bundle. Let $|D| = (R,q)$ be an elliptic divisor. If $D$ is cooriented, using \autoref{rem:ellipticomplexlog} there is up to diffeomorphism a unique complex log divisor structure on $D$ such that $|D|$ is extracted via realification. For the associated Lie algebroids this results in the following (see the discussion below \autoref{prop:complexlogtransverse}).
\begin{prop}\label{prop:ellfiberprod} Let $|D| = (R,q)$ be an elliptic divisor such that $D$ is cooriented, and let $D = (U,\sigma)$ be a complex log divisor such that $(R,q) = ((U\otimes \overline{U})_\R, \sigma \otimes \overline{\sigma})$. Then $\mc{A}_{|D|} \otimes \C \cong \mc{A}_D \times_{TX} \mc{A}_{\overline{D}}$.
\end{prop}
In other words, not only do all cooriented elliptic divisors come from complex log divisors, the same holds for their associated ideal Lie algebroids. Using the discussion below \autoref{def:fiberproduct}, we thus obtain the following.
\begin{cor} Under the previous hypotheses, there is a short exact sequence
\be
		TX(-\log |D|) \otimes \C \stackrel{\iota}{\hookrightarrow} TX(-\log D) \to TX_\C \to 0.
\ee
\end{cor}
In other words, $TX(-\log |D|) \otimes \C$ is a complex $\mc{A}_D$-Lie algebroid.
\subsection{Edge tangent bundle}
\label{sec:edgetangentbundle}
In this section we describe another type of Lie algebroid that can be obtained using the rescaling procedure of Section \ref{sec:rescaling}. This Lie algebroid has appeared under various guises, either due to its relevance as the Poisson algebroid of a log-Poisson structure (see Section \ref{sec:logpoisson}, or e.g.\ \cite{GualtieriLi14}), or in the study of analysis on singular spaces \cite{Mazzeo91}.

Let $(X,Z)$ be a log pair and let $F \subseteq TZ$ be a corank-one involutive distribution defining a foliation of $Z$. Consider the sheaf $\mc{V}_X(I_Z,F) = \{v \in \mc{V}_X \, v|_Z \in F\}$. As $F$ is involutive, we see that this sheaf falls in the general framework of rescaling discussed in Section \ref{sec:rescaling}. Consequently, by \autoref{prop:rescaling} there exists a Lie algebroid $TX(-\log(Z,F)) \to X$ for which $\Gamma(TX(-\log (Z,F))) = \mc{V}_X(I_Z,\mc{F})$.
\begin{rem} By \autoref{rem:bbprimerescaling} it is clear that $\mc{A}_{Z,F} := TX(-\log(Z,F))$ is an $\mc{A}_Z$-Lie algebroid, and that $\mc{A}_{Z,F}$ has isomorphism locus equal to $X \setminus Z$.
\end{rem}
In terms of a coordinate system $(z,y,x_3,\dots,x_n)$ where $z \in I_Z$ is a local defining function for $Z$ and $y$ is a coordinate adapted to the foliation, we have
\be
	\Gamma(TX(-\log (Z,F))) = \mc{V}_X(I_Z,F) = \langle z \partial_z, z \partial_y, \partial_{x_i} \rangle.
\ee
In a different vein, instead of the existence of $F$ we can assume that there is a fibration structure $p\colon Z \to B$ onto some base $B$. Consider the sheaf $\mc{V}_X(\ker p) = \{v \in \mc{V}_X \, | \, v|_Z \in \ker p\}$ of vector fields that are tangent to the fibers of $p$. It is immediate that $\ker p \subseteq TZ$ is an involutive subbundle, so that it is a Lie subalgebroid of $TX$ supported on $Z$. Again using the method of rescaling, we conclude there exists a Lie algebroid ${}^{\ker p} TX \to X$ with $\Gamma({}^{\ker p} TX) = \mc{V}_X(\ker p)$.
\begin{rem} This Lie algebroid is usually called the \emph{edge tangent bundle} when $Z = \partial X$, and is then denoted by ${}^e TX$, where the $e$ stands for `edge'.
\end{rem}
Another way to obtain such a Lie algebroid is to consider the situation where there is a nowhere vanishing closed one-form $\alpha \in \Omega^1_{\rm cl}(Z)$. We define a sheaf $\mc{V}_X(\ker \alpha) = \{v \in \mc{V}_X \, | \, v|_Z \in \ker \alpha\}$. We allow ourselves a slightly more general situation, where we use Lie algebroid one-forms instead.
\begin{prop}\label{prop:closed1formsubalgd} Let $\mc{A} \to X$ be a Lie algebroid and $\alpha \in \Omega^1_{\rm cl}(\mc{A})$ be nowhere zero. Then $\ker \alpha \subseteq \mc{A}$ is a Lie subalgebroid.
\end{prop}
\bp Certainly $\ker \alpha$ forms a vector subbundle, and we check it is closed under $[\cdot,\cdot]_\mc{A}$. Let $v,w \in \ker \alpha$. Then we have $d_\mc{A} \alpha(v,w) = \mc{L}_v \alpha(w) - \mc{L}_w \alpha(v) - \alpha([v,w]_\mc{A})$. As $d_\mc{A} \alpha = 0$, we see that $\alpha([v,w]_\mc{A}) = 0$. But then $[v,w]_\mc{A} \in \ker \alpha$.
\ep
Going back to our previous discussion, by an application of \autoref{prop:closed1formsubalgd} we conlude that $\ker \alpha \subseteq TZ$ is a Lie subalgebroid. As $TZ \subseteq TX|_Z$ is in turn a Lie subalgebroid of $TX$ supported on $Z$, we conclude that the same holds for $\ker \alpha$. Consequently, by \autoref{prop:rescaling} there exists a Lie algebroid $\mc{A}_{\ker \alpha} := {}^{\ker\alpha} TX \to X$ for which $\Gamma(\mc{A}_{\ker \alpha}) = \mc{V}_X(\ker \alpha)$. As this is a rescaling of $TX$ at $Z$, it immediately follows that $\mc{A}_{\ker \alpha}$ is a $\mc{A}_Z$-Lie algebroid.

This type of Lie algebroid is present in the work of Lanius \cite[Section 3.2]{Lanius16two}, who computes its Lie algebroid cohomology using the realization that it is an $\mc{A}_Z$-Lie algebroid. The quotient complex is shown to have no cohomology, so that the $\mc{A}_Z$-anchor is a quasi-isomorphism. This is also noted in \cite{GualtieriLiPelayoRatiu17}.
\begin{prop}\label{prop:edgelacohomology} Let $(X,Z)$ be a log pair with corank-one involutive distribution $F \subseteq TZ$. Then $H^k(\mc{A}_{Z,F}) \cong H^k(\mc{A}_Z) \cong H^k(X) \oplus H^{k-1}(Z)$.
\end{prop}
\subsection{Zero tangent bundle}
\label{sec:zerotangentbundle}
We next discuss the zero tangent bundle asscociated to a log pair $(X,Z)$. This is an extreme case of the edge tangent bundle, where $p\colon Z \to B$ is the trivial fibration $p = {\rm id}_Z$, and $\ker p$ is trivial. In other words, its sections are all vector fields on $X$ which vanish at $Z$. This results in the following Lie algebroid (see \autoref{defn:secondaryidealla}).
\begin{defn} Let $(X,Z)$ be a log pair. The \emph{zero tangent bundle} is the secondary ideal Lie algebroid $\mc{B}_Z$ associated to $L_Z$. In other words, $\Gamma(\mc{B}_Z) = I_Z \cdot \Gamma(TX)$.
\end{defn}
The isomorphism locus of $\mc{B}_Z$ is again $X \setminus Z$. The usual notation for the zero tangent bundle, especially when $Z = \partial X$, is ${}^0 TX$. It is immediate that $\mc{B}_Z$ is the $(0,Z)$-rescaling of $TX$. In terms of local generators, for adapted coordinates $(z,x_1,\dots,x_n)$ around $Z$ we have
\be
	\Gamma(\mc{B}_Z) = \langle z \partial_z, z \partial_{x_i} \rangle.
\ee
\begin{rem} Using either \autoref{rem:bbprimerescaling}, \autoref{rem:isoloci}, or their local descriptions, we see that the zero tangent bundle $\mc{B}_Z$ is an $\mc{A}_Z$-Lie algebroid.
\end{rem}

Note that we determined the isomorphism type of $\mc{B}_Z$ with \autoref{prop:seclaisotype}.
\begin{prop}\label{prop:zeroisotangent} Let $(X,Z)$ be a log pair. Then $\mc{B}_Z \cong TX \otimes L_Z^*$.
\end{prop}
We will use this in Section \ref{sec:compcharclasses} to compute characteristic classes of $\mc{B}_Z$. While for the log-tangent bundle $\mc{A}_Z$ its degeneracy locus $Z$ was a transitive invariant submanifold, for $\mc{B}_Z$ the degeneracy locus is as far away from transitive as possible. Indeed, the restriction is totally intransitive, with $\rho_{\mc{B}_Z}|_Z$ being the trivial map.
\subsection{Scattering tangent bundle}
\label{sec:scatteringtangentbundle}
Our next Lie algebroid is the scattering tangent bundle, again associated to a log pair $(X,Z)$. It is obtained from the log-tangent bundle $\mc{A}_Z$ by performing $0$-rescaling.
\begin{defn} Let $(X,Z)$ be a log pair. The \emph{scattering tangent bundle} is given by $\mc{C}_Z := {}^0 \mc{A}_Z$, the $0$-rescaling of $\mc{A}_Z$.
\end{defn}
As $\mc{C}_Z$ is a rescaling of $\mc{A}_Z$, it immediately follows that $\mc{C}_Z$ is a $\mc{A}_Z$-Lie algebroid. Moreover, the isomorphism locus of $\mc{C}_Z$ is given by $X \setminus Z$.
\begin{rem} The scattering tangent bundle is usually denoted by ${}^{sc} TX$ \cite{Melrose95}, which is more consistent with the notations ${}^b TX$, ${}^e TX$ and ${}^0 TX$. However, we wish to emphasize the hypersurface in the notation.
\end{rem}
In terms of the usual coordinate system $(z,x_i)$ with $\{z = 0\} = Z$, the scattering tangent bundle is locally generated by $\Gamma(\mc{C}_Z) = \langle z^2 \partial_z, z \partial_{x_i} \rangle$. Similarly to \autoref{prop:zeroisotangent}, using \autoref{prop:seclaisotype} we determine the isomorphism type of $\mc{C}_Z$.
\begin{prop}\label{prop:scatteringisolog} Let $(X,Z)$ be a log pair. Then $\mc{C}_Z \cong \mc{A}_Z \otimes L_Z^*$.
\end{prop}
We will compute characteristic classes of the scattering tangent bundle in Section \ref{sec:compcharclasses}. See also \cite{Lanius16} and Section \ref{sec:scatteringsymp} where $\mc{C}_Z$-symplectic structures are studied. We finish by mentioning the Lie algebroid cohomology of the scattering tangent bundle, computed in \cite{Lanius16} using the strategy mentioned at the end of Section \ref{sec:aliealgebroids}.
\begin{thm}[{\cite[Theorem 2.15]{Lanius16}}] Let $(X,Z)$ be a log pair. Then we have $H^k(\mc{C}_Z) \cong H^k(\mc{A}_Z) \oplus \Omega^{k-1}(Z;E_{-k}) \cong H^k(X) \oplus H^{k-1}(Z) \oplus \Omega^{k-1}(Z;E_{-k})$, where $E_s = |N^*Z|^{s} \to Z$ is the bundle of $s$-densities on $N^*Z$.
\end{thm}
In particular, note that unlike for the log-tangent bundle $\mc{A}_Z$ (see \autoref{thm:loglacohomology}), the Lie algebroid cohomology of $\mc{C}_Z$ is very far from being finite-dimensional.
\subsection{$b^k$-tangent bundles}
\label{sec:bktangentbundle}
In this section we discuss in our language a variant of the log-tangent bundle of Section \ref{sec:logtangentbundle} defined and studied by Scott \cite{Scott16}. All results in this section are his.

While the log-tangent bundle $TX(-\log Z)$ has a single ``singular'' generator with linear vanishing along $Z$, these Lie algebroids are meant to capture $k$th order tangency along $Z$ for any $k \geq 1$. In order for this to make invariant sense, more data is required.

Let $(X,Z)$ be a log pair, with ideal sheaf $I_Z$ and $j_Z\colon Z \hookrightarrow X$ the inclusion.
%
\begin{defn} Let $k \geq 1$. The \emph{sheaf of $k$-jets} at $Z$ is $\mc{J}_Z^k := j_Z^{-1}(C^\infty(X) / I_Z^{k+1})$. A \emph{$k$-jet} at $Z$ is a global section of $\mc{J}_Z^k$. Denote the set of $k$-jets at $Z$ by $J_Z^k$.
\end{defn}
Given a function $f$ defined in a neighbourhood of $Z$, let $[f]^k_Z \in J_Z^k$ denote the $k$-jet at $Z$ represented by $f$. For a given $k$-jet at $Z$, $j \in J_Z^k$, write $f \in j$ if $[f]_Z^k = j$. Let $j_{k-1} \in J_Z^{k-1}$ be a choice of $(k-1)$-jet and define the sheaf of vector fields
\be
	\mc{V}_X(I^k_Z, j_{k-1}) = \{v \in \mc{V}_X \, | \, \mc{L}_v f \in I_Z^k \text{ for all } f \in j_{k-1}\}.
\ee
Scott shows this is a locally free sheaf, and that given any $v \in \mc{V}_X(I_Z)$ and $f \in C^\infty(X)$, the jet $[\mc{L}_v f]_Z^{k-1}$ depends only on $[f]_Z^{k-1}$, so that its definition makes sense.
\begin{defn}[{\cite[Definition 2.11]{Scott16}}] Let $(X,Z)$ be a log pair and $j_{k-1} \in J_Z^{k-1}$. The \emph{$b^k$-tangent bundle} is the Lie algebroid $\mc{A}_Z^k \to X$ with $\Gamma(\mc{A}_Z^k) = \mc{V}_X(I^k_Z, j_{k-1})$.
\end{defn}
The notation that is used by Scott is that $\mc{A}_Z^k = {}^{b^k} TX$. We suppress the choice of $(k-1)$-jet $j_{k-1}$ in our notation. Note that for $k = 1$ we obtain that $\mc{A}_Z^k = {}^{b^k} TX \cong {}^b TX = \mc{A}_Z$, as any local generator $f$ of $I_Z$ represents the trivial $0$-jet. Using a local coordinate $z \in j_{k-1}$ for $Z$ and subsequent coordinate system $(z,x_2,\dots,x_n)$, we have that $\Gamma(\mc{A}_Z^k) = \langle z^k \partial_z, \partial_{x_i} \rangle$. This shows these Lie algebroids indeed capture the notion of $k$th order tangency along $Z$.
\begin{rem} Defining a sheaf by demanding its sections $v \in \mc{V}_X$ satisfy the condition $\mc{L}_v I_Z \subseteq I_Z^k$ is too naive, as Scott explains (e.g.\ \cite[Example 2.7]{Scott16}). The local description of $\mc{A}^k_Z$ is as one would expect, but one has to use local defining functions $z$ taken from the chosen $(k-1)$-jet.
\end{rem}
\begin{rem} Given $[z] = j_{k-1} \in J_Z^{k-1}$, the sheaf $\mc{V}_X(I^k_Z,j_{k-1})$ is unchanged if one modifies $z$ to $f \circ z$ for some $f\colon \R \to \R$ for which $f(0) = 0$ and $f'(0) > 0$. In other words, for $j_{k-1}' := [f \circ z]$ we have $\mc{V}_X(I^k_Z,j_{k-1}) = \mc{V}_X(I^k_Z,j_{k-1}')$. This means the corresponding Lie algebroids $\mc{A}_Z^k$ will be identical.
\end{rem}
The natural maps $J_Z^{k-1} \to J_Z^{k-1-i}$ for $0 < i \leq k-1$ give inclusions of sheaves
\be
	\mc{V}_X(I^k_Z, j_{k-1}) \subseteq \mc{V}_X(I^{k-1}_Z, j_{k-1}) \subseteq \dots \subseteq \mc{V}_X(I_Z).
\ee
Consequently, given a $(k-1)$-jet, the associated $\mc{A}_Z^{k-i}$ are compatibly $\mc{A}_Z^{k-j}$-Lie algebroids with maps $\varphi_{k-i,k-j}$ for any $0 \leq i \leq j \leq k$ in the sense of \autoref{defn:aliealgebroid}.
\begin{rem} It is an interesting question whether a similar definition can be made for other types of divisors, providing more examples of ideal Lie algebroids.
\end{rem}
\begin{que} Given a log divisor $(X,Z)$, can one obtain the bundles $\mc{A}_Z^k$ out of $TX$ using an iterated rescaling procedure using a jet filtration as in \cite[Chapter 8]{Melrose93}?
\end{que}
As for the log-tangent bundle, the isomorphism locus of $\mc{A}_Z^k$ is $X \setminus Z$. Moreover, its restriction to $Z$ is transitive, and the isotropy $\mathbb{L}_{Z,k} := \ker \rho_{\mc{A}_Z^k}|_Z$ is again canonically trivial, locally generated as $\langle z^k \partial_z \rangle$ for $z \in j_{k-1}$. There is again a residue map ${\rm Res}_Z^k\colon \Omega^\bullet(\mc{A}_Z^k) \to \Omega^{\bullet-1}(Z)$ which can be described as contraction with the canonical section of $\mathbb{L}_{Z,k}$. Moreover, similar arguments as for the log-tangent bundle allow one to compute the Lie algebroid cohomology of $b^k$-tangent bundles.
\begin{thm}[{\cite[Proposition 4.3]{Scott16}}] Let $(X,Z)$ be a log pair and $j_{k-1} \in J_Z^{k-1}$. Then $H^i(\mc{A}_Z^k) \cong H^i(X) \oplus \left(H^{i-1}(Z)\right)^k$.
\end{thm}
The proof of the above theorem uses the exactness of the sequence of complexes
\be
	0 \to \Omega^\bullet(\mc{A}_Z^{k-1}) \stackrel{\varphi_{k,k-1}^*}{\to} \Omega^\bullet(\mc{A}_Z^k) \stackrel{{\rm Res}_Z^k}{\to} \Omega^{\bullet-1}(Z) \to 0,
\ee
the fact that this sequence splits, and induction on $k$. Finally, in anticipation of the next section we remark that Scott defines a \emph{$b^k$-map} to be a $b$-map $f\colon (X,Z_X) \to (Y,Z_Y)$ for which $f^*(j_{Z_Y,k-1}) = j_{Z_X,k-1}$, and invites the reader to verify that these induce Lie algebroid morphisms $(\varphi,f)\colon \mc{A}_{Z_X}^k \to \mc{A}_{Z_Y}^k$ for which $\varphi = Tf$ on sections. We hope the similarity with \autoref{prop:bmapideal} and \autoref{prop:logdivmorph} is clear.
\section{Morphisms between ideal Lie algebroids}
\label{sec:morphilas}
We close this chapter by discussing the relation between morphisms of divisors and between the associated ideal Lie algebroids they give rise to. We show that in certain cases, given a morphism of divisors $f\colon (X,I_X) \to (Y,I_Y)$ (written in terms of the divisor ideals), there is an induced Lie algebroid morphism $(\varphi,f)\colon \mc{A}_{I_X} \to \mc{A}_{I_Y}$ between their ideal Lie algebroids. However, there is no converse statement. 

Let $f\colon X \to Y$ with $f^* I_Y = I_X$ be a morphism of divisors. Then by \autoref{prop:divmorphstrong} we know that $f$ is a strong map of pairs between zero sets. By taking complements, $f^{-1}(Z_{I_Y}) = Z_{I_X}$ implies that $f^{-1}(Y \setminus Z_{I_Y}) = X \setminus Z_{I_X}$. In other words, if $\mc{A}_{I_X}$ and $\mc{A}_{I_Y}$ exist, $f\colon (X, X_{\mc{A}_{I_X}}) \to (Y, Y_{\mc{A}_{I_Y}})$ is a strong map of pairs. Consider now \autoref{prop:denseisolocusstrongmap} and \autoref{prop:denselam}.

We cannot in general decide the following, posed as a question to the reader.
\begin{que}\label{que:ilamorph} Let $f$ be a morphism of divisors giving rise to ideal Lie algebroids. When does $Tf$ induce a Lie algebroid morphism $\varphi$ for which $\varphi = Tf$ on sections?
\end{que}

\subsection{Between log divisors}
\label{sec:morphlogdivisors}
In this section we discuss the case of morphisms between log divisors. In this case, \autoref{que:ilamorph} has an affirmative answer. Consider log pairs $(X,Z_X)$ and $(Y,Z_Y)$. Recall from Section \ref{sec:logdivisor} a morphism of log divisors is the same thing as a $b$-map, i.e.\ a transverse strong map of pairs $f\colon (X,Z_X) \to (Y,Z_Y)$.
\begin{prop}\label{prop:logdivmorph} Let $f\colon (X,Z_X) \to (Y,Z_Y)$ be a morphism of divisors. Then $Tf$ induces a Lie algebroid morphism $(\varphi,f)\colon TX(-\log Z_X) \to TY(-\log Z_Y)$ such that $\varphi = Tf$ on sections.
\end{prop}
\bp It is immediate that $f^{-1}(Z_Y) = Z_X$. Note that the isomorphism loci of log-tangent bundles are dense. By \autoref{prop:denselam} it thus suffices to show that $Tf$ induces a vector bundle morphism, which in turn is equivalent to showing that $f^*$ extends to an algebra morphism from $\Omega^\bullet(Y;\log Z_Y)$ to $\Omega^\bullet(X; \log Z_X)$. Let $x \in Z_X$ and denote $y = f(x) \in Z_Y$. Consider suitable tubular neighbourhood coordinates $(z,x_2,\dots, x_n)$ in a neighbourhood $U$ of $x$ such that $U \cap Z_X = \{z = 0\}$ with $I_{Z_X} = \langle z \rangle$, and $(z', y_2, \dots, y_m)$ in a neighbourhood $V$ of $y$ such that $V \cap Z_Y = \{z' = 0\}$ and $I_{Z_Y} = \langle z' \rangle$. In these coordinates we have $\Omega^1(U; - \log Z_X) = \langle d \log z, dx_2, \dots, d x_n\rangle$ and $\Omega^1(V; -\log Z_Y) = \langle d\log z', dy_2, \dots, dy_m\rangle$.

The Lie algebroid one-forms $dy_i$ can be pulled back using $f^*$ as these are smooth. These smooth one-form inject into $\Omega^1(U;- \log Z_X)$ using the anchor. We are left with checking that $d \log z'$ is pulled back to a Lie algebroid one-form. As $f$ is a morphism of divisors we have $f^* I_{Z_Y} = I_{Z_X}$, so that $f^*(z') = e^h z$ for $h$ a smooth function on $U$. Consequently $f^* d \log z' = d \log f^*(z') = d\log(e^h z) = dh + d\log z \in \Omega^1(U;\log Z_X)$. We conclude that $Tf$ induces a Lie algebroid morphism as desired.
\ep
We finish by discussing the relation between the above Lie algebroid morphisms and the residue maps ${\rm Res}_Z$ that were introduced in Section \ref{sec:logtangentbundle} (see also Section \ref{sec:residuemaps}). A simple consequence of \autoref{cor:residuemapsdenseiso} is the following.
\begin{prop} Let $(\varphi,f)\colon TX(-\log Z_X) \to TY(-\log Z_Y)$ be a Lie algebroid morphism. Then ${\rm Res}_{Z_X} \circ \varphi^* = f^* \circ {\rm Res}_{Z_Y}$.
\end{prop}
\bp A Lie algebroid morphism between these two Lie algebroids necessarily satifies $f^{-1}(Z_Y) = Z_X$. Moreover, their degeneracy loci are indeed transitive as mentioned before, and their codimensions agree.
\ep
\subsection{From elliptic to log divisors}
We next discuss the case of an elliptic divisor mapping to a log divisor. We will see that again there is a positive answer to \autoref{que:ilamorph}. Let $(X, |D|)$ be an elliptic pair and $(Y,Z)$ a log pair, coming with divisor ideals $I_{|D|} \subset C^\infty(X)$ and $I_Z \subset C^\infty(Y)$. A morphism of the corresponding divisors gives rise to a Lie algebroid morphism.
\begin{prop}\label{prop:elltologmorphism} Let $f\colon (X,|D|) \to (Y,Z)$ be a morphism of divisors. Then $Tf$ induces a Lie algebroid morphism $(\varphi,f)\colon TX(-\log|D|) \to TY(-\log Z)$ such that $\varphi = Tf$ on sections.
\end{prop}
\bp This is almost identical to \autoref{prop:logdivmorph}. It is immediate that $f^{-1}(Z) = D$. Note that the isomorphism loci of the elliptic tangent bundle and the log-tangent bundle are dense. By \autoref{prop:denselam} it thus suffices to show that $Tf$ induces a vector bundle morphism, which in turn is equivalent to showing that $f^*$ extends to an algebra morphism from $\Omega^\bullet(Y;\log Z)$ to $\Omega^\bullet(X; \log |D|)$. Let $x \in D$ and denote $y = f(x) \in Z$. Consider suitable tubular neighbourhood coordinates $(r, \theta, x_3,\dots, x_n)$ in a neighbourhood $U$ of $x$ such that $U \cap D = \{r = 0\}$ with $I_{|D|} = \langle r^2 \rangle$, and $(z, y_2, \dots, y_m)$ in a neighbourhood $V$ of $y$ such that $V \cap Z = \{z = 0\}$ and $I_Z = \langle z \rangle$. In these coordinates we have $\Omega^1(U; - \log |D|) = \langle d \log r, d \theta, dx_3, \dots, d x_n\rangle$ and $\Omega^1(V; -\log Z) = \langle d\log z, dy_2, \dots, dy_m\rangle$.

The Lie algebroid one-forms $dy_i$ can be pulled back using $f^*$ as these are smooth. Moreover, the smooth one-form inject into $\Omega^1(U;- \log D)$ using the anchor. We are left with checking that $d \log z$ is pulled back to a Lie algebroid one-form. As $f$ is a morphism of divisors we have $f^* I_Z = I_{|D|}$, so that $f^*(z) = e^h r^2$ for $h$ a smooth function on $U$. Consequently $f^* d \log z = d \log f^*(z) = d \log (e^h r^2) = d h + 2 d\log r \in \Omega^1(U;-\log |D|)$. We conclude that $Tf$ induces a Lie algebroid morphism as desired.
\ep
\begin{rem} In light of the proofs of \autoref{prop:logdivmorph} and \autoref{prop:elltologmorphism}, we see it is relatively straightforward to obtain a morphism to a log-tangent bundle. This is because the log-tangent bundle has only one ``singular'' generator, i.e.\ the kernel of the restriction of the anchor to its degeneracy locus is one-dimensional.
\end{rem}

A Lie algebroid morphism from the elliptic tangent bundle $TX(-\log |D|)$ to the log tangent bundle $TY(-\log Z)$ intertwines the residue maps that were discussed in Section \ref{sec:logtangentbundle} and Section \ref{sec:elltangentbundle}.
\begin{prop}\label{prop:residuemaps} Let $(\varphi,f)\colon TX(-\log |D|) \to TY(-\log Z)$ be a Lie algebroid morphism. Then ${\rm Res}_q \circ \varphi^* = 0$. Moreover, $({\rm Res}_r+{\rm Res}_\theta) \circ \varphi^* = f^* \circ {\rm Res}_Z$.
\end{prop}
\bp By assumption we have $f^{-1}(Z) = D$ so that $df\colon TD \to TZ$. Restricting $TX(-\log |D|)$ to $D$ and $TY(-\log Z)$ to $Z$ gives the following commutative diagram.
\begin{center}
	\begin{tikzpicture}
	
	\matrix (m) [matrix of math nodes, row sep=2.5em, column sep=2.5em,text height=1.5ex, text depth=0.25ex]
	{	\wt{\mc{S}}\colon 0 & \underline{\R} \oplus \mf{k} & TX(-\log |D|)|_D & TD & 0 \\ \mc{S}\colon 0 & \mathbb{L}_Z & TY(-\log Z)|_Z & TZ & 0 \\};
	\path[-stealth]
	(m-1-1) edge (m-1-2)
	(m-1-2) edge (m-1-3)
	(m-1-4) edge (m-1-5)
	(m-2-1) edge (m-2-2)
	(m-2-2) edge (m-2-3)
	(m-2-4) edge (m-2-5)
	(m-1-3) edge node [left] {$\varphi$} (m-2-3)
	(m-1-2) edge node [left] {$\varphi$} (m-2-2)
	(m-1-4) edge node [right] {$df$} (m-2-4)
	(m-1-3) edge node [above] {$\rho_X$} (m-1-4)
	(m-2-3) edge node [above] {$\rho_Y$} (m-2-4);	
	\end{tikzpicture}
\end{center}
Consequently, we obtain a map $\varphi^*\colon \mc{S}^* \to \wt{\mc{S}}^*$ between dual sequences, and also between spaces of sections. Using the notation preceding this proof we have $E = \mathbb{L}_Z$ so that $\ell = \dim(E) = 1$, and $\wt{E} = \underline{\R} \oplus \mf{k}$ so that $\wt{\ell} = \dim(\wt{E}) = 2$. Recall that $\mathbb{L}_Z$ carries a canonical trivialization. Given a form $\alpha \in \Omega^k(\log Z)$, we can identify ${\rm Res}(\alpha) \in \Gamma(Z;\wedge^{k-1} T^*Z \otimes \mathbb{L}_Z^*)$ with ${\rm Res}_Z(\alpha) \in \Omega^{k-1}(Z)$. Similarly, a choice of coorientation for $ND$ trivializes $\mf{k} = \wedge^2 N^*D \otimes R$. Given $\beta \in \Omega^k(\log |D|)$ this identifies $\wt{\rm Res}(\beta) \in \Gamma(D; \wedge^{k-2} T^*D \otimes \mf{k}^*)$ with ${\rm Res}_q(\beta) \in \Omega^{k-2}(D)$, using that $\wedge^2(\underline{\mathbb{R}} \oplus \mf{k}) \cong \mf{k}$. Moreover, for $\beta$ with ${\rm Res}_q(\beta) = 0$, the radial residue ${\rm Res}_r(\beta) \in \Omega^{k-1}(D)$ together with the $\theta$-residue ${\rm Res}_\theta(\beta) \in \Omega^{k-1}(D)$ is identified with $\wt{\rm Res}_{-1}(\beta) \in \Omega^{k-1}(D;\underline{\R} \oplus \mf{t})$. As $2 = \wt{\ell} \geq \ell = 1$ with $\ell - \wt{\ell} = -1$ we obtain immediately from \autoref{lem:residues} and \autoref{lem:residuecommute} that ${\rm Res}_q(\varphi^* \alpha) = 0$ and that $({\rm Res}_r+{\rm Res}_\theta)(\varphi^* \alpha) = f^*({\rm Res}_Z(\alpha))$.
\ep
\chapter{\texorpdfstring{$\mc{A}$}{A}-symplectic structures}
\label{chap:asymplecticstructures}
\fancypagestyle{empty}{%
	\fancyhf{}%
	\renewcommand\headrulewidth{0pt}%
	\fancyhead[RO,LE]{\thepage}%
}
In this chapter we discuss Lie algebroid symplectic structures. We assume the reader is familiar with the notion of a \symp, which in itself will not play a huge role in this thesis. We focus our attention on their analogue defined for Lie algebroids. These have garnered recent interest especially for the types of Lie algebroids that we are emphasizing in this thesis. The reason for this is that there is often a powerful link to Poisson geometry, as an $\mc{A}$-symplectic structure is dual to that of a nondegenerate $\mc{A}$-Poisson structure, which we will discuss in Chapter \ref{chap:apoissonstructures}.

This chapter should serve to convince the reader that various basic techniques from symplectic geometry carry over to $\mc{A}$-symplectic structures. In particular, the Moser argument carries over to Lie algebroids with dense isomorphism locus. Note that most results in this chapter are for arbitrary Lie algebroids $\mc{A}$, while for some we need to assume a closer link to the underlying manifold. We often take the assumption that $\mc{A}$ has dense isomorphism locus.
\subsection*{Organization of the chapter}
This chapter is built up as follows. In Section \ref{sec:definitionsasymp} we discuss the definition of an $\mc{A}$-symplectic structure and several simple consequences. In Section \ref{sec:anambu} we discuss $\mc{A}$-Nambu structures (of top degree), which to some extent can be treated on equal footing. This is made more precise in Section \ref{sec:amosertechniques}, where we develop Moser-type techniques simultaneously for both $\mc{A}$-symplectic and $\mc{A}$-Nambu structures. Then in Section \ref{sec:examplessymp} we discuss examples of $\mc{A}$-symplectic structures for some of the concrete Lie algebroids discussed in Section \ref{sec:laexamples}, among which log-symplectic and elliptic symplectic structures are the most important ones for us. Finally, Section \ref{sec:blogsandfsymps} discusses the relation between log-symplectic and \fsymp{}s.
\section{Definitions}
\label{sec:definitionsasymp}
In this section we start by defining the basic conceps relevant to the study of $\mc{A}$-symplectic geometry. Let $\mc{A} \to X$ be a Lie algebroid.
\begin{defn} An \emph{$\mc{A}$-symplectic structure}\index{$\mc{A}$-symplectic structure} is a closed and nondegenerate $\mc{A}$-two-form ${\omega_{\mc{A}} \in \Omega^2_{\mc{A}}(X)}$. Denote the space of $\mc{A}$-symplectic forms by ${\rm Symp}(\mc{A})$.
\end{defn}
In other words, $\omega_\mc{A}$ satisfies $d_\mc{A} \omega_\mc{A} = 0$ and $\omega_\mc{A}^n \neq 0$, where ${\rm rank}(\mc{A}) = 2n$. Any $\mc{A}$-two-form $\omega_\mc{A}$ defines a map $\omega_\mc{A}^\flat\colon \mc{A}^* \to \mc{A}$ given by $v \mapsto \iota_v \omega_\mc{A}$ for $v \in \mc{A}$, and nondegeneracy of $\omega_\mc{A}$ is equivalent to $\omega_\mc{A}^\flat$ being an isomorphism. Due to the nondegeneracy condition, $\mc{A}$-symplectic structures can only exist if the rank of $\mc{A}$ is even (but note that $\dim X$ need not necessarily be even). Note that any $\mc{A}$-symplectic structure $\omega_\mc{A}$ defines an $\mc{A}$-cohomology class $[\omega_\mc{A}] \in H^2(\mc{A})$.
\begin{prop} Let $\omega_\mc{A}$ be an $\mc{A}$-symplectic structure. Then $\mc{A}$ is orientable.
\end{prop}
\bp This follows as $\omega_\mc{A}^n \in \Gamma(\det(\mc{A}^*))$ is nonvanishing, with ${\rm rank}(\mc{A}) = 2n$.
\ep
Consequently, any Lie algebroid $\mc{A}$ for which there exists an $\mc{A}$-symplectic structure is in particular orientable as a vector bundle, i.e.\ $w_1(\mc{A}) = 0$. The consequences of this for obstructing the existence of $\mc{A}$-\symp{}s for specific Lie algebroids $\mc{A}$ are explored in Chapter \ref{chap:homotopicalobstrs}.
\begin{rem}\label{rem:asympranktwo} If the rank of a Lie algebroid $\mc{A} \to X$ is two, the manifold $X$ admits an $\mc{A}$-\symp if and only if $\mc{A}$ is orientable, i.e.\ $w_1(\mc{A}) = 0$. In particular, this settles when ideal Lie algebroids admit $\mc{A}$-symplectic structures on surfaces.
\end{rem}
There is a standard notion of morphism between $\mc{A}$-symplectic manifolds.
\begin{defn} Let $(X, \mc{A}, \omega_\mc{A})$ and $(X', \mc{A}', \omega_{\mc{A}'})$ be Lie algebroids with symplectic structures. An \emph{$\mc{A}$-symplectomorphism} from $\mc{A}$ to $\mc{A}'$ is a Lie algebroid morphism $(\varphi,f)\colon \mc{A} \to \mc{A}'$ such that $\varphi^* \omega_{\mc{A}'} = \omega_\mc{A}$.
\end{defn}
\begin{rem} Note that $\varphi^* \omega_{\mc{A}'} = \omega_\mc{A}$ is equivalent to $\omega_\mc{A}^\flat = \varphi \circ \omega_\mc{B}^\flat \circ \varphi^*$ as maps. As both $\omega_\mc{A}$ and $\omega_{\mc{A}'}$ are nondegenerate, any $\mc{A}$-symplectomorphism $(\varphi,f)$ must in particular be an isomorphism between $\mc{A}$ and $\mc{A}'$. The base map $f\colon X \to X'$ is not necessarily a diffeomorphism.
\end{rem}
We next define the accompanying notion of an $\mc{A}$-almost-complex structure.
\begin{defn} Let $\mc{A} \to X$ be a Lie algebroid. An \emph{$\mc{A}$-\acs} is a vector bundle complex structure $J$ for $\mc{A}$, i.e.\ $J\colon \mc{A} \to \mc{A}$ with $J^2 = -{\rm id}_\mc{A}$.
\end{defn}
The usual argument using an auxiliary bundle metric for $\mc{A}$ shows that any $\mc{A}$-\symp{} has an accompanying compatible $\mc{A}$-\acs (see e.g.\ \cite{McDuffSalamon98}). Note that their definition and existence is independent of the Lie algebroid structure on $\mc{A}$. We will use \acs{}s in Chapter \ref{chap:constrasymp} as a means of constructing $\mc{A}$-\symp{}s, and in Chapter \ref{chap:homotopicalobstrs} to obstruct their existence in dimension four using this simple consequence.
\begin{prop}\label{prop:asympaacs} Let $\mc{A} \to X$ be a Lie algebroid such that there exists an $\mc{A}$-\symp. Then $X$ admits an $\mc{A}$-\acs.
\end{prop}
Let $\mc{A} \to X$ be a Lie algebroid with $\mc{A}$-\symp $\omega_\mc{A}$, and let $D \subseteq X$ be a transitive $\mc{A}$-invariant submanifold. This means there is a residue map (see Section \ref{sec:residuemaps}), and we obtain a form ${\rm Res}_D(\omega_\mc{A}) \in \Omega^{2-\ell}(D; \wedge^\ell E^*)$, where $E = \ker \rho_{\mc{A}}|_D$. Such residues capture important information about $\omega_{\mc{A}}$ near $D$.
\begin{rem} If $\mc{A} \to X$ is a Lie algebroid whose isomorphism locus $X_\mc{A}$ is nonempty (noting it is always open), one can perform any operation available to symplectic structures to a given $\mc{A}$-\symp. For example, the symplectic fiber sum \cite{Gompf95} or symplectic blow-up \cite{McDuffSalamon98} procedures can be performed in $X_\mc{A}$, which then naturally result in new Lie algebroid symplectic structures.
\end{rem}
We finish by posing several questions to the reader.
\begin{que} When does the top power of an $\mc{A}$-symplectic structure define a nonzero class in Lie algebroid cohomology (possibly twisted by the $\mc{A}$-module $Q_\mc{A}$)?
\end{que}
\begin{que} Can the approximately-holomorphic methods of Donaldson \cite{Donaldson99} be adapted to $\mc{A}$-symplectic geometry? The existence of compact $\mc{A}$-symplectic submanifolds results in compact Poisson transversals \cite{FrejlichMarcut17} for their dual Poisson structures.
\end{que}
\begin{que} When can one define analogues of $J$-holomorphic curve invariants for $\mc{A}$-symplectic structures? This relies on normal form results around $Z_\mc{A}$.
\end{que}
\section{\texorpdfstring{$\mc{A}$}{A}-Nambu structures}
\label{sec:anambu}
In this section we briefly discuss the notion of an $\mc{A}$-Nambu structure. These are the Lie algebroid analogues of Nambu structures as studied in e.g.\ \cite{Nambu73, Takhtajan94, Dufour00, MartinezTorres04}. This is done in this chapter as they can to some extent be treated on similar footing with $\mc{A}$-symplectic structures. This will become apparent further in Section \ref{sec:amosertechniques} when we develop Moser techniques simultaneously for $\mc{A}$-symplectic structures and $\mc{A}$-Nambu structures. In Chapter \ref{chap:homotopicalobstrs} we will discuss criteria for their existence for certain Lie algebroids $\mc{A}$ we have encountered before. Let $\mc{A} \to X$ be a Lie algebroid.

\begin{defn} An \emph{$\mc{A}$-Nambu structure} is a nondegenerate section $\Pi$ of $\det(\mc{A})$.
\end{defn}
Assuming that ${\rm rank}(\mc{A}) = n$, the nondegeneracy condition is the fact that the map $\Pi^\sharp\colon \mc{A}^* \to \wedge^{n-1} \mc{A}$ given by $v \mapsto \iota_v \Pi$ for $v \in \mc{A}^*$ is nondegenerate. In other words, an $\mc{A}$-Nambu structure is just a nonzero section of $\det(\mc{A})$. This notion makes sense on any vector bundle, as the Lie algebroid structure of $\mc{A}$ is not used. This is because as an $\mc{A}$-Nambu structure is a section of the line bundle $\det(\mc{A})$, any natural integrability condition one would impose is automatically satisfied.
\begin{rem} To be more consistent with other authors we should call these $\mc{A}$-Nambu structure \emph{of top degree}. However, as we will only consider these types of $\mc{A}$-Nambu structures in this thesis, we have chosen to use this shorter name for brevity.
\end{rem}
An $\mc{A}$-Nambu structure is dual to a nondegenerate $\mc{A}$-$n$-form $\Pi^{-1} \in \Omega^n(\mc{A})$, hence defines an $\mc{A}$-cohomology class $[\Pi^{-1}] \in H^n(\mc{A})$. By definition, without using the Lie algebroid structure of $\mc{A}$, we have the following, as $w_1(\det(\mc{A})) = w_1(\mc{A})$.

\begin{prop}\label{prop:anambustructure} Let $\mc{A} \to X$ be a Lie algebroid. Then $\mc{A}$ admits an $\mc{A}$-Nambu structure if and only if $\mc{A}$ is orientable, i.e.\ if and only if $w_1(\mc{A}) = 0$.
\end{prop}
If $\mc{A}$ has rank equal to two, an $\mc{A}$-Nambu structure is a nondegenerate $\mc{A}$-Poisson structure (see Chapter \ref{chap:apoissonstructures}), as the Poisson condition is automatic. By \autoref{prop:asympapoisson} in the next chapter (or directly) we see that $\mc{A}$-Nambu structures on rank two Lie algebroids are in one-to-one correspondence with $\mc{A}$-\symp{}s. In Chapter \ref{chap:homotopicalobstrs} we will use the previous proposition to obstruct the existence of $\mc{A}$-Nambu structures.
\section{\texorpdfstring{$\mc{A}$}{A}-Moser techniques}
\label{sec:amosertechniques}
In this section we discuss how the Moser argument \cite{Moser65, Weinstein71} (see also \cite{McDuffSalamon98}) can be adapted to $\mc{A}$-forms for a Lie algebroid $\mc{A}$ with dense isomorphism locus. Moreover, we discuss a version of the Poincar\'e lemma for such Lie algebroids. Extensions of Moser techniques to Lie algebroid forms have been considered before, notably for the log-tangent bundles. Indeed, the results in this section encapsulate several results in the literature, including \cite{GuilleminMirandaPires14,Lanius16,MarcutOsornoTorres14,MirandaPlanas16,NestTsygan96,Radko02,Scott16}. The contents of this section are joint with Melinda Lanius and will appear in \cite{KlaasseLanius17}.

Let $\mc{A} \to X$ be a Lie algebroid and let $\tau \in \Omega^k(\mc{A})$ be a \emph{nondegenerate $\mc{A}$-k-form} (after \cite{MarcutOsornoTorres14}). This means that $\tau^\flat\colon \mc{A} \to \wedge^{k-1} \mc{A}^*$ given by $\iota_v \mapsto \iota_v \tau$ for $v \in \mc{A}$ is surjective. By comparing dimensions, we see that $k$ must equal $1$, $2$, or ${\rm rank}(\mc{A})$. Note that ${\rm rank}(\mc{A}) = \dim X$ if $\mc{A}$ has dense isomorphism locus $X_\mc{A} \subseteq X$. Moreover, recall that $Z_\mc{A} = X \setminus X_\mc{A}$ is the degeneracy locus of $\mc{A}$.

We first cover the following $\mc{A}$-analogue of the relative Poincar\'e lemma.
\begin{lem}\label{lem:arelpoincare} Let $\mc{A} \to X$ be a Lie algebroid with dense isomorphism locus for which $Z_\mc{A}$ is smooth. Let $\tau \in \Omega^k(\mc{A})$ be $d_\mc{A}$-closed such that $(\rho_{\mc{A}}^{-1})^*(\tau)$ extends smoothly by $0$ over $Z_\mc{A}$. Then there exists a neighbourhood $U$ of $Z_\mc{A}$ and $\sigma \in \Omega^{k-1}(U;\mc{A})$ vanishing on $Z_\mc{A}$ such that $\tau = d_\mc{A} \sigma$ on $U$.
\end{lem}
In other words, under certain conditions we can find $\mc{A}$-primitives around $Z_\mc{A}$.
\bp Denote by $\tau' \in \Omega^k(M)$ the smooth extension of $(\rho_{\mc{A}}^{-1})^*(\tau)$. Then $\tau'|_{Z_\mc{A}} = 0$, hence by the standard relative Poincar\'e lemma there exists a neighbourhood $U$ of $Z_\mc{A}$ which retracts to $Z_\mc{A}$ and a smooth form $\sigma' \in \Omega^{k-1}(U)$ vanishing on $Z_\mc{A}$ such that $\tau' = d \sigma'$. Define $\sigma := \rho_\mc{A}^* \sigma'$, which hence also vanishes on $Z_\mc{A}$. Then $d_\mc{A} \sigma = d_\mc{A}(\rho_\mc{A}^*\sigma') = \rho_\mc{A}^* (d \sigma') = \rho_\mc{A}^* \tau'$, and hence on $U \cap X_\mc{A}$ we conclude that $\tau = d_\mc{A} \sigma$. But then as $X_\mc{A}$ is dense we have $\tau = d_\mc{A} \sigma$ on all of $U$ by continuity.
\ep

The following is the relative $\mc{A}$-Moser theorem for Lie algebroids $\mc{A}$ with dense isomorphism loci. Recall that the isomorphism locus $X_\mc{A}$ is open, so that we can restrict $\mc{A}$ to $X_\mc{A}$ as a Lie algebroid. We can then consider the dual of the inverse of the anchor, i.e.\ $(\rho_\mc{A}^{-1})^*\colon \Omega^\bullet(X_\mc{A};\mc{A}) \to \Omega^\bullet(X_\mc{A})$, allowing us to view $\mc{A}$-forms as smooth forms on $X_\mc{A}$.
\begin{thm}\label{thm:amoser} Let $\mc{A} \to X$ be a Lie algebroid with dense isomorphism locus for which $Z_\mc{A}$ is smooth, and $k \in \{1, 2, \dim X\}$. Let $\omega, \omega' \in \Omega^k(\mc{A})$ be $d_\mc{A}$-closed $\mc{A}$-$k$-forms which are nondegenerate on $Z_\mc{A}$. Assume that either
	\bi
	\item $[\omega] = [\omega'] \in H^k(\mc{A})$, or
	\item $(\rho_\mc{A}^{-1})^*(\omega' - \omega)$ extends smoothly by $0$ over $Z_\mc{A}$.
	\ei
	Then there exists a Lie algebroid isomorphism $(\varphi,f)\colon (\mc{A}|_{U},\omega) \to (\mc{A}|_{U'},\omega')$ on neighbourhoods $U$ and $U'$ of $Z_\mc{A}$ for which $\varphi^* \omega' = \omega$, which in the second case can be chosen such that $f|_{Z_\mc{A}} = {\rm id}$.
\end{thm}
\begin{rem} As a consequence of \autoref{thm:amoser}, to establish an $\mc{A}$-Darboux theorem providing a pointwise normal form for $\mc{A}$-\symp{}s, one need only establish what an $\mc{A}$-\symp must look like locally near a point in $Z_\mc{A}$.
\end{rem}
\bp[ of \autoref{thm:amoser}] For $t \in [0,1]$, define $\omega_t := \omega + t(\omega' - \omega)$. In the first case, by definition there exists $\sigma \in \Omega^{k-1}(\mc{A})$ such that $\omega' = \omega + d_\mc{A} \sigma$. In the second case, instead by \autoref{lem:arelpoincare} there exists a neighbourhood $U$ of $Z_\mc{A}$ and $\sigma \in \Omega^1(U;\mc{A})$ vanishing on $Z_\mc{A}$ for which the same conclusion holds. Hence in both cases $\omega_t = \omega + t d_\mc{A} \sigma$ on some neighbourhood $U$ of $Z_\mc{A}$. By openness of the nondegeneracy condition and compactness of $[0,1]$, by shrinking $U$ if necessary we can ensure that $\omega_t$ is nondegenerate for all $t \in [0,1]$.

Nondegeneracy of $\omega_t$ then implies the existence of a (unique unless $k = 1$) smooth family of sections $v_t \in \Gamma(\mc{A})$ satisfying $\sigma + \iota_{v_t} \omega_t = 0$. As $\mc{A}$ has dense isomorphism locus, the anchor $\rho_\mc{A}\colon \Gamma(\mc{A}) \to \Gamma(TX)$ is injective, so that $V_t := \rho_{\mc{A}}(v_t)$ is a smooth family of vector fields on $X$. Again shrinking $U$ if necessary, ensure the existence of a family of diffeomorphisms $\psi_t$ for $t \in [0,1]$ such that $\psi_0 = {\rm id}$ and $\frac{d}{dt} \psi_t = V_t \circ \psi_t$. Then $T \psi_t$ determines a smooth family of Lie algebroid isomorphisms $\varphi_t$ of $\mc{A}$ such that $\varphi_t \equiv T \psi_t$ on sections over $X_\mc{A}$.

Note that over $U \cap X_\mc{A}$ we have $(\rho_{\mc{A}}^{-1})^* \varphi_t^* = (\varphi_t  \circ \rho_{\mc{A}}^{-1})^* = (\rho_\mc{A}^{-1} \circ \psi_t)^* = \psi_t^* (\rho_\mc{A}^{-1})^*$. Consider the $\mc{A}$-$k$-forms $\eta_t := \varphi_t^* \omega_t$ and introduce the notation $\omega_t' := (\rho_\mc{A}^{-1})^* \omega_t$ and $\sigma' := (\rho_\mc{A}^{-1})^* \sigma$. Then $\eta_0 = \omega_0 = \omega$. We compute over $U \cap X_\mc{A}$ that $(\rho_{\mc{A}}^{-1})^* \frac{d}{dt} \eta_t = (\rho_{\mc{A}}^{-1})^* \frac{d}{dt} \varphi_t^* \omega_t = \frac{d}{dt} (\rho_{\mc{A}}^{-1})^* \varphi_t^* \omega_t = \frac{d}{dt} \psi_t^* (\rho_\mc{A}^{-1})^* \omega_t = \frac{d}{dt} \psi_t^* \omega_t' = \psi_t^* \left(\frac{d}{dt}\omega_t' + \iota_{V_t} d\omega_t' + d \iota_{V_t} \omega_t'\right)$. We have that $\frac{d}{dt} \omega_t' = \frac{d}{dt} (\rho_\mc{A}^{-1})^* \omega_t = (\rho_\mc{A}^{-1})^* \frac{d}{dt} \omega_t$\\ ${= (\rho_\mc{A}^{-1})^* d_\mc{A} \sigma = d (\rho_\mc{A}^{-1})^* \sigma = d\sigma'}$. As $\omega_t$ is $d_\mc{A}$-closed so that $d\omega_t' = 0$, we thus see that $(\rho_{\mc{A}}^{-1})^* \frac{d}{dt} \eta_t = 0$ for all $t \in [0,1]$, as $\sigma' + \iota_{V_t} \omega_t' = (\rho_{\mc{A}}^{-1})^* \sigma + \iota_{\rho_\mc{A}(v_t)} (\rho_{\mc{A}}^{-1})^* \omega_t = (\rho_{\mc{A}}^{-1})^*(\sigma + \iota_{v_t} \omega_t) = 0$. Consequently, by continuity and density of $X_\mc{A}$, we conclude that $\frac{d}{dt} \eta_t = 0$ for all $t \in [0,1]$ on all of $U$. Hence $\eta_t = \omega$ for all $t \in [0,1]$, which in particular implies that for $\varphi := \varphi_1$ we have $\varphi^* \omega = \omega'$ as desired. In the second case $\sigma$ vanishes on $Z_\mc{A}$, so that the families $v_t$ and $V_t$ do as well. But then each $\psi_t$ is the identity on $Z_\mc{A}$, so that in particular this is true for $f = \psi_1$.
\ep
Inspecting the proof, we see that a similar statement can be made globally, as the size of the open $U$ is only dictated by the existence of a primitive and nondegeneracy.
\begin{cor}\label{cor:amoserglobal} Let $\mc{A} \to X$ be a Lie algebroid with dense isomorphism locus over a compact manifold for which $Z_\mc{A}$ is smooth, and $k \in \{1, 2, \dim X\}$. Let $\omega, \omega' \in \Omega^k(\mc{A})$ be $d_\mc{A}$-closed $\mc{A}$-$k$-forms such that $(1-t) \omega + t \omega'$ is nondegenerate for all $t \in [0,1]$ and $[\omega] = [\omega'] \in H^k(\mc{A})$. Then there exists a Lie algebroid isomorphism $(\varphi,f)\colon (\mc{A},X) \to (\mc{A},X)$ for which $\varphi^* \omega' = \omega$.
\end{cor}
Using this result, we see that $\mc{A}$-Nambu structures (hence $\mc{A}$-symplectic structures on surfaces) specifying the same $\mc{A}$-orientation are classified by their $\mc{A}$-cohomology class. Namely, letting $n = \dim X$, we note that $\det(\mc{A}^*)$ is a line bundle. Consequently, given cohomologous forms $\omega,\omega' \in \Omega^n(\mc{A}) = \Gamma(\det(\mc{A}^*))$ we have $\omega = f \omega'$ for some nonvanishing function $f \in C^\infty(X)$. This function must be strictly positive as $\omega$ and $\omega'$ give rise to the same $\mc{A}$-orientation, so that $(1-t)\omega + t \omega' = ((1-t) + t f) \omega$ is nondegenerate for all $t \in [0,1]$. If one has more knowledge about the Lie algebroid $\mc{A}$, the assumption on induced $\mc{A}$-orientations can sometimes be dropped. We summarize this discussion in the following proposition.
\begin{prop} Let $\mc{A} \to X$ be a Lie algebroid with dense isomorphism locus for which $Z_\mc{A}$ is smooth. An $\mc{A}$-Nambu structures $\Pi$ inducing a given $\mc{A}$-orientation is classified up to $\mc{A}$-orientation-preserving isomorphism by its $\mc{A}$-cohomology class.
\end{prop}
Let us reiterate this also classifies $\mc{A}$-symplectic structures on surfaces.

\section{Examples}
\label{sec:examplessymp}
In this section we discuss examples of $\mc{A}$-\symp{}s, using the concrete Lie algebroids discussed in Section \ref{sec:laexamples}. We will mostly focus our attention on those carried by the log-tangent and elliptic tangent bundles, as it is these that we will construct in Chapters \ref{chap:constructingblogs} and \ref{chap:constructingsgcs}. In all cases, it is useful to keep in mind that $\mc{A}$-\symp{}s are dual to nondegenerate $\mc{A}$-Poisson structures that are to be defined in the next chapter.
\subsection{Log-symplectic structures}
\label{sec:logsympstr}
We start with the symplectic structures for the log-tangent bundle of Section \ref{sec:logtangentbundle}. These have received much attention recently, see for example \cite{Cavalcanti17,FrejlichMartinezTorresMiranda15,GualtieriLi14,GuilleminMirandaPires14,MarcutOsornoTorres14,MarcutOsornoTorres14two}. Throughout this section, assume that $X$ is a $2n$-dimensional manifold.
\begin{defn}\label{defn:logsymplectic} Let $(X,Z)$ be a log pair with log-tangent bundle $\mc{A}_Z$. A \emph{log-symplectic structure} is an $\mc{A}_Z$-symplectic structure. It is called \emph{bona fide} if $Z \neq \emptyset$.
\end{defn}
\begin{rem} Such forms are also called \emph{$b$-symplectic structures} \cite{GuilleminMirandaPires14}, noting in Section \ref{sec:logtangentbundle} that the log-tangent bundle is also called the $b$-tangent bundle.
\end{rem}
Let $\omega_{\mc{A}_Z} \in {\rm Symp}(\mc{A}_Z)$ be a log-symplectic structure. We establish that by applying the residue map ${\rm Res}_Z\colon \Omega^2(\mc{A}_Z) \to \Omega^1(Z)$, we obtain a cosymplectic structure.
\begin{defn} Let $Y^{2n-1}$ be a manifold. A \emph{cosymplectic structure} on $Y$ is pair $(\alpha,\beta) \in \Omega^1_{\rm cl}(Y) \times \Omega^2_{\rm cl}(Y)$ of closed one- and two-forms such that $\alpha \wedge \beta^{n-1} \neq 0$.
\end{defn}
It is immediate that a cosymplectic structure induces an orientation. Let $z \in I_Z$ be a generator and locally write $\omega_{\mc{A}_Z} = d \log z \wedge \alpha_0 + \alpha_1$ with $\alpha_i$ smooth and closed. Then ${\rm Res}_Z(\omega_{\mc{A}_Z}) = j_Z^* \alpha_0$ for $j_Z \colon Z \hookrightarrow X$. Alternatively, it is given by the interior contraction $\iota_{\mathbb{L}_Z} \omega_{\mc{A}_Z}$, with $z \partial_z$ providing a trivialization of $\mathbb{L}_Z$. Consider the pair $(\alpha,\beta) = (j_Z^* \alpha_0, j_Z^* \alpha_1)$ on $Z$. Then as $\omega_{\mc{A}_Z}$ is nondegenerate we get that $0 \neq \omega_{\mc{A}_Z}^n = d \log_z \wedge \alpha_0 \wedge \alpha_1^{n-1}$. This shows that in particular that $(\alpha,\beta)$ is a cosymplectic structure on $Z$. While $\alpha$ is invariantly defined (being a residue), $\beta$ is not and depends on the choice of generator $z \in I_Z$. Changing $z$ to $z' = e^h z$ for some smooth function $h$, we see that $d\log z' = d\log z + dh$. Consequently, $\alpha_1$ gets changed to $\alpha_1' = \alpha_1 + \alpha_0 \wedge dh$. In summary, what is well-defined is the following.
\begin{prop}[\cite{GuilleminMirandaPires14,MarcutOsornoTorres14}]\label{prop:blogcosymp} Let $(X,Z)$ be a log pair. Then a log-symplectic structure $\omega_{\mc{A}_Z}$ determines a \emph{log equivalence class} of cosymplectic structures on $Z$, i.e.\ a cosymplectic structure $(\alpha,\beta)$ where $\beta$ is defined up to forms $d(h \alpha)$ for $h \in C^\infty(Z)$.
\end{prop}
Using the Moser argument (see Section \ref{sec:amosertechniques}), we obtain a Darboux theorem for \blog{}s around points in $Z$.
\begin{thm}[{\cite[Theorem 37]{GuilleminMirandaPires14}}] Let $(X,Z)$ be a log pair with \blog{} $\omega_{\mc{A}_Z}$, and let $x \in Z$. Then there exists an open neighbourhood $U$ of $x$ which is $\mc{A}_Z$-symplectomorphic to $d\log z \wedge dy + \omega_0$, where $\{z = 0\} = U \cap Z$ and $\omega_0$ is the standard symplectic structure on $\R^{2n-2}$.
\end{thm}
Using the relative Moser argument to $Z$ we can further obtain the following result extracted from \cite{GuilleminMirandaPires14} describing a \blog around its singular locus, see also \cite[Theorem 3.2]{Cavalcanti17} and the discussion before \cite[Theorem 1.12]{GualtieriLi14}.
\begin{prop}\label{prop:bloglocalform} Let $(X^{2n},Z_X,\omega_{\mc{A}_{Z_X}})$ be a compact log-symplectic manifold with induced log equivalence class of cosymplectic structures $[(\alpha,\beta)]$. Then around any connected component $Z$ of $Z_X$, the two-form $\omega_{\mc{A}_{Z_X}}$ is equivalent to $d \log |x| \wedge \alpha + p^*(\beta)$ in a neighbourhood of the zero section of the normal bundle $p\colon NZ \to Z$, where $|x|$ is the distance to the zero section with respect to a fixed metric on $NZ$.
\end{prop}
\begin{rem} In the above proposition, the use of a distance function is required because $Z$ need not be coorientable. This will be the case if $X$ is orientable, however.
\end{rem}
\begin{rem} The fact that a result such as \autoref{prop:bloglocalform} holds is essentially due to the fact that the residue sequence of $\mc{A}_Z$ splits, as was used to compute its Lie algebroid cohomology (see \autoref{thm:loglacohomology}).
\end{rem}
Any compact cosymplectic manifold, such as $Z_\pi$, fibers over $S^1$ \cite{GualtieriLi14, Li08, MarcutOsornoTorres14, OsornoTorres15}. This is due to the fact that it has a closed nonvanishing one-form, namely $\alpha$. Call a cosymplectic manifold $(Y,\alpha,\beta)$ \emph{proper} if the distribution given by $\ker \alpha$ has a compact leaf. Given a \blog $\omega_{\mc{A}_Z}$ on $(X,Z)$, call a connected component of $Z$ \emph{proper} if its induced cosymplectic structure proper. Finally, call $\omega_{\mc{A}_Z}$ \emph{proper} if all connected components of $Z$ are proper. We then have the following (see also \cite{OsornoTorres15}).
\begin{prop}[{\cite[Theorem 3.6]{Cavalcanti17}}]\label{prop:blogproper} Let $(X,Z)$ be a log pair equipped with a \blog. Then if the connected components of $Z$ are compact, the \blog can be deformed into one which is proper.
\end{prop}
As we will always be working with compact log pairs, we can thus without loss of generality assume any \blog is proper, so that $Z$ fibers over $S^1$ using $\alpha = {\rm Res}_Z(\omega_{\mc{A}_Z})$. For later use we record the following.
\begin{lem}\label{lem:blogsurface} Let $\Sigma^2$ be a compact oriented surface. Then $(\Sigma,\partial \Sigma)$ admits a \blog. For a hypersurface $Z$, $(\Sigma,Z)$ admits a \blog if and only if $[Z] = 0 \in H_1(\Sigma;\Z_2)$.
\end{lem}
\bp The log pair $(\Sigma, Z)$ carries a \blog if and only if $\mc{A}_Z$ is orientable. Using \autoref{prop:logbundleiso} we have $w_1(\mc{A}_Z) = w_1(\Sigma) + w_1(L_Z)$. As an orientable manifold with boundary has orientable boundary, the result follows by \autoref{prop:logdivstiefelwhitney}, which gives that $w_1(L_Z) = {\rm PD}_{\Z_2}([Z]) \in H^1(\Sigma;\Z_2)$.
\ep
Moreover, let us mention that there are cohomological obstructions for a compact manifold to admit a \blog.
\begin{thm}[\cite{MarcutOsornoTorres14two}]\label{thm:blogobstr} Let $X^{2n}$ be a compact log-symplectic manifold. Then there exists a class $a \in H^2(X;\R)$ such that $a^{n-1} \neq 0$.
\end{thm}
\begin{thm}[\cite{Cavalcanti17}]\label{thm:bonafideobstr} Let $X^{2n}$ be a compact oriented bona fide log-symplectic manifold. Then there exists a nonzero class $b \in H^2(X;\R)$ such that $b^2 = 0$. Moreover, if $n > 1$ then $b_2(X) \geq 2$.
\end{thm}
As a further consequence, if $X$ is four-dimensional we see that it must have indefinite intersection form in order for it to admit a bona fide \blog. A log-symplectic structure can be modified in a tubular neighbourhood of a specific type of codimension-two submanifold to increase its singular locus. We phrase the following using the dual bivectors (see Section \ref{sec:logpoisson}) to avoid having to introduce more notation.
\begin{thm}[{\cite[Theorem 5.1]{Cavalcanti17}}]\label{thm:addsinglocus} Let $(X^{2n},Z_X,\pi)$ be a compact log-Poisson manifold and let $k > 0$ be an integer. Assume that $X$ has a compact symplectic submanifold $F^{2n-2} \subset X \setminus Z_X$ with trivial normal bundle. Then $(X,Z_X')$ admits a log-Poisson structure $\pi'$ agreeing with $\pi$ away from $Z_X' \setminus Z_X$, where $Z_X'$ is the disjoint union of $Z_X$ with $k$ copies of $F \times S^1$.
\end{thm}
\subsection{Elliptic symplectic structures}
\label{sec:ellipticsymp}
In this section we discuss symplectic structures \cite{CavalcantiGualtieri15} for the elliptic tangent bundle of Section \ref{sec:elltangentbundle}. We focus mainly on those elliptic symplectic structures which have zero elliptic residue. This is because these link to a geometric structure of independent interest, namely that of a \sgcs (see Section \ref{sec:scgs}).
\begin{defn} Let $(X,|D|)$ be an elliptic pair with elliptic tangent bundle $\mc{A}_{|D|}$. An \emph{elliptic symplectic structure} is an $\mc{A}_{|D|}$-symplectic structure.
\end{defn}
Recall that given an elliptic divisor we can find coordinates $(r,\theta,x_i)$ for which $I_{|D|} = \langle r^2 \rangle$. Then $\omega_{\mc{A}_{|D|}}$ is given by $\alpha_0 d\log r \wedge d\theta + d \log r \wedge \alpha_1 + d \theta \wedge \alpha_2 + \alpha_3$ with $\alpha_i$ smooth. If $\omega_{\mc{A}_{|D|}}$ has zero elliptic residue, this means that $\alpha_0 = 0$. The Moser argument (see Section \ref{sec:amosertechniques}) implies the following Darboux theorem.
\begin{thm}[\cite{CavalcantiGualtieri15}] Let $(X,|D|)$ be an elliptic pair with elliptic symplectic structure $\omega_{\mc{A}_{|D|}}$ having zero elliptic residue, and let $x \in D$. Then there exists an open neighbourhood $U$ of $x$ which is $\mc{A}_{|D|}$-symplectomorphic to $d \log r \wedge d x + d\theta \wedge dy + \omega_0$, where $\{r = 0\} = U \cap D$ and $\omega_0$ is the standard symplectic structure on $\R^{2n-4}$.
\end{thm}
\subsection{Zero symplectic structures}
\label{sec:zerosymp}
In this section we discuss symplectic structures in the zero tangent bundle $\mc{B}_Z$ associated to a log pair $(X,Z)$.
Sadly, this Lie algebroid cannot admit symplectic structures on manifolds of dimension higher than two, due to the following reuslt.
\begin{prop}[{\cite[Proposition 2.21]{Lanius16}}]\label{prop:zerotangentnosymp} Let $(X,Z)$ be a log pair such that $\dim X > 2$. Then $X$ does not admit a $\mc{B}_Z$-symplectic structure.
\end{prop}
The same proposition further says there cannot exist Lie algebroid symplectic structures for any of the Lie algebroids $\mc{A}_Z^{k,m}$ if $m > 0$ and $k \neq 1$, where $\mc{A}_Z^{k,m}$ is obtained from the $b^k$-tangent bundle $\mc{A}^l_Z$ by performing $(0,Z)$-rescaling a total of $m$ times. The above case is when $m = 1$ and $k = 0$, where we agree that $\mc{A}_Z^0 = TX$.

On a log surface $(\Sigma,Z)$, by an application of the Moser argument (see Section \ref{sec:amosertechniques}) it is immediate that any $\mc{B}_Z$-symplectic structure can be written around points $x \in Z$ as $d\log z \wedge \frac{dy}{z}$ for a generator $z \in I_Z$.
\subsection{Scattering symplectic structures}
\label{sec:scatteringsymp}
In this section we briefly discuss symplectic structures in the scattering tangent bundle $\mc{C}_Z$ of Section \ref{sec:scatteringtangentbundle} associated to a log pair $(X,Z)$. For more information, see \cite{Lanius16}. In terms of the Lie algebroids $\mc{A}_Z^{k,m}$ introduced in the previous sections, we see that $\mc{C}_Z = \mc{A}_Z^{1,1}$, so that the scattering tangent bundle is not excluded from having symplectic structures by the full statement of \autoref{prop:zerotangentnosymp}.

Contrary to the log-symplectic case, where $Z$ gets an induced cosymplectic structure, in the scattering case there is instead an induced \emph{contact structure}, i.e.\ a one-form $\alpha$ for which $\alpha \wedge d \alpha^{n-1} \neq 0$ (see \cite[Proposition 3.1]{Lanius16}). Using \autoref{thm:amoser} and the normal form for contact structures (see e.g.\ \cite{Geiges08}) we obtain the Darboux theorem for scattering-symplectic structures.
\begin{prop}[{\cite[Proposition 3.2]{Lanius16}}] Let $(X,Z,\omega_{\mc{C}_Z})$ be a $\mc{C}_Z$-symplectic log pair and $x \in Z$. Then there exists an open neighbourhood $U$ of $x$ for which $\omega_{\mc{C}_Z} = d (\frac1{z^2} \alpha_0) = \frac{d z}{z^3} \wedge \alpha_0 + \frac1{z^2} \omega_0$, where $\{z = 0\} = U \cap Z$ and $\alpha_0$ is the standard contact structure on $\R^{2n-1}$ and $\omega_0$ is the standard symplectic structure on $\R^{2n-2}$.
\end{prop}
\subsection{$b^k$-symplectic structures}
\label{sec:bksymp}
In this section we discuss $b^k$-symplectic structures, which are symplectic structures on the $b^k$-tangent bundles $\mc{A}_Z^k$ of Section \ref{sec:bktangentbundle}. These are studied by Scott in \cite{Scott16}.
\begin{defn} Let $(X,Z)$ be a log pair and $j_{k-1} \in J_Z^{k-1}$. A \emph{$b^k$-symplectic structure} is an $\mc{A}_Z^k$-symplectic structure. It is called \emph{bona fide} if $Z \neq \emptyset$.
\end{defn}
Given a $b^k$-symplectic structure $\omega_{\mc{A}_Z^k}$, we can apply ${\rm Res}_Z^k$ to $\omega_{\mc{A}_Z^k}$ to obtain a nonzero closed one-form. This form can also be viewed as the contraction $\iota_{\mathbb{L}_{z,k}} \omega$. Locally using $z \in j_{k-1}$ so that we can write $\mathbb{L}_{z,k} = \langle z^k \partial_z \rangle$ and $\omega_{\mc{A}_Z^k} = \frac{dz}{z^k} \wedge \alpha + \beta$. Note that $\alpha, \beta$ are not uniquely determined by $\omega_{\mc{A}_Z^k}$ but depend on $z$. However, we have
\be
	\iota_{\mathbb{L}_{z,k}} \omega = \iota_{z^k \partial_z} \left(\frac{dz}{z^k} \wedge \alpha + \beta\right) = j_Z^* \alpha.
\ee
Consequently, there is an induced cosymplectic structure on $Z$, as was true for log-symplectic structures (see \autoref{prop:blogcosymp}). Moreover, this cosymplectic structure can again be used to construct a local form for $b^k$-symplectic structures (see \cite{Scott16}).
\begin{prop}\label{prop:bksymplocalform} Let $(X,Z)$ be a compact $b^k$-symplectic manifold for $j_{k-1} \in J_{Z_X}^{k-1}$, with induced log equivalence class of cosymplectic structures $[(\alpha,\beta)]$. Then around $Z$, the two-form $\omega_{\mc{A}^k_{Z_X}}$ is equivalent to $\frac1{|x|^{k-1}} d\log |x| \wedge \alpha + p^*(\beta)$ in a neighbourhood of the zero section of the normal bundle $p\colon NZ \to Z$, where $|x|$ is the distance to the zero section with respect to a fixed metric on $NZ$.
\end{prop}
\subsection{Complex log-symplectic structures}
\label{sec:complexlogsymp}
In this section we discuss a type of geometric structure that strictly speaking does not fall under the name $\mc{A}$-symplectic structure. Recall from Section \ref{sec:elltangentbundle} that given a complex log divisor, there is a natural Lie algebroid morphism $\iota\colon \mc{A}_{|D|} \otimes \C \to \mc{A}_D$ given by an inclusion of sheaves.

Let $(X,H)$ be a manifold with closed three-form and let $D = (U,\sigma)$ be a complex log divisor with complex log-tangent bundle $\mc{A}_D$.
\begin{defn}[\cite{CavalcantiGualtieri15}]  A \emph{complex log-symplectic structure} is a $\mc{A}_D$-two-form $\omega_{\mc{A}_D}$ such that $d_{\mc{A}_D} \omega_{\mc{A}_D} = \rho_{\mc{A}_D}^* H$, and the complex elliptic form $\iota^* \omega_{\mc{A}} = b + i \omega_{\mc{A}_{|D|}}$ has nondegenerate imaginary part.
\end{defn}
As the three-form $H$ is real, the condition $d_{\mc{A}_D} \omega_{\mc{A}_D} = \rho_{\mc{A}_D}^* H$ implies that $d_{\mc{A}_{|D|}} \omega_{\mc{A}_|D|} = 0$. Together with the nondegeneracy assumption this means that $\omega_{\mc{A}_{|D|}}$ is an $\mc{A}_{|D|}$-symplectic structure, which necessarily has zero elliptic residue. Using results from \cite{CavalcantiGualtieri15} (Theorem 3.7 loc.\ cit.), one can see that this map is surjective if $D$ is coorientable: given an elliptic symplectic structure (with zero elliptic residue) associated to a coorientable elliptic divisor $|D|$, there exists a complex log divisor $D = (U,\sigma)$ (see \autoref{rem:ellipticomplexlog}) inducing $|D|$. Moreover, there exists a closed three-form $H$ and a complex log-symplectic structure $\omega_{\mc{A}_D}$ with respect to $H$ for which ${\rm im}^*(\omega_{\mc{A}_{D}}) = \omega_{\mc{A}_{|D|}}$.

\section{Folded-\symp{}s}
\label{sec:blogsandfsymps}
\renewcommand{\b}[1][]{{^b}{#1}}
We close this chapter by taking a slight detour to discuss \fsymp{}s, which are closed two-forms that are in a controlled sense almost-everywhere symplectic. We show that \blog{}s naturally give rise to \fsymp{}s. \Fsymp{}s are studied amongst others in \cite{Baykur06, CannasDaSilvaGuilleminPires11, CannasDaSilva10, CannasDaSilvaGuilleminWoodward00}. The results in this section appeared before in \cite{CavalcantiKlaasse16}.

\begin{defn} A \emph{\fsymp} on a compact $2n$-dimensional manifold $X$ is a closed two-form $\omega$ such that $\bigwedge^n \omega$ is transverse to the zero section in $\bigwedge^{2n} T^*X$, and such that $\omega^{n-1}|_{Z_\omega} \neq 0$, where $Z_\omega = (\bigwedge^n \omega)^{-1}(0)$. The hypersurface $Z_\omega$ is called the \emph{folding locus} of $\omega$, while its complement $X \setminus Z_\omega$ is called the \emph{symplectic locus}. A \fsymp is called \emph{bona fide} if $Z_\omega \neq \emptyset$.
\end{defn}
This definition should be compared with \autoref{defn:logsymplectic} (and \autoref{defn:logpoisson} below). We will say the pair $(X,Z_X)$ \emph{admits a \fsymp} if $X$ admits a \fsymp $\omega$ for which $Z_\omega = Z_X$. According to the Darboux model, a \fsymp $\omega$ is locally given by $\omega = x_1 d x_1 \wedge dx_2 + \dots + dx_{2n-1} \wedge dx_{2n}$, using coordinates $x_i$ in a neighbourhood $U$ such that $Z_\omega \cap U = \{x_1 = 0\}$.

Cannas da Silva gave a homotopical characterization for an orientable manifold to admit a \fsymp.
\begin{thm}[\cite{CannasDaSilva10}]\label{thm:fsympstableacs} Let $X$ be an orientable manifold. Then $X$ admits a \fsymp if and only if $X$ admits a stable \acs. In particular, every orientable four-manifold admits a \fsymp.
\end{thm}
Further, Baykur has given a construction showing there is a relation between \alf{}s on four-manifolds and \fsymp{}s. This should be compared with \autoref{thm:alflogsymp}.
\begin{thm}[{\cite[Proposition 3.2]{Baykur06}}]\label{thm:alffsymp} Let $f\colon X^4 \to \Sigma^2$ be an \alf between compact connected manifolds. Assume that the generic fiber $F$ is orientable and $[F] \neq 0 \in H_2(X;\R)$. Then $X$ admits a \fsymp.
\end{thm}
From our point of view this theorem does not come as a surprise. Indeed, every \blog on $(X,Z_X)$ gives rise to a \fsymp on $(X,Z_X)$, as is known in the Poisson community. We learned the following proof of this result from M{\u{a}}rcu{\c{t}} and Frejlich. A similar proof in the context of desingularizing orientable $b^k$-manifolds can be found in \cite{GuilleminMirandaWeitsman15} (note that a $b^1$-manifold is essentially a $b$-manifold).
\begin{thm}\label{thm:blogfolded} Let $(X^{2n},Z_X,\pi)$ be a compact log-symplectic manifold. Then $(X,Z_X)$ admits a \fsymp $\omega$ for which $\omega = \pi^{-1}$ outside a neighbourhood of $Z_X$.
\end{thm}
\bp By \autoref{prop:bloglocalform} a neighbourhood of each connected component $Z$ of $Z_X$ is equivalent to a neighbourhood $U$ of the zero section of the normal bundle $NZ$ equipped with a distance function $|x|$, so that $\pi^{-1} = d \log |x| \wedge \theta + \sigma$ for closed one- and two-forms $\theta$ and $\sigma$ on $Z$ satisfying $\theta \wedge \sigma^{n-1} \neq 0$. By rescaling $|x|$ we can assume that $U$ contains all points of distance at most $e^2 + 1$ away from the zero section. Denote $\omega_Z = d |x|^2 \wedge \theta + \sigma$ and let $f\colon \R_+ \to \R$ be a smooth monotone interpolation between the functions $f_0\colon [0,1] \to \R$, $f_0(x) = x^2$ and $f_1\colon [e^2,\infty) \to \R$, $f_1(x) = \log x$. Consider the closed two-form $\omega_f = df(|x|) \wedge \theta + \sigma$, extended by $\pi^{-1}$ outside of $U$. Then $\omega_f = \pi^{-1}$ away from $Z$, while near $Z$ we have $\omega_f = \omega_Z$. Moreover, $\omega_f$ is symplectic on $X \setminus Z$ by monotonicity of $f$. Perform this procedure for all connected components of $Z_X$ to obtain a closed two-form $\omega$ on $X$ for which $\omega = \pi^{-1}$ away from $Z_X$. By the local description near $Z_X$ it follows that $\omega^n$ vanishes transversally with $Z_\omega = Z_X$. Further, the restriction of $\omega^{n-1}$ to $Z_\omega$ is equal to $\sigma^{n-1}$, hence is nonvanishing as $\theta \wedge \sigma^{n-1} \neq 0$. We conclude that $\omega$ is a \fsymp for $(X,Z_X)$.
\ep
The previous theorem, together with \autoref{thm:alflogsymp}, implies \autoref{thm:alffsymp}. The \fsymp of \autoref{thm:alffsymp} agrees with the one obtained through our methods, as is hinted at by the fact that in Baykur's construction the folding locus fibers over the circle.

The converse to \autoref{thm:blogfolded} does not hold. For example, $S^4$ does not admit a \blog by \autoref{thm:blogobstr}, while it does admit \fsymp{}s by \autoref{thm:fsympstableacs}. Similarly, by \autoref{thm:bonafideobstr} and results from Seiberg-Witten theory due to Taubes \cite{Taubes94}, the four-manifold $\C P^2 \# \C P^2$ does not admit \blog{}s, bona fide or not. However, by \autoref{thm:fsympstableacs} it does admit a \fsymp.

The reason the converse to \autoref{thm:blogfolded} is false is essentially because the information contained in the one-form $\theta$ determined by the \blog (see \autoref{prop:bloglocalform}) is lost when passing to the folded-symplectic world. Note here that a \fsymp $\omega$ restricts to $Z_\omega$ to define a one-dimensional foliation $\ker(\left.\omega\right|_{Z_\omega})$ called the \emph{null foliation}, and $\left.\omega\right|_{Z_\omega}$ is a \emph{pre-symplectic structure} on $Z_\omega$, i.e.\ a closed two-form of maximal rank. On the other hand, a \blog $\pi$ induces a cosymplectic structure on $Z_\pi$ by \autoref{prop:blogcosymp}, so that the associated nowhere-vanishing closed one-form gives a codimension-one foliation on $Z_\pi$. Further, there is a symplectic structure on the leaves.

The following result makes precise that it is exactly the existence of a suitable closed one-form $\theta$ on $Z_X$ that makes the converse to \autoref{thm:blogfolded} hold.
\begin{thm}\label{thm:foldedblog} Let $(X^{2n},Z_X,\omega)$ be a compact folded-symplectic manifold. Assume that there exists a closed one-form $\theta \in \Omega^1(Z_X)$ such that $\theta \wedge \omega^{n-1}|_{Z_X} \neq 0$. Then $(X,Z_X)$ admits a \blog $\pi$ for which $\pi = \omega^{-1}$ outside a neighbourhood of $Z_X$.
\end{thm}
In other words, when the folded-symplectic form $\omega$ can be complemented to give a cosymplectic structure on $Z_X$, one can turn $\omega$ into a \blog.
\bp Consider the normal bundle $NZ_X$ and let $U$ be a neighbourhood of the zero section. Choose a distance function $|x|$ for $Z_X$ which is constant outside of $\overline{U}$. Note that $d \log |x| \wedge \theta \wedge \omega^{n-1}$ is nonzero at $Z_X$ as $\theta \wedge \omega^{n-1}|_{Z_X} \neq 0$ and $|x|$ is transverse to $Z_X$. By continuity it is still nonzero in a neighbourhood $V \subset U$ of $Z_X$. Let $f = f(|x|)$ be a smooth function on $NZ_X$ so that $f \equiv 1$ on $V$ and $f \equiv 0$ near $\partial U$, which is then extended to $X$ by being identically $0$ outside $U$. Define a closed $b$-two-form $\omega_f = t \, d (f \log |x|) \wedge \theta + \omega \in \b{\Omega}^2(X)$ for the $b$-manifold $(X,Z_X)$, where $t \neq 0$ is a real parameter. Choose the sign of $t$ so that the forms $t \, d(f \log|x|) \wedge \theta$ and $\omega$ give the same orientation on $V \setminus Z_X$. We have $\omega_f^n = t n \, d (f \log |x|) \wedge \theta \wedge \omega^{n-1} + \omega^n$, so by choosing $t$ small enough we conclude that $\omega_f$ is a $b$-symplectic form for which $\omega_f = \omega$ outside $\overline{U}$. By \autoref{prop:blogbsymp} the dual bivector $\pi$ to $\omega_f$ is a \blog for $(X,Z_X)$, and $\pi = \omega^{-1}$ away from $Z_X$.
\ep
\begin{rem} A result that is close in spirit to our \autoref{thm:blogfolded} and \autoref{thm:foldedblog} combined regarding the relation between log-symplectic and \fsymp{}s was obtained in \cite{FrejlichMartinezTorresMiranda15}. There the authors show that an orientable open manifold $X$ admits a \blog if and only if $X$ admits a \fsymp. Contrastive with our results, their proof relies on an $h$-principle for open manifolds hence is purely existential and the different structures produced have no relation in general. In particular the corresponding loci may be completely different.
\end{rem}
\chapter{\texorpdfstring{$\mc{A}$}{A}-Poisson structures}
\label{chap:apoissonstructures}
\fancypagestyle{empty}{%
	\fancyhf{}%
	\renewcommand\headrulewidth{0pt}%
	\fancyhead[RO,LE]{\thepage}%
}
In this chapter we define Poisson structures for general Lie algebroids. These structures are the main object of study in this thesis. Their definition is identical to that of usual Poisson structures, except one uses Lie algebroid objects. Their behavior heavily depends on which Lie algebroid is used. Of particular importance is the notion of lifting Poisson structures to other Lie algebroids. This allows us to incorporate the degenerate behavior of a Poisson structure into that of the Lie algebroid it is lifted to. This procedure can be repeated and the hope is to end up with an $\mc{A}$-Poisson structure that is nondegenerate, so that it can be studied as an $\mc{A}$-symplectic structure of the previous chapter. To this end, we consider a specific class of Poisson structures which give rise to the divisors of Chapter \ref{chap:divisors} and discuss when they can be lifted.
\subsection*{Organization of the chapter}
This chapter is built up as follows. In Section \ref{sec:poissonstructures} we start with a discussion the usual notion of a Poisson structure on a manifold. Then in Section \ref{sec:apoissonstrs} we define $\mc{A}$-Poisson structures and discuss $\mc{A}$-analogues of concepts from Poisson geometry. After this, in Section \ref{sec:poissondivisors} we relate Poisson structures to divisors. We note in Section \ref{sec:riggedalgebroids} the relation between such Poisson structures and the concepts of rigged algebroids \cite{Lanius16} and almost-regularity \cite{AndroulidakisZambon17}. We continue in Section \ref{sec:poissonexamples} by discussing examples of Poisson structures of divisor-type and/or $\mc{A}$-Poisson structures. Finally, in Section \ref{sec:liftingpoisson} we introduce the concept of lifting Poisson structures between Lie algebroids, and discuss this procedure for some concrete divisors and ideal Lie algebroids.
\section{Poisson structures} 
\label{sec:poissonstructures}
In this section we give a brief introduction to Poisson geometry. For more information on Poisson structures, see e.g.\ \cite{LaurentGengouxPichereauVanhaecke13,Vaisman94, FernandesMarcut15}. Let $X$ be a manifold and denote the space of bivectors on $X$ by $\mf{X}^2(X) = \Gamma(\wedge^2 TX)$. Recall that $\mf{X}^\bullet(X) = \Gamma(\wedge^\bullet TX)$ carries an extension of the Lie bracket, called the \emph{Schouten bracket}, again denoted by $[\cdot,\cdot]$.
\begin{defn} A bivector $\pi \in \mf{X}^2(X)$ is \emph{Poisson} if $[\pi,\pi] = 0$. Denote by ${\rm Poiss}(X)$ the space of all Poisson structures on $X$.
\end{defn}
Such bivectors are related to a bracket on the space of functions of $X$.
\begin{defn} A \emph{Poisson bracket} on a manifold $X$ is a skew-symmetric billinear bracket $\{\cdot,\cdot\}\colon C^\infty(X) \times C^\infty(X) \to C^\infty(X)$ satisfying for all $f,g,h \in C^\infty(X)$
	\bi
	\item The Leibniz rule $\{h, f g\} = \{h, f\}g + f \{h, g\}$;
	\item The Jacobi identity $\{f, \{g,h\}\} + \{g, \{h,f\}\} + \{h,\{f,g\}\} = 0$.
	\ei
\end{defn}
A Poisson structure $\pi$ on $X$ is in one-to-one correspondence with a Poisson bracket $\{\cdot,\cdot\}_\pi$ on $C^\infty(X)$ via $\{f,g\}_\pi = \pi(df,dg)$ for $f,g \in C^\infty(X)$. The condition that $[\pi,\pi] = 0$ is equivalent to $\{\cdot,\cdot\}_\pi$ satisfying the Jacobi identity.

Given a bivector $\pi$, we obtain a map $\pi^\sharp\colon T^*X \to TX$ using $\pi^\sharp(\xi) = \iota_{\xi} \pi$ for $\xi \in T^*X$. We call the rank of $\pi^\sharp$ the \emph{rank} of the Poisson structure $\pi$, which by skew-symmetry must be even. A point $x \in X$ is \emph{regular} if the rank of $\pi$ is constant in some open neighbourhood of $x$, and singular otherwise. We denote by $X_{\pi, {\rm reg}} \subseteq X$ the space of regular points of $\pi$, also called the \emph{regular locus} of $\pi$. This is an open dense subspace of $X$, with its complement $X_{\pi, {\rm sing}}$, the \emph{singular locus} of $\pi$, being closed and nowhere dense. Finally, we call $\pi$ \emph{regular} if $X_{\pi, {\rm reg}} = X$. This means $\pi^\sharp$ has constant rank. Finally, we say $\pi$ is \emph{nondegenerate} if $\pi^\sharp$ is an isomorphism. We briefly discuss some examples of Poisson structures.
\begin{exa} Any manifold carries the trivial Poisson structure $\pi \equiv 0$.
\end{exa}
\begin{exa}\label{exa:symppoisson} Let $(X,\omega)$ be a symplectic manifold. Then $\pi = \omega^{-1}$ defined by $\pi^\sharp = (\omega^\flat)^{-1}$ is a Poisson structure, with $d\omega = 0$ being equivalent to $\pi$ being Poisson.
\end{exa}
\begin{exa}\label{exa:standardnondeg} Let $X = \R^{2n}$ with standard coordinate system $(x_1,\dots,x_{2n})$. Then $\pi_0 = - \sum_{i=1}^n \partial_{x_{2i-1}} \wedge \partial_{x_{2i}}$ is the \emph{standard nondegenerate Poisson structure} on $\R^{2n}$.
\end{exa}
The dual of a Lie algebra carries a canonical Poisson structure. See also \autoref{rem:lalgdpoisson}. We will not emphasize these Poisson structures much in this thesis.
\begin{exa}\label{exa:liepoisson} Let $\mf{g}$ be a finite-dimensional Lie algebra. Then $\mf{g}^*$ carries a natural linear Poisson structure defined by $\{f,g\}(\xi) = \langle [d_\xi f, d_\xi g]_\mf{g}, \xi \rangle$, where $f, g\in C^\infty(\mf{g}^*)$ and $\xi \in \mf{g}^*$. Here we view the differential $d_\xi f\colon T_\xi \mf{g}^* \to \R$ as an element of $\mf{g}$ using the relation $\langle d_\xi f, v \rangle = \frac{d}{dt}|_{t = 0} f(\xi + t v)$ for $v \in \mf{g}^*$.
\end{exa}
Note moreover that if $X$ is two-dimensional, any bivector $\pi$ on $X$ is automatically Poisson, as we have $[\pi,\pi] \in \mf{X}^3(X) = 0$. The existence of a nondegenerate Poisson structure on a surface is equivalent to orientability, as is true for symplectic structures.
\begin{defn} We say a Poisson structure $\pi$ is \emph{generically symplectic} if the space $X_{\pi, {\rm symp}}$ where $\pi^\sharp$ has full rank is open and dense. In this case $X_{\pi, {\rm symp}} = X_{\pi, {\rm reg}}$.
\end{defn}
We call $X_{\pi, {\rm symp}}$ the \emph{symplectic locus} of $\pi$, which in general may be empty. In case it is nonempty, its complement is the singular locus $X_{\pi, {\rm sing}}$ discussed before.
\begin{defn} Let $\pi \in {\rm Poiss}(X)$ be a Poisson structure. Given $f \in C^\infty(X)$, we call $V_f := \pi^\sharp(df)$ the \emph{$\pi$-Hamiltonian vector field} associated to $f$ by $\pi$.
\end{defn}
Let ${\rm Ham}(\pi)$ be the space of $\pi$-Hamiltonian vector fields. The property $V_{\{f,g\}_\pi} = [V_f,V_g]$ for all $f,g \in C^\infty(X)$ is equivalent to the Jacobi identity for $\{\cdot,\cdot\}_\pi$, as
\be
\{f,g\}_\pi = \mc{L}_{V_f} g = - \mc{L}_{V_g} f = d g(V_f) = - df(V_g).
\ee
The \emph{$\pi$-Poisson vector fields} are those $V \in \mf{X}(X)$ for which $\mc{L}_V \pi = 0$. Denote the space of $\pi$-Poisson vector fields by ${\rm Poiss}(\pi)$.
\begin{prop}\label{prop:hamispoiss} Let $\pi  \in {\rm Poiss}(X)$. Then ${\rm Ham}(\pi) \subseteq {\rm Poiss}(\pi)$.
\end{prop}
\bp Let $V \in \mf{X}(X)$ and $f,g \in C^\infty(X)$. Then by a short computation using the Lie derivative we have $\mc{L}_V \pi(df,dg) = \mc{L}_V \{f,g\}_\pi - \{\mc{L}_V f,g\}_\pi - \{f, \mc{L}_V g\}_\pi$. If $V = V_h \in {\rm Ham}(\pi)$ for $h \in C^\infty(X)$, then by definition $\mc{L}_{V_h} = \{h, \cdot\}_\pi$, so that $\mc{L}_{V_h} \pi \equiv 0$ by the Jacobi identity for $\{\cdot,\cdot\}_\pi$. We conclude that $V_h \in {\rm Poiss}(\pi)$.
\ep
\begin{defn} Given two Poisson manifold $(X,\pi)$ and $(X',\pi')$, a \emph{Poisson map} is a map $f\colon X \to X'$ such that $f_*(\pi) = \pi'$. In terms of Poisson brackets, we have $\{g,h\}_{\pi'} \circ f = \{g \circ f, h \circ f\}_\pi$ for all $g,h \in C^\infty(X')$.
\end{defn}
\begin{exa} Let $(X,\pi)$ and $(X',\pi')$ be Poisson manifolds. Then their product $(X \times X', \pi + \pi')$ is Poisson, and the projections $X \times X' \to X$ and $X \times X' \to X'$ are Poisson maps.
\end{exa}
We can describe Poisson structures locally using the splitting theorem due to Weinstein \cite{Weinstein83}, combined with the symplectic Darboux theorem \cite{McDuffSalamon98}.
\begin{thm}\label{thm:poissonsplitting} Let $(X,\pi)$ be a Poisson manifold and $x \in X$ with ${\rm rank}(\pi) = 2k$. Then there exists an open neighbourhood $U$ of $x$ which is Poisson diffeomorphic to a product $(N, \pi_N) \times (\R^{2k},\pi_0)$, where $\pi_N$ vanishes at $x$ and $\pi_0$ is the standard nondegenerate Poisson structure of \autoref{exa:standardnondeg}.
\end{thm}
We define a \emph{Poisson submanifold} of $\pi$ as a Poisson manifold $(N,\pi_N)$ equipped with an injective immersive Poisson map $i\colon N \hookrightarrow X$. We can also consider ideals which are compatible with the Poisson bracket.
\begin{defn} Let $(X,\pi)$ be a Poisson manifold. An ideal $I \subseteq C^\infty(X)$ is a \emph{Poisson ideal} for $\pi$ if $\{I,C^\infty(X)\}_\pi \subseteq I$.
\end{defn}
There are several equivalent characterizations of when a submanifold is Poisson. See for example \cite[Proposition 2.2]{FernandesMarcut15}. Recall from Section \ref{sec:poissonalgebroids} that a Poisson structure $\pi \in {\rm Poiss}(X)$ gives rise to a Lie algebroid called the Poisson algebroid $T^*_\pi X$ associated to $\pi$.
\begin{lem}\label{lem:poissonsubmanifold} Let $(X,\pi)$ be a Poisson manifold and $N \subseteq X$ a closed submanifold. Then the following are equivalent:
	\bi
	\item $N$ is a $\pi$-Poisson submanifold;
	\item $N$ is a $T^*_\pi X$-invariant submanifold;
	\item $I_N$ is a Poisson ideal for $\pi$;
	\item $\pi^\sharp(T_x X) \subseteq T_x N$ for all $x \in N$;
	\item $\mc{L}_{V_f} I_N \subseteq I_N$ for all Hamiltonian vector fields $V_f \in {\rm Ham}(\pi)$.
	\ei
\end{lem}
We see it makes sense to more generally define a subspace $N \subseteq X$ to be Poisson if it is $T^*_\pi X$-invariant (so if $N$ is closed, by demanding $I_N$ to be Poisson ideal for $\pi$).
\begin{rem} Other classes of submanifolds of Poisson manifolds exist, such as isotropic and coisotropic submanifolds, and Poisson transversals. We will not directly make use of these in this thesis.
\end{rem}
Using \autoref{thm:poissonsplitting}, one can show that any Poisson structure gives rise to a singular foliation by symplectic leaves. Let $\pi \in {\rm Poiss}(X)$ be given.
\begin{defn} A \emph{$\pi$-symplectic leaf} $L$ is a maximal integral submanifold of $\pi^\sharp$, i.e.\ a maximal path-connected submanifold of $X$ such that $T_x L = {\rm im}\, \pi^\sharp_x$ for all $x \in L$.
\end{defn}
The inclusion $i\colon L \hookrightarrow X$ of a symplectic leaf is Poisson, and the Poisson structure on $L$ is nondegenerate, explaining the name. We say that a subset $M \subseteq X$ is \emph{$\pi$-saturated} if it is a union of $\pi$-symplectic leaves. Then we have the following.
\begin{thm} Let $(X,\pi)$ be a Poisson manifold. Then each point $x \in X$ is contained in a unique $\pi$-symplectic leaf $L_x$ which contains all $\pi^\sharp$-integral submanifolds containing $x$. Moreover, both $X_{\pi,{\rm reg}}$ and $X_{\pi, {\rm sing}}$ are $\pi$-saturated subsets.
\end{thm}
The collection of $\pi$-symplectic leaves will be referred to as the symplectic foliation induced by $\pi$. In general this is a singular foliation as the leaves may vary in dimension. As $T^*X$ is locally generated by exact one-forms, the tangent spaces of $\pi$-symplectic leaves are spanned by the evaluation of all $\pi$-Hamiltonian vector fields.
\begin{exa} The trivial Poisson structure $\pi$ on $X$ has rank $0$ everywhere, and its symplectic leaves are the points of $X$.
\end{exa}
\begin{exa} Let $(X,\omega)$ be a symplectic manifold. Then $\pi$ has one symplectic leaf, namely $X$ itself.
\end{exa}
\begin{exa} Let $\mf{g}^*$ be the dual of finite dimensional Lie algebra, equipped with the linear Poisson structure from $\mf{g}$ as in \autoref{exa:liepoisson}. Its symplectic leaves are exactly the orbits of the coadjoint action of a Lie group $G$ with ${\rm Lie}(G) = \mf{g}$ on $\mf{g}^*$.
\end{exa}
\begin{exa} Let $\pi$ be a regular Poisson structure of rank $2m$. Then its symplectic foliation is regular, and all symplectic leaves have constant dimension $2m$.
\end{exa}
Most Poisson structures $\pi$ we will encounter in this thesis are generically symplectic. This means that $X_{\pi, {\rm reg}}$ is open and dense, and its connected components are $\pi$-symplectic leaves. Using the Poisson algebroid we see that in this case, $X_{\pi, {\rm reg}}$ is also the isomorphism locus of the Lie algebroid $T^*_\pi X$, as its anchor is given by $\pi^\sharp$.
\begin{defn}\label{defn:picasimir} Let $(X,\pi)$ be a Poisson manifold. A function $f \in C^\infty(X)$ is called a \emph{$\pi$-Casimir} if $\{f,\cdot\}_\pi = 0$.
\end{defn}
Hence, $f$ is a $\pi$-Casimir if and only if its $\pi$-Hamiltonian vector field is zero, $V_f \equiv 0$. These functions are constant along the leaves of $\pi$, as for any $\pi$-Hamiltonian vectof field $V_g$ tangent to a symplectic leaf $L$ we have $\mc{L}_{V_g} f = -\{f,g\}_\pi = 0$. A quick computation using the Schouten bracket shows that given a $\pi$-Casimir $f$, the bivector $\pi' := f \cdot \pi$ is a new Poisson structure on $X$.
\subsection{Poisson cohomology}
\label{sec:poissoncohomology}
In this section we briefly discuss Poisson(-Lichnerowicz) cohomology. For more information the reader is again directed to e.g.\ \cite{LaurentGengouxPichereauVanhaecke13,DufourZung05}. This takes place on $\mf{X}^\bullet(X)$ with differential $\delta_\pi\colon \mf{X}^\bullet(X) \to \mf{X}^{\bullet+1}(X)$ given by $\delta_\pi = [\pi, \cdot]$, which squares to zero if and only if $[\pi,\pi] = 0$, as follows from the Jacobi identity for the Schouten bracket.
\begin{defn} Let $\pi \in {\rm Poiss}(X)$. Then the \emph{Poisson cohomology} of $\pi$ is given by $H^k_\pi(X) := H^k(\mf{X}^\bullet(X), \delta_\pi)$ for $k \in \N \cup \{0\}$.
\end{defn}
\begin{rem} Straight from the definition, we see that $H^0_\pi(X)$ is the space of $\pi$-Casimirs, while $H^1_\pi(X)$ is the quotient ${\rm Poiss}(\pi) / {\rm Ham}(\pi)$.
\end{rem}
\begin{rem} As $[\pi,\pi] = 0$, we see that $\pi$ defines a class $\pi \in H^2_\pi(X)$.
\end{rem}
As with Lie algebroid cohomology (section \ref{sec:lacohomology}), Poisson cohomology is in general very hard to compute. When $\pi$ is either maximally degenerate or nondegenerate, this is not the case.
\begin{exa} Let $\pi \equiv 0$ be trivial. Then $H^\bullet_\pi(X) = \mf{X}^\bullet(X)$.
\end{exa}
\begin{exa}\label{exa:nondegpoissoncoh} Let $\pi$ be nondegenerate. Then $H_\pi^\bullet(X) \cong H^\bullet_{\rm dR}(X)$.
\end{exa}
As was mentioned as \autoref{prop:poissoncohlacoh}, Poisson cohomology of $\pi$ is isomorphic to the Lie algebroid cohomology of its associated Poisson algebroid $T^*_\pi X$: we have $H^\bullet_\pi(X) \cong H^\bullet(T^*_\pi X)$. There is always a map $H_{\rm dR}^\bullet(X) \to H^\bullet_\pi(X)$ from the anchor of $T^*_\pi X$ using \autoref{prop:lalgdmorphcoh}. As $\pi^\sharp\colon T^*_\pi X \to TX$ is an isomorphism if $\pi$ is nondegenerate and moreover satisfies $\pi^\sharp(d\alpha) = - [\pi, \pi^\sharp \alpha] = - \delta_\pi(\pi^\sharp \alpha)$ for $\alpha \in \Omega^\bullet(X)$, this is one way to see \autoref{exa:nondegpoissoncoh}. We will see in Section \ref{sec:poissonexamples} several interesting examples of Poisson structures whose Poisson cohomology can be computed.
\section{\texorpdfstring{$\mc{A}$}{A}-Poisson structures}
\label{sec:apoissonstrs}
In this section we define the notion of a Poisson structure in a Lie algebroid. We will refer to these as $\mc{A}$-Poisson structure.

Let $\mc{A} \to X$ be a Lie algebroid. The bracket $[\cdot,\cdot]_\mc{A}$ on $\Gamma(\mc{A})$ extends to an $\mathcal{A}$-Schouten bracket on the exterior algebra $\Gamma(\wedge^\bullet \mathcal{A})$, again denoted by $[\cdot,\cdot]_\mc{A}$.
Denote the space of $\mc{A}$-bivectors by $\mf{X}^2(\mc{A}) = \Gamma(\wedge^2 \mc{A})$ (\emph{not} $\Gamma(\wedge^2 T\mc{A})$), so that somewhat confusingly we have $\mf{X}^2(X) = \mf{X}^2(TX)$.
\begin{defn} An \emph{$\mathcal{A}$-Poisson structure} is an $\mathcal{A}$-bivector $\pi_\mathcal{A}\in \Gamma(\wedge^2\mathcal{A})$ with $[\pi_\mathcal{A},\pi_\mathcal{A}]_{\mc{A}} = 0$. Denote the space of such $\mc{A}$-bivectors by ${\rm Poiss}(\mc{A})$.
\end{defn}
We see that ${\rm Poiss}(X) = {\rm Poiss}(TX)$. Any $\mc{A}$-bivector $\pi_\mc{A}$ gives rise to a map $\pi_\mc{A}^\sharp\colon \mc{A}^* \to \mc{A}$. An $\mc{A}$-Poisson structure $\pi_{\mc{A}}$ is called \emph{nondegenerate} if $\pi_{\mc{A}}^\sharp\colon \mc{A}^* \to \mc{A}$ is an isomorphism, whose existence forces the rank of $\mc{A}$ to be even. Alternatively, one can demand the Pfaffian $\wedge^n \pi_\mc{A}$ to be nowhere vanishing. The notions of regular and singular points carry over from the case where $\mc{A} = TX$, as do their respective regular and singular loci.
\begin{rem} $\mc{A}$-Poisson structures are also considered in \cite{Lanius16}, especially for those Lie algebroids whose isomorphism locus is the complement of a hypersurface.
\end{rem}
\begin{rem} An $\mc{A}$-Poisson structure also goes under the name of a \emph{(strong) Hamiltonian operator} for the Lie bialgebroid $(\mc{A},\mc{A}^*,0)$, see \cite{LiuWeinsteinXu97}.
\end{rem}
As when $\mc{A} = TX$, for any Lie algebroid of even rank there is a bijection between $\mc{A}$-symplectic forms and nondegenerate $\mc{A}$-Poisson structures (compare with \autoref{exa:symppoisson}). Namely, given an $\mc{A}$-symplectic structure $\omega_{\mc{A}}$, nondegeneracy implies we can invert the map $\omega_{\mc{A}}^\flat\colon \mc{A} \to \mc{A}^*$ to $(\omega_{\mc{A}}^\flat)^{-1} = \pi_{\mc{A}}^\sharp\colon \mc{A}^* \to \mc{A}$ for an $\mc{A}$-Poisson structure $\pi_{\mc{A}}$. The conditions $d_{\mc{A}} \omega_{\mc{A}} = 0$ and $[\pi_{\mc{A}}, \pi_{\mc{A}}]_{\mc{A}} = 0$ are equivalent.
\begin{prop}\label{prop:asympapoisson} Let $\mc{A} \to X$ be a Lie algebroid of even rank. Then there is a bijection between $\mc{A}$-symplectic forms and nondegenerate $\mc{A}$-Poisson structures.
\end{prop}
For an $\mc{A}$-\symp{} $\omega_\mc{A}$, we call $\pi = \rho_{\mc{A}}(\pi_\mc{A})$ the \emph{dual bivector} to $\omega_{\mc{A}}$. In later chapters we are mainly interested in nondegenerate $\mc{A}$-Poisson structures, as we wish to use symplectic techniques. Indeed, the goal of the lifting process described in Section \ref{sec:liftingpoisson} is to find a Lie algebroid $\mc{A} \to X$ for which a given Poisson structure $\pi \in {\rm Poiss}(X)$ can be viewed as a nondegenerate $\mc{A}$-Poisson structure.
\begin{rem} If $\mc{A}$ has rank two, any $\mc{A}$-bivector is $\mc{A}$-Poisson, as $[\pi_\mc{A}, \pi_{\mc{A}}]_\mc{A} \in \mf{X}^3(\mc{A}) = 0$. In particular, a Lie algebroid $\mc{A}$ of rank two admits a nondegenerate $\mc{A}$-Poisson structure if and only if $\mc{A}$ is orientable (compare with \autoref{rem:asympranktwo}).
\end{rem}
An $\mc{A}$-Poisson structure $\pi_{\mc{A}}$ is $\emph{m-regular}$ for an integer $m \geq 0$ if $\pi_\mc{A}^\sharp$ has constant rank equal to $2m$. Alternatively, one can demand the existence of a line subbundle $K \subseteq \wedge^{2m} \mc{A}$ such that $\wedge^m \pi_\mc{A}$ is a nowhere vanishing section of $K$. Note that it is implied that $\wedge^{m+1} \pi_\mc{A} = 0$. Moreover, note that being nondegenerate is the same thing as being $n$-regular, where ${\rm rank}(\mc{A}) = 2n$. An $\mc{A}$-Poisson structure $\pi_\mc{A}$ gives rise to a (singular) $\mc{A}$-distribution $D_{\pi_\mc{A}} = \pi_\mc{A}^\sharp(\mc{A}^*)$. Regularity of $\pi_\mc{A}$ amounts to regularity of $D_{\pi_\mc{A}}$. Note that then $\det(D_\mc{A}) = K$, with $K \subseteq \wedge^{2m} \mc{A}$ as before.
\begin{defn} An \emph{$\mc{A}$-Poisson map} from $(X,\mc{A},\pi_{\mc{A}})$ to $(X',\mc{A}',\pi_{\mc{A}'})$ is a Lie algebroid morphism $(\varphi,f)\colon (X,\mc{A},\pi_{\mc{A}}) \to (X',\mc{A}',\pi_{\mc{A}'})$ such that $\varphi(\pi_{\mc{A}}) = \pi_{\mc{A}'}$.
\end{defn}
As for Poisson maps, it is not true that any $\mc{A}$-Poisson structure can be pushed forward to an $\mc{A}'$-Poisson structure along a given Lie algebroid morphism $(\varphi,f)\colon \mc{A} \to \mc{A}'$.
Given an $\mc{A}$-Poisson structure, we obtain its \emph{$\mc{A}$-Hamiltonian vector fields} as the image of the map $f \mapsto \pi_\mc{A}^\sharp(d_\mc{A} f) =: v_{\mc{A}, f} \in \mf{X}(\mc{A})$. Similarly we have the \emph{$\mc{A}$-Poisson vector fields}, i.e.\ those $v_\mc{A} \in \mf{X}(\mc{A})$ such that $\mc{L}_{v_\mc{A}} \pi_{\mc{A}} = 0$, using the $\mc{A}$-Lie derivative.

Define the spaces ${\rm Ham}(\pi_\mc{A}) \subseteq {\rm Poiss}(\pi_\mc{A}) \subseteq \mf{X}(\mc{A})$ as for Poisson structures. On functions we have that $d_\mc{A} = \rho_{\mc{A}}^* \circ d$. Hence $\pi_{\mc{A}}^\sharp(d_{\mc{A}} f) = \pi_{\mc{A}}^\sharp(\rho_{\mc{A}}^*(df))$. Setting $\pi := \rho_\mc{A}(\pi_\mc{A})$, we see from the relation $\pi^\sharp = \rho_{\mc{A}} \circ \pi_{\mc{A}}^\sharp \circ \rho_{\mc{A}}^*$ that $\rho_\mc{A}(V_{\mc{A},f}) = V_f$, i.e. an $\mc{A}$-Hamiltonian vector field $\pi_\mc{A}$ maps onto the Hamiltonian vector field of $\pi$.
\begin{prop} There is a surjection $\rho_\mc{A}\colon {\rm Ham}(\pi_\mc{A}) \to {\rm Ham}(\pi)$.
\end{prop}
%
%
\begin{rem}\label{rem:productapoisson} Let $(X,\mc{A},\pi_\mc{A})$ and $(X',\mc{A}',\pi_{\mc{A}'})$ be two $\mc{A}$-Poisson manifolds. Then their product $(X \times X', \mc{A} \oplus \mc{A}', \pi_{\mc{A}} + \pi_{\mc{A}'})$ is also $\mc{A}$-Poisson, where $\mc{A} \oplus \mc{A}' \to X \times X'$ is the direct product of Lie algebroids (see \cite{HigginsMackenzie90}).
\end{rem}
\subsection{$\mc{A}$-Poisson cohomology}
\label{sec:apoissoncoh}
We can repeat the definition of Poisson cohomology of Section \ref{sec:poissoncohomology} for $\mc{A}$-Poisson structures $\pi_\mc{A}$. Consider $\mf{X}^\bullet(\mc{A})$, equipped with the differential $\delta_{\pi_{\mc{A}}} = [\pi_\mc{A},\cdot]_\mc{A}$. This again squares to zero as $[\pi_\mc{A},\pi_\mc{A}]_\mc{A} = 0$.
\begin{defn} Let $\pi_\mc{A} \in {\rm Poiss}(\pi_\mc{A})$. Then the \emph{$\mc{A}$-Poisson cohomology} of $\pi_\mc{A}$ is given by $H^k_{\pi_\mc{A}}(\mc{A}) := H^k(\mf{X}^\bullet(\mc{A}), \delta_{\pi_\mc{A}})$ for $k \in \N \cup \{0\}$.
\end{defn}
As for usual Poisson structures (see \autoref{prop:poissoncohlacoh}), $\mc{A}$-Poisson cohomology is isomorphic to the Lie algebroid cohomology of the associated $\mc{A}$-Poisson algebroid.
\begin{prop}\label{prop:apoissoncoh} Let $\pi_\mc{A}$ be an $\mc{A}$-Poisson structure. Then $H^k_{\pi_\mc{A}}(\mc{A}) \cong H^\bullet(\mc{A}^*_{\pi_{\mc{A}}})$.
\end{prop}
If $\pi_\mc{A}$ is a nondegenerate $\mc{A}$-Poisson structure, we have that the $\mc{A}$-anchor of its $\mc{A}$-Poisson algebroid is an isomorphism, $\pi_\mc{A}^\sharp\colon \mc{A}^*_{\pi_\mc{A}} \stackrel{\cong}{\to} \mc{A}$. Consequently, we obtain as for Poisson structures that in this case we have $H^k_{\pi_{\mc{A}}}(\mc{A}) \cong H^k(\mc{A})$.
\subsection{Degeneraci loci}
In this section we discuss degeneraci loci of $\mc{A}$-Poisson structures $\pi_\mc{A}$, inspired by \cite{Pym13}. Considering what was done for Lie algebroids in Section \ref{sec:degenloci}, we could define these using the map $\pi_\mc{A}^\sharp$. However, $\pi^\sharp$ is skew-symmetric, so its degree of vanishing is too high (see \autoref{lem:detofabivector} below). Instead we use the Pfaffians $\wedge^k\pi_\mc{A} \in \Gamma(\wedge^{2k} \mc{A})$ for $k \geq 0$. These give maps $\wedge^k \pi_\mc{A} \Gamma(\wedge^{2k} \mc{A}^*) \to C^\infty(X)$, which define ideals $I_{\pi_\mc{A},2k}$.
\begin{defn} Let $\pi_\mc{A} \in {\rm Poiss}(\mc{A})$. The \emph{$2k$th degeneracy locus} of $\pi_\mc{A}$ is the subspace $X_{\pi_\mc{A},2k} \subseteq X$ defined by the ideal $I_{\pi_\mc{A},2k+2} \subseteq C^\infty(X)$.
\end{defn}
These degeneracy loci are the subspaces of $X$ where the rank of $\pi_{\mc{A}}$ is $2k$ or less. Tying things in with the $\mc{A}$-Poisson algebroid $\mc{A}^*_{\pi_\mc{A}}$, we see that the degeneraci loci of $\mc{A}_{\pi_\mc{A}}^*$ and $\pi_\mc{A}$ agree as subspaces of $X$, but their ideals do not.
\begin{prop}\label{prop:pidegenpipoisson} Let $\pi \in {\rm Poiss}(X)$. Each degeneracy locus $X_{\pi,2k}$ is a $\pi$-Poisson submanifold if it is smooth.
\end{prop}
\bp This follows from \autoref{lem:poissonsubmanifold} together with \autoref{prop:degenlocus}.
\ep
Consequently, each degeneracy locus has an induced Poisson structure.
\section{Poisson structures of divisor-type}
\label{sec:poissondivisors}
Divisors as defined in Chapter \ref{chap:divisors} provide a convenient way to define and study specific classes of Poisson structures. We first describe the general way in which divisors relate to Poisson geometry. We will immediately do so for the general situation of $\mc{A}$-Poisson structures. Let $\mc{A} \to X$ be a Lie algebroid with ${\rm rank}(\mc{A}) = 2n$.
\begin{defn} The \emph{Pfaffian} of $\pi_\mc{A}\in \mf{X}^2(\mc{A})$ is the section $\wedge^n \pi_\mc{A} \in \Gamma(\det(\mc{A}))$.
\end{defn}
There is a general way in which divisors give rise to interesting $\mc{A}$-Poisson structures by weakening the nondegeneracy or $m$-regularity conditions above.
\begin{defn}\label{defn:adivisortype} An $\mc{A}$-bivector $\pi_\mc{A}$ is said to be of \emph{divisor-type} if $(\det \mc{A}, \wedge^n \pi_\mc{A})$ is a divisor. It is of \emph{$m$-divisor-type} if $m \geq 0$ is such that $\wedge^{m+1} \pi_\mc{A} = 0$ and $\wedge^m \pi_\mc{A} \neq 0$, and there exists a subbundle $K \subseteq \wedge^{2m} \mc{A}$ such that $(K, \wedge^m \pi_{\mc{A}})$ is a divisor.
\end{defn}
We denote the divisor of $\pi_\mc{A}$ by ${\rm div}(\pi_\mc{A})$, and the associated divisor ideal by $I_{\pi_\mc{A}}$. Being of $n$-divisor-type is the same as being of divisor-type. Given a divisor ideal $I$, we denote by $\mf{X}^2_I(\mc{A}) \subseteq \mf{X}^2(\mc{A})$ the space of divisor $\mc{A}$-bivectors $\pi_\mc{A}$ such that $I_{\pi_\mc{A}} = I$. We denote by ${\rm Poiss}_I(\mc{A})$ the space of $\mc{A}$-Poisson structures of $I$-divisor-type, and by ${\rm Poiss}_{I,m}(\mc{A})$ the space of $\mc{A}$-Poisson structures of $m$-$I$-divisor-type.
\begin{rem} Interesting examples of Poisson structures of divisor-type can be found in Section \ref{sec:poissonexamples}.
\end{rem}
\begin{rem}\label{rem:mloglogf} For a Poisson structure $\pi \in {\rm Poiss}(X)$, being of $m$-log-divisor-type is the same as being a \emph{log-f Poisson structure} as in \cite{AndroulidakisZambon17}.
\end{rem}
\begin{rem} The class of $\mc{A}$-Poisson structures of divisor-type is closed under products (\autoref{rem:productapoisson}), using the external tensor product of divisors (see \autoref{defn:divdirectstum}). This is noted in \cite{AndroulidakisZambon17} for $\mc{A} = TX$ using \autoref{cor:divtypealmostreg} in the next section.
\end{rem}
Let $\pi_\mc{A}$ be an $\mc{A}$-bivector and let $\pi_\mc{A}^\sharp\colon \mc{A}^* \to \mc{A}$ be the associated bundle map. Then $\det(\pi_\mc{A}^\sharp)\colon \det(\mc{A}^*) \to \det(\mc{A})$, i.e.\ $\det(\pi_\mc{A}^\sharp) \in \Gamma(\det(\mc{A}) \otimes \det(\mc{A}))$, using that $\det(E)^* \cong (\det(E))^*$. On the other hand, we have $\wedge^n \pi_\mc{A} \in \Gamma(\det(\mc{A}))$.
\begin{lem}\label{lem:detofabivector} Let $\pi_\mc{A}$ be a $\mc{A}$-bivector. Then $\det(\pi_\mc{A}^\sharp) = \wedge^n \pi_\mc{A} \otimes \wedge^n \pi_\mc{A}$.
\end{lem}
In particular, if $\wedge^n \pi$ vanishes transversely, i.e.\ has a zeroes of order one, then $\det(\pi^\sharp)$ will have zeroes of order two. Given a Poisson structure $\pi$ of divisor-type and a $\pi$-Casimir function $f$, we can form a new Poisson structure $\pi' = f \cdot \pi$ (see the discussion below \autoref{defn:picasimir}). Assuming that the zero set of $f$ is nowhere dense and viewing $f$ as a section of the trivial bundle $\underline{\R}$, the realization that $\wedge^n \pi' = f^n \cdot \wedge^n \pi$ allows us to conclude that $\pi'$ is again of divisor-type. If $(\det(TX), \wedge^n \pi) \cong (U,\sigma)$ with divisor ideal $I$, we have that $\pi'$ is of $I'$-divisor-type, with $I' := \langle f^n \rangle \cdot I$ the ideal of the product $(U \otimes \underline{\R}^n, \sigma \otimes f^n)$. A similar discussion holds for Poisson structures of $m$-divisor-type. See also \cite[Lemma 3.3]{AndroulidakisZambon17} and the next section.
\section{Rigged algebroids and almost-regularity}
\label{sec:riggedalgebroids}
In this section we describe relations between Poisson structures of divisor-type and the notions of rigged algebroids \cite{Lanius16} and almost-regular Poisson structures \cite{AndroulidakisZambon17}.

For several Poisson structures that are almost-everywhere of maximal rank, it has been possible to compute their Poisson cohomology. This otherwise very complicated computation is made tractable by using the lifting procedure we describe in the next section. In particular, in \cite{Lanius16} Lanius uses a specific type of Lie algebroid called rigged algebroids for this purpose. We now define their $\mc{A}$-analogue.
\begin{defn} Let $\pi_\mc{A}$ be an $\mc{A}$-Poisson structure. The \emph{$\mc{A}$-rigged algebroid} $\mc{R}_{\pi_\mc{A}}$ of $\pi_\mc{A}$ is the unique Lie algebroid such that $\Gamma(\mc{R}_{\pi_\mc{A}}) = \pi^\#_\mc{A}(\Gamma(\mc{A}^*))$, if it exists.
\end{defn}
The $\mc{A}$-rigged algebroid will be an $\mc{A}$-Lie algebroid if it exists, as it is immediate from the definition that $\Gamma(\mc{R}_{\pi_\mc{A}}) \subseteq \Gamma(\mc{A})$, providing a natural $\mc{A}$-anchor. In particular, the rigged algebroid $\mc{R}_\pi$ of a Poisson structure $\pi \in {\rm Poiss}(X)$ makes sense, and it is this Lie algebroid that is used in \cite{Lanius16}. Indeed, realize that $H^\bullet(\mc{A}^*_{\pi_{\mc{A}}}) \cong H^\bullet(\mc{R}_{\pi_{\mc{A}}})$ where $\mc{A}^*_{\pi_\mc{A}}$ is the $\mc{A}$-Poisson algebroid associated to $\pi_\mc{A}$ (see \autoref{rem:apoissonalgebroids}). This together with the strategy outlined at the end of Section \ref{sec:aliealgebroids} can be used to compute $\mc{A}$-Poisson cohomology, and thus Poisson cohomology, in favorable cases.

In \cite{AndroulidakisZambon17}, Androulidakis--Zambon study a specific class of Poisson structures, inspired by trying to pinpoint which Poisson structures give rise to a smooth holonomy groupoid. Recall that given a Poisson structure $\pi$ on $X$, it gives rise to a sheaf $\mc{F} = \pi^\sharp(\Gamma(T^*X))$ of Hamiltonian vector fields of $\pi$, which is an involutive $C^\infty(X)$-submodule of $\mf{X}(X)$. It is not always a locally free sheaf, however.
\begin{defn}[{\cite[Definition 2.5]{AndroulidakisZambon17}}] A Poisson structure $\pi$ on $X$ is \emph{almost-regular} if there exists a vector bundle $E$ on $X$ such that $\Gamma(E) = \mc{F}$.
\end{defn}
The vector bundle $E$ above will immediately be a Lie algebroid if it exists. Using our language, by definition a Poisson structure is almost-regular if and only if its rigged algebroid exists. There is another characterization of almost-regularity. Denote by $X_{\pi, {\rm reg}} \subseteq X$ the open subspace where $\pi$ has maximal rank. It immediately follows that $X_{\pi, {\rm reg}} = (\wedge^k \pi)^{-1}(0)$, where $2k$ is the maximal rank of $\pi$ on $X$.
\begin{prop}[{\cite[Theorem 2.8]{AndroulidakisZambon17}}]\label{prop:almostregulardistr} Let $\pi$ be a Poisson structure on $X$. Then $\pi$ is almost-regular if and only if $X_{\pi, {\rm reg}}$ is dense and there exists a distribution $D \subseteq TX$ such that $D_x = T_x L$ for all $x \in X_{\pi, {\rm reg}}$, where $L$ is the $\pi$-symplectic leaf through $x$.
\end{prop}
The distribution $D$ in the above proposition is automatically unique and involutive, and integrates to a regular foliation by $\pi$-Poisson submanifolds. These three properties all follow by continuity, as they hold in $X_{\pi, {\rm reg}}$. Finally, there is the following characterization.
\begin{prop}[{\cite[Proposition 2.11]{AndroulidakisZambon17}}]\label{prop:almostregularline} Let $\pi$ be a Poisson structure on $X$, and let $2k$ be the maximal rank of $\pi$. Then $\pi$ is almost-regular if and only if $X_{\pi, {\rm reg}}$ is dense and there exists a real line bundle $K \subseteq \wedge^{2k} TX$ such that $\wedge^k \pi \in \Gamma(K)$.
\end{prop}
The distribution $D$ of \autoref{prop:almostregulardistr} is related to the line bundle $K$ of \autoref{prop:almostregularline} as $\wedge^k D = K$. Together with \autoref{defn:adivisortype} we obtain the following.
\begin{cor}\label{cor:divtypealmostreg} Let $\pi \in {\rm Poiss}(X)$ be given. Then the following are equivalent:
\bi
	\item $\pi$ is of $m$-divisor-type for some $m \geq 0$;
	\item $\pi$ is almost-regular;
	\item The rigged algebroid $\mc{R}_\pi$ exists.
\ei
\end{cor}
In Section \ref{sec:liftingpoisson} we will discuss the process of lifting Poisson structures to a Lie algebroid $\mc{A}$. Combining \autoref{cor:divtypealmostreg} with \autoref{prop:liftofdivtype} we see that if such a Poisson structure $\pi$ admits an $\mc{A}$-lift $\pi_\mc{A}$, the $\mc{A}$-rigged algebroid of $\pi_\mc{A}$ will also exist. Moreover, note that $\mc{R}_{\pi_{\mc{A}}} \cong \mc{A}$ if $\pi_{\mc{A}}$ is nondegenerate.
\begin{rem} Using \autoref{lem:detofabivector} it is not hard to show that if $\pi_\mc{A}$ is of $I$-divisor-type for some divisor ideal $I$, then $\mc{R}_{\pi_\mc{A}}$ is of $I^2$-divisor-type (using its $\mc{A}$-anchor).
\end{rem}
In joint work with Lanius \cite{KlaasseLanius17two} we explore the consequences of this in developing a general scheme to compute the Poisson cohomology of such Poisson structures using rigged algebroids and the process of lifting.

\section{Examples}
\label{sec:poissonexamples}
In this section we discuss some examples of Poisson structures of divisor-type.
\subsection{Log-Poisson structures}
\label{sec:logpoisson}
In this section we discuss Poisson structures built out of log divisors (Section \ref{sec:logdivisor}).
\begin{defn}\label{defn:logpoisson} Let $X$ be a $2n$-dimensional manifold. A \emph{log-Poisson structure} is a Poisson structure $\pi$ on $X$ that is of log divisor-type.
\end{defn}
While the name ``log-Poisson structure'' is most consistent with other Poisson structures built out of divisors, these Poisson structures also go under the name of $b$-Poisson, $b$-symplectic, and \emph{log-symplectic structures} \cite{Cavalcanti17,FrejlichMartinezTorresMiranda15,GualtieriLi14,GuilleminMirandaPires14,MarcutOsornoTorres14,MarcutOsornoTorres14two}. The latter two names emphasize the ability to view log-Poisson structures as suitably degenerate symplectic forms (namely as what we called log-symplectic structures in Section \ref{sec:logsympstr}, see \autoref{prop:blogbsymp}). This means the results of Secion \ref{sec:logsympstr} all apply to log-Poisson structures, and in particular one can consider the examples presented there. In this section however we will stick to using Poisson-theoretic techniques. There is an analogous theory in the holomorphic setting \cite{Goto02, Pym13, Pym15, Dongho12}. We say a log pair $(X,Z)$ \emph{admits a log-Poisson structure} if there exists a log-Poisson structure $\pi$ on $X$ such that $Z_{\wedge^n \pi} = Z$. We will also denote $Z_{\wedge^n \pi}$ as $Z_\pi$.
\begin{exa}\label{exa:logonr2} Let $X = \R^2$ with coordinates $(x,y)$ and consider $\pi = x \partial_x \wedge \partial_y$ with $Z_\pi = \{x = 0\}$. This is a log-Poisson structure on $(X,Z)$. Incidentally, this is also the Lie-Poisson structure given to the dual of the two-dimensional Lie algebra $\mf{g} = \langle e_1, e_2 \rangle$ with Lie bracket $[e_1,e_2] = e_1$ (see \autoref{exa:liepoisson}).
\end{exa}
\begin{exa} Let $X = S^2$ with standard coordinates $(h,\theta)$ and the Poisson bivector $\pi = h \partial_h \wedge \partial_\theta$, for which $Z_\pi$ is the equator $\{h = 0\}$.
\end{exa}
\begin{exa} More generally, Let $(X,\pi)$ be a nondegenerate oriented Poisson surface and let $f \in C^\infty(X)$ be a function with only transverse zeroes. Then $f \cdot \pi$ is a log-Poisson structure for the log pair $(X,f^{-1}(0))$. For example, one can take as singular locus any curve $\gamma \subseteq X$ which admits a global defining function.
\end{exa}
Log-Poisson structures can also exist on non-orientable manifolds. Indeed using the correspondence with log-symplectic structures (\autoref{prop:blogbsymp}), a two-dimensional log pair $(X,Z)$ admits a log-Poisson structure if and only if the associated log-tangent bundle is orientable. This in turn happens if and only if $w_1(TX) + {\rm PD}_{\Z_2}[Z] = 0$ (see \autoref{prop:logcharclass}). If $X$ is not orientable so that $w_1(TX) \neq 0$, this allows for connected components of $Z$ which have nontrivial normal bundle..
\begin{exa}[{\cite[Example 1.13]{GualtieriLi14}}] Let $X = \R P^2$ with $\pi = (g(x) - y^2) \partial_x \wedge \partial_y$, with $g(x) = x(x-1)(x-t)$ for $0 < t < 1$ a cubic polynomial. This extends smoothly from $\R^2$ to $\R P^2$ with singular locus $Z_\pi = \{y^2 = g(x)\}$ a real elliptic curve. This consists of two connected components, $Z_1$ and $Z_0$, with $Z_0$ containing $\{(0,0), (t,0)\}$ having trivial normal bundle, and $Z_1$ containing $\{(1,0), (\infty,0)\}$ with nontrivial normal bundle. Note that $X$ is not orientable.
\end{exa}
Note that any surface $\Sigma$, compact or not, orientable or not, admits a log-Poisson structure. This is because any vector bundle admits transverse sections, hence so does $\wedge^2 T\Sigma$; the Poisson condition is immediate in dimension two. However, not every log pair $(\Sigma,Z_\Sigma)$ admits a log-Poisson structure as we mentioned before, with an easy counterexample being $(\R P^2, \emptyset)$.

The rank of a log-Poisson structure $\pi$ is equal to $2n$ on $X \setminus Z$, and $2n-2$ on $Z$. Indeed, $\pi^n$ is nonvanishing on $X \setminus Z$, while on $Z$ we make use of the relation $d\pi^n = n d\pi \otimes \pi^{n-1}$ to conclude by transversality that $\pi^{n-1}|_Z$ is nonzero. We thus see that $X \setminus Z$ is the symplectic locus of $\pi$, while on its singular locus $Z$ there is a corank-one symplectic foliation given by the induced Poisson structure $\pi|_Z$. We can readily obtain the following local description of a log-Poisson structure.
\begin{prop}[\cite{GuilleminMirandaPires14}]\label{prop:loglocalform} Let $\pi \in {\rm Poiss}(X^{2n})$ be a log-Poisson structure. Then $\pi$ is locally equivalent around points in $Z_\pi$ to $z \partial_z \wedge \partial_y + \pi_0$ on $\R^2 \times \R^{2n-2}$, with $\pi_0$ the standard nondegenerate Poisson structure, and $Z_\pi = \{z = 0\}$
\end{prop}
\bp Using \autoref{thm:poissonsplitting} we need only show that given a log-Poisson structure $\pi'$ on $\R^2$, we can find coordinates expressing it as the one of \autoref{exa:logonr2}. On $\R^2$ we have $\pi' = f(x,y) \partial_x \wedge \partial_y$, where $f(x,y)$ vanishes transversally. We can then find new coordinates $(x' = f(x,y) ,y')$, so that $\pi' = x' \partial_{x'} \wedge \partial_{y'}$ as desired.
\ep
In terms of the local coordinates given by the above proposition, the induced Poisson structure on $Z$ is given by $\pi_Z = \pi_0$, which indeed has rank $2n-2$. With respect to the standard measure $\mu = dz \wedge dy \wedge dx_i$, the modular vector field of $\pi$ is given by $V_\mu = \partial_y$. Its restriction $V_Z := V_\mu|_Z$ is nonvanishing on $Z$ (and is determined up to Hamiltonian vector fields), and satisfies $V_Z \wedge \pi_Z^{n-1} \neq 0$, i.e.\ is a Poisson vector field transverse to the symplectic leaves of $\pi_Z$. We see further that $(Z,\pi_Z)$ is so-called \emph{unimodular}. Note that $X_\mu$ is also given by the ``singular Hamiltonian vector field'' $X_{d \log z} = \pi^\sharp(d \log z)$ for $z \in I_Z$.
\begin{rem} \autoref{prop:loglocalform} alludes to the lifting result (\autoref{prop:blogbsymp}) mentioned before. Namely, inverting the local form of $\pi$ we get $\pi^{-1} = d \log x \wedge dy + \omega_0$, which is the Darboux normal form of a \blog{} (see Section \ref{sec:logsympstr}). 
\end{rem}                                                                                       
Following up on this remark, $(\pi_Z, V_Z)$ corresponds to the induced cosymplectic structure $(\alpha,\beta)$ of the log-symplectic structure $\omega = \pi^{-1}$: the symplectic foliation of $\pi_Z$ is given by the kernel of $\alpha$ with the pullback of $\beta$, and $V_Z$ satisfies $\iota_{V_Z} \beta = 0$ and $\iota_{V_Z} \alpha = 1$ (see also \cite{GuilleminMirandaPires11}).
The following will be important later in Section \ref{sec:boundarymaps}.
\begin{prop}\label{prop:orlogseparating} Let $(X,Z)$ be an orientable log pair which admits a log-Poisson structure. Then $Z$ is separating.
\end{prop}
\bp Choose an orientation volume form ${\rm vol}_X$ and let $\pi$ be a log-Poisson structure on $(X,Z)$. Then $\wedge^n \pi = h {\rm vol}_X^{-1}$ for some $h \in C^\infty(X)$. It follows that $Z = h^{-1}(0)$ is coorientable, and that $Z$ separates $X \setminus Z$ into $X_+ \sqcup X_-$ using the sign of $h$.
\ep
\begin{exa} Consider the torus $X = T^2$ and let $Z$ be a single transverse arc. Then $X \setminus Z$ is connected, so that by \autoref{prop:orlogseparating} we see that $(X,Z)$ cannot admit a log-Poisson structure. This can also be explained using the dual picture, as the log-tangent bundle $\mc{A}_Z$ is not orientable (see \autoref{cor:orlogpairseparating}).
\end{exa}
We finish by discussing the Poisson cohomology of a log-Poisson structure, using the strategy briefly outlined in Section \ref{sec:riggedalgebroids} (see also \cite{Lanius16}). As the local description of a log-Poisson structure shows, its rigged algebroid is locally given by $\langle z \partial_z, z \partial_y, \partial_{x_i} \rangle$. This is the local description of the Lie algebroid $TX(-\log (Z,F))$ of Section \ref{sec:edgetangentbundle}, with $F$ the corank-one symplectic foliation on $Z$ induced by $\pi$. Using \autoref{prop:edgelacohomology} we conclude the following result.
\begin{prop}[\cite{GuilleminMirandaPires14, MarcutOsornoTorres14}] Let $\pi \in {\rm Poiss}(X)$ be a log-Poisson structure. Then
	\be
		H_\pi^k(X) \cong H^k(\mc{A}_{Z_\pi}) \cong H^k(X) \oplus H^{k-1}(Z_\pi).
	\ee
\end{prop}
\subsection{$m$-log Poisson structures}
Next we discuss Poisson structures that are of $m$-divisor type and define a log divisor. As mentioned in \autoref{rem:mloglogf}, these are called log-f Poisson structures in \cite{AndroulidakisZambon17}.
\begin{defn} An \emph{$m$-log Poisson structure} is a Poisson structure $\pi \in {\rm Poiss}(X)$ that is of $m$-log-divisor type.
\end{defn}
In other words, $\pi$ is regular of rank $2m$ outside of $Z_\pi = (\wedge^m \pi)^{-1}(0)$, along which it vanishes transversally, so that it has rank $2m-2$ on $Z_\pi$.
\begin{exa}\label{exa:mlogpoisson} Consider $X = \R^3$ with coordinates $(x,y,z)$. Then both $\pi = x \partial_x \wedge \partial_y$ and $\pi' = z \partial_x \wedge \partial_y$ are $1$-log Poisson. Indeed, given any function $f\colon \R^3 \to \R$ for which $0$ is a regular value, the bivector $\pi_f = f \partial_x \wedge \partial_y$ is $m$-log-Poisson (\cite[Example 3.8]{AndroulidakisZambon17}).
\end{exa}
We will see in Section \ref{sec:liftingmlogpoisson} that the behavior of $\pi$ depends on its modular vector field, which in the above examples (with respect to $\mu = dx \wedge dy \wedge dz$) is given by $V_\mu =\partial_y$ respectively $V_\mu = 0$. In future work we hope to study these Poisson structures in more detail.
\subsection{Elliptic Poisson structures}
\label{sec:ellpoissonstr}
We can also construct Poisson structures out of elliptic divisors (Section \ref{sec:elldivisor}), obtaining the notion of an elliptic Poisson structure. These are called Poisson structures of elliptic log-symplectic type in \cite{CavalcantiGualtieri15}.
\begin{defn}[{\cite[Definition 3.3]{CavalcantiGualtieri15}}] Let $X$ be a $2n$-dimensional manifold. An \emph{elliptic Poisson structure} is a Poisson structure $\pi$ on $X$ that is of elliptic divisor-type.
\end{defn}
Let $(X,|D|)$ be an \emph{elliptic pair}, i.e.\ a manifold equipped with an elliptic divisor $|D|$. An elliptic pair $(X,|D|)$ \emph{admits an elliptic Poisson structure} if there exists an elliptic Poisson structure $\pi$ on $X$ such that $|D_{\wedge^n \pi}| = |D|$. When there is no elliptic divisor structure on $D$, we say that $(X,D)$ admits an elliptic Poisson structure if there exists some elliptic divisor structure $|D|$ on $D$ such that $(X,|D|)$ admits an elliptic Poisson structure.
\begin{exa} Let $X = \R^2$ with coordinates $(x,y)$ and consider $\pi = (x^2 + y^2) \partial_x \wedge \partial_y$ with $D_\pi = \{x = y = 0\}$. This is an elliptic Poisson structure on $(X,D)$.
\end{exa}
We are mainly interested elliptic Poisson structures that are of a specific type. Namely, in \autoref{prop:elllogcorrespondence} we discuss that elliptic Poisson structures are the dual bivectors of the elliptic symplectic structures of Section \ref{sec:ellipticsymp}. In light of the link to \sgcs{}s mentioned there, we care mostly about elliptic Poisson structures dual to $\mc{A}_{|D|}$-symplectic structures with zero elliptic residue.

A thorough study of elliptic Poisson structures for their own sake has not been undertaken. In particular, their Poisson cohomology has not yet been computed. We address this in work with Melinda Lanius \cite{KlaasseLanius17two}.
\subsection{Other examples}
We can also consider Poisson structures of divisor-type for any of the divisors of Section \ref{sec:divexamples}. Mostly one has focused on (nondegenerate) $\mc{A}$-Poisson structures for some of the Lie algebroids from Section \ref{sec:laexamples}. Namely, by \autoref{prop:asympapoisson} these admit a description as $\mc{A}$-symplectic structures (see also Section \ref{sec:examplessymp}). It would be nice to have a more intrinsic description of these Poisson structures without direct use of the Lie algebroids they can be lifted to. 
\subsubsection{Nondegenerate scattering-Poisson structures}
These are considered by Lanius \cite{Lanius16}, see also Section \ref{sec:scatteringsymp}. Let $(X,Z)$ be a log pair and consider a nondegenerate $\mc{C}_Z$-Poisson structure $\pi_{\mc{C}_Z}$ with underlying Poisson structure $\pi$, where $\mc{C}_Z$ is the scattering tangent bundle of Section \ref{sec:scatteringtangentbundle}. Lanius \cite[Theorem 1.3]{Lanius16} computes the Poisson cohomology of $\pi$ by computing the Lie algebroid cohomology of its Poisson algebroid (see Section \ref{sec:poissonalgebroids}), or more precisely, of its rigged algebroid (see Section \ref{sec:riggedalgebroids}).

The symplectic foliation of a nondegenerate scattering-Poisson structure consists of $X \setminus Z$, and the points of $Z$. Indeed, in Section \ref{sec:scatteringsymp} we saw that there is an induced contact structure on $Z$.
\subsubsection{Nondegenerate $b^k$-Poisson structures}
These are considered by Scott \cite{Scott16}, see also Section \ref{sec:bksymp}. Using similar methods as for nondegenerate scattering-Poisson structures, Lanius \cite[Theorem 1.4]{Lanius16} further computes the Poisson cohomology of Poisson structures $\pi$ underlying nondegenerate $\mc{A}_Z^k$-Poisson structures, where $\mc{A}_Z^k$ is one of the $b^k$-tangent bundles of Section \ref{sec:bktangentbundle}.

The symplectic foliation of (the bivector underlying) a nondegenerate $b^k$-Poisson structure consists of $X \setminus Z$, and a codimension-one foliation on $Z$. This is the same as for log-Poisson structures, and is due to the induced cosymplectic structure on $Z$ (see Section \ref{sec:bksymp}). Their local form can be deduced from \autoref{prop:bksymplocalform}.
\begin{prop}[\cite{Scott16}]\label{prop:bklocalform} Let $\pi \in {\rm Poiss}(X^{2n})$ be the underlying bivector of a nondegenerate $b^k$-Poisson structure with respect to a $(k-1)$-jet $j_{k-1} \in J^{k-1}_Z$. Then $\pi$ is locally equivalent around points in $Z_\pi$ to $z^k \partial_z \wedge \partial_y + \pi_0$ on $\R^2 \times \R^{2n-2}$, with $\pi_0$ the standard nondegenerate Poisson structure, and $Z_\pi = \{z = 0\}$ with $z \in j_{k-1}$.
\end{prop}

\section{Lifting Poisson structures}
\label{sec:liftingpoisson}
In this section we discuss the process of lifting $\mc{A}$-Poisson structures through Lie algebroid morphisms. For now, fix a Lie algebroid $\mc{A} \to X$.
\begin{defn} A bivector $\pi \in \mf{X}^2(X)$ is of \emph{$\mc{A}$-type} if there exists $\pi_\mc{A} \in \mf{X}^2(\mc{A})$ such that $\rho_\mc{A}(\pi_\mc{A}) = \pi$, and $\pi_\mc{A}$ is called an \emph{$\mc{A}$-lift} of $\pi$.
\end{defn}
The space of $\mc{A}$-liftable bivectors is given by $\mf{X}_\mc{A}^2(X) := \rho_\mc{A}(\mf{X}^2(\mc{A})) \subset \mf{X}^2(X)$. The space of $\mc{A}$-liftable Poisson structures is denoted by ${\rm Poiss}_{\mc{A}}(X) \subseteq {\rm Poiss}(X)$. We say $\pi$ is \emph{of nondegenerate $\mc{A}$-type} if it admits a nondegenerate $\mc{A}$-lift $\pi_{\mc{A}}$, and similarly for $m$-regularity, or even for being of divisor-type. For example, a Poisson structure is of log-$\mc{A}$-type if it admits an $\mc{A}$-lift which is of log divisor-type.
\begin{rem} In \cite[Definition 2.16]{Lanius16} this notion is considered using slightly different terminology: there $\pi \in {\rm Poiss}(X)$ is called $\mc{A}$-Poisson if it is of $\mc{A}$-type.
\end{rem}
\begin{rem} Our notion of nondegenerate $b^k$-type is almost the same as the notion of $b^k$-type in \cite{Scott16} due to Scott. There is one difference, namely that there is a class of (non-canonically) isomorphic Lie algebroids, namely $\mc{A}_Z^k$ with respect to different jet data, and Scott's $b^k$-type demands liftability to any one of these Lie algebroids.
\end{rem}
An $\mc{A}$-lift is in general not unique if it exists (consider $\pi = 0$ and $\mc{A}$ with trivial anchor). However, it is unique if the isomorphism locus $X_\mc{A}$ is dense, as will often be the case for us. Moreover, when $X_{\mc{A}}$ is dense, any $\mc{A}$-lift of a Poisson structure is automatically $\mc{A}$-Poisson. All ideal Lie algebroids built using divisor ideals are of this type (see Section \ref{sec:idealliealgebroids}). More generally this holds for almost-injective Lie algebroids $\mc{A}$, which also includes all involutive subbundles of $TX$. Note that we have $\pi^\sharp = \rho_{\mc{A}} \circ \pi_{\mc{A}}^\sharp \circ \rho_{\mc{A}}^*$ as maps, summarized in the following commutative diagram.
\begin{center}
	\begin{tikzpicture}
	\matrix (m) [matrix of math nodes, row sep=2.5em, column sep=2.5em,text height=1.5ex, text depth=0.25ex]
	{	\mc{A}^* & \mc{A} \\ T^*X & TX \\};
	\path[-stealth]
	(m-1-1) edge node [above] {$\pi_\mc{A}^\sharp$} (m-1-2)
	(m-2-1) edge node [left] {$\rho_\mc{A}^*$} (m-1-1)
	(m-2-1) edge node [above] {$\pi^\sharp$} (m-2-2)
	(m-1-2) edge node [right] {$\rho_\mc{A}$} (m-2-2);
	\end{tikzpicture}
\end{center}
\begin{prop}\label{prop:liftofdivtype} Let $\pi$ be a bivector of $m$-divisor-type that is also of $\mc{A}$-type. Then any $\mc{A}$-lift of $\pi$ is of $m$-divisor-type.
\end{prop}
\bp Let $\pi_\mc{A}$ be an $\mc{A}$-lift of $\pi$, so that it satisfies $\rho_\mc{A}(\pi_\mc{A}) = \pi$, or more precisely $(\wedge^2 \rho_\mc{A})(\pi_\mc{A}) = \pi$. The divisor associated to $\pi$ is given by $(K, \wedge^m \pi)$ for some line bundle $K \subseteq \wedge^{2m} TX$. As $\wedge^{2m} \rho_\mc{A}\colon \wedge^{2m} \mc{A} \to \wedge^{2m} TX$, we find a line bundle $K_\mc{A} \subseteq \wedge^{2m} \mc{A}$ for which $\wedge^{2m}\rho_\mc{A}\colon K_\mc{A} \to K$ and $\wedge^m \pi_\mc{A} \in \Gamma(K_\mc{A})$. If $\pi_\mc{A}$ vanishes then so must $\pi$, so that it is immediate that $(K_\mc{A}, \wedge^m \pi_\mc{A})$ defines a divisor, as any subset of a nowhere dense set is nowhere dense.
\ep
\begin{rem}\label{rem:poissonliftdivisors} More precisely, we have that ${\rm div}(\pi) = {\rm div}(\pi_\mc{A}) \otimes {\rm div}(\mc{A})$ using the product of divisors. In terms of divisor ideals, we thus have $I_\pi = I_{\pi_\mc{A}} \cdot I_\mc{A}$.
\end{rem}
Next we generalize to discuss the process of lifting an $\mc{A}'$-Poisson structure $\pi_{\mc{A}'}$ to a Lie algebroid $\mc{A} \to X$ which has a Lie algebroid morphism $(\varphi,{\rm id}_X)$ to another Lie algebroid $\mc{A}' \to X$. Typically, $\mc{A}$ will be a rescaling of $\mc{A}'$. As was mentioned at the start of this section, the aim of lifting a Poisson structure is to make it closer to being nondegenerate. This is done by incorporating the degeneracy of $\pi_{\mc{A}'}$ in the degeneracy of the Lie algebroid morphism $\varphi$.
\begin{defn}\label{defn:poissonalift} Let $(\varphi,{\rm id}_X)\colon \mathcal{A} \to \mathcal{A}'$ be a Lie algebroid morphism. An $\mc{A}'$-bivector $\pi_{\mathcal{A}'} \in \Gamma(\wedge^2 \mc{A}')$ is \emph{$\mc{A}$-liftable} if there exists an $\mc{A}$-bivector $\pi_{\mc{A}} \in \Gamma(\wedge^2 \mc{A})$ such that $\varphi(\pi_{\mc{A}}) = \pi_{\mc{A}'}$.
\end{defn}
In the above situation we call $\pi_{\mc{A}}$ the \emph{$\mc{A}$-lift} of $\pi_{\mc{A}'}$, and say that $\pi_{\mc{A}'}$ is of \emph{$\mc{A}$-type}. Note that we have $\pi_{\mc{A}'}^\sharp = \varphi \circ \pi_\mc{A}^\sharp \circ \varphi^*$ as maps, and again a commutative diagram.
\begin{center}
	\begin{tikzpicture}
	\matrix (m) [matrix of math nodes, row sep=2.5em, column sep=2.5em,text height=1.5ex, text depth=0.25ex]
	{	\mc{A}^* & \mc{A} \\ \mc{A}'{}^* & \mc{A}' \\};
	\path[-stealth]
	(m-1-1) edge node [above] {$\pi_\mc{A}^\sharp$} (m-1-2)
	(m-2-1) edge node [left] {$\varphi^*$} (m-1-1)
	(m-2-1) edge node [above] {$\pi_{\mc{A}'}^\sharp$} (m-2-2)
	(m-1-2) edge node [right] {$\varphi$} (m-2-2);
	\end{tikzpicture}
\end{center}
\begin{rem} For $\mathcal{A}' = TX$ and $\varphi = \rho_{\mathcal{A}}$, we get the notion of $\mathcal{A}$-type used in \cite{Lanius16}.
\end{rem}
It is not immediate that an $\mc{A}$-lift is itself Poisson. However, when it is, the lifting condition is exactly that $(\varphi,{\rm id}_X)\colon (\mc{A},\pi_\mc{A}) \to (\mc{A}',\pi_{\mc{A}'})$ is an $\mc{A}$-Poisson map. If $\pi_{\mc{A}'}$ has a nondegenerate $\mc{A}$-Poisson lift $\pi_\mc{A}$, we say $\pi_{\mc{A}'}$ is of \emph{nondegenerate $\mc{A}$-type}, and similarly for other types of $\mc{A}$-Poisson structures. We often omit the Lie algebroid morphism $\varphi$ in our notation. We will always use the anchor map $\rho_{\mc{A}}\colon \mc{A} \to TX$ when considering $\mc{A}$-lifts of bivectors on the Lie algebroid $TX$. The following is a consequence of density of isomorphism loci and continuity.
\begin{prop}\label{prop:uniquelifts} Let $(\varphi,{\rm id}_X)\colon \mathcal{A} \to \mathcal{A}'$ be a Lie algebroid morphism with dense isomorphism locus. Then $\mc{A}$-lifts are unique if they exist. Moreover, if $\pi_{\mc{A}'}$ is $\mc{A}'$-Poisson, then any $\mc{A}$-lift will be $\mc{A}$-Poisson.
\end{prop}
This also holds if $\varphi$ is injective on sections, i.e.\ $\varphi$ is almost-injective. Certainly if no assumptions on $\varphi$ are made, the above proposition is false. As before, consider for example $\mc{A}' = TX$ and $\mc{A}$ having trivial anchor. Then any section of $\wedge^2 \mc{A}$ is an $\mc{A}$-lift of $\pi_{\mc{A}'} = 0$. Essentially, what is used is the fact that the map on sheaves of sections is injective, which is equivalent to having dense isomorphism locus (when the ranks of $\mc{A}$ and $\mc{A}'$ agree).
\begin{rem}\label{rem:apoissonliftdivisors} Inspecting \autoref{prop:liftofdivtype} we see immediately that the same conclusion is true when lifting bivectors between Lie algebroids. Moreover, consider a Lie algebroid morphism $(\varphi,{\rm id}_X)\colon \mathcal{A} \to \mathcal{A}'$ with dense isomorphism locus and an $\mc{A}'$-Poisson structure $\pi_{\mc{A}'}$ of $\mc{A}$-type. Then as in \autoref{rem:poissonliftdivisors} we have ${\rm div}(\pi_{\mc{A}'}) = {\rm div}(\pi_\mc{A}) \otimes {\rm div}(\varphi)$ and $I_{\pi_{\mc{A}'}} = I_{\pi_\mc{A}} \cdot I_\varphi$. In this sense the degeneracy of $\pi_{\mc{A}'}$ is partially absorbed by $\varphi$.
\end{rem}
\begin{rem} Using \autoref{prop:uniquelifts}, Lie algebroids $\mc{A}$ with dense isomorphism loci provide unique $\mc{A}$-Poisson lifts of Poisson structures on $TX$, if they exist.
\end{rem}
We can iterate the lifting procedure as follows, whose proof is immediate.
\begin{prop}\label{prop:liftiterate} Let $(\varphi,{\rm id}_X)\colon \mc{A} \to \mc{A}'$ and $(\varphi',{\rm id}_X)\colon \mc{A}' \to \mc{A}''$ be Lie algebroid morphisms, and $\pi_{\mc{A}''} \in \rm{Poiss}(\mc{A}'')$. Then if $\pi_{\mc{A}''}$ is of $\mc{A}$-type, it is of $\mc{A}'$-type. If $\pi_\mc{A}$ is an $\mc{A}$-Poisson lift of $\pi_{\mc{A}''}$, then $\pi_{\mc{A}'} := \varphi(\pi_\mc{A})$ is an $\mc{A}'$-Poisson lift of $\pi_{\mc{A}''}$ which itself is of $\mc{A}$-type.
\end{prop}
Consequently, an $\mc{A}'$-Poisson structure $\pi_{\mathcal{A}'}$ being of $\mathcal{A}$-type is the same thing as its underlying Poisson structure $\pi := \rho_{\mathcal{A}'}(\pi_{\mathcal{A}'})$ being of $\mathcal{A}$-type. Let $\pi_{\mc{A}'}$ be a nondegenerate $\mc{A}'$-Poisson structure. Then any $\mc{A}$-Poisson map $(\varphi,f)\colon (\mc{A},\pi_\mc{A}) \to (\mc{A}',\pi_{\mc{A}'})$ has to be a Lie algebroid submersion, i.e.\ $\varphi$ must be fiberwise surjective. When $\mc{A}$ and $\mc{A}'$ are Lie algebroids of the same rank over $X$, this implies the following.
\begin{prop} Let $\pi \in {\rm Poiss}(X)$ be of nondegenerate $\mc{A}' $-type and consider  $(\varphi,{\rm id}_X)\colon \mc{A} \to \mc{A}'$ a Lie algebroid morphism. Assume that ${\rm rank}(\mc{A}) = {\rm rank}(\mc{A}')$. Then the following are equivalent:
	\bi
	\item $\varphi$ is an isomorphism;
	\item $\pi$ is of $\mc{A}$-type;
	\item $\pi$ is of nondegenerate $\mc{A}$-type.
	\ei
\end{prop}
As a consequence, once one has lifted $\pi$ to being nondegenerate, one cannot (meaningfully) lift further. This makes sense, as the lifting process is meant to desingularize the Poisson structure, and nondegenerate ones are maximally nonsingular. This is consistent with \autoref{rem:apoissonliftdivisors}, as isomorphisms specify trivial divisors. Moreover, we see that $\pi_{\mc{A}'}$ is of nondegenerate $\mc{A}$-type if and only if $\pi = \rho_{\mc{A}'}(\pi_{\mc{A}'})$ is.
\begin{rem} It is possible for a Poisson structure to nondegenerately lift to multiple non-isomorphic Lie algebroids, as it can for example happen that different Lie algebroids specify the same divisor. The consequences of this are explored in \cite{KlaasseLanius17two}.
\end{rem}
\begin{exa}\label{exa:mregularapoisson} Let $\pi_\mc{A}$ be an $m$-regular $\mc{A}$-Poisson structure. Then its image $D_{\pi_\mc{A}} = \pi_\mc{A}^\sharp(\mc{A}^*)$ is a regular $\mc{A}$-distribution. By \autoref{exa:adistribution}, $D_{\pi_{\mc{A}}}$ is a Lie subalgebroid of $\mc{A}$ and thus carries an $\mc{A}$-Lie algebroid structure $\varphi\colon D_{\pi_\mc{A}} \to \mc{A}$. It is immediate that $\pi_{\mc{A}}$ is of nondegenerate $D_{\pi_\mc{A}}$-type.
\end{exa}
We next pose a question to the reader, which we can possibly settle upon establishing the splitting results for $\mc{A}$-Poisson structures mentioned in Chapter \ref{chap:splittingtheorems}, or by using the language of jets of sections.
\begin{que} Let $\pi \in {\rm Poiss}(X)$ be of nondegenerate $\mc{A}_Z^2$-type for a certain choice of $1$-jet $j_1 \in J_Z^1$ at $Z$. In other words, let $\pi$ be the dual bivector of an $\mc{A}_Z^2$-symplectic structure. Using the natural morphisms $\mc{A}_Z^2 \to \mc{A}_Z \to TX$ (see Section \ref{sec:bktangentbundle}), we know that $\pi$ is of $\mc{A}_Z$-type by \autoref{prop:liftiterate} with $\mc{A}_Z$-lift $\pi_{\mc{A}_Z}$. Is $\pi$ in fact of log-$\mc{A}_Z$-type, with $Z(\pi_{\mc{A}_Z}) = Z$? Similarly for the higher $b^k$-bundles $\mc{A}_Z^k$.
\end{que}
\subsection{Lifting Poisson structures of divisor-type}
In this section we discuss the process of lifting Poisson structures of divisor-type. Let $(U,\sigma)$ be a divisor with divisor ideal $I$, and assume that the associated ideal Lie algebroid $\mc{A}_I$ exists. Denote the space of $\mc{A}_I$-liftable divisor-bivectors by $\mf{X}^2_{\mc{A}_I,I}(X) = \mf{X}^2_{\mc{A}_I}(X) \cap \mf{X}^2_{I}(X)$. Similarly, define ${\rm Poiss}_{\mc{A}_I, I}(X) = {\rm Poiss}_{\mc{A}_I}(X) \cap {\rm Poiss}_{I}(X)$.

In general, divisor bivectors are not liftable to the ideal Lie algebroid they define, i.e.\ $\mf{X}^2_I(X) \not\subseteq \mf{X}^2_{\mc{A}_I}(X)$. The other inclusion is also false. For example, a bivector $\pi$ of log divisor-type, i.e. such that $\wedge^n \pi$ is transverse, with $Z_\pi = (\wedge^n \pi)^{-1}(0)$, does not necessarily lift to the log-tangent bundle $\mc{A}_{Z_\pi}$ it defines.
\begin{exa} Consider the bivector $\pi = x \partial_x \wedge \partial_y + \partial_z \wedge \partial_w + \partial_x \wedge \partial_w$ on $\R^4$ with coordinates $(x,y,z,w)$. Note that $\wedge^2 \pi = x \partial_x \wedge \partial_y \wedge \partial_z \wedge \partial_w$, which vanishes transversally on $Z_\pi = \{x = 0\}$. Hence $\pi$ is of log divisor-type, but does not lift due to the presence of $\partial_x \wedge \partial_w$. However, note that $\pi$ is not Poisson.
\end{exa}
\begin{exa} Consider the bivector $\pi = x^2 \partial_x \wedge \partial_y$ on $\R^2$ with coordinates $(x,y)$. Then $\pi$ lifts to $\mc{A}_Z$ where $Z = \{x = 0\}$ with $\mc{A}_Z$-lift $\pi_{\mc{A}_Z} = x (x \partial_x) \wedge \partial_y$, but it does not lift nondegenerately, nor does $\pi$ specify a log divisor structure on $Z$. In other words, we see that $\mf{X}^2_{\mc{A}_I}(X) \not\subseteq \mf{X}^2_I(X)$ in this case.
\end{exa}
The Poisson condition is important for liftability, as one sees by comparing the results in Section \ref{sec:liftinglogpoisson} and beyond with the above example. Namely, we will show that ${\rm Poiss}_I(X)$ is contained in ${\rm Poiss}_{\mc{A}_I}(X)$ when $I$ is a log divisor ideal.
We first discuss which vector fields are liftable to $\mc{A}_I$, when $\pi \in \mf{X}^2_I(X)$.
\begin{prop} A vector field $V \in \mf{X}(X)$ lies in $\mf{X}_{\mc{A}_I}(X)$ if $\mc{L}_V (\wedge^n \pi) = 0$.
\end{prop}
\bp We denote $\wedge^n \pi$ by $\pi^n$. By definition of $\Gamma(\mc{A}_I) = \mc{V}_X(I)$, we must show that $\mc{L}_V I \subset I$. Let ${\rm vol}_X$ be a local volume form. Define $f := \pi^n({\rm vol}_X)$ so that $I = \langle f \rangle$. Then $\mc{L}_V f = \mc{L}_V \pi^n({\rm vol}_X) = (\mc{L}_V \pi^n)({\rm vol}_X) + \pi^n \mc{L}_V \rm{vol}_X$. Note that $\mc{L}_V {\rm vol}_X$ is again a volume form on a small enough open where $V$ is nonzero (if $V$ vanishes it clearly lifts). Hence $\mc{L}_V {\rm vol}_X = g {\rm vol}_X$ for some function $g$. We conclude that $\mc{L}_V f= \mc{L}_V \pi^n ({\rm vol}_X) + \pi^n(g {\rm vol}_X) = \mc{L}_V \pi^n ({\rm vol}_X) + g \pi^n({\rm vol}_X) = \mc{L}_V \pi^n ({\rm vol}_X) + g f$. The last term lies in $I$, so that if $\mc{L}_V \pi^n = 0$ we see that $V$ lifts.
\ep
The assumption in the previous proposition is quite strong, but does allow for the following conclusion.
\begin{cor}\label{cor:poissonhamiltonian} Let $\pi \in {\rm Poiss}_I(X)$. Then ${\rm Ham}_\pi(X) \subseteq \mf{X}_{\mc{A}_I}(X)$.
\end{cor}
\bp By the previous proposition it suffices to show that all $\pi$-Hamiltonian vector fields $V_f$ preserve $\pi^n$. Using \autoref{prop:hamispoiss} we know that $V$ is $\pi$-Poisson, i.e. $\mc{L}_{V_f} \pi = 0$. But then $\mc{L}_{V_f} \pi^n = n (\mc{L}_V \pi) \wedge \pi^{n-1} = 0$, so that $V_f \in \mf{X}_{\mc{A}_I}(X)$.
\ep
Due to this, we see that for $\pi \in {\rm Poiss}_I(X)$ there exists a map $\wt{\pi}^\sharp\colon T^*X \to \mc{A}_I$ fitting in the following diagram, which is not yet the existence of an $\mc{A}_I$-lift.
\begin{center}
	\begin{tikzpicture}
	\matrix (m) [matrix of math nodes, row sep=2.5em, column sep=2.5em,text height=1.5ex, text depth=0.25ex]
	{	\mc{A}_I^* & \mc{A}_I \\ T^*X & TX \\};
	\path[-stealth]
	(m-2-1) edge node [left] {$\rho_{\mc{A}}^*$} (m-1-1)
	(m-2-1) edge node [above] {$\pi^\sharp$} (m-2-2)
	(m-2-1) edge node [above] {$\wt{\pi}^\sharp$} (m-1-2)
	(m-1-2) edge node [right] {$\rho_\mc{A}$} (m-2-2);
	\draw[dotted, ->] (m-1-1) to (m-1-2);
	\end{tikzpicture}
\end{center}
However, we can now dualize the bundle morphism $\wt{\pi}^\sharp$ to a map $(\wt{\pi}^\sharp)^*\colon \mc{A}_I^* \to TX$. To test whether this lifts to a map $\pi_{\mc{A}_I}^\sharp\colon \mc{A}_I^* \to \mc{A}_I$, we can proceed as in \autoref{cor:poissonhamiltonian} by checking whether $(\wt{\pi}^\sharp)^*$ maps to $\mc{V}_I(X) \subseteq \mc{V}_X = \Gamma(TX)$ on sections. Assuming that $\Gamma(\mc{A}_I^*)$ admits local bases of closed sections, we can by continuity answer this in the isomorphism locus $X \setminus Z_I$, where $(\wt{\pi}^\sharp)^* = (\pi^\sharp)^*$ using the isomorphism given by $\rho_{\mc{A}_I}$. Here the lifting property follows as in \autoref{cor:poissonhamiltonian}, using that $(\pi^\sharp)^* = - \pi$ by skew-symmetry. In conclusion we obtain the following.
\begin{prop}\label{prop:ipoissonaitype} Let $I$ be a divisor ideal on $X$ for which $\mc{A}_I$ exists and admits local bases of closed sections. Then ${\rm Poiss}_I(X) \subseteq {\rm Poiss}_{\mc{A}_I}(X)$.
\end{prop}
Hence, in this case, any Poisson structure of $I$-divisor-type is of $\mc{A}_I$-type. There is a converse to this proposition which we address in the nondegenerate case. In light of \autoref{rem:poissonliftdivisors} we must consider to what extent ${\rm div}(\mc{A}_I)$ is related to $I$. Call a divisor ideal $I$ associated to a divisor $(U,\sigma)$ \emph{standard} if its ideal Lie algebroid $\mc{A}_I$ exists and satisfies ${\rm div}(\mc{A}_I) = (U,\sigma)$, i.e.\ $I_{\mc{A}_I} = I$. For example, by \autoref{prop:logtangentlogdiv} and \autoref{prop:elltangentlogdiv}, both log and elliptic divisors are standard.
\begin{rem}\label{rem:nonstandard} There are many divisor ideals which are not standard. For example, using \autoref{prop:locprincidealmodule} we obtain that $\mc{V}_X(I^k) = \mc{V}_X(I)$ for all $k \geq 0$. Consequently, if $I$ is a divisor ideal for which $\mc{A}_I$ exists, we have that $\mc{A}_{I^k} = \mc{A}_I$. Hence, if $I$ is standard, we get that $I_{\mc{A}_{I^k}} = I_{\mc{A}_I} = I$, which is not equal to $I^k$ unless $I$ is trivial.
\end{rem}
\begin{prop}\label{prop:standarddivideal} Let $I$ be a standard divisor ideal with divisor $(U,\sigma)$. Let $\pi \in {\rm Poiss}_{\mc{A}_I}(X)$ be of nondegenerate $\mc{A}_I$-type. Then $\pi \in {\rm Poiss}_I(X)$.
\end{prop}
\bp Let $\pi_{\mc{A}_I}$ be a nondegenerate $\mc{A}_I$-lift of $\pi$ and let $2n = \dim X = {\rm rank} \mc{A}_I$. Then $\wedge^n \pi_{\mc{A}_I}$ is a nonzero, hence $(\det(\mc{A}_I), \wedge^n \pi_{\mc{A}_I})$ is isomorphic to the trivial divisor $(\underline{\R}, \underline{1})$. But the map $\det \rho_{\mc{A}_I}$ sends $(\det(\mc{A}_I), \wedge^n \pi_{\mc{A}_I})$ to $(\det(TX), \wedge^n \pi)$. As $I$ is standard, we conclude that $(\det(TX), \wedge^n \pi) = (U,\sigma)$ as desired.
\ep
\begin{rem} An equivalent proof of \autoref{prop:standarddivideal} using the language of divisors goes as follows: as $\pi_{\mc{A}_I}$ is nondegenerate, it specifies the trivial divisor. Using \autoref{rem:poissonliftdivisors} we obtain that ${\rm div}(\mc{A}_I) = {\rm div}(\pi)$. However, by assumption $I_{\mc{A}_I} = I$, so that we conclude that $\pi \in {\rm Poiss}_I(X)$.
\end{rem}
\begin{rem} Following up on \autoref{rem:nonstandard}, let $I$ be a nontrivial divisor ideal for which $\mc{A}_I$ exists. A given $\pi \in {\rm Poiss}_{I^k}(X)$ for some $k > 1$, while it can be of $\mc{A}_I$-type using \autoref{prop:ipoissonaitype}, can never be of nondegenerate $\mc{A}_I$-type.
\end{rem}
In other words, Poisson structures liftable to standard divisor ideal Lie algebroids need not be of that divisor-type, but those which are nondegenerately liftable are. Moreover, if a Poisson structure of standard divisor-type lifts to that primary ideal Lie algebroid, it does so nondegenerately.
We turn to a related question. Note that if $\pi$ is of $\mc{A}_I$-type, the zero set $Z_I$ must be a $\pi$-Poisson subset, i.e.\ $I_{Z_I}$ is a $\pi$-Poisson ideal. There is a similar statement for $I$ (recall that $I \subseteq I_{Z_I}$ by \autoref{prop:divvanishingideal}).
\begin{prop} Let $\pi \in {\rm Poiss}_I(X)$. Then $I$ is a $\pi$-Poisson ideal.
\end{prop}
\bp We show this on generators. Let ${\rm vol}_X$ be a local volume form so that $I = \langle f \rangle$ for $f = \pi^n({\rm vol}_X)$. Let $g \in C^\infty(X)$. Then we have $\{g, f\}_\pi = \{g, \pi^n({\rm vol}_X)\}_\pi = \pi(dg, d(\pi^n{\rm vol}_X)) = \pm \mc{L}_{V_g} (d \pi^n({\rm vol}_X)) = \pi^n \mc{L}_{V_g} {\rm vol}_X$, using that as $V_g$ is Hamiltonian, it preserves $\pi^n$. On a small enough open we have $\mc{L}_{V_g} {\rm vol}_X = h {\rm vol}_X$ for some function $h$, so that $\pi^n \mc{L}_{V_g} {\rm vol}_X = h f$, which lies in $I$.
\ep
\subsection{Lifting log-Poisson structures}
\label{sec:liftinglogpoisson}
In this section we discuss the process of lifting log-Poisson structures, i.e.\ Poisson structures which define log divisors. Let $X$ be a $2n$-dimensional manifold. Poisson hypersurfaces $Z \subseteq X$ for a given Poisson structure $\pi \in {\rm Poiss}(X)$ are related to the log-tangent bundle $\mc{A}_Z = TX(-\log Z)$ constructed out of them as follows.
\begin{prop}[{\cite[Proposition 4.4.1]{Pym13}}]\label{prop:poissonhypersurface} Let $\pi \in {\rm Poiss}(X)$ be given. Then a hypersurface $Z \subseteq X$ is $\pi$-Poisson if and only if $\pi$ is of $\mc{A}_Z$-type.
\end{prop}
\bp The converse statement follows as if $\pi$ is of $\mc{A}_Z$-type, the degeneracy locus of $\mc{A}_Z$ is $\pi$-Poisson. For the direct implication, let $z \in I_Z$ be a local generator at $Z$, so that $\Omega^2(\log Z)$ is generated by $d\log z \wedge \Omega^1(X)$ and $\Omega^2(X)$. As $\Omega^1(X)$ is generated by exact forms, it suffices to check that $\pi(d\log z, dg) = z^{-1} \{z,g\}$ is smooth for all $g \in C^\infty(X)$. As $Z$ is $\pi$-Poisson, the ideal $I_Z$ is Poisson for $\pi$, so that $\{z,g\}_\pi = h z$ for some $h \in C^\infty(X)$. Hence $\pi(d\log z,dg) = h$ which is smooth.
\ep
We can determine when the $\mc{A}_Z$-lift will be nondegenerate, i.e.\ dual to a log-symplectic structure. Let $\pi$ be generically symplectic, so that $\pi^\sharp$ is generically an isomorphism, i.e.\ $T^*_\pi X$ has dense isomorphism locus. Then $X_{\pi,2n-2} = X_{T^*_\pi X, n-1}$ is a Poisson subspace, and equal to the zero set of $\wedge^n \pi$. Assume that this is a smooth submanifold of $X$. To state the following result, let a \emph{Poisson divisor} for $\pi$ be an element of the free abelian group generated by all smooth $\pi$-Poisson hypersurfaces.
\begin{prop}[{\cite[Proposition 4.4.2]{Pym13}}]\label{prop:nondegaztype} Let $\pi \in {\rm Poiss}(X)$ and $Z \subseteq X$ a smooth Poisson hypersurface for $\pi$. Assume that $\pi$ is generically symplectic. Then $\pi$ is of nondegenerate $\mc{A}_Z$-type if and only if $Z = X_{\pi,2n-2}$ as Poisson divisors.
\end{prop}
In other words, for $\pi$ to nondegenerately lift to $\mc{A}_Z$, $Z$ must be a Poisson hypersurface for $\pi$, and the precise Poisson divisor where $\pi$ has rank equal to $2n-2$.
\bp As $Z$ is a $2n-1$-dimensional Poisson subspace of $X$, the rank of $\pi$ on $Z$ can be at most $2n-2$, as it must be even. As the rank of a Poisson bivector can only change in even amounts and is generically equal to $2n$, we conclude that $Z \subseteq X_{\pi,2n-2}$. By \autoref{prop:poissonhypersurface}, $\pi$ admits an $\mc{A}_Z$-lift $\pi_{\mc{A}_Z}$. Note that $\pi_{\mc{A}_Z}^\sharp\colon \mc{A}_Z^* \to \mc{A}_Z$ is an isomorphism if and only if $\pi_{\mc{A}_Z}^n \in \Gamma(\det \mc{A}_Z)$ is nowhere zero. Using \autoref{prop:logtangentlogdiv} we have $\det(\mc{A}_Z) \cong \det(TX) \otimes L^*_Z$, so that the subset of $X$ where $\pi_\mc{A}^\sharp$ is not an isomorphism is exactly the zero set of $\pi_\mc{A}^n \in \Gamma(\det TX \otimes L^*_Z)$. In fact we have $\pi_\mc{A}^n = \pi^n \otimes z^{-1}$ where $z$ is a local generator of $I_Z$, which thus vanishes at the Poisson subspace $X_{\pi,2n-2} - Z$ (with subtracting making sense as $Z \subseteq X_{\pi,2n-2}$). We conclude that $\pi_\mc{A}$ is nondegenerate if and only if $X_{\pi,2n-2} = Z$.
\ep
\begin{rem}\label{rem:x2notlog} In the statement of \autoref{prop:nondegaztype}, it does not suffice to demand that $Z = X_{\pi,2n-2}$ as subspaces. Namely, consider $\R^2$ with coordinates $(x,y)$ and $\pi = x^2 \partial_x \wedge \partial_y$ with $Z_\pi = \{x = 0\}$. Then $\pi$ admits an $\mc{A}_{Z_\pi}$-lift $\pi_{\mc{A}_Z} = x (x \partial_x \wedge \partial_y)$, which is clearly not nondegenerate. The divisor ideal of $\pi$ is given by $\langle x^2 \rangle$, which is not equal to $\langle x \rangle = I_{Z_\pi}$.
\end{rem}
From the above two results we conclude the following bijective correspondence between log-Poisson structures, and \blog{}s. This will be used in Chapter \ref{chap:constructingblogs} to construct log-Poisson structures.
\begin{cor}[{\cite[Proposition 20]{GuilleminMirandaPires14}}]\label{prop:blogbsymp} A Poisson structure $\pi$ on $X^{2n}$ is log if and only if it is of nondegenerate $\mc{A}_Z$-type with log divisor $Z = (\wedge^{2n} TX, \wedge^n \pi)$.
\end{cor}
\bp In light of \autoref{prop:standarddivideal} and \autoref{prop:nondegaztype}, we need only establish that if $\pi$ is a log-Poisson structure for the log pair $(X,Z)$, then $\pi$ has rank $2n-2$ on $Z$. However, this follows immediately from the Leibniz rule $d \pi^n = n d\pi \otimes \pi^{n-1}$ together with the fact that $\pi^n$ vanishes transversally, so that $\pi^{n-1}|_Z$ must be nonzero. This was also remarked in Section \ref{sec:logpoisson}.
\ep
\begin{rem} The Poisson structure of \autoref{rem:x2notlog} is not a log Poisson structure; instead its associated divisor is the square of a log divisor.
\end{rem}
\begin{rem} An alternate proof of \autoref{prop:blogbsymp} uses the determinant of the lifting diagram (see \cite[Proposition 1.4]{GualtieriLi14})
	\begin{center}
		\begin{tikzpicture}
		\matrix (m) [matrix of math nodes, row sep=2.5em, column sep=2.5em,text height=1.5ex, text depth=0.25ex]
		{	\mc{A}_Z^* & \mc{A}_Z \\ T^*X & TX \\};
		\path[-stealth]
		(m-1-1) edge node [above] {$\pi_{\mc{A}_Z}^\sharp$} (m-1-2)
		(m-2-1) edge node [left] {$\rho_{\mc{A}_Z}^*$} (m-1-1)
		(m-2-1) edge node [above] {$\pi^\sharp$} (m-2-2)
		(m-1-2) edge node [right] {$\rho_{\mc{A}_Z}$} (m-2-2);
		\end{tikzpicture}
	\end{center}
Using \autoref{prop:logtangentlogdiv} we see that both $\det \rho_{\mc{A}_Z}$ and $\det \rho_{\mc{A}_Z}^*$ vanish to first order at $Z$, while using \autoref{lem:detofabivector} we see that $\det(\pi) = \pi \otimes \pi$ vanishes to second order. This implies that $\pi_{\mc{A}_Z}^\sharp$ is an isomorphism, so that $\pi_{\mc{A}_Z}$ is nondegenerate as desired.
\end{rem}
Similar statements can be made for normal-crossing log divisors (see \cite{Pym13}).
\begin{rem} A consequence of the above corollary is ${\rm Poiss}_{I}(X) \subseteq {\rm Poiss}_{\mc{A}_I}(X)$ for log divisors ideals $I = I_Z$.
\end{rem}
A similar characterization can be given of when a Poisson structure can be lifted to the scattering tangent bundle $\mc{C}_Z$ associated to a hypersurface. This can be found in \cite{Lanius16}, and is translated into our language. Recall that $\mc{C}_Z$ is an $\mc{A}_Z$-Lie algebroid.
\begin{prop}[{\cite[Lemma 5.2]{Lanius16}}] Let $\pi$ be a Poisson structure and $Z \subseteq X$ a smooth Poisson hypersurface for $\pi$ with scattering tangent bundle $\mc{C}_Z$. Then $\pi$ is of $\mc{C}_Z$-type if the first jet of its $\mc{A}_Z$-lift $\pi_{\mc{A}_Z}$ vanishes at $Z$. Moreover, it is of nondegenerate $\mc{C}_Z$-type if and only if $Z = X_{\pi_{\mc{A}_Z}, 2n-2}$ as $\pi_{\mc{A}_Z}$-Poisson divisors.
\end{prop}
\subsubsection{Lifting $m$-log-Poisson structures}
\label{sec:liftingmlogpoisson}
We next discuss lifting of $m$-log-Poisson structures. Let $\pi \in {\rm Poiss}(X)$ be an $m$-log-Poisson structure. Thus there exists a line subbundle $K \subseteq \wedge^{2m} \mc{A}$ such that $(K, \wedge^m \pi)$ is a log divisor. Let $Z = Z_\pi = (\wedge^m \pi)^{-1}(0)$ be the zero set of this divisor.
\begin{prop} The zero set $Z_\pi$ is a Poisson hypersurface for $\pi$.
\end{prop}
\bp As $Z_\pi$ is the zero set of $\wedge^m \pi$, we see immediately that $Z_\pi = X_{\pi,2m-2}$. By \autoref{prop:pidegenpipoisson}, this means that $Z_\pi$ is a Poisson subspace for $\pi$, hence is a Poisson hypersurface for $\pi$.
\ep
Consequently, by \autoref{prop:poissonhypersurface}, we know that $\pi$ is of $\mc{A}_{Z}$-type. As $\mc{A}_{Z}$ has isomorphism locus $X \setminus Z_\pi$ which is dense, by \autoref{prop:uniquelifts} we know that $\pi$ hence can be uniquely lifted to an $\mc{A}_{Z}$-Poisson structure $\pi_{\mc{A}_Z}$. Using the discussion around \autoref{exa:mlogpoisson}, it is not hard to verify the following.
\begin{prop} Let $\pi$ be an $m$-log-Poisson structure such that $V_\pi$ is zero on $Z_\pi$. Then $\pi$ is of $m$-regular $\mc{A}_{Z}$-type.
\end{prop}
In other words, the condition on the vanishing of the modular vector field implies its unique $\mc{A}_{Z}$-lift will be $m$-regular. Then, as we now have an $m$-regular $\mc{A}_Z$-Poisson structure, as in \autoref{exa:mregularapoisson} we know we can lift $\pi_{\mc{A}_Z}$ further to its $\mc{A}_Z$-distribution $D_{\pi_{\mc{A}_Z}}$, where it will be nondegenerate. On the other hand, using \autoref{cor:divtypealmostreg} and \autoref{prop:almostregulardistr} there exists a distribution $D \subseteq TX$ on $X$ for which $\det D = K$. It is immediate that $\rho_{\mc{A}_Z}$ maps $D_{\pi_{\mc{A}_Z}}$ onto $D$ over $X \setminus Z$. Moreover, $\pi$ is of log-$D$-type, i.e.\ admits a lift $\pi_D$ to $D$ which defines a log divisor. We summarize this in the following diagram (where $\mc{A} = \mc{A}_Z$).
\begin{center}
	\begin{tikzpicture}
	\matrix (m) [matrix of math nodes, row sep=2.5em, column sep=2.5em,text height=1.5ex, text depth=0.25ex]
	{	(D_{\pi_{\mc{A}}}, \pi_{D_\mc{A}} \,\text{nondeg.}) & (D, \pi_D \,\text{log}) \\ (\mc{A}, \pi_\mc{A} \, m\text{-reg.}) & (TX, \pi \, m\text{-log}) \\};
	\path[-stealth]
	(m-2-1) edge node [above] {$\rho_\mc{A}$} (m-2-2);
	\draw[right hook-latex]
	(m-1-1) edge (m-2-1)
	(m-1-2) edge (m-2-2);
	\draw[dotted, ->] (m-1-1) edge node [above] {$\rho_\mc{A}$} (m-1-2);
	\end{tikzpicture}
\end{center}
\subsection{Lifting elliptic Poisson structures}

As for log-Poisson structures, the elliptic Poisson structures introduced in Section \ref{sec:ellpoissonstr} can be lifted to the elliptic tangent bundle they give rise to.
\begin{prop}[{\cite[Lemma 3.4]{CavalcantiGualtieri15}}]\label{prop:elllogcorrespondence} A Poisson structure $\pi$ on $X^{2n}$ is elliptic if and only if it is of nondegenerate $\mc{A}_{|D|}$-type with elliptic divisor $|D| = (\wedge^{2n} TX, \wedge^n \pi)$.
\end{prop}
\bp	The fact that $\pi$ lifts to $\mc{A}_{|D|}$ follows from \autoref{prop:ipoissonaitype}. The rest follows from \autoref{prop:standarddivideal} as elliptic divisors are standard (see \autoref{prop:elltangentlogdiv}). Alternatively (or equivalently), consider the lifting diagram
\begin{center}
	\begin{tikzpicture}
	\matrix (m) [matrix of math nodes, row sep=2.5em, column sep=2.5em,text height=1.5ex, text depth=0.25ex]
	{	\mc{A}_{|D|}^* & \mc{A}_{|D|} \\ T^*X & TX \\};
	\path[-stealth]
	(m-1-1) edge node [above] {$\pi_{\mc{A}_{|D|}}^\sharp$} (m-1-2)
	(m-2-1) edge node [left] {$\rho_{\mc{A}_{|D|}}^*$} (m-1-1)
	(m-2-1) edge node [above] {$\pi^\sharp$} (m-2-2)
	(m-1-2) edge node [right] {$\rho_{\mc{A}_{|D|}}$} (m-2-2);
	\end{tikzpicture}
\end{center}
from which by \autoref{lem:detofabivector} we see that $\pi^\sharp_{\mc{A}_{|D|}}$ is an isomorphism.
\ep
\chapter{Dirac structures}
\label{chap:diracgeometry}
\fancypagestyle{empty}{%
	\fancyhf{}%
	\renewcommand\headrulewidth{0pt}%
	\fancyhead[RO,LE]{\thepage}%
}
In this chapter we give an overview of the language of Dirac geometry as is used in this thesis. Dirac structures form a convenient setting to simultaneously treat involutive distributions, closed two-forms and Poisson bivectors. Moreover, they admit a description using spinors, which we will make ample use of. As will become clear, it is possible to define the $\mc{A}$-analogue of Dirac structures, so that this chapter could also be called ``$\mc{A}$-Dirac structures'' in analogy with Chapters \ref{chap:asymplecticstructures} and \ref{chap:apoissonstructures}.

We can moreover relate Dirac structures to divisors as was done in Section \ref{sec:poissondivisors} for Poisson structures, and there is also a lifting concept for Dirac structures. Note that this chapter is a prerequisite for the next chapter on \gcs{}s, as a \gcs{} can be described as a specific type of complex Dirac structure. General references on Dirac structures include \cite{Bursztyn13, Courant90, Gualtieri11, LiBlandMeinrenken09}. This chapter reviews the existing literature, except for Sections \ref{sec:adiracstructures}, \ref{sec:diracdivisors} and \ref{sec:liftingdirac}.
\subsection*{Organization of the chapter}
This chapter is built up as follows. In Section \ref{sec:diraclinearalgebra} we consider linear Dirac structures and spinors on vector spaces, emphasizing the concept of a double of a vector space. In Section \ref{sec:diracstructures} we then globalize to manifolds, defining Courant algebroids and Dirac structures along with their spinor bundle description. In Section \ref{sec:adiracstructures} we then consider the $\mc{A}$-analogues of both concepts. We continue in Section \ref{sec:diracdivisors} by relating Dirac structures to divisors using their spinor bundle, and finish this chapter with Section \ref{sec:liftingdirac} on lifting Dirac structures to Lie algebroids.
\section{Linear algebra and spinors}
\label{sec:diraclinearalgebra}
We start with linear algebra and develop the required concepts on vector spaces, before globalizing to manifolds. We follow in part the treatment of \cite{Cavalcanti07}, yet consider also the references mentioned at the start of this chapter. Throughout, let $V$ be a real vector space of dimension $m$.
\subsection{Doubles of a vector space}
Recall that given a nondegenerate pairing $\langle \cdot, \cdot \rangle_W$ on a vector space $W$, a subspace $E < W$ is \emph{isotropic} if $E \subseteq E^\perp$, i.e.\ $\langle E, E \rangle_W = 0$. If the pairing $\langle \cdot, \cdot \rangle_W$ has signature $(k,\ell)$, then the maximal dimension of an isotropic subspace of $W$ is ${\rm min}(k,\ell)$.
\begin{defn}\label{defn:double} A \emph{double} of $V$ is a $2m$-dimensional vector space $\mc{D}V$ equipped with a nondegenerate pairing $\langle \cdot, \cdot \rangle_{\mc{D}V}$ and a surjective linear map $\rho_V\colon \mc{D}V \to V$ such that $\ker \rho_V$ is isotropic with respect to $\langle \cdot, \cdot \rangle_{\mc{D}V}$.
\end{defn}
Let $\mc{D}V$ be a double of $V$. Note that $\ker \rho_V$ has dimension $m$, so that $\langle \cdot, \cdot \rangle_{\mc{D}V}$ has signature $(m,m)$. The pairing gives an isomorphism $(\mc{D}V)^* \cong \mc{D}V$ given by $v \mapsto \langle v, \cdot \rangle_{\mc{D}V}$ for $v \in \mc{D}V$. As $\rho_V$ is surjective, its dual $\rho_V^*$ is injective, so that $\rho_V^*\colon V \to (\mc{D}V)^* \cong \mc{D}V$ allows us to view $V^*$ as a subset of $\mc{D}V$. We will sometimes write $V^* \subseteq \mc{D}V$, omitting the map $\rho_V^*$. By definition we have $\langle \rho_V^*(V^*), \ker \rho_V \rangle_{\mc{D}V} = 0$ so that $\rho_V^*(V^*) \subseteq (\ker \rho_V)^\perp$. However, as $\ker \rho_V$ is maximally isotropic we have $(\ker \rho_V)^\perp = \ker \rho_V$, whence $\rho_V^*(V^*) = \ker \rho_V$. From this discussion we obtain a short exact sequence
\be
0 \to V^* \stackrel{\frac12 \rho_V^*}{\to} \mc{D}V \stackrel{\rho_V}{\to} V \to0.
\ee
It will be very useful to relate doubles $\mc{D}V$ to the following standard one.
\begin{defn} The \emph{standard double} of $V$ is the vector space $\mathbb{V} = V \oplus V^*$ with pairing $\langle X + \xi, Y + \eta \rangle_{\mathbb{V}} = \frac12(\xi(Y) + \eta(X))$ and $\rho_V\colon \mathbb{V} \to V$ the projection.
\end{defn}
The identification between $\mc{D}V$ and $\mathbb{V}$ is obtained by splitting the above sequence.
\begin{defn} Let $\mc{D}V$ be a double of $V$. An \emph{isotropic splitting} for $\mc{D}V$ is a linear map $\nabla\colon V \to D\mc{V}$ satisfying $\rho_V \circ \nabla = {\rm Id}$ such that $\nabla(V)$ is isotropic.
\end{defn}
Let $\nabla$ be an isotropic splitting for a double $\mc{D}V$. Then $p_\nabla = (\rho_V, \nabla^*)\colon \mc{D}V \to V \oplus V^*$ gives an isomorphism $\mc{D}V \cong \mathbb{V}$. This is an isomorphism of doubles of $V$.
\begin{lem} The isomorphism $p_\nabla\colon \mc{D}V \to \mathbb{V}$ is an isometry for $\langle \cdot, \cdot \rangle_{\mc{D}V}$ and $\langle \cdot, \cdot \rangle_{\mathbb{V}}$.
\end{lem}
\subsection{Linear Dirac structures}
\label{sec:lineardirac}
In this section we introduce linear Dirac structures. Let $\mc{D}V$ be a double of $V$.
\begin{defn} A \emph{linear Dirac structure} is a maximal isotropic subspace $E$ of $\mc{D}V$.
\end{defn}
We thus have $E = E^\perp$ and $\dim E = m$. Such subspaces are also called \emph{Lagrangian}. There is a complex analogue as follows. Note that it is immediate that the complexification $\mc{D} V \otimes \C$ is a double of $V \otimes \C$.
\begin{defn} A \emph{linear complex Dirac structure} is a maximal isotropic subspace $E_\C$ of $\mc{D}V \otimes \C$.
\end{defn}
We will denote complexifications by a subscript, i.e.\ $\mc{D}V_\C = \mc{D}V \otimes \C$.
\begin{exa} By definition $\rho_V^*(V^*)$ is isotropic, but it is also of maximal dimension, so that it is a linear Dirac structure.
\end{exa}
\begin{exa} Given a splitting $\mc{D}V \cong \mathbb{V}$, it is immediate that $V$ is a linear Dirac structure. Note one cannot view $V$ as naturally sitting inside $\mc{D}V$ without choosing an isotropic splitting.
\end{exa}
\begin{exa} Choose a splitting $\mc{D}V \cong \mathbb{V}$ and let $\omega \in \wedge^2 V^*$ be a two-form. Then ${\rm Gr}(\omega) = \{X + \omega^\flat(X) \, | \, X \in V\} = e^B V^*$ is a linear Dirac structure.
\end{exa}
\begin{exa} Choose a splitting $\mc{D}V \cong \mathbb{V}$ and let $\pi \in \wedge^2 V$ be a bivector. Then ${\rm Gr}(\pi) = \{\pi^\sharp(\xi) + \xi \, | \, \xi \in V^*\} = e^\beta V$ is a linear Dirac structure.
\end{exa}
\begin{exa} Choose a splitting $\mc{D}V \cong \mathbb{V}$ and let $D \subseteq V$ be a subspace. Then $D \oplus {\rm Ann}(D)$ is a linear Dirac structure, where ${\rm Ann}$ denotes the annihilator.
\end{exa}
We will discuss examples of linear complex Dirac structures in the next chapter, after we define the notion of a linear \gcs{}.
\begin{defn} Let $E$ be a linear Dirac structure in $\mc{D}V$. The \emph{type} ${\rm type}(E)$ of $E$ is the codimension of $\rho_V(E) \subseteq V$, hence is an integer ranging from $0$ to $m$.
\end{defn}
The type can be viewed as measuring the size of the ``intersection'' of $E$ with $V$. In the extreme cases of ${\rm type}(E) \in \{0,m\}$, the Dirac structure $E$ is the graph of a two-form or bivector (after a choice of splitting).
\begin{prop} Choose a splitting $\mc{D}V \cong \mathbb{V}$ and let $E$ be a linear Dirac structure such that $E \cap V^* = 0$. Then there exists a two-form $\omega \in \wedge^2 V^*$ such that $E = {\rm Gr}(\omega)$.
\end{prop}
\bp As $E \cap V^* = 0$, the map $\rho_V|_E\colon E \to V$ is injective. By counting dimensions it is an isomorphism. Define a map $F\colon V \to V^*$ as the composition $F = {\rm pr}_{V^*} \circ (\rho_V|_E)^{-1}$. Then $E = {\rm Gr}(F)$. As $E$ is isotropic, the map $\omega\colon V \times V \to \R$ given by $\omega(X,Y) = \iota_Y F(X)$ is skew-symmetric, so that $\omega \in \wedge^2 V^*$ and $F = \omega^\flat$.
\ep
A similar statement is true for bivectors, whose proof we omit.
\begin{prop} Choose a splitting $\mc{D}V \cong \mathbb{V}$ and let $E$ be a linear Dirac structure such that $E \cap V = 0$. Then there exists a bivector $\pi \in \wedge^2 V$ such that $E = {\rm Gr}(\pi)$.
\end{prop}
We finish this section by noting that given a linear complex Dirac structure $E_\C$, its complex conjugate $\overline{E}_\C$ is another linear complex Dirac structures.
\subsection{Spinors}
\label{sec:spinors}
There is a useful alternative description of linear Dirac structures in terms of spinors. As before, fix a vector space $V$ and let $\mc{D}V$ be a double of $V$.
\begin{defn} The \emph{Clifford algebra} ${\rm Cl}(\mc{D}V)$ is defined as the free algebra generated by elements $v \in \mc{D}V \subseteq {\rm Cl}(\mc{D}V)$ subject to the relation $v^2 = \langle v, v \rangle$.
\end{defn}
Note that $\rho_V^*(V^*)$ is maximally isotropic, so that its exterior algebra $\wedge^\bullet V^*$ is a subalgebra of ${\rm Cl}(\mc{D}V)$. Elements of this subalgebra will be called \emph{spinors}. In particular, the line $\wedge^m V^*$ generates a left ideal $\mc{I}$. Let $\mc{D}V \cong \mathbb{V}$ be an isotropic splitting which we will fix for the remainder of this section. Then we obtain an isomorphism $\mc{I} \cong \wedge^\bullet V \cdot \wedge^m V^* \cong \wedge^\bullet V^*$. The induced Clifford action on $\wedge^\bullet V^*$ is given by $(X + \xi) \cdot \alpha = \iota_X \alpha + \xi \wedge \alpha$ for $X + \xi \in \mathbb{V}$ and $\alpha \in \wedge^\bullet V^*$.
\begin{prop} This is an action of the Clifford algebra, i.e.\ $(X + \xi) \cdot (X + \xi) \cdot \alpha = \langle X + \xi, X + \xi \rangle \alpha$.
\end{prop}
\bp This follows from $\iota_X(\iota_X \alpha + \xi \wedge \alpha) + \xi \wedge (\iota_X \alpha + \xi \wedge \alpha) = 2 \iota_X \xi \cdot \alpha$. 
\ep
Let $\rho \in \wedge^\bullet V^*$ be a nonzero spinor. Then $\rho$ gives rise to its \emph{Clifford annihilator} $E_\rho := {\rm Ann}(\rho) = \{v \in \mathbb{V} \, | \, v \cdot \rho = 0\}$. We can quickly check the following.
\begin{lem} Let $\rho$ be a spinor. Then $E_\rho$ is isotropic.
\end{lem}
\bp Let $v \in E_\rho$. Then $0 = v^2 \cdot \rho = \langle v, v \rangle \rho$, so that $\langle v, v \rangle = 0$.
\ep
There is a complex analogue of the above. Given a complex spinor $\rho \in \wedge^\bullet V^* \otimes \C$, its Clifford annihilator is given by $E_\rho = {\rm Ann}(\rho) = \{v \in \mathbb{V} \otimes \C \, | \, v \cdot \rho = 0\}$. Again, $E_\rho$ is an isotropic subspace.
One wonders which spinors have maximally isotropic Clifford annihilators. These spinors deserve a special name.
\begin{defn} A (complex) spinor $\rho$ is \emph{pure} if $E_\rho$ is maximally isotropic.
\end{defn}
As one can show, a complex spinor being pure implies it is of the form $\rho = e^{B + i \omega} \wedge \Omega$, where $B, \omega \in \wedge^2 V^*$ are real two-forms and $\Omega$ is a decomposable complex form. The following proposition is crucial.
\begin{prop} Let $E$ ($E_\C$) be a maximally isotropic (complex) subspace. Then there exists a (complex) pure spinor $\rho$ such that $E_\rho = E$ respectively $E_\C$. Moreover, if $E_\rho = E_{\rho'}$ for two (complex) pure spinors $\rho, \rho'$, then $\rho = \lambda \rho'$ for some nonzero scalar $\lambda$.
\end{prop}
As a consequence, maximally isotropic subspaces are in one-to-one correspondence with \emph{pure spinor lines} $K_\rho = \langle \rho \rangle \subset \wedge^\bullet V^*$ (or $\wedge^\bullet V^* \otimes \C$). It is possible to relate intersections of maximal isotropic subspaces with the spinors defining them by introducing a pairing on spinors.  Let $\tau$ be the antiautomorphism of ${\rm Cl}(\mathbb{V})$ defined on decomposable elements by $\tau(v_1 \cdot \dots \cdot v_k) = v_k \cdot \ldots \cdot v_1$ for $v_i \in \mathbb{V}$.
\begin{defn}\label{defn:chevalley} The \emph{Chevalley pairing} on $\wedge^\bullet V^* \subset {\rm Cl}(\mathbb{V})$ is defined by $(\rho_1, \rho_2) \mapsto (\tau(\rho_1) \wedge \rho_2)_{\rm top}$, where ${\rm top}$ denotes the top degree component (i.e.\ $m$) of the form.
\end{defn}
Using this pairing, intersections are related to the pairing as follows.
\begin{prop}\label{prop:pairingintersection} Let $\rho, \tau$ be pure spinors. Then $E_\rho \cap E_\tau = 0$ if and only if $(\rho,\tau) \neq 0$.
\end{prop}
\subsection{Transformations}
In this section we discuss automorphisms of doubles $\mc{D}V$. It will often be useful to choose an isotropic splitting $\nabla$ giving an isomorphism $\mc{D}V \cong \mathbb{V}$ as discussed above.

Let $\mc{D}V$ be a double of $V$, and let $B \in \wedge^2 V^*$ be a real two-form.
\begin{defn} The \emph{$B$-field transform} of $B$ is the map $v \mapsto v - B^\flat(\rho_V(v))$.
\end{defn}
A $B$-field transform should be thought of as a shearing transformation in the direction of $V^*$. It can be shown that any two isotropic splittings of $\mc{D}V$ differ by a $B$-field transform. Because of this, to study Lagrangian subspaces of arbitrary doubles $\mc{D}V$ it suffices to study Lagrangian subspaces of $\mathbb{V}$ and how they transform under $B$-field transformations.

There is an analogous transformation associated to bivectors, but for this we need to choose a splitting. Let $\mc{D}V \cong V \oplus V^*$ be an isotropic splitting for $\mc{D}V$ and $\beta \in \wedge^2 V$ be a real bivector.
\begin{defn} The \emph{$\beta$-transform} associated to $\beta$ is the map $V + \xi \mapsto V + \xi + \beta^\sharp(\xi)$.
\end{defn}
Given a linear Dirac structure $E$ in $\mc{D}V$, we can apply the $B$-field transform to $E$ to obtain a new linear subspace, $E_B = \{v - B^\flat(\rho_V(v)) \, | \, v \in E\}$, and similarly for $E_\beta$. The following is a result of skew-symmetry of $B$ and $\beta$.
\begin{prop}\label{prop:bfieldtransform} The space $E_B$ and $E_\beta$ is maximally isotropic, i.e.\ a are linear Dirac structure in $\mc{D}V$ (the latter after a choice of splitting).
\end{prop}
The $B$-field and $\beta$-transforms act on spinors in the obvious way, namely $\rho \mapsto e^B \wedge \rho$ and $\rho \mapsto e^\beta \rho$. It is immediate that $(e^{-B} \rho, e^{-B} \rho') = (\rho, \rho')$. In terms of a splitting $\mc{D}V \cong \mathbb{V}$, we can write $e^B$ and $e^\beta$ on $\mathbb{V}$ in matrix form as:
\be
e^B = \begin{pmatrix} 1 & 0 \\ B & 1 \end{pmatrix} \qquad \text{and} \qquad e^\beta = \begin{pmatrix} 1 & \beta \\ 0 & 1 \end{pmatrix}.
\ee
Moreover, as a $B$-field transform preserves the projection of $E$ onto $V$, we have ${\rm type}(E_B) = {\rm type}(E)$. However, the $\beta$-transform can change the type.
\section{Dirac structures}
\label{sec:diracstructures}
We are now ready to globalize to manifolds, where we will work with vector bundles and perform the previous linear algebra discussion fiberwise.
\subsection{Courant algebroids}
Let $X$ be a $2n$-dimensional manifold equipped with a closed three-form $H \in \Omega^3_{\rm cl}(X)$. We start with the analogue of the standard double of a vector space.
\begin{defn} The \emph{double tangent bundle} is $\mathbb{T}X := TX \oplus T^*X$ with anchor the projection $p\colon \mathbb{T}X \to TX$.  It carries a natural pairing $\langle V + \xi, W + \eta\rangle = \frac12(\eta(V) + \xi(W))$ of split signature and an $H$-twisted Courant bracket $\llbracket V + \xi, W + \eta\rrbracket_H = [V, W] + \mc{L}_V \eta - \iota_W d\xi + \iota_V \iota_W H$ for $V, W \in \Gamma(TX)$ and $\xi, \eta \in \Gamma(T^*X)$.
\end{defn}
The properties of the double tangent bundle are axiomatized by that of a Courant algebroid \cite{LiuWeinsteinXu97}, as we now recall.
\begin{defn}\label{defn:courantalgebroid} A \emph{Courant algebroid} is an anchored vector bundle $(\mc{E},\rho_\mc{E})$ equipped with a nondegenerate pairing $\langle \cdot, \cdot \rangle_\mc{E}$, and a bracket $\llbracket\cdot,\cdot\rrbracket_{\mc{E}}$ on $\Gamma(\mc{E})$ satisfying for all $v_1,v_2,v_3 \in \Gamma(\mc{E})$ and $f \in C^\infty(X)$:
	\bi
		\item $\llbracket v_1,\llbracket v_2,v_3\rrbracket_\mc{E}\rrbracket_\mc{E} = \llbracket \llbracket v_1,v_2\rrbracket_\mc{E},v_3\rrbracket_\mc{E} + \llbracket v_2,\llbracket v_1,v_3\rrbracket_\mc{E}\rrbracket_\mc{E}$;
		\item $\mc{L}_{\rho_{\mc{E}}(v_1)} \langle v_2, v_3 \rangle_\mc{E} = \langle \llbracket v_1,v_2 \rrbracket_{\mc{E}}, v_3 \rangle_\mc{E} + \langle v_2, \llbracket v_1, v_3 \rrbracket_{\mc{E}} \rangle_\mc{E}$;
		\item $\llbracket v_1, v_2\rrbracket_{\mc{E}} + \llbracket v_2, v_1\rrbracket_{\mc{E}} = \rho_{\mc{E}}^*(d \langle v_1, v_2 \rangle_\mc{E})$.
	\ei
\end{defn}
The first two axioms state that the operation $\llbracket v_1, \cdot \rrbracket_{\mc{E}}$ is a derivation of both the bracket and the pairing. The third axiom states that the lack of skew-symmetry of the bracket is governed by the pairing.
\begin{rem} Given a Courant algebroid $\mc{E} \to X$, a consequence of the definition is that for all $v_1, v_2 \in \Gamma(\mc{E})$ and $f \in C^\infty(X)$ we have
\bi
	\item $\llbracket v_1, f v_2 \rrbracket_\mc{E} = f \llbracket v_1, v_2 \rrbracket_\mc{E} + \mc{L}_{\rho_{\mc{E}}(v_1)} f \cdot v_2$;
	\item $\rho_\mc{E}(\llbracket v_1, v_2 \rrbracket_\mc{E}) = \llbracket \rho_\mc{E}(v_1), \rho_\mc{E}(v_2) \rrbracket_\mc{E}$;
	\item $\rho_\mc{E} \circ \rho_\mc{E}^* = 0$, where $\rho_\mc{E}^*\colon T^*X \to \mc{E}^* \cong \mc{E}$, using the pairing $\langle \cdot, \cdot \rangle_\mc{E}$.
\ei
These properties should be compared with that of Lie algebroid (see \autoref{defn:liealgebroid} and \autoref{prop:anchorliealgmorph}), and with that of a double of a vector space (\autoref{defn:double}).
\end{rem}
The global study of (complex) Dirac structures and hence of \gcs{} mainly takes place on certain Courant algebroids that are called exact.
\begin{defn} A Courant algebroid $\mc{E}$ is \emph{exact} if the following sequence is exact:
	\be
	0 \to T^*X \stackrel{\rho_{\mc{E}}^*}{\to} \mc{E} \stackrel{\rho_{\mc{E}}}{\to} TX \to 0.
	\ee
\end{defn}
We see that the fibers $\mc{E}_x$ of an exact Courant algebroid are doubles of $T_x X$. In particular, exactness at $\mc{E}$ implies that $\langle \cdot, \cdot \rangle_\mc{E}$ is of split signature. Note that $\mathbb{T}X$ is an example of an exact Courant algebroid. Let $\mc{E}$ be an exact Courant algebroid. Then there always exists an isotropic splitting $\nabla\colon TX \to \mc{E}$ of the above sequence.
\begin{defn}\label{defn:curvature} The \emph{curvature} of an isotropic splitting $\nabla$ for $\mc{E}$ is the closed three-form $H \in \Omega^3_{\rm cl}(X)$ defined by $H(V_1,V_2,V_3) = \frac12 \langle \llbracket\nabla(V_1), \nabla(V_2)\rrbracket_{\mc{E}}, \nabla(V_3) \rangle_{\mc{E}}$ for $V_1, V_2, V_3 \in \Gamma(TX)$.
\end{defn}
An isotropic splitting determines an isomorphism of Courant algebroids $\mc{E} \cong (\mathbb{T}X, H)$ between $\mc{E}$ and the double tangent bundle with $H$-twisted Courant bracket, where $H$ is the curvature of the splitting. As before, isotropic splittings differ by a two-form $B \in \Omega^2(X)$, which modifies the curvature by the addition of $d B$. We see that an exact Courant algebroid has a well-defined cohomology class $[H] \in H^3(X;\R)$ associated to it, which is called the \emph{\v{S}evera class} of $\mc{E}$ \cite{Severa}. This class completely determines the exact Courant algebroid structure on $\mc{E}$ up to isomorphism.
\subsection{Dirac structures}
We are now ready to define Dirac structures, which are the global analogues of the linear Dirac structures of Section \ref{sec:lineardirac}. Let $\mc{E} \to X$ be a Courant algebroid.
\begin{defn} A \emph{almost-Dirac structure} in $\mc{E}$ is a Lagrangian subbundle $E \subseteq \mc{E}$.
\end{defn}
We see each fiber of an almost-Dirac structure is a linear Dirac structure.
\begin{defn} A \emph{Dirac structure} in $\mc{E}$ is an almost-Dirac structure $E$ in $\mc{E}$ that is involutive with respect to $\llbracket\cdot,\cdot\rrbracket_\mc{E}$, i.e.\ $\llbracket\Gamma(E),\Gamma(E)\rrbracket_\mc{E} \subseteq \Gamma(E)$. A \emph{complex Dirac structure} in $\mc{E}$ is an involutive complex almost-Dirac structure $E_\C$.
\end{defn}
Recalling that the lack of skew-symmetry of the Courant bracket $\llbracket \cdot, \cdot \rrbracket_\mc{E}$ is governed by the pairing $\langle \cdot, \cdot \rangle_\mc{E}$ (see \autoref{defn:courantalgebroid}), we immediately see the following. 
\begin{prop} Let $E$ be a Dirac structure in $\mc{E}$. Then $\rho_E = \rho_{\mc{E}}|_E$ together with $[\cdot,\cdot]_E = \llbracket\cdot,\cdot\rrbracket_\mc{E}$ defines a Lie algebroid structure on $E$.
\end{prop}
Recall from Section \ref{sec:spinors} that there is a one-to-one correspondence between linear Dirac structures and pure spinor lines. This has a global analogue as follows. Asssume that $\mc{E} \to X$ is an exact Courant algebroid, and fix an isotropic splitting $\mc{E} \cong (\mathbb{T}X,H)$.
\begin{defn} A \emph{spinor line} is a line subbundle $K \subseteq \wedge^\bullet T^*X$. A spinor line is \emph{pure} if it can be locally generated by a spinor of the form $\rho = e^B \wedge \Omega$ with $\Omega$ a decomposable form. A \emph{complex spinor line} is a complex line subbundle $K \subseteq \wedge^\bullet V^* \otimes \C$. A complex spinor line is \emph{pure} if it can be locally generated by a complex spinor of the form $\rho = e^{B + i \omega} \wedge \Omega$ with $\Omega$ a decomposable complex form.
\end{defn}
Any almost-Dirac structure $E$ corresponds to a pure spinor line $K_E$ with $E = {\rm Ann}(K_E)$. We can characterize the integrability condition of $E$ in terms of $K_E$. Let $d^H = d + H \wedge$ be the $H$-twisted de Rham differential.
\begin{prop} Let $K$ be a pure spinor line. Then the almost-Dirac structure $E_K = {\rm Ann}(K)$ in $\mc{E}$ is a Dirac structure if and only if $d^H \Gamma(K) \subseteq \Gamma(E^* \cdot K)$.
\end{prop}
Hence for any local trivialization $\rho$ of $K$, there must exist a local section $v = X + \xi \in \Gamma(\mathbb{T}X)$ such that $d^H \rho = v \cdot \rho$. 
Let $E$ be an almost-Dirac structure. The \emph{type} of $E$ at a point $x \in X$ is the codimension of the projection $\rho_E(E_x)$ in $T_x X$. In terms of spinors, it is given by the degree of $\Omega$, where $\rho = e^B \wedge \Omega$ is a local generator of the pure spinor line $K_E$.

There are three main examples of Dirac structures, with $H = 0$.
\begin{exa} Given $\omega \in \Omega^2(X)$, its graph ${\rm Gr}(\omega) := \{V + \iota_V \omega \, | \, V \in TX\} \subset \mathbb{T}X$ is an almost-Dirac structure, which is a Dirac structure if and only if $d\omega = 0$.
\end{exa}
\begin{exa} Given $\pi \in \mf{X}^2(X)$, its graph ${\rm Gr}(\pi) := \{\iota_\alpha \pi + \alpha \, | \, \alpha \in T^*X\} \subset \mathbb{T}X$ is an almost-Dirac structure, which is a Dirac structure if and only if $[\pi,\pi] = 0$.
\end{exa}
\begin{exa} Let $\mc{F}$ be a foliation of dimension $2m$ on $X$. Then $W := T \mc{F}$ is an involutive distribution, and $E_{\mc{F}} := W \oplus W^\circ$ is a Dirac structure. In fact, a distribution $W$ defines an almost-Dirac structure which is Dirac if and only if $W$ is involutive.
\end{exa}
Let us discuss the spinor description of two of the above examples.
\begin{exa} Given a two-form $\omega$, we have ${\rm Gr}(\omega) = {\rm Ann}(e^{-\omega})$, as for $\varphi := e^{-\omega}$ we have $0 = (V + \alpha) \cdot e^{-\omega} = (-\iota_V \omega + \alpha) \wedge e^{-\omega}$ implying $\alpha = \iota_V \omega$, so that $V+\alpha \in {\rm Ann}(e^{-\omega})$ implies $V + \alpha = V + \iota_V \omega \in {\rm Gr}(\omega)$.
\end{exa}
The spinor line of a Poisson structure $\pi$ is generated locally by $e^{-\pi} \cdot {\rm vol}_X$. More invariantly, we can rephrase the definition of a Poisson structure as follows.
\begin{defn} A \emph{Poisson structure} is a spinor line $K \subset \wedge^\bullet T^* X$ locally generated by a pure spinor $\varphi \in \Gamma(K)$, pointwise of the form $\varphi = e^\omega \wedge \Omega$, and $\varphi_{2n} = \omega^{\frac{n-k}{2}} \wedge \Omega \neq 0$, and there exists $V \in \mf{X}(X)$ such that $d\rho = \iota_V \rho$.
\end{defn}
We have the following useful proposition, saying when a Dirac structure is a graph of a closed two-form or Poisson bivector near a point.
\begin{prop}\label{prop:diraclocal} Let $E \subseteq \mathbb{T}X$ be a Dirac structure and $x \in X$. Then:
	\bi
	\item[i)] $E = {\rm Gr}(\omega)$ for some $\omega$ near $x$ if and only if $\rho_E$ is surjective at $x$;
	\item[ii)] $E = {\rm Gr}(\pi)$ for some $\pi$ near $x$ if and only if $E \cap T_x X = 0$.
	\ei
\end{prop}
There is further the following global statement, due to the fact that Lagrangian subspaces are half the dimension of $\mathbb{T}X$.
\begin{prop} Let $E_1, E_2 \subseteq \mathbb{T}X$ be Lagrangians. Then the following are equivalent:
	\bi
	\item $E_1$ is transverse to $E_2$;
	\item $E_1 + E_2 = \mathbb{T}X$;
	\item $E_1 \oplus E_2 = \mathbb{T}X$;
	\item $E_1 \cap E_2 = 0$.
	\ei
\end{prop}
Note that $TX$ and $T^*X$ are Lagrangian inside $\mathbb{T}X$. Consequently, $E = {\rm Gr}(\omega)$ globally if and only if $E \cap T^*X = 0$, and $E = {\rm Gr}(\pi)$ if and only if $E \cap TX = 0$.
\section{\texorpdfstring{$\mc{A}$}{A}-Dirac structures}
\label{sec:adiracstructures}

In this section we introduce the notion of an $\mc{A}$-Dirac structure given a fixed Lie algebroid $\mc{A} \to X$. This is done by replacing the tangent bundle $TX$ by $\mc{A}$ in all of the definitions. The reason for introducing this concept, as in Chapter \ref{chap:apoissonstructures}, is to develop the language necessary to lift Dirac structures to Lie algebroids, as we discuss in Section \ref{sec:diracdivisors}. Further, we formulate $\mc{A}$-Poisson structures in terms of $\mc{A}$-Dirac geometry.

\subsection{\texorpdfstring{$\mc{A}$}{A}-Courant algebroids}

In this section we discuss the natural Courant algebroid associated to a Lie algebroid. Let $\mc{A} \to X$ be a Lie algebroid of rank $n$ and fix a closed $\mc{A}$-three-form $H_\mc{A} \in \Omega^3_{\rm cl}(\mc{A})$.
\begin{defn} The \emph{$\mc{A}$-Courant algebroid} is the direct sum $\mathbb{A} = \mc{A} \oplus \mc{A}^*$, equipped with the natural pairing $\langle v + \alpha, w + \beta \rangle = \frac12(\alpha(w) + \beta(v))$ and the $H_\mc{A}$-twisted Dorfman bracket $\llbracket v + \alpha, w + \beta\rrbracket_{\mathbb{A},H_\mc{A}} := [v,w]_{\mc{A}} + \mc{L}_v \beta - \iota_w d_{\mc{A}} \alpha + \iota_v \iota_w H_\mc{A}$ for $v,w, \in \Gamma(\mc{A})$, $\alpha,\beta \in \Gamma(\mc{A}^*)$.
\end{defn}
As before, the natural pairing is a nondegenerate symmetric bilinear form of signature $(n,n)$. The Courant algebroid $\mathbb{A}$ is also called the \emph{standard double} of $\mc{A}$, and is the Courant algebroid associated to the Lie bialgebroid $(\mc{A},\mc{A}^*)$ where $\mc{A}^*$ carries the trivial Lie algebroid structure \cite{LiuWeinsteinXu97}.
\begin{rem} It stands to reason that there is a general definition of an $\mc{A}$-Courant algebroid $E_{\mathbb{A}}$ with an $\mc{A}$-anchor $\varphi_{E_\mathbb{A}}\colon E_\mathbb{A} \to \mc{A}$ which is \emph{$\mc{A}$-exact} if it fits in the short exact sequence $0 \to \mc{A}^* \to E_{\mathbb{A}} \to \mc{A} \to 0$. Moreover, probably then any $\mc{A}$-exact $\mc{A}$-Courant algebroid satisfies $E_\mathbb{A} \cong (\mathbb{A},\llbracket\cdot,\cdot\rrbracket_{\mathbb{A},H_\mc{A}})$ via the map induced by $\varphi_{E_\mathbb{A}}$ and an isotropic splitting, in full analogy with the discussion below \autoref{defn:curvature}.
\end{rem}
An $\mc{A}$-two-form $B_\mc{A} \in \Omega^2(\mc{A})$ defines a \emph{$B_\mc{A}$-transformation} $\mc{R}_{B_\mc{A}}\colon (\mb{A},H_\mc{A}) \to (\mb{A},H_\mc{A} + d_\mc{A} B_\mc{A})$, given by $v + \alpha \mapsto v + \alpha + \iota_v B_\mc{A}$ for $v + w \in \mb{A}$.
\subsection{\texorpdfstring{$\mc{A}$}{A}-Dirac structures}
 With the $\mc{A}$-Courant algebroid in hand we can turn to introducing $\mc{A}$-Dirac structures.
\begin{defn} An \emph{$\mc{A}$-Dirac structure} relative to $H_\mc{A}$ is a subbundle $E_\mc{A} \subseteq \mathbb{A}$ that is Lagrangian for $\langle \cdot, \cdot \rangle$, i.e.\ $E_\mc{A} = E_\mc{A}^\perp$, and is involutive with respect to $\llbracket\cdot,\cdot\rrbracket_{\mathbb{A},H_\mc{A}}$.
\end{defn}
There is also the notion of an \emph{almost-$\mc{A}$-Dirac structure}, which is a Lagrangian subbundle $E_\mc{A} \subseteq \mathbb{A}$. Given an $\mc{A}$-Dirac structure $E_\mc{A}$ and an $\mc{A}$-two-form $B_\mc{A}$, the image $E^{B_\mc{A}}_\mc{A} := \mc{R}_{B_\mc{A}}(E_\mc{A})$ is an $\mc{A}$-Dirac structure relative to the $\mc{A}$-three-form $H_\mc{A} + d_\mc{A} B_\mc{A}$. Any $\mc{A}$-Dirac structure $E_\mc{A}$ is an $\mc{A}$-Lie algebroid with bracket induced from $\llbracket\cdot,\cdot\rrbracket_{\mb{A}, H_\mc{A}}$ and projection to $\mc{A}$ as $\varphi_{E_\mc{A}}$. 
\begin{exa} An $\mc{A}$-bivector $\pi_\mc{A} \in \mf{X}^2(\mc{A})$ gives rise to an almost-$\mc{A}$-Dirac structure $E_{\pi_{\mc{A}}} := {\rm Gr}(\pi_{\mc{A}}) = \{\pi_{\mc{A}}^\sharp(v) + v \in \mathbb{A} \, | \, v \in \mc{A}^*\} \subseteq \mathbb{A}$. Moreover, $E_{\pi_\mc{A}}$ is an $\mc{A}$-Dirac structure relative to $H_\mc{A} = 0$ if and only $\pi_\mc{A}$ is $\mc{A}$-Poisson.
\end{exa}
If $H_\mc{A} = 0$, any $\mc{A}$-Dirac structure $E_\mc{A}$ such that $E_\mc{A} \cap \mc{A} = 0$ is given by an $\mc{A}$-Poisson structure as in the above example. Similarly for $\mc{A}$-symplectic structures.
\begin{exa} Let $D_\mc{A} \subseteq \mc{A}$ be an $\mc{A}$-distribution. Then $E_{D_\mc{A}} = D_\mc{A} \oplus D_\mc{A}^\circ \subseteq \mathbb{A}$ is an almost-$\mc{A}$-Dirac structure. It is an $\mc{A}$-Dirac structure relative to $H_\mc{A}$ if and only if $D_\mc{A}$ is involutive with respect to the bracket $[v,w]_{\mc{A},H_\mc{A}} = [v,w]_{\mc{A}} - \iota_v \iota_w H_\mc{A}$ for $v,w \in \Gamma(\mc{A})$.
\end{exa}
While $\mc{A}$-anchored vector bundles and $\mc{A}$-Lie algebroids are also anchored vector bundles and Lie algebroids, it is not true that an $\mc{A}$-Dirac structure is a Dirac structure. To study this we discuss morphisms between Dirac structures.
\subsection{Forward and backward images}
In this section we discuss morphisms between Dirac structures. This we will immediately do for $\mc{A}$-Dirac structures, as there is no essential difference for $\mc{A} \neq TX$. Let $(\varphi,f)\colon (\mc{A},X) \to (\mc{B},N)$ be a Lie algebroid morphism, $B_\mc{A}$ a closed $\mc{A}$-two-form, and consider Courant algebroids $(\mb{A},H_\mc{A})$ and $(\mb{B},H_\mc{B})$.
\begin{defn} A pair $\Phi = (\varphi,B_\mc{A})$ is a \emph{generalized map} if $\varphi^* H_\mc{B} = H_\mc{A} + d_\mc{A} B_\mc{A}$.
\end{defn}
We will work with $(\varphi,0)$ abbreviated to $\varphi$ and introduce the $B_\mc{A}$-transformation at the end. While a generalized map does not induce a morphism between $\mb{A}$ and $\mb{B}$ unless $\varphi$ is an isomorphism, it does define a relation. Define the bundle $\mb{D} := \mb{A} \oplus \overline{f^*(\mb{B})}$ over $M$, where the line denotes equipping $\mb{B}$ with the pairing $- (\cdot,\cdot)_{\mb{B}}$.
\begin{defn} A \emph{Lagrangian distribution} from $\mc{A}$ to $\mc{B}$ is a subspace $L \subseteq \mathbb{D}$ along the graph of $f$ consisting of Lagrangian subspaces. A Lagrangian distribution is a \emph{Lagrangian relation} if it also forms a smooth subbundle of $\mathbb{D}$. The \emph{graph} of $\varphi$ is defined by $\Gamma_\varphi = \{(X + \xi, Y + \eta) \in \mathbb{A} \times f^*\mathbb{B} \, | \, \varphi(X) = Y, \varphi^*(\eta) = \xi\} \subseteq \mathbb{D}$.
\end{defn}
One checks that $\Gamma_\varphi$ is a Lagrangian relation from $\mc{A}$ to $\mc{B}$, and that an $\mc{A}$-Dirac structure $E_\mc{A} \subseteq \mathbb{A} \times \{0\}$ is a Lagrangian relation from $\mc{A}$ to the trivial vector space $\{0\}$, while any $\mc{B}$-Dirac structure $E_\mc{B}$ gives a Lagrangian relation $\varphi^* E_\mc{B} \subseteq \{0\} \times \overline{f^*(\mathbb{B})}$ from $\{0\}$ to $\mc{B}$. Lagrangian distributions can be composed to form new Lagrangian distributions, but Lagrangian relations can not always be composed as Lagrangian relations without further assumptions (see \cite{Bursztyn13} for more information). Note that $\mb{D}$ is equipped with two projection maps ${\rm pr}_{\mb{A}}\colon \mb{D} \to \mb{A}$ and ${\rm pr}_{\mb{B}}\colon \mb{D} \to f^*(\mb{B})$.
\begin{defn} Let $E_\mc{A}$ be an $\mc{A}$-Dirac structure and $E_{\mc{B}}$ be a $\mc{B}$-Dirac structure. The \emph{forward image} of $\mc{E}_{\mc{A}}$ is given by $\mf{F}_\varphi(E_\mc{A}) = {\rm pr}_{\mb{B}}(\Gamma_\varphi \cap (E_\mc{A} \times f^*(\mb{B}))) \subseteq f^*(\mb{B})$. Moreover, the \emph{backward image} of $\mc{E}_{\mc{B}}$ is given by $\mf{B}_\varphi(E_\mc{B}) = {\rm pr}_{\mb{A}}(\Gamma_\varphi \cap (\mb{A} \times f^*(E_\mc{B}))) \subseteq \mb{A}$.
\end{defn}
In other words, we have for the forward image that
\be
\mf{F}_\varphi(E_\mc{A}) = \{Y +\eta \in \mc{B} \oplus \mc{B}^* \, | \, \exists X \in \mc{A} \text{ with } \varphi(X) = Y \text{ and } X + \varphi^*\eta \in E_\mc{A}\},
\ee
and similarly for the backward image
\be
\mf{B}_\varphi(E_\mc{B}) = \{X + \xi \in \mc{A} \oplus \mc{A}^* \, | \, \exists \eta \in \mc{B}^* \text{ with } \varphi^* \eta = \xi \text{ and } \varphi(X) + \eta \in E_\mc{B}\}.
\ee
The forward and backward image define Lagrangian distributions as $\mf{F}(E_\mc{A}) = \Gamma_\varphi \circ E_\mc{A}$ and $\mf{B}(E_\mc{B}) = f^*(E_\mc{B}) \circ \Gamma_\varphi$. However, it is not immediate that they are smooth subbundles. When they are, they are Dirac structures as they are automatically involutive relative to the appropriate closed three-forms. Moreover, for the forward image in $f^*(\mathbb{B})$ to define a subbundle of $\mathbb{B}$, it must be invariant under the map $f$. In general there is no reason to expect either of these conditions to hold. We will not discuss criteria for smoothness here (see e.g.\ \cite{Bursztyn13,VanderLeerDuran16} for when $\mc{A} = TX$ and $\mc{B} = TN$).
\begin{defn} Let $E_\mc{A}$ be an $\mc{A}$-Dirac structure relative to $H_\mc{A}$ and $E_{\mc{B}}$ be a $\mc{B}$-Dirac structure relative to $\eta_{\mc{B}}$. Then a generalized map $\Phi = (\varphi,B_\mc{A})$ is a \emph{forward Dirac map} if $\mf{F}_\varphi(\mc{R}_{- B_\mc{A}}(E_\mc{A})) = f^* E_\mc{B}$. Similarly, $\Phi$ is a \emph{backward Dirac map} if $\mc{R}_{B_\mc{A}}(\mf{B}_\varphi(E_\mc{B})) = E_\mc{A}$.
\end{defn}
Given an $\mc{A}$-Dirac structure $E_\mc{A} \subseteq \mathbb{A}$, we can consider its forward Dirac image $E = \mf{F}_{\rho_\mc{A}}(E_{\mc{A}}) \subseteq \mathbb{T}X$ under the anchor of $\mc{A}$, but this need not be a Dirac structure.
\begin{exa} Let $\mc{A} \to X$ be a Lie algebroid such that $\rho_\mc{A}(\mc{A}) \subseteq TX$ is not a subbundle. Then $E_\mc{A} := \mc{A} \subseteq \mathbb{A}$ is an $\mc{A}$-Dirac structure, but its forward Dirac image $\mf{F}_{\rho_\mc{A}}(E_\mc{A}) = \rho_\mc{A}(\mc{A}) \oplus \ker \rho_\mc{A}^*$ is not a Dirac structure on $X$.
\end{exa}
\section{Dirac structures of divisor-type}
\label{sec:diracdivisors}
In this section we discuss how the divisors of Chapter \ref{chap:divisors} interact with Dirac structures. This is similar to Section \ref{sec:poissondivisors}, and is partially inspired by the work of Blohmann \cite{Blohmann14}.

Let $X$ be a manifold with split exact Courant algebroid $\mc{E} \cong (\mathbb{T}X, H)$ and let $K$ be a pure spinor line. Then $K \subseteq \wedge^\bullet T^*X$ so that the dual $K^*$ has a canonical section $s \in \Gamma(K^*)$, defined by $s(\rho) := \rho_0$ for $\rho \in \Gamma(K)$, where $\rho_0 \in \wedge^0 T^*X = C^\infty(X)$ is the degree-zero part of $\rho$. As a Dirac structure $E$ is described by a pure spinor line $K_E$ satisfying $E = {\rm Ann}(K_E)$, we relate Dirac structures to divisors as follows.
\begin{defn} A Dirac structure $E$ is of \emph{divisor-type} if $(K_E^*, s)$ is a divisor.
\end{defn}
Let $E$ be a Dirac structure of divisor-type. Then $s$ has nowhere dense zero set $Z_s$, so that $E$ has type zero on $X \setminus Z_s$.
\begin{rem} The notion of a Dirac structure of log divisor-type has been considered in \cite{Blohmann14} under the name log-Dirac structure.
\end{rem}
Recall that a Poisson structure $\pi \in {\rm Poiss}(X)$ defines an associated Dirac structure $E_\pi = {\rm Gr}(\pi)$. Its associated pure spinor line $K_\pi$ is locally generated by the spinor $\varphi_\pi = e^{-\pi} \cdot {\rm vol}_X$. A consequence of this is the following.
\begin{prop} Let $\pi \in {\rm Poiss}(X)$, and let $n = \dim X$. Then $\pi$ is of divisor-type if and only if the Dirac structure ${\rm Gr}(\pi)$ is of divisor-type.
\end{prop}
\bp Let $\varphi_\pi$ be as above. Then $s(\varphi_\pi) = (e^{-\pi} \cdot {\rm vol}_X)_0 = (-\pi)^n \cdot {\rm vol}_X$.
\ep
\begin{rem} The previous proposition is proven in the log divisor case in \cite{Blohmann14}.
\end{rem}
We hope to explore Dirac structures of divisor-type more in future work, especially in relation with the notion of lifting as described in the next section.
\section{Lifting Dirac structures}
\label{sec:liftingdirac}
In this section we briefly discuss the notion of lifting Dirac structures. The discussion in the previous section motivates the following definition, analogous to \autoref{defn:poissonalift}. Fix a Lie algebroid morphism $(\varphi,{\rm id}_X)\colon \mc{A} \to \mc{A}'$.
\begin{defn} An $\mc{A}'$-Dirac structure $E_{\mc{A}'} \subseteq \mathbb{A}'$ relative to $H_{\mc{A}'}$ is \emph{$\mc{A}$-liftable} if there exists a $\mc{A}$-two-form $B_\mc{A}$ and an $\mc{A}$-Dirac structure $E_\mc{A} \subseteq \mathbb{A}$ relative to the three-form $H_\mc{A} := \varphi^* H_{\mc{A}'} - d_\mc{A} B_\mc{A}$ such that $\mf{F}_\varphi(E_\mc{A}) = E_{\mc{A}'}$.
\end{defn}
As for $\mc{A}$-lifting of Poisson structures, we call $E_{\mc{A}}$ the \emph{$\mc{A}$-lift} of $E_{\mc{A}'}$, and say that $E_{\mc{A}'}$ is of \emph{$\mc{A}$-type}. For the above reason, we will often restrict to $\mc{A}$-Dirac structures which are obtained as $\mc{A}$-lifts of Dirac structures in $\mathbb{T}X$. Note that then $H_\mc{A}$ is `essentially smooth', as it is the pullback of a smooth three-form $H \in \Omega^3_{\rm cl}(M)$ via the anchor, up to addition by an exact $\mc{A}$-three-form.

The forward images of $\mc{A}$-Dirac structures which are graphs of $\mc{A}$-Poisson bivectors behave nicely with respect to liftings of the Poisson bivectors. In fact, the following also holds for almost-Dirac structures and bivectors. In the proposition below, we set $H_{\mc{A}'} = H_\mc{A} = 0$.
\begin{prop} Let $\pi_{\mc{A}'} \in {\rm Poiss}(\mc{A}')$. Then $\pi_{\mc{A}'}$ is of $\mc{A}$-type if and only if $E_{\pi_{\mc{A}'}}$ is of $\mc{A}$-type, and if $\pi_\mc{A}$ is an $\mc{A}$-lift of $\pi_{\mc{A}'}$, then $\mf{F}_\varphi(E_{\pi_\mc{A}}) = E_{\pi_{\mc{A}'}}$.
\end{prop}
\bp Let $W + \eta \in \mf{F}_\varphi(E_{\pi_\mc{A}})$. Then there exists $V + \xi \in \mb{A}$ such that $\varphi(V) = W$, and $V = \pi_\mc{A}^\sharp(\xi)$. Moreover we have $\pi_{\mc{A}'}^\sharp(\eta) = \varphi (\pi_\mc{A}^\sharp (\varphi^* \eta)) = \varphi(\pi_{\mc{A}}^\sharp(\xi)) = \varphi(X) = W$. We conclude that $W + \eta \in E_{\pi_{\mc{A}'}}$ hence $\mf{F}_\varphi(E_{\pi_\mc{A}}) \subseteq E_{\pi_{\mc{A}'}}$. By counting dimensions we in fact have equality.
\ep
In future work we will study how the lifting procedure interacts with Dirac structures of divisor-type, as we pursued in Section \ref{sec:liftingpoisson} for Poisson structures.
\chapter{Generalized complex structures}
\label{chap:gencomplexstrs}
In this chapter we discuss the theory of generalized complex structures \cite{Hitchin03, Gualtieri04, Gualtieri11}. These are geometric structures on exact Courant algebroids which can be characterized as particular kinds of complex Dirac structures. Generalized complex structures capture aspects of both symplectic and complex geometry. While slightly misleading, a generalized complex structure can be seen as a Poisson structure together with a suitably compatible complex structure normal to the (singular) symplectic leaves. In general the symplectic leaves may have varying dimension, leading to the notion of type change, where complex and symplectic behavior is mixed.

Our main interest in this thesis is a class of generalized complex manifolds for which the type need not be constantly equal to zero, but differs from zero only in the mildest way possible. These are called \sgcs{}s \cite{CavalcantiGualtieri15,Goto16}, and are those \gcs{}s whose defining anticanonical section vanishes transversally. Whenever a \gcs is of type zero, it is isomorphic after a $B$-field transformation to a symplectic structure. In this sense a \sgcs is a symplectic-like structure that fails to be symplectic on at most a codimension-two submanifold.

Particular emphasis is put on the spinor description of \gcs{}s, due to our study of \sgcs{}s. Associated to any \sgcs{} is a generalized Calabi--Yau structure (of type one), so that we are led to consider such structures as well. This chapter contains no new results apart from the short discussion in Section \ref{sec:agcs}, and follows the references \cite{Cavalcanti07, CavalcantiGualtieri15, Gualtieri11}.
\subsection*{Organization of the chapter}
This chapter is built up as follows. In Section \ref{sec:lingcs} we discuss the \gcs{}s on vector spaces, in analogy with Section \ref{sec:diraclinearalgebra} for Dirac structures. In Section \ref{sec:gcss} we globalize to manifolds, and discuss the standard examples and pointwise splitting results. Moreover, we briefly discuss generalized Calabi--Yau structures. In Section \ref{sec:scgs} we then discuss \sgcs{}s and recall how they can be viewed as Lie algebroid symplectic structures for the complex log- or elliptic tangent bundle. In Section \ref{sec:agcs} we very briefly note how results from \cite{CavalcantiGualtieri15} can possibly be extended to \gcs{}s with similar behavior.
\section{Linear generalized complex structures}
\label{sec:lingcs}
In this section we discuss the linear algebra of generalized complex geometry. We should call these structures \emph{linear} \gcs{}s, but we will speak of just \gcs{}s for brevity. Let $\mc{D}V$ be a fixed double of $V$.
\begin{defn} A \emph{\gcs} on $V$ is a linear complex structure $\mc{J}$ on $\mc{D}V$ that is orthogonal with respect to $\langle \cdot, \cdot \rangle$, i.e.\ such that $\langle J \cdot, J \cdot \rangle = \langle \cdot, \cdot \rangle$.
\end{defn}
In other words, we have $\mc{J}\colon \mc{D}V \to \mc{D}V$ such that $\mc{J}^2 = - {\rm Id}_{\mc{D}V}$.
\begin{rem} As $\mc{J}$ is orthogonal, nondegeneracy of the pairing shows that $\langle \mc{J} \cdot, \cdot \rangle$ is a linear symplectic structure on $\mc{D}V$. Hence an alternative name would be a \emph{generalized symplectic} structure.
\end{rem}
Note that $\mc{J}$ admits an extension to the complexification $\mc{D}V \otimes \C$ of $\mc{D}V$ by setting $\mc{J}(v \otimes z) = \mc{J} v \otimes z$. As usual, this means $\mc{J}$ splits $\mc{D}V \otimes \C$ into a direct sum of $\pm i$-eigenspaces $E_\mc{J}$ and $E_\mc{J}'$. Explicitly, we have $E_\mc{J} = \{v - i \mc{J} v \, | \, v \in \mc{D}V\}$. Complex conjugation shows that $\overline{E}_\mc{J} = E_\mc{J}'$. Given a splitting $\mc{D}V \cong \mathbb{V}$, we can write $\mc{J} \in {\rm End}(\mathbb{V})$ in matrix form as
\be
\mc{J} = \begin{pmatrix} A & \pi_\mc{J}^\sharp \\ B & -A^* \end{pmatrix}.
\ee
Here $\pi_\mc{J} \in \wedge^2 V$ is a bivector, which we view as a map $\pi^\sharp_\mc{J}\colon V^* \to V$. Note that the image $D = \rho_V(\mc{J} V^*)$ is the image of $\pi_\mc{J}^\sharp$, which turns out to be a symplectic distribution. Its annihilator ${\rm Ann}(D) = V^* \cap \mc{J} V^*$ is a complex subspace, and $\mc{E} / {\rm Ann}(D) \cong V / D$, so that normal to $D$ there is a complex structure. One quickly realizes that these eigenspace are maximally isotropic.
\begin{lem} Let $\mc{J}$ be a \gcs on $V$. Then $E_\mc{J}$ is maximally isotropic for $\langle \cdot, \cdot \rangle_{\mc{D}V}$.
\end{lem}
\bp Let $v, w \in E_\mc{J}$. Then $\langle v, w \rangle = \langle \mc{J} v, \mc{J} w \rangle = \langle i v, i w \rangle = - \langle v, w \rangle$, so that $\langle v ,w \rangle = 0$. We conclude that $E_\mc{J} \subseteq E_\mc{J}^\perp$, so that $E_\mc{J}$ is isotropic. Because $\overline{E}_\mc{J} = E_\mc{J}'$ and $\mc{D}V \otimes \C \cong E_\mc{J} \oplus E_\mc{J}'$, we see that $\dim_\C(E_\mc{J}) = m$, i.e.\ $E_\mc{J}$ is maximal.
\ep
Hence $E_\mc{J}$ is a linear complex Dirac structure in $\mc{D}V$. We can determine which linear complex Dirac structures are eigenspaces of \gcs{}s.
\begin{prop} Let $E_\mathbb{C}$ be a linear complex Dirac structure in $\mc{D}V$. Then $E_\C = E_\mc{J}$ for a \gcs $\mc{J}$ if and only if $E_\C \cap \overline{E}_\C = 0$.
\end{prop}
\bp Certainly if $E_\C = E_\mc{J}$, then $E_\C \cap \overline{E}_\C = E_\C \cap E_\C' = 0$, as eigenspaces are disjoint. Conversely, given $E_\C$ such that $E_\C \cap \overline{E}_\C = 0$, we have a direct sum decomposition $\mc{D}V \otimes \C \cong E_\C \oplus \overline{E}_\C$. Using this decomposition, define $\mc{J} \in {\rm End}(\mc{D}V)$ by $\mc{J}|_{E_\C} = i$ and $\mc{J}|_{\overline{E}_\C} = - i$. Then $\mc{J}^2 = - {\rm Id}$, and moreover $\langle \mc{J} E_\C, \mc{J} E_\C \rangle = 0$ and similarly for $\overline{E}_\C$ because $\mc{J}$ preserves these isotropic subspaces. Finally, for $v \in E_\C$ and $w \in \overline{E}_\C$ we have $\langle \mc{J} v, \mc{J} w \rangle = \langle i v, -i w \rangle = - i^2 \langle v, w \rangle = \langle v, w \rangle$, so that $\mc{J}$ is orthogonal with respect to $\langle \cdot, \cdot \rangle$. But then $\mc{J}$ is a \gcs{} on $V$ with $E_\C = E_\mc{J}$.
\ep
We discussed before that linear Dirac structures admit a description in terms of spinors. We explore which pure spinors give rise to \gcs{}s.
\begin{lem} Let $\rho$ be a complex pure spinor. Then $\overline{E}_\rho = E_{\overline{\rho}}$.
\end{lem}
Combining this with \autoref{prop:pairingintersection} we obtain the following.
\begin{prop} Let $\rho = e^{B + i \omega} \wedge \Omega$ be a complex pure spinor. Then $E_\rho$ defines a \gcs{}s $\mc{J}_\rho$ if and only if $(\rho, \overline{\rho}) \neq 0$.
\end{prop}
Hence for $V$ to admit a \gcs{}, its dimension $m$ must be even, say $m = 2n$ (as only then can $(\rho, \overline{\rho})$ be nonzero). Then $(\rho, \overline{\rho}) = (e^{2 i \omega} \wedge \Omega, \overline{\Omega})$, which is equal to $\Omega \wedge \overline{\Omega} \wedge \omega^{n-k}$ up to constants, where $k$ is the degree of $\Omega$.
\begin{defn} Let $\mc{J}$ be a \gcs{} on $V$ where $\dim_\R V = 2n$ and let $E_\mc{J} = E_\rho$ with $\rho = e^{B + i \omega} \wedge \Omega$ with respect to a splitting $\mc{D}V \cong \mathbb{V}$. The \emph{type} of $\mc{J}$ is $k$, the degree of $\Omega$. The line $K_\mc{J} = K_\rho$ is the \emph{canonical line} of $\mc{J}$.
\end{defn}
The type is thus an integer ranging from $0$ to $n$, which is half the dimension of $V$. We also have ${\rm type}(\mc{J}) = \frac12 \dim_\R (V^* \cap \mc{J} V^*) = {\rm type}(E_\mc{J})$. We discuss how the $B$-field  transform and $\beta$-transform interact with a \gcs.
\begin{prop} Let $\mc{J}$ be a \gcs{} with $+i$-eigenspace $E_\mc{J}$, and let $B \in \wedge^2 V^*$ be a two-form. Then $E_{\mc{J},B}$ satisfies $E_{\mc{J},B} \cap \overline{E}_{\mc{J},B} = 0$, hence defines a \gcs{} $\mc{J}_B$ such that $E_{\mc{J}_B} = E_{\mc{J},B}$. Given a splitting $\mc{D}V \cong \mathbb{V}$, we have
	\be
	\mc{J}_B = \begin{pmatrix} 1 & 0 \\ -B & 1 \end{pmatrix} \mc{J} \begin{pmatrix} 1 & 0 \\ B & 1 \end{pmatrix}.
	\ee
\end{prop}
\bp As $B$ is real, we have $E_{\mc{J},B} \cap \overline{E}_{\mc{J},B} = ({\rm Id} - B) (E_\mc{J} \cap \overline{E}_\mc{J}) = 0$. That $E_{\mc{J},B}$ is maximally isotropic follows from \autoref{prop:bfieldtransform}.
\ep
A similar statement is true for bivectors.
\begin{prop} Let $\mc{J}$ be a \gcs{} with $+i$-eigenspace $E_\mc{J}$, and let $\beta \in \wedge^2 V$ be a bivector. Then $E_{\mc{J},\beta}$ satisfies $E_{\mc{J},\beta} \cap \overline{E}_{\mc{J},\beta} = 0$, hence defines a \gcs{} $\mc{J}_\beta$ such that $E_{\mc{J}_{\beta}} = E_{\mc{J},\beta}$. Given a splitting $\mc{D}V \cong \mathbb{V}$, we have
	\be
	\mc{J}_\beta = \begin{pmatrix} 1 & -\beta \\ 0 & 1 \end{pmatrix} \mc{J} \begin{pmatrix} 1 & \beta \\ 0 & 1 \end{pmatrix}.
	\ee
\end{prop}
%
In conclusion, a \gcs on a split double $\mc{D}V \cong \mathbb{V}$ of $V$ is defined by one of the following three equivalent definitions:
\bi
\item An orthogonal linear complex structure $\mc{J}$ on $\mc{D}V$;
\item An complex linear Dirac structure $E_\mc{J}$ such that $E_\mc{J} \cap \overline{E}_\mc{J} = 0$;
\item A complex pure spinor line $K_\mc{J}$ such that $(K_{\mc{J}}, \overline{K}_{\mc{J}}) \neq 0$.
\ei
Next we turn to describing some examples of \gcs{}s.
\begin{exa} Let $\mc{D}V \cong V \oplus V^*$ be an isotropic splitting for $\mc{D}V$ and $J$ a complex structure on $V$, i.e.\ an endomorphism $J\colon V \to V$ such that $J^2 = - {\rm Id}_V$. Then $J$ gives rise to a \gcs{} $\mc{J}_J$ on $V$ which using the splitting can be written in matrix form as
	\be
	\mc{J}_J = \begin{pmatrix} - J & 0 \\ 0 & J^* \end{pmatrix}
	\ee
	The $+i$-eigenspace of $\mc{J}_J$ is given by $E_{\mc{J}_J} = V^{0,1} \oplus V^{*1,0} \subseteq \mathbb{V}_\C$. It is immediate that $E_{\mc{J}_J}$ is maximally isotropic and satisfies $E_{\mc{J}_J} \cap \overline{E}_{\mc{J}_J} = 0$. The canonical line of $\mc{J}_J$ is given by $K_{\mc{J}_J} = \wedge^{n,0} V^*$, the canonical line of $I$ in the usual sense. The type of $\mc{J}_J$ is $n$. The action of a $B$-field transformation results in
	\be
	e^{-B} \mc{J}_J e^{B} = \begin{pmatrix} -J & 0 \\ B J + J^* B & I^* \end{pmatrix},
	\ee
	and $e^{-B} E_{\mc{J}_J} = \{X + \xi - \iota_X B \, | \, X + \xi \in V^{0,1} \oplus V^{*1,0}\}$. Only the $(0,2)$-part of $B$ has any effect.
\end{exa}
Given a type $n$ \gcs{}, its pure spinor is of the form $\rho = e^{B + i \omega} \wedge \Omega$ with $\Omega$ of degree $n$. As $\Omega \wedge \overline{\Omega} \neq 0$, the form $\Omega$ determines a complex structure on $V$ for which it is of type $(n,0)$. Note that $e^{B + i \omega} \wedge \Omega = e^c \wedge \Omega$, where $c = (B + i \omega)^{0,2}$, the type being determined by this complex structure. Let $B' := c + \overline{c}$. Then $\rho = e^{B'} \wedge \Omega$, so that we conclude that any generalized complex structure of type $n$ is the $B$-field transform of a complex structure.
\begin{exa} Let $\mc{D}V \cong V \oplus V^*$ be an isotropic splitting for $\mc{D}V$ and $\omega$ a symplectic structure on $V$, i.e.\ a two-form $\omega \in \wedge^2 V^*$ such that $\omega^\flat\colon V \to V^*$ is an isomorphism. Then $\omega$ gives rise to a \gcs{} $\mc{J}_\omega$ on $V$ which using the splitting can be written in matrix form as
	\be
	\mc{J}_\omega = \begin{pmatrix} 0 & - \omega^{-1} \\ \omega & 0 \end{pmatrix}
	\ee
	The $+i$-eigenspace of $\mc{J}_\omega$ is given by $E_{\mc{J}_\omega} = \{X - i \omega^\flat(X) \, | \, X \in V \otimes \C\} = e^{-i\omega} V$. As $\omega$ is nondegenerate we obtain $E_{\mc{J}_\omega} \cap \overline{E}_{\mc{J}_\omega} = 0$. Moreover, $\omega$ being skew-symmetric implies that $E_{\mc{J}_\omega}$ is isotropic. The canonical line of $\mc{J}_\omega$ is given by $K_{\mc{J}_\omega} = \langle e^{i \omega} \rangle$. The type of $\mc{J}_\omega$ is zero. The action of a $B$-field transformation results in
	\be
	e^{-B} \mc{J}_\omega e^B = \begin{pmatrix} - \omega^{-1} B & - \omega^{-1} \\ \omega + B \omega^{-1} B & B \omega^{-1} \end{pmatrix},
	\ee
	and $e^{-B} E_{\mc{J}_\omega} = \{X - (B + i \omega)(X) \, | \, X \in V \otimes \C\} = e^{-(B + i \omega)} V$ and $e^{-B} \rho = e^{B + i \omega}$, which again describe a \gcs{} of type $0$.
\end{exa}
Note that the description in terms of spinors shows that a \gcs of type $0$ is a $B$-field transform of a symplectic structure.
\section{\Gcs{}s}
\label{sec:gcss}
We can now consider \gcs{}s on manifolds proper, in analogy with our treatment of Dirac structures in Chapter \ref{chap:diracgeometry}. Let $X$ be a $2n$-dimensional manifold and $\mc{E} \to X$ an exact Courant algebroid which we will consider to be split as $\mc{E} \cong (\mathbb{T}X, H)$ for some closed three-form $H \in \Omega^3_{\rm cl}(X)$.
\begin{defn} A \emph{\gcs} on $(X,H)$ is a complex structure $\mc{J}$ on $\mathbb{T}X$ that is orthogonal with respect to $\langle\cdot,\cdot\rangle$, and whose $+i$-eigenbundle is involutive under $[[\cdot,\cdot]]_H$.
\end{defn}
There is an alternative definition of a \gcs using spinors. Recall that sections $v = V + \alpha \in \Gamma(\mathbb{T}X)$ of $\mathbb{T}X$ act on differential forms via Clifford multiplication, given by $v \cdot \rho = \iota_V \rho + \alpha \wedge \rho$ for $\rho \in \Omega^\bullet(X)$.
\begin{defn} A \emph{\gcs} on $(X,H)$ is a pure spinor line $K_{\mc{J}} \subset \Omega^\bullet_\C(X)$ pointwise generated by a pure spinor $\rho = e^{B + i \omega} \wedge \Omega$ for which $\Omega \wedge \overline{\Omega} \wedge \omega^{n-k} \neq 0$, and such that $d \rho + H \wedge \rho = v \cdot \rho$ for any local section $\rho \in \Gamma(K_{\mc{J}}) \setminus \{0\}$ and some $v \in \Gamma(\mathbb{T}X)$.
\end{defn}
Both definitions are related using that $K_\mc{J} = {\rm Ann}(E_\mc{J})$ is the annihilator under the Clifford action of the complex Dirac structure $E_\mc{J}$, the $+i$-eigenbundle of $\mc{J}$, which satisfies $E_\mc{J} \cap \overline{E}_{\mc{J}} = 0$. The bundle $K_{\mc{J}}$ is called the \emph{canonical bundle} of $\mc{J}$. 
\begin{exa}\label{exa:gcss} The following are \gcs{}s on $(X,H = 0)$.
	\bi
	\item Let $\omega$ be a symplectic structure on $X$. Then $K_{\mc{J}_\omega} := \langle e^{i \omega} \rangle$ defines a \gcs $\mc{J}_\omega$.
	\item Let $J$ be a complex structure on $X$ with canonical bundle $K_J = \wedge^{n,0} T^* X$. Then $K_{\mc{J}_J} := K_J$ defines a \gcs $\mc{J}_J$.
	\item Let $P \in \Gamma(\wedge^{2,0} TX)$ a holomorphic Poisson structure with respect to a complex structure $J$. Then $K_{\mc{J}_{P,J}} := e^{P} K_J$ defines a \gcs $\mc{J}_{P,J}$.
	\ei
	The automorphisms $\mc{J}_\omega$, $\mc{J}_J$ and $\mc{J}_{P,J}$ are given by, with $\pi = {\rm Im}(P)$:
	\begin{align*}
		\mc{J}_\omega = \begin{pmatrix} 0 & -\omega^{-1} \\ \omega & 0\end{pmatrix},  \qquad\qquad &	\mc{J}_J = \begin{pmatrix} -J & 0 \\ 0 & J^* \end{pmatrix}, \qquad \qquad \mc{J}_{P,J} = \begin{pmatrix} -J & \pi \\ 0 & J^* \end{pmatrix}.
	\end{align*}
	
\end{exa}
A \gcs constitutes a reduction of the structure group of $\mathbb{T}X$ to $U(n,n) \subseteq {\rm SO}(2n,2n)$. Realizing this, it is not hard to prove the following.
\begin{prop}[c.f.\ {\cite[Proposition 3.4]{Gualtieri11}}] Let $X$ be a generalized complex manifold. Then $X$ admits an almost-complex structure.
\end{prop}
It is unknown whether the converse also holds. More precisely, as of this writing there is no example of an almost-complex manifold which is not generalized complex.
\begin{exa} For the manifolds $X_{m,n} = m \C P^2 \# n \overline{\C P}^2$ with $m,n \geq 0$, this is the only obstruction \cite{CavalcantiGualtieri07}, i.e.\ there exists a \gcs on $X_{m,n}$ if and only if $m$ is odd (which is precisely when it is almost-complex, see \autoref{rem:acs}). Note that $X_{m,n}$ is symplectic or complex if and only if $m = 1$.
\end{exa}
We next introduce the type of a \gcs$\mc{J}$, which colloquially provides a measure for how many complex directions there are. The type is an integer-valued upper semicontinuous function on $X$ whose parity is locally constant.
\begin{defn} Let $\mc{J}$ be a \gcs on $X$. The \emph{type} of $\mc{J}$ is a map ${\rm type}(\mc{J})\colon X \to \Z$ whose value at a point $x \in X$ is the integer $k$ above, the degree of $\Omega$. Alternatively, it is given by ${\rm type}_x(E_\mc{J})$. The \emph{type change locus} $D_{\mc{J}}$ of $\mc{J}$ is the subset of $X$ where ${\rm type}(\mc{J})$ is not locally constant.
\end{defn}
Any \gcs $\mc{J}$ determines a Poisson structure $\pi_\mc{J}$ instead of just a bivector as follows, as was realized early on in the study of \gcs{}s (see e.g.\ \cite[Proposition 3.21]{Gualtieri11} for a Dirac-geometric proof).
\begin{defn} Let $\mc{J}$ be a \gcs on $X$. The \emph{underlying Poisson structure} $\pi_\mc{J}$ of $\mc{J}$ the bivector defined by $\pi_{\mc{J}}^\sharp := p \circ \mc{J}|_{T^*X}$.
\end{defn}
The type of $\mc{J}$ is related to the rank of $\pi_{\mc{J}}$ through the formula ${\rm rank}(\pi_{\mc{J}}) = 2n - 2 \, {\rm type}(\mc{J})$. Using $\pi_{\mc{J}}$ one can view a \gcs $\mc{J}$ as specifying a singular foliation on $X$ with symplectic leaves, and a suitably compatible complex structure normal to the leaves. Call a point $x \in X$ \emph{regular} for $\mc{J}$ if $\pi_{\mc{J}}$ is regular at $x$, i.e.\ ${\rm type}(\mc{J})$ is locally constant at $x$.

The complement of $D_{\mc{J}}$ is an open dense set where the type is minimal. Using the type, \gcs{}s are seen to interpolate between symplectic and complex structures. Recall that given $B \in \Omega^2_{\rm cl}(X)$ we have its $B$-field transformation, given by $e^B\colon \mathbb{T}X \to \mathbb{T}X$, $e^B\colon V + \alpha \mapsto V + \alpha + \iota_V B$.
\begin{lem}\label{lem:gcstype0} Let $(X,H,\mc{J})$ be a generalized complex manifold and $x \in X$ such that ${\rm type}_x(\mc{J}) = 0$. Then there is an open neighbourhood $U \subseteq X$ of $x$ and a closed two-form $B \in \Omega^2_{\rm cl}(U)$ such that $e^B \mc{J} = \mc{J}_{\omega}$ for $\omega$ a symplectic structure on $U$.
\end{lem}
A similar statement holds at points of maximal type.
\begin{lem}\label{lem:gcstypemax} Let $(X^{2n},H,\mc{J})$ be a generalized complex manifold and $x \in X$ such that ${\rm type}_x(\mc{J}) = n$. Then there exists an open neighbourhood $U \subseteq X$ of $x$ and a closed two-form $B \in \Omega^2_{\rm cl}(U)$ such that $e^B \mc{J} = \mc{J}_{J}$ for $J$ a complex structure on $U$.
\end{lem}
Both results have global consequences, relating \gcs{} of constant type $0$ or $n$ to the symplectic respectively complex examples found in \autoref{exa:gcss}. Both \autoref{lem:gcstype0} and \autoref{lem:gcstypemax} follow from the following more general statements describing the local behavior of \gcs{}s. In the statements of \autoref{thm:gcssplittingreg} and \autoref{thm:gcssplitting}, by `equivalence' we mean equivalence up to $B$-field transformation.
\begin{thm}[Gualtieri, \cite{Gualtieri11}]\label{thm:gcssplittingreg} Let $\mc{J}$ be a \gcs{} on $X$ and $x \in X$ a regular point of type $k$ for $\mc{J}$. Then there exists an open neighbourhood $U$ of $x$ and an equivalence $(U,\mc{J}) \cong (\C^k, \mc{J}_J) \times (\R^{2n-2k}, \mc{J}_{\omega_0})$, where $\omega_0$ is the standard symplectic structure.
\end{thm}
The proof uses the Weinstein splitting theorem for $\pi_\mc{J}$ (see \autoref{thm:poissonsplitting}), and an understanding of the induced complex structure on the leaf space of $\pi_\mc{J}$ near $x$. A large strengthening of this result in which regularity is dropped, encompassing in particular the Newlander--Nirenberg theorem for complex structures, is the following.
\begin{thm}[Bailey, \cite{Bailey13}]\label{thm:gcssplitting} Let $\mc{J}$ be a \gcs{} on $X$ and $x \in X$ a point of type $k$ for $\mc{J}$. Then there exists an open neighbourhood $U$ of $x$ and an equivalence $(U,\mc{J}) \cong (\C^k, P, \mc{J}_{P,J}) \times (\R^{2n-2k}, \mc{J}_{\omega_0})$, where $P$ is a holomorphic Poisson structure which vanishes at $0$, and $\omega_0$ is the standard symplectic structure.
\end{thm}
We should also mention a partial result due to Abouzaid--Boyarchenko \cite{AbouzaidBoyarchenko06}.
\subsection{Generalized Calabi--Yau structures}
For later use, we further introduce the analogue of a Calabi--Yau manifold in generalized geometry. Denote by $d^H = d + H \wedge$ the $H$-twisted de Rham differential.
\begin{defn}[\cite{Hitchin03}] A generalized complex structure $\mc{J}$ on $(X,H)$ is \emph{generalized Calabi--Yau} if $K_\mc{J}$ has a global nowhere vanishing $d^H$-closed pure spinor.
\end{defn}
The canonical bundle of a generalized Calabi--Yau structure is not only trivial, but further holomorphically trivial, in that the induced generalized holomorphic structure it carries is the standard one. We refer the reader to \cite{CavalcantiGualtieri15} for further discussion.
\begin{exa} Let $\rho = e^{i \omega}$ be the spinor of a symplectic structure, and $H = 0$. This is a globally defined spinor, so $K_{\mc{J}_\omega}$ is trivial. Then $d \rho = d e^{i \omega} = i d \omega e^{i \omega} = 0$. We see thus that any symplectic structure is generalized Calabi--Yau.
\end{exa}
\begin{exa} Let $J$ be a complex structure on $X$. Then $(X,0)$ is generalized complex. As $K_{\mc{J}_J} = K_J$ and any (local) generator is automatically closed, we see that $\mc{J}_J$ is generalized Calabi--Yau if and only if $J$ is Calabi--Yau.
\end{exa}
The main result we want to mention here is that any compact type-1 generalized Calabi--Yau manifold fibers over the torus $T^2$. This is similar to the case of compact cosymplectic manifolds mentioned in Section \ref{sec:logsympstr}.
\begin{thm}[\cite{CavalcantiGualtieri15,BaileyCavalcantiGualtieri16}]\label{thm:type1gency} Let $(X,H,\mc{J})$ be a compact type-1 generalized Calabi--Yau manifold. Then there exists a surjective submersion $p\colon X \to T^2$, and if $X$ has a compact leaf the map $p$ can be chosen such that its fibers are $\pi_{\mc{J}}$-symplectic leaves.
\end{thm}
The proof relies on extracting two independent closed one-forms out of $\mc{J}$ and applying the Tischler argument.
\begin{rem} In \cite{Hitchin03}, Hitchin mentions that compact four-dimensional generalized Calabi--Yau manifolds of odd type fiber over $T^2$.
\end{rem}
\section{Stable generalized complex structures}
\label{sec:scgs}
Generalized complex structures which are stable were introduced in \cite{CavalcantiGualtieri15, Goto16}. Their defining property is a natural condition and since stable generalized complex structures are close to being symplectic, one can use symplectic techniques to study them.

In \cite{CavalcantiGualtieri09, GotoHayano16, Torres12, TorresYazinski14}, many examples of \sgcs{}s on four-manifolds were constructed, in particular on manifolds without a symplectic nor complex structure. A thorough study of stable generalized complex manifolds was initiated in \cite{CavalcantiGualtieri15}. \Sgcs{}s can be seen as the generalized geometric analogue of a \blog (or log-Poisson structure) discussed in Section \ref{sec:logsympstr} and Section \ref{sec:logpoisson}.

Let $\mc{J}$ be a \gcs on $(X,H)$. The anticanonical bundle $K_{\mc{J}}^*$ has a section $s \in \Gamma(K_{\mc{J}}^*)$, given by $s(\rho) := \rho_0$ for $\rho \in \Gamma(K_{\mc{J}})$, with $\rho_0$ the degree-zero part of $\rho$.
\begin{defn} A \gcs $\mc{J}$ on $(X,H)$ is \emph{stable} if $s$ is transverse to the zero section in $K_{\mc{J}}^*$. The set $D_{\mc{J}} := s^{-1}(0)$ is a codimension-two smooth submanifold of $X$ called the \emph{anticanonical divisor} of $\mc{J}$.
\end{defn}
\begin{rem} Using the language of Chapter \ref{chap:divisors} and Section \ref{sec:diracdivisors}, a \gcs{} is stable if its complex Dirac structure is of complex log divisor-type.
\end{rem}
Outside of $D_{\mc{J}}$, the section $s$ is nonvanishing hence the type of $\mc{J}$ is equal to zero, while over $D_{\mc{J}}$ it is equal to two. Consequently, \sgcs{} can be seen as \gcs{}s which are close to being symplectic.

\begin{exa}\label{exa:sgcsc2} Consider $(\C^2,0)$ with coordinates $(z,w)$ and with holomorphic Poisson structure $\pi = z \partial_z \wedge \partial_w$. This gives a \sgcs with $K_{\mc{J}} = \langle z + dz \wedge dw \rangle$ and $D_\mc{J} = \{z = 0\}$. As the spinor $\rho = z + dz \wedge dw$ is equal to $dz \wedge dw$ along $D_{\mc{J}}$, the type of $\mc{J}$ is equal to two there. Away from $D_\mc{J}$, we can write $\rho = z \exp(d \log z \wedge dw)$, so that there the type is seen to be zero. Note that $d \rho = dz = - \iota_{\partial w} \rho$, so that $\mc{J}$ is indeed integrable.
\end{exa}
By \autoref{exa:gcss}, any holomorphic Poisson structure $(J,P)$ defines a \gcs $\mc{J}_{P,J}$ on $(X,0)$ with $K_{\mc{J}_{P,J}} = e^P K_J$. Thus $K_{\mc{J}_{P,J}}$ is locally generated by $e^P \Omega$, with $\Omega$ a local trivialization of $K_J$. This proves the following.
\begin{prop}[{\cite[Example 2.17]{CavalcantiGualtieri15}}] Let $(X,0,J,P)$ be a holomorphic Poisson structure. Then $\mc{J}_{P,J}$ is a \sgcs if and only if the Pfaffian $\wedge^n P$ is transverse to the zero section in $\Gamma(\wedge^{2n,0} TX)$.
\end{prop}
Let $(X,H,\mc{J})$ be a stable generalized complex manifold with type change locus $D_\mc{J}$. Then $D_\mc{J}$ is the degeneracy locus of $\pi_{\mc{J}}$ and is in fact a \emph{strong generalized Poisson submanifold} for $\mc{J}$. This means it inherits a generalized complex structure.
\begin{prop}[{\cite[Theorem 2.19]{CavalcantiGualtieri15}}] Let $(X,H,\mc{J})$ be a stable generalized complex manifold. Then $D_{\mc{J}}$ inherits a type-1 generalized Calabi--Yau structure.
\end{prop}
This structure is obtained via a residue procedure. Using \autoref{thm:type1gency}, we see that $D_\mc{J}$ fibers over $T^2$. There is further information induced on the normal bundle $ND_{\mc{J}}$, but we will not discuss this here. Let us describe the local form of a \sgcs.
\begin{prop}[\cite{CavalcantiGualtieri15}] Let $(X,H,\mc{J})$ be a stable generalized complex manifold and $x \in D_{\mc{J}}$. Then there exists a neighbourhood $U \subseteq X$ of $x$ and a generalized holomorphic map $\psi\colon U \to (\C^2 \times \R^{n-4}, 0,\mc{J}_{\rm st})$, where $K_{\mc{J}_{\rm st}} = K_{\mc{J}_{P_0,J}} \wedge K_{\mc{J}_{\omega_0}}$, where $P_0$ is the holomorphic Poisson structure of \autoref{exa:sgcsc2}.
\end{prop}
\bp Let $n = {\rm type}_x(\mc{J})$. By \autoref{thm:gcssplitting} we have locally around $x$ that $\mc{J} = \mc{J}_{P,J} \oplus \mc{J}_{\omega_0}$, where $P$ is a $J$-holomorphic Poisson structure on $\C^{2n}$ vanishing at $x$, and $\omega_0$ is the standard symplectic structure on $\R^{2k}$. The spinor for $\mc{J}$ is thus given by $\rho = e^P \Omega^{2n,0} \cdot e^{i\omega}$, and $\rho_0 = P^n \Omega^{2n,0}$. As $P(0) = 0$, we see that $P^n$ has a zero of order $n$. But $\rho_0$ is transverse to zero as $\mc{J}$ is stable, so that we must have $n = 1$. But then $P$ is a holomorphic Poisson structure on standard $\C^2$ which vanishes transversally. In standard coordinates we have $P = f(z,w) \partial_z \wedge \partial_w$, but then we can find new coordinates in which $P = z \partial_z \wedge \partial_w$ and $\{z = 0\} = D_\mc{J}$.
\ep
Hence any \sgcs $\mc{J}$ is locally equivalent around points in $D_\mc{J}$ to $\langle e^{i \omega_0}(z + dz \wedge dw) \rangle$ on $\C^2 \times \R^{2n-4}$, with $\omega_0$ the standard symplectic form on $\R^{2n-4}$, and $D_{\mc{J}} = \{z = 0\}$. This also shows the type of $\mc{J}$ over $D_{\mc{J}}$ is two.

If one incorporates the induced structure on the normal bundle $ND_{\mc{J}}$, there is a semilocal form for a \sgcs around its type change locus is given by its \emph{linearization} along $D_{\mc{J}}$, which is the \sgcs naturally present on the normal bundle to this type-1 generalized Calabi--Yau manifold (see \cite[Theorem 3.32]{CavalcantiGualtieri15}). This result is proven by an application of Moser-type arguments for the complex log tangent bundle associated to the complex log divisor $D_{\mc{J}}$. Indeed, let us consider how $\mc{J}$ relates to the Lie algebroids it gives rise to. In terms of our language we have the following, which locally is the statement that $z + dz \wedge dw = z \exp(d \log z \wedge dw)$ as in \autoref{exa:sgcsc2}.
\begin{prop}[c.f.\ {\cite[Theorem 2.19]{CavalcantiGualtieri15}}] Let $(X,H,\mc{J})$ be a \sgcs. Then its canonical spinor line $K_\mc{J}$ can be seen as described by $e^\sigma$ for $\sigma \in \Omega^2(\mc{A}_{D_\mc{J}})$ a complex log-symplectic form on $\mc{A}_{D_\mc{J}}$ with respect to $H$.
\end{prop}
Following \cite{CavalcantiGualtieri15}, recall the Chevalley pairing on spinors from \autoref{defn:chevalley}. As we know that $K_\mc{J}$ is generated by $e^\sigma$ with $\sigma$ a complex log-symplectic form associated to the complex log divisor $D_\mc{J}$, the Chevalley pairing restricts to a map
\be
	(\cdot,\cdot)\colon K_\mc{J} \otimes \overline{K}_\mc{J} \to \wedge^{2n} T^*_\C X
\ee
which is an isomorphism, sending $e^\sigma \otimes e^{\overline{\sigma}}$ to $(e^{\sigma - \overline{\sigma}})_{\rm top} = (2i)^n \omega^n$ (up to a possible factor $n !$), where $\omega$ is the imaginary part of $\sigma$. This is an elliptic symplectic structure for the associated elliptic divisor determined by $D_\mc{J}$.

By dualizing we can phrase things in terms of the anticanonical divisor $(K_\mc{J}^*, s)$. The above implies $s \otimes \overline{s}$ gets sent to $(2i)^{-n} \wedge^n\pi_{\mc{J}}$, where $\pi_{\mc{J}}$ is the dual bivector to $\omega$. In other words, the Chevalley pairing gives an isomorphism of divisors $(K_\mc{J}^* \otimes \overline{K}^*_\mc{J}, s \otimes \overline{s}) \cong (\wedge^{2n} T^* X, \wedge^n \pi_{\mc{J}})$. This shows that $\pi_\mc{J}$ is an elliptic Poisson structure. In short, the underlying Poisson structure $\pi_\mc{J}$ of a \sgcs $\mc{J}$ is elliptic, and this characterizes when $\mc{J}$ is stable.
\begin{prop}[{\cite[Theorem 3.7]{CavalcantiGualtieri15}}]\label{prop:stableellpoiss} Let $(X,H,\mc{J})$ be a generalized complex manifold. Then $\mc{J}$ is a \sgcs if and only if $\pi_{\mc{J}}$ is an elliptic Poisson structure.
\end{prop}
An elliptic pair $(X,|D|)$ \emph{admits a \sgcs} if there exists a closed three-form $H \in \Omega^3_{\rm cl}(X)$ and a \sgcs $\mc{J}$ on $(X,H)$ such that $|D_{\mc{J}}| = |D|$. When no elliptic divisor structure is specified on $D$, we say that $(X,D)$ admits a \sgcs if there exists some elliptic divisor structure $|D|$ on $D$ such that $(X,|D|)$ admits a \sgcs.

In fact, \sgcs{}s $\mc{J}$ on $X$ are in one-to-one correspondence with certain types of elliptic Poisson structures $\pi$ via the map $\mc{J} \mapsto \pi_{\mc{J}}$. Moreover, the closed three-form $H$ required to state integrability of $\mc{J}$ is determined by $\pi$. We can now state the extension of \autoref{prop:stableellpoiss}, giving a characterization of \sgcs{} purely in terms of elliptic symplectic structures, through the use of complex log-symplectic structures (see Section \ref{sec:complexlogsymp}).
\begin{thm}[{\cite[Theorem 3.7]{CavalcantiGualtieri15}}]\label{thm:sgcscorrespondence} Let $X$ be a compact manifold. There is a bijection $(\mc{J}, H) \to (\pi_{\mc{J}}^{-1}, \mf{o})$ between \sgcs{}s up to gauge equivalence and elliptic symplectic structures with vanishing elliptic residue and cooriented degeneracy locus.
\end{thm}
The associated closed three-form $H$ in the definition of a \gcs can be determined via $[H] = {\rm Res}_r([\omega_{\mc{A}}]) \wedge {\rm PD}[D]$, where $\omega_{\mc{A}}$ is the Lie algebroid symplectic structure for $\mc{A} = TX(-\log |D_{\mc{J}}|)$ whose dual bivector is $\pi_{\mc{J}}$, and $D = (\wedge^n \pi_\mc{J})^{-1}(0)$. The Poincar\'e dual of $D$ requires a choice of coorientation. Moreover, the closed one-forms on $D_{\mc{J}}$ which exhibit $D_{\mc{J}}$ as fibering over $T^2$ (see \autoref{thm:type1gency}), are after perturbation given by the residue pair $({\rm Res}_r,{\rm Res}_\theta)(\omega_{\mc{A}})$.
\section{\texorpdfstring{$\mc{A}_\C$}{AC}-generalized complex structures}
\label{sec:agcs}
In this section we very briefly discuss a natural question, namely to what extent it makes sense to define generalized complex structures using Lie algebroid objects. Such a concept has been defined \cite{Barton07}, but to the author's knowledge, no use is made there of a possible link between the Lie algebroid $\mc{A}$ and the tangent bundle $TX$.

Recall that a Lie algebroid $\mc{A}$ gives rise to a Courant algebroid $\mathbb{A} = \mc{A} \oplus \mc{A}^*$.
\begin{defn} Let $\mc{A} \to X$ be a Lie algebroid and $H_\mc{A} \in \Omega^3_{\rm cl}(\mc{A})$. An \emph{$\mc{A}$-\gcs{}} is a complex structure $\mc{J}_{\mc{A}}$ on $\mathbb{A}$ whose $+i$-eigenbundle is involutive under $[[\cdot,\cdot]]_{\mc{A},H_\mc{A}}$.
\end{defn}
As of this writing, it felt premature to devote significant attention to this concept. In the future we plan to consider the merits of developing such a theory in combination with the concept of lifting. Indeed, it is immediate that one can define \gcs{}s of general complex divisor-type.
\begin{defn} A generalized complex structure $\mc{J}$ on $(X,H)$ is said to be of \emph{divisor-type} if $(K^*_\mc{J},s)$ is a divisor.
\end{defn}
A similar definition can be made for $\mc{A}$-generalized complex structures. Given a generalized complex structure of divisor-type, we can wonder how they interact with the associated complex ideal Lie algebroids by lifting. The results of Cavalcanti--Gualtieri \cite{CavalcantiGualtieri15} on \sgcs{}s can be seen as a first instance of such a procedure, for complex log divisors. As for \sgcs{}s, the underlying of Poisson structure of a \gcs{} of divisor-type is itself of divisor-type, as one sees by the same procedure explained above \autoref{prop:stableellpoiss}.
\chapter{Constructing \texorpdfstring{$\mc{A}$}{A}-symplectic structures}
\label{chap:constrasymp}
As discussed in the introduction to this thesis, the link between symplectic geometry and the study of Lefschetz fibrations has proven very fruitful. In this chapter we focus on one side of this link, namely the process of equipping the total space of a Lefschetz fibration with a symplectic structure. The techniques used are principally developed by Thurston \cite{Thurston76} and Gompf \cite{GompfStipsicz99,Gompf01,Gompf04,Gompf04two,Gompf05}. In this chapter we extend such Gompf--Thurston symplectic techniques to arbitrary Lie algebroids. We introduce Lie algebroid fibrations and Lie algebroid Lefschetz fibrations, and give criteria when these can be equipped with compatible Lie algebroid symplectic structures.

The following is the main result of this chapter, as was mentioned in the introduction to this thesis. A more precise version exists (\autoref{thm:thurstonlfibration}), as well as an analogous result for Lie algebroid fibrations (\autoref{thm:thurstonfibration}).
\newtheorem*{thm:introasymp}{Theorem \ref{thm:introasymp}}
\begin{thm:introasymp} Let $(\varphi,f)\colon \mc{A}^4_{X} \to \mc{A}^2_{\Sigma}$ be a Lie algebroid \lf with connected fibers. Assume that $\mc{A}_\Sigma$ admits a Lie algebroid symplectic structure and there exists a closed $\mc{A}_X$-two-form $\eta$ such that the restriction $\eta|_{\ker \varphi}$ is nondegenerate. Then $X$ admits an $\mc{A}_X$-\symp.
\end{thm:introasymp}
In Chapters \ref{chap:constructingblogs} and \ref{chap:constructingsgcs}, we will use this result to establish construction results for \blog{}s and \sgcs{} respectively. Apart from the initial section on Lefschetz fibrations, the contents of this chapter have appeared before in \cite{CavalcantiKlaasse17} and are joint with Gil Cavalcanti.
\subsubsection*{Organization of the chapter}
The structure of this chapter is as follows. We start in Section \ref{sec:lalfs} with a brisk introduction to Lefschetz fibrations, before introducing its Lie algebroid version and the notion of a Lie algebroid fibration. In Section \ref{sec:gompfthurston} we then prove precise versions of the main result above, by extending the techniques of Gompf to Lie algebroid fibrations and Lie algebroid Lefschetz fibrations.
\section{Lie algebroid Lefschetz fibrations}
\label{sec:lalfs}
In this section we introduce the appropriate notions of fibration and Lefschetz fibration in the context of Lie algebroid morphisms. Before doing this, we define Lefschetz fibrations and pencils, and quickly review their relation with symplectic geometry.

\subsection{Lefschetz fibrations}
This section discusses the concept of a \lf. We will be somewhat brief with regard to details, and refer the reader to e.g.\ \cite{GompfStipsicz99} for more information.
\begin{defn}\label{defn:lefschetzfibration} A \emph{\lf} is a proper map $f\colon X^{2n} \to \Sigma^2$ between oriented connected manifolds which is injective on critical points and such that for each critical point $x \in X$ there exist orientation-preserving complex coordinate charts centered at $x$ and $f(x)$ in which $f$ takes the form $f\colon \C^n \to \C$, $f(z_1, \dots, z_n) = z_1^2 + \dots + z_n^2$.
\end{defn}
A Lefschetz fibration can be thought of as a complex analogue of a real Morse function. While there are critical points, they are very mild and well-understood.
\begin{rem} Requiring that \lf{}s are injective on critical points is convenient but not essential, as it can always be ensured by a small perturbation.
\end{rem}
We will focus mostly on Lefschetz fibrations on four-dimensional manifolds. In this dimension one can alternatively take the local form of a Lefschetz singularity to be $f(z_1,z_2) = z_1 z_2$ by applying a linear change of coordinates. In dimension four there is the accompanying notion of a Lefschetz pencil when $\Sigma = \C P^1 \cong S^2$.
\begin{defn} Given a closed connected oriented four-manifold $X^4$, a \emph{Lefschetz pencil} consists of a nonempty finite subset $B \subseteq X$ called the \emph{base locus}, and smooth map $f\colon X \setminus B \to \C P^1$ with only Lefschetz critical points, such that for each point $b \in B$ there is an orientation-preserving coordinate chart where $f$ is given by the projectivization $f\colon \C^2 \setminus \{0\} \to \C P^1$, more concretely $f(z_1,z_2) = z_1 / z_2$.
\end{defn}
A Lefschetz pencil can be viewed as a Lefschetz fibration outside of $B$. Any Lefschetz pencil $f\colon X^4 \setminus B \to \C P^1$ can be blown-up to give a Lefschetz fibration $f\colon X \#_n \overline{\C P}^2 \to \C P^1$, where $n = |B|$ is the amount of base points. Lefschetz pencils can also be defined in higher dimensions, but we will not have use for this here (see e.g.\ \cite[Definition 1.4]{Gompf04}).
\begin{rem} On $\C^n$, complex conjugation provides an orientation-reversing diffeomorphism if and only if $n$ is odd. These preserve the local description of the Lefschetz singularities. Consequently we can always assume that the charts are orientation-preserving on $\Sigma^2$, but not necessarily on $X^{2n}$. See also \autoref{defn:achirallefschetzfibration}.
\end{rem}
From the local description of the Lefschetz singularities it immediately follows that the set $\Delta := {\rm Crit}(f)$ is discrete. Assuming that $\Sigma$ is compact (as will always happen for us) we see that $\Delta$ is finite, hence so is the set $\Delta' := f(\Delta)$ of critical values. This is in contrast with a general map $f\colon X \to \Sigma$, whose regular values are generic by Sard's theorem, but its critical values certainly need not be discrete.
\begin{exa} Any fibration $f\colon X^{2n} \to \Sigma^2$ is a \lf with $\Delta = \emptyset$.
\end{exa}
Given a Lefschetz fibration $f\colon X \to \Sigma$, the restriction $f\colon X \setminus f^{-1}(\Delta') \to \Sigma \setminus \Delta'$ is a fibration.  Because $\Sigma$ is connected, we see that all but finitely many fibers of $f$ are smooth submanifolds of $X$ with the same diffeomorphism type. When $\dim X = 2n = 4$, these fibers are smooth surfaces of a certain genus $g$. This number is called the \emph{genus} of the Lefschetz fibration. We denote the fibers of $f$ by $F_y = f^{-1}(y)$ for $y \in \Sigma$, and the generic fibers by $F$. Moreover, for $y \in \Sigma$ we let $\Delta_y = F_y \cap \Delta$ be set of the critical points on $F_y$, which has cardinality at most one if $f$ is injective on critical points. Using a handlebody description one can calculate the following, starting with the trivial fibration whose Euler characteristic is given by the product.
\begin{prop}\label{prop:lfeulerchar} Let $f\colon X^4 \to \Sigma^2$ be a genus $g$ Lefschetz fibration over a base of genus $h$ with $\mu$ singular fibers. Then $\chi(X) = (2-2g)(2-2h) + \mu$.
\end{prop}
\begin{rem}\label{rem:lfacs} Given a Lefschetz fibration or pencil on a four-manifold $X^4$, one can construct an \acs as follows: near $\Delta$ and $B$, the defining charts give a complex structure. At other points the map $f$ is a fibration, hence splits $TX$ as a sum of two oriented subbundles of rank two (normal and tangent to the fibers), which one declares to be complex line bundles. This is compatible with the local models near $\Delta$ and $B$, as the fibers are holomorphic submanifolds in those charts.
\end{rem}
\begin{rem}\label{rem:acs} It is a standard fact in four-manifold theory (\cite[Theorem 1.4.13]{GompfStipsicz99}) that if $X^4$ is almost-complex, $b_2^+(X) + b_1(X)$ must be odd. This is an if and only if for simply-connected four-manifolds. Consequently, the manifolds $S^4$, $\# m (S^2 \times S^2)$ and $m \C P^2 \# n \overline{\C P}^2$ for $m$ even do not admit \acs{}s.
\end{rem}
We consider the singular fibers in more detail in the four-dimensional case. Each point in $\Delta$ has a chart in which $f(z_1,z_2) = z_1^2 +  z_2^2$, or equivalently one for which $f(z_1,z_2) = z_1 z_2$. This shows that the unique critical value is the origin $0$, with $f^{-1}(0) = \{(z_1,z_2) \, | \, z_1 = 0 \text{ or } z_2 = 0\}$. This is a pair of intersecting planes also referred to as a \emph{nodal singularity}. Thus each singular fiber is a smoothly immersed surface, with each critical point corresponding to a positive transverse self-intersection.

Let us now consider coordinates for which $f(z_1,z_2) = z_1^2 + z_2^2$. Consider a nearby regular value, which for simplicity we shall put on the real axis at some $\varepsilon > 0$. Consider the set $\gamma = [0,\varepsilon] \subseteq \C$ and its inverse image $f^{-1}(\gamma)$. If we now intersect with $\R^2 \subseteq \C^2$, we obtain the circle $\{x_1^2 + x_2^2 = t\}$ for the regular fiber $F_t \cap \R^2$ with $t > 0$ (here $z_i = x_i + i y_i$ for $i = 1,2$). This circle bounds a disk $D_t$ in $f^{-1}([0,t]) \cap \R^2$, called a \emph{Lefschetz thimble}. This disk shrinks to a point in $F_0$ as $t$ approaches $0$. From this the singular fiber $F_0$ is seen to be obtained from the regular fibers $F_t$ by collapsing the circle $\partial D_t = F_t \cap \R^2$, called a \emph{vanishing cycle}, to a point. These vanishing cycles allow one to capture Lefschetz fibrations combinatorially, as the presence of singular fibers amounts to the attachment of a two-handle along the vanishing cycles. First however, let us consider connectedness of the fibers.
\begin{prop}[{\cite[Proposition 8.1.9]{GompfStipsicz99}}]\label{prop:lfexactsequence} Let $f\colon X^4 \to \Sigma^2$ be a Lefschetz fibration with generic fiber $F$. Then there is an exact sequence $\pi_1(F) \to \pi_1(X) \to \pi_1(Y) \to \pi_0(F) \to 0$.
\end{prop}
As a consequence, if $\Sigma$ is simply-connected, each fiber of $f$ is connected. Moreover, if the fibers of $f$ are connected and $b_1(X) = 0$, then $\Sigma$ must be either $S^2$ or $D^2$ (depending on whether $\partial X = \emptyset$). More interestingly, if $f$ is not connected, $\pi_1(X)$ maps to a finite-index subgroup of $\pi_1(\Sigma)$ as $X$ is assumed to be compact. This specifies a cover $h\colon \wt{\Sigma} \to \Sigma$, and $f$ factors through $h$ giving a Lefschetz fibration $\wt{f}\colon X^4 \to \wt{\Sigma}$. However, this map is surjective on $\pi_1$ by definition, so that by \autoref{prop:lfexactsequence} the map $\wt{f}$ has connected fibers.

The description of a neighbourhood of a singular fiber using vanishing cycles lead to describing a genus $g$ \lf $f$ using its monodromy representation. For a gentle overview, we recommend also \cite{Fuller03}.  Namely, the boundary $\partial U$ of neighbourhood of a singular fiber is a fibration (using $f$) over some circle surrounding its singular value. This means it is a mapping torus $\partial U = \Sigma_g \times I / \sim$, with $(x,0) \sim (\psi(x),1)$ for some homeomorphism $\psi\colon \Sigma_g \to \Sigma_g$. This map is the \emph{monodromy} of the singular value. It turns out this map is given by a right-handed \emph{Dehn twist} around the aforementioned vanishing cycle for the singular fiber.

To describe the entire \lf{}, consider the set of singular values, $\Delta' \subseteq \Sigma$. For each $p_i \in \Delta'$, consider a small disk around $p_i$ for which we have the description above. Given a smooth reference fiber $F_0 = f^{-1}(p_0)$ for $p_0 \in \Sigma \setminus \Delta'$, one can draw arcs from $p_0$ to each $p_i$. This allows us to describe $X$ as the trivial $\Sigma_g$-bundle together with an attachment of two-handles for each $p_i$ according to the vanishing cycles. What is independent of choices is the total monodromy of the \lf, which is the composition of each of the Dehn twists (up to isotopy, and so-called \emph{elementary transformations} and conjugation). This is a word in the \emph{mapping class group} of $\Sigma_g$.

We continue by discussing how four-dimensional \lf{}s relate to symplectic geometry. The following is the result we will extend to Lie algebroids (\autoref{thm:thurstonlfibration}).
\begin{thm}[Gompf, {\cite[Theorem 10.2.18]{GompfStipsicz99}}]\label{thm:gompflfsymp} Let $f\colon X^4 \to \Sigma^2$ be a \lf whose generic fiber $F$ satisfies $[F] \neq 0 \in H_2(X;\R)$. Then $X$ admits a \symp such that the fibers of $f$ are symplectic.
\end{thm}
\begin{rem} The condition on the homology class $[F]$ of the generic fiber cannot be dropped, as if $X$ admits such a \symp, it pairs nontrivially with $[F]$. Moreover, $S^1 \times S^3$ admits a \lf without critical points (using the Hopf fibration $S^3 \to S^2$), but $H_2(S^1 \times S^3; \R) = 0$ so cannot be symplectic. Note that this fibration has torus fibers (c.f.\ the discussion below \autoref{rem:entyrefulleralf}).
\end{rem}
There is a partial converse, proven using approximately-holomorphic techniques.
\begin{thm}[Donaldson, {\cite{Donaldson99}}]\label{thm:donaldsonlp} Let $X^4$ be a compact symplectic four-manifold. Then $X$ admits a Lefschetz pencil with symplectic fibers.
\end{thm}
We continue our discussion by relaxing the definition of a \lf, as we can drop the assumption that the charts are compatible with orientations.
\begin{defn}\label{defn:achirallefschetzfibration} An \emph{\alf} is a map $f\colon X^{2n} \to \Sigma^2$ between compact connected manifolds so that for each critical point $x \in X$ there exist complex coordinate charts centered at $x$ and $f(x)$ in which $f$ takes the form $f\colon \C^n \to \C$, $f(z_1, \dots, z_n) = z_1^2 + \dots + z_n^2$.
\end{defn}
Note that we do not require $X$ nor $\Sigma$ to be orientable. If they are however, after choosing orientations one can assign a sign to each critical point of $f$. A given critical point $x \in X$ obtains a sign by demanding that the complex structure of the chart on $\Sigma$ is compatible with its orientation; we then say $x$ is \emph{positive} if the complex structure on $X$ is compatible with its orientation, and \emph{negative} otherwise. Note that any \lf is also an \alf.
\begin{rem}\label{rem:entyrefulleralf} Etnyre--Fuller \cite{EtnyreFuller06} have shown that any four-manifold $X$ admits an \alf{} with a section after performing surgery on a certain framed circle. This implies in particular that if $X$ is simply-connected, both $X \# \C P^2 \# \overline{\C P}^2$ and $X \# S^2 \times S^2$ admit \alf{}s with sections.
\end{rem}
We consider the assumption on the homology class of the generic fiber $F$. If its homology class satisfies $[F] \neq 0 \in H_2(X;\R)$, we say $f$ is \emph{homologically essential}. Given an \alf{} $f\colon X \to \Sigma$, the kernel of $Tf$ defines a complex line bundle $L \to X$ (see also \autoref{rem:lfacs}), whose first Chern class restricts to the Euler class of $F$ and satisfies $\langle c_1(L)|_F, [F] \rangle = \chi(F)$. This is nonzero if $g(F) \neq 1$, in which case the (achiral) Lefschetz fibration will be homologically essential.
\begin{rem} Let $f\colon X \to \Sigma$ be an (achiral) \lf over an oriented surface which admits a section $s$, i.e.\ a map $s\colon \Sigma \to X$ such that $f \circ s = {\rm id}_\Sigma$. Then the fibers of $f$ are homologically essential, as they intersect nontrivially with the image of the section.
\end{rem}
As a consequence of the previous remark, if $X^4$ admits a Lefschetz pencil, it can be equipped with a \symp. Namely, one first blows-up to a \lf $f\colon X \# n \overline{\C P}^2 \to \C P^1$ with $n = |B|$, which admits a section (the exceptional spheres) so that it is homologically essential. One then applies \autoref{thm:gompflfsymp} to obtain a \symp on $X \# n \overline{\C P}^2$, and then blows down each of the symplectic exceptional spheres, resulting in a symplectic structure on $X$. In this sense \autoref{thm:gompflfsymp} and \autoref{thm:donaldsonlp} are inverses to each other: a four-manifold admits a \symp{} if and only if it admits a Lefschetz pencil.
\begin{rem}\label{rem:alfs4} Unlike for \lf{}s (see \autoref{rem:lfacs}), a four-manifold equipped with an \alf need not be almost-complex. For example, $X = S^4$ admits an \alf $f$ due to Matsumoto \cite[Example 8.4.7]{GompfStipsicz99}. As $b_2(S^4) = 0$, the fibers of $f$ are homologically trivial, so must be tori.
\end{rem}

\subsection{Lie algebroid Lefschetz fibrations}
\label{sec:lalf}
We are now ready to define the notion of a Lefschetz fibration between Lie algebroids. To do this, we first need the accompanying notion of a Lie algebroid fibration.
\begin{defn} A \emph{Lie algebroid fibration} $(\varphi,f)\colon \mc{A} \to \mc{A}'$ is a Lie algebroid morphism for which $\varphi\colon \mc{A} \to f^* \mc{A}'$ is surjective. Equivalently, $\varphi$ is fiberwise surjective.
\end{defn}
Note that if $f\colon X \to Y$ is a fibration, then $(Tf,f)\colon TX \to TY$ is a Lie algebroid fibration using \autoref{exa:tgntlamorph} to conclude it is a Lie algebroid morphism.
%
\begin{rem} Our notion of Lie algebroid fibration differs from the one used by other authors, notably Mackenzie \cite{Mackenzie05}. Our Lie algebroid fibrations are not required to cover a surjective submersion. In other words, only $\varphi$ is fiberwise surjective, not both $\varphi$ and $Tf$. This is in line with viewing $\mc{A}$ as the replacement of $TX$.
\end{rem}
We next introduce Lie algebroid \lf{}s, which are a simultaneous generalization of \lf{}s, as well as of Lie algebroid fibrations. We essentially separate two types of singular behavior, namely that of the anchors of the Lie algebroids, and that of the morphism between them. In the isomorphism locus of the Lie algebroid, the condition of being a Lie algebroid fibration is just that it be a fibration. We can weaken this condition here, and only here, to allow for Lefschetz-type singularities. This results in the following definition.
\begin{defn}\label{defn:lalf} A \emph{Lie algebroid \lf} $(\varphi,f)\colon \mc{A}^{2n}_{X} \to \mc{A}^2_{\Sigma}$ is a Lie algebroid morphism which is a Lie algebroid fibration outside a discrete set $\Delta \subset X_{{\mc{A}}_X}$ with $f(\Delta) \subset \Sigma_{\mc{A}_\Sigma}$ such that $f|_{X_{\mc{A}_X}}\colon X_{{\mc{A}}_X} \to \Sigma_{{\mc{A}}_\Sigma}$ is a \lf.
\end{defn}
\begin{rem} Note that if $X_{\mc{A}_X}$ is empty, the notions of Lie algebroid \lf and fibration coincide. There are no Lefschetz singularities outside of $X_{\mc{A}_X}$. Moreover, whenever $\Delta$ is nonempty, both $X_{\mc{A}_X}$ and $\Sigma_{\mc{A}_\Sigma}$ are nonempty, and hence $\dim(X) = 2n$ and $\dim(\Sigma) = 2$. Finally, when $X_{\mc{A}_X} = X$ and $\Sigma_{\mc{A}_\Sigma} = \Sigma$, a Lie algebroid \lf is just a \lf $f\colon X^{2n} \to \Sigma^2$.
\end{rem}
In Section \ref{sec:gompfthurston} we will use Lie algebroid \lf{}s whose generic fibers in $X_{\mc{A}_X}$ are connected. Unlike for usual Lefschetz fibrations (see \autoref{prop:lfexactsequence}), for Lie algebroid Lefschetz fibrations there is in general no exact sequence in homotopy by which we can ensure connected fibers.
\section{Constructing Lie algebroid symplectic structures}
\label{sec:gompfthurston}
In this section we consider the Thurston argument \cite{Thurston76} for constructing \symp{}s extended by Gompf \cite{Gompf04}, and adapt it to the context of Lie algebroid symplectic forms. The guiding principle is to combine suitably symplectic-type structures from the base of a fibration-like map with a form that is symplectic on the tangent spaces of the fibers of that map. In the Lie algebroid case one uses Lie algebroid morphisms $(\varphi,f)\colon \mc{A}_X \to \mc{A}_Y$. Special attention is required because $\rho_{\mc{A}_X}\colon \ker \varphi \to \ker Tf$ need not be an isomorphism (nor injective or surjective), hence one should interpret the tangent space to the fibers suitably.

We will use Lie algebroid \acs{}s as certificates for nondegeneracy of forms, by using the notion of tameness. Let $\mc{A}_X \to X$ be a Lie algebroid and $\omega \in {\rm Symp}(\mc{A}_X)$. Recall from Section \ref{sec:definitionsasymp} that an $\mc{A}_X$-\acs is a vector bundle complex structure $J$ for $\mc{A}_X$.
\begin{defn}\label{defn:taming} An $\mc{A}_X$-\acs $J$ is \emph{$\omega$-tame} if $\omega(v, J v) > 0$ for all $v \in \Gamma(\mc{A}_X)$. Given a Lie algebroid morphism $(\varphi,f)\colon \mc{A}_X \to \mc{A}_Y$ and $\omega_Y \in {\rm Symp}(\mc{A}_Y)$, $J$ is \emph{$(\omega_Y,\varphi)$-tame} if $(\varphi_x ^* \omega_Y)(v, J v) > 0$ for all $v \in \mc{A}_{X,x} \setminus \ker \varphi_x$ and $x \in X$.
\end{defn}
As usual, the space of taming $\mc{A}_X$-\acs{}s for $\omega$ is convex and nonempty and is denoted by $\mc{J}(\omega)$. Note that any $\omega \in \Omega_{\mc{A}_X}^2(X)$ taming a $\mc{A}_X$-\acs $J$ is necessarily nondegenerate. Hence if $\omega$ is a closed $\mc{A}_X$-two-form taming some $J$, then $\omega$ is $\mc{A}_X$-symplectic and $J$ induces the same $\mc{A}_X$-orientation as $\omega$. Note moreover that if $J$ is $(\omega_Y, \varphi)$-tame, then $\ker \varphi$ is a $J$-complex subspace of $\mc{A}_X$. Indeed, if $v \in \ker \varphi$ and $J v \not\in \ker \varphi$, we would have $0 = \varphi^*\omega_Y(v, Jv) = \varphi^*\omega_Y(Jv, J(Jv)) > 0$, which is a contradiction.
\begin{prop}\label{prop:tamingopen} The taming condition is open, i.e.\ it is preserved under sufficiently small perturbations of $\omega$ and $J$, and of varying the point in $X$.
\end{prop}
\bp The taming condition $\omega(v, J v) > 0$ for the pair $(\omega, J)$ holds provided it holds for all $v \in \Sigma\mc{A}_X \subset \mc{A}_X$, the unit sphere bundle with respect to some preassigned metric, as $X$ is compact. As $\Sigma\mc{A}_X$ is also compact, the continuous map $\wt{\omega}\colon \Sigma\mc{A}_X \to \R$ given by $\wt{\omega}(v) := \omega(v, J v)$ for $v \in \Sigma\mc{A}_X$ is bounded from below by a positive constant on $\Sigma\mc{A}_X$. But then $\wt{\omega}$ will remain positive under small perturbations of $\omega$ and $J$. Similarly the condition of $\omega$ taming $J$ on $\ker \varphi$ is open. Consider $x \in X$ so that $\omega(v, J v) > 0$ for all $v \in \ker \varphi_x$. As $\wt{\omega}$ is continuous and $\Sigma\mc{A}_X$ is compact, there exists a neighbourhood $U$ of $\ker \varphi_x \cap \Sigma\mc{A}_X$ in $\Sigma\mc{A}_X$ on which $\wt{\omega}$ is positive. Points $x' \in X$ close to $x$ will then have $\ker \varphi_{x'} \subset U$ because $\ker \varphi$ is closed.
\ep
We will not use the associated notion of compatibility, where $J$ further leaves $\omega$ invariant, as we use \acs{}s as auxiliary structures to show non-degeneracy, and make use of the openness of this condition. For this reason, all \acs{}s will only be required to be continuous as this avoids arguments to ensure smoothness.

The following is the Lie algebroid version of \cite[Theorem 3.1]{Gompf04} and is our main tool to construct $\mc{A}_X$-\symp{}s.
\begin{thm}\label{thm:lathurstontrick} Let $(\varphi,f)\colon (X,\mc{A}_X) \to (Y, \mc{A}_Y)$ be a Lie algebroid morphism between compact connected manifolds, $J$ an $\mc{A}_X$-\acs, $\omega_Y$ a $\mc{A}_Y$-symplectic form and $\eta$ a closed $\mc{A}_X$-two-form. Assume that
	\bi
	\item[i)] $J$ is $(\omega_Y,\varphi)$-tame;
	\item[ii)] $\eta$ tames $J|_{\ker \varphi_x}$ for all $x \in X$.
	\ei
	Then $X$ admits an $\mc{A}_X$-\symp.
\end{thm}
The proof of this result is modelled on that by Gompf of \cite[Theorem 3.1]{Gompf04}.
\bp Let $t>0$ and consider the form $\omega_t := \varphi^* \omega_Y + t \eta$. We show $J$ is $\omega_t$-tame for $t$ small enough. By \autoref{prop:tamingopen} it is enough to show that there exists a $t_0 > 0$ so that $\omega_t(v, J v) > 0$ for every $t \in (0,t_0)$ and $v$ in the unit sphere bundle $\Sigma \mc{A}_X \subset \mc{A}_X$ with respect to some metric. For $v \in \mc{A}_X$ we have $\omega_t(v, J v) = \varphi^* \omega_Y(v,J v) + t \eta(v, J v)$. As $J$ is $(\omega_Y,\varphi)$-tame, the first term is positive for $v \in \mc{A}_X \setminus \ker \varphi$ and is zero otherwise. The second term $\eta(v, J v)$ is positive on $\ker \varphi$ because $J|_{\ker \varphi}$ is $\eta$-tame, hence is also positive for all $v$ in some neighbourhood $U$ of $\ker \varphi \cap \Sigma\mc{A}_X$ in $\mc{A}_X$ by openness of the taming condition. We conclude that $\omega_t(v, Jv) > 0$ for all $t > 0$ when $v \in U$. The function $v \mapsto \eta(v, J v)$ is bounded on the compact set $\Sigma\mc{A}_X \setminus U$. Furthermore, $\varphi^*\omega_Y(v, J v)$ is also bounded from below there by a positive constant, as it is positive away from $\ker \varphi$, and thus also away from $\ker \varphi \cap \Sigma\mc{A}_X \subset U$. But then $\omega_t(v, J v) > 0$ for all $0 < t < t_0$ for $t_0$ sufficiently small, so that $\omega_t$ is $\mc{A}_X$-symplectic for $t$ small enough.
\ep
Given a map $f\colon X \to Y$ we denote by $F_y = f^{-1}(y)$ for $y \in Y$ the level set, or \emph{fiber}, of $f$ over $y$. In order to meaningfully apply \autoref{thm:lathurstontrick} we must be able to construct closed $\mc{A}_X$-forms $\eta$ as in ii) of the statement. Note firstly that it suffices to have such closed forms in neighbourhoods of fibers, all lying in the same global Lie algebroid cohomology class.
\begin{prop}\label{prop:gluingfiberform} Assume that there exists a class $c \in H^2_{\mc{A}_X}(X;\R)$ and that for each $y \in Y$, $F_y$ has a neighbourhood $W_y$ with a closed $\mc{A}_X$-two-form $\eta_y \in \Omega^2_{\mc{A}_X}(W_y)$ such that $[\eta_y] = c|_{W_y} \in H^2_{\mc{A}_X}(W_y;\R)$, and $\eta_y$ tames $J|_{\ker \varphi_x}$ for all $x \in W_y$. Then there exists a closed $\mc{A}_X$-two-form $\eta$ such that $[\eta] = c$ and $\eta$ tames $J|_{\ker \varphi_x}$ for all $x \in X$.
\end{prop}
\bp Let $\xi \in \Omega^2_{\mc{A}_X}(X)$ be closed and such that $[\xi] = c$. Then for each $y \in Y$ we have $[\eta_y] = c|_{W_y} = [\xi]|_{W_y}$, so on $W_y$ we have $\eta_y = \xi + d_{\mc{A}_X} \alpha_y$ for some $\alpha_y \in \Omega^1_{{\mc{A}}_X}(W_y)$. As each $X \setminus W_y$ and hence $f(X \setminus W_y)$ is compact, each $y \in Y$ has a neighbourhood disjoint from $f(X\setminus W_y)$. Cover $Y$ by a finite amount of such open sets $U_i$ so that each $f^{-1}(U_i)$ is contained in some $W_{y_i}$. Let $\{\psi_i\}$ be a partition of unity subordinate to the cover $\{U_i\}$ of $Y$, so that $\{\psi_i \circ f\}$ is a partition of unity of $X$. Define an $\mc{A}_X$-two-form $\eta \in \Omega^2_{\mc{A}_X}(X)$ on $X$ via
\be
\eta := \xi + d_{\mc{A}_X}(\sum_i (\psi_i \circ f) \alpha_{y_i}) = \xi + \sum_i (\psi_i \circ f) d_{\mc{A}_X}\alpha_{y_i} + \sum_i \rho_{\mc{A}_X}^*(d\psi_i \circ Tf) \wedge \alpha_{y_i}.
\ee
Then $d_{\mc{A}_X}\eta = 0$ and $[\eta] = c$. As $\rho_{\mc{A}_X}^*(d\psi_i \circ Tf) = d\psi_i \circ Tf\circ \rho_{\mc{A}_X} = dg \circ \rho_{\mc{A}_Y} \circ \varphi$ using that $\varphi$ is a Lie algebroid morphism, the last of the above three terms vanishes when applied to a pair of vectors in $\ker \varphi_x$ for any $x \in X$, so on each $\ker \varphi_x$ we have
\be
\eta = \xi + \sum_i (\psi_i \circ f) d_{\mc{A}_X}\alpha_{y_i} = \sum_i (\psi_i \circ f)(\xi + d_{\mc{A}_X}\alpha_{y_i}) = \sum_i (\psi_i \circ f) \eta_{y_i}.
\ee
From the above we see that on $\ker \varphi$, the $\mc{A}_X$-form $\eta$ is a convex combination of $J$-taming $\mc{A}_X$-forms, so $J|_{\ker \varphi}$ is $\eta$-tame.
\ep
One can further look for (local or global) closed two-forms $\wt{\eta} \in \Omega^2(X)$ so that $\eta = \rho_{\mc{A}_X}^* \wt{\eta}$ satisfies hypothesis ii) of \autoref{thm:lathurstontrick}. When using such $\mc{A}_X$-forms which are pullbacks of regular forms, the behavior of the map $\rho_{\mc{A}_X}\colon \ker \varphi_x \to \ker Tf_x$ is important. We will see in Chapter \ref{chap:constructingsgcs} an example where there cannot be an $\eta$ of the form $\eta = \rho_{\mc{A}_X}^* \wt{\eta}$ making $\ker \varphi$ symplectic. This is unlike the situation for log-symplectic forms, where the anchor of the log-tangent bundle provides an isomorphism between $\ker \varphi$ and $\ker Tf$ (see \autoref{prop:kerbdf}).

Using \autoref{thm:lathurstontrick} we can prove the Lie algebroid version of Thurston's result for symplectic fiber bundles with two-dimensional fibers \cite{Thurston76}, adapting the proof by Gompf in \cite{Gompf01}.
\begin{thm}\label{thm:thurstonfibration} Let $(\varphi,f)\colon (X,\mc{A}_X^{2n}) \to (Y, \mc{A}_Y^{2n-2})$ be a Lie algebroid fibration between compact connected manifolds. Assume that $Y$ is $\mc{A}_Y$-symplectic and there exists a closed $\mc{A}_X$-two-form $\eta$ which is nondegenerate on $\ker \varphi$. Then $X$ admits an $\mc{A}_X$-\symp.
\end{thm}
\bp Let $\omega_Y \in {\rm Symp}(\mc{A}_Y)$ and choose $J_Y \in J(\omega_Y)$. Fix the orientation for $\ker \varphi$ so that $\eta$ is positive. Let $g$ be a metric on $\mc{A}_X$ and let $H \subset \mc{A}_X$ be the subbundle of orthogonal complements to $\ker \varphi$, so that $\varphi\colon H \to \mc{A}_Y$ is a fiberwise isomorphism. Define an $\mc{A}_X$-\acs $J$ by letting $J|_H = \varphi^* J_Y$, and on $\ker \varphi$, use the metric and define $J$ by $\frac{\pi}2$-counterclockwise rotation, demanding $\varphi$ is orientation preserving via the fiber-first convention. This determines $J$ uniquely on $\mc{A}_X$ by linearity. Moreover, $J$ is $(\omega_Y,\varphi)$-tame as $\varphi^* \omega_Y(v, J v) = \omega_Y(\varphi v, J_Y \varphi v) > 0$, for all $v \in \mc{A}_X \setminus \ker \varphi \cong H$. Further, $\eta$ tames $J$ on $\ker \varphi$ as $J$ is compatible with the orientation on $\ker \varphi$ determined by $\eta$. By \autoref{thm:lathurstontrick} we obtain an $\mc{A}_X$-\symp.
\ep
\begin{rem} One can combine \autoref{thm:thurstonfibration} with \autoref{prop:gluingfiberform} to obtain a statement requiring only the existence of local forms $\eta_y$ governed by a global cohomology class.
\end{rem}

Sadly for general Lie algebroid morphisms there is no direct analogue of homological essentialness of generic fibers to replace the hypothesis of the existence of $\eta$ as in \autoref{thm:lathurstontrick} and \autoref{thm:thurstonfibration}. In other words, we cannot replace the condition on the existence of $\eta$ by demanding that the generic fiber $F$ satisfies $[F] \neq 0 \in H_2(X;R)$. This is a well-known necessary condition when constructing symplectic-like structures out of Lefschetz-type fibrations using Gompf--Thurston techniques. The reason for the lack of such an analogue is again the behavior of the map $\rho_{\mc{A}_X}\colon \ker \varphi \to \ker Tf$. While the codomain can be seen as the tangent space to the fiber at regular points, the domain cannot. For Lie algebroid submersions, surjectivity is only demanded of $\varphi$, not of $Tf$, and there is no duality pairing between homology and Lie algebroid cohomology in general.
\begin{rem} In the special cases of the log-tangent bundle and the elliptic tangent bundle one can state such an analogue (see Chapters \ref{chap:constructingblogs} and \ref{chap:constructingsgcs}), precisely because we understand their respective Lie algebroid cohomologies (see \autoref{thm:loglacohomology} and \autoref{thm:elllacohomology}).
\end{rem}

Recall that we defined a Lie algebroid \lf to be a Lie algebroid fibration away from the set of Lefschetz singularities, which must be a subset of the isomorphism locus of $\mc{A}_X$. Note that $X_{\mc{A}_X}$ is either empty or of codimension zero, so that we can restrict $\mc{A}_X$ to $X_{\mc{A}_X}$ and again obtain a Lie algebroid. In this way the inclusion $i\colon X_{\mc{A}_X} \hookrightarrow X$ covers a Lie algebroid morphism $\iota$, and we obtain a restriction map in cohomology $\iota^*\colon H_{\mc{A}_X}(X) \to H_{\mc{A}_X}(X_{\mc{A}_X})$. The anchor gives a bijection $\rho_{\mc{A}_X}^*\colon \Omega^\bullet(X_{\mc{A}_X}) \to \Omega^\bullet_{\mc{A}_X}(X_{\mc{A}_X})$ and hence an isomorphism $\rho_{\mc{A}_X}^*\colon H_{\rm dR}(X_{\mc{A}_X}) \to H_{\mc{A}_X}(X_{\mc{A}_X})$ in cohomology over $X_{\mc{A}_X}$. 

Assuming ${\rm rank}(\mc{A}_X) = {\rm rank}(\mc{A}_Y) + 2$, the second hypothesis in \autoref{thm:lathurstontrick} demands the existence of a class $c \in H_{\mc{A}_X}(X)$ for which the pullback de Rham class $(\rho_{\mc{A}_X}^*)^{-1} \iota^*(c)$ evaluates nonzero on each generic fiber. Even if such a class exists, a similar statement must hold over $X \setminus X_{\mc{A}_X}$. Hence, we need the existence of a two-form $\eta \in \Omega^2_{\rm cl}(X;\mc{A}_X)$ such that $(\rho_{\mc{A}_X}^*)^{-1} \iota^* [\eta]$ evaluates nonzero on the fibers, and $\eta|_{\ker \varphi}$ is nowhere zero over $X \setminus X_{\mc{A}_X}$.

The notion of a Lie algebroid \lf is such that the hypotheses of \autoref{thm:lathurstontrick} are still satisfied, when in dimension four.
\begin{thm}\label{thm:thurstonlfibration} Let $(\varphi,f)\colon \mc{A}^4_{X} \to \mc{A}^2_{\Sigma}$ be a Lie algebroid \lf with connected fibers between compact connected manifolds. Assume that $\mc{A}_\Sigma$ admits a symplectic structure and there exists a closed $\mc{A}_X$-two-form $\eta$ such that $(\rho_{\mc{A}_X}^*)^{-1} \iota^* [\eta]$ evaluates nonzero on the fibers, and $\eta|_{\ker \varphi}$ is nondegenerate over $X \setminus X_{\mc{A}_X}$. Then $X$ admits an $\mc{A}_X$-\symp.
\end{thm}
The proof of this result is based on Gompf's proof of \cite[Theorem 10.2.18]{GompfStipsicz99} using almost-complex structures. We will use \autoref{prop:gluingfiberform} to glue together the given form $\eta$ with a suitable adaptation in the isomorphism locus of $\mc{A}_X$.
\begin{rem} Note that if $\mc{A}_X = TX$, then $X_{\mc{A}_X} = X$ and we are demanding that the generic fiber is homologically essential, recovering the result by Gompf \cite[Theorem 10.2.18]{GompfStipsicz99}.
\end{rem}
\bp If $X_{\mc{A}_X}$ is empty, $\eta|_{\ker \varphi}$ is nondegenerate everywhere and $\varphi$ is a Lie algebroid fibration, so that the result follows from \autoref{thm:thurstonfibration}. If not, then $X$ is four-dimensional. Denote $\xi := (\rho_{\mc{A}_X}^*)^{-1} \iota^* \eta \in \Omega^2(X_{\mc{A}_X})$ and $c = [\xi] \in H^2_{dR}(X_{\mc{A}_X})$. Note that $\xi$ orients the generic fibers $F$, as these are two dimensional. Recall that $\Delta \subset X_{\mc{A}_X}$ is the set of Lefschetz singularities, and let $\Delta' := f(\Delta)$ be the set of singular values of $f|_{X_{\mc{A}_X}}$.

Let $\omega_\Sigma$ be an $\mc{A}_\Sigma$-\symp and use the proof of \autoref{thm:thurstonfibration} to obtain a $(\varphi,\omega_\Sigma)$-tame almost-complex structure $J$ on $\mc{A}_X$ over $X \setminus f^{-1}(\Delta')$ compatible with the orientation on $\ker \varphi$, noting that $f\colon X \setminus f^{-1}(\Delta') \to \Sigma \setminus \Delta'$ is a Lie algebroid fibration. As $\Delta'$ is contained in $\Sigma_{\mc{A}_\Sigma}$, let $V \subseteq \Sigma$ be the disjoint union of open balls $V_y$ disjoint from $Z_\Sigma$ and centered at each point $y \in \Delta'$. Let $W := f^{-1}(V) \subset X$ be the union of the neighbourhoods $W_y := f^{-1}(V_y)$ of singular fibers $F_y$. Let $C \subset X$ be the disjoint union of open balls $C_y \subseteq \Sigma_{\mc{A}_\Sigma}$ centered at each point $f^{-1}(y) = x$ for all $y \in \Delta'$ so that on each ball $f$ takes on the local form in \autoref{defn:lefschetzfibration}. Possibly shrink $C$ so that $\overline{C_y} \subset W_y$. The local description of $f$ gives an \acs on $C$ with the fibers being holomorphic, and we glue this to the existing almost-complex structure $J$ on $X \setminus C$. This gives a global $(\omega_\Sigma, \varphi)$-tame $\mc{A}_X$-\acs $J$.

As $\eta|_{\ker \varphi}$ is nondegenerate over $X \setminus X_{\mc{A}_X}$, the same is true in a neighbourhood $S$ around $X \setminus X_{\mc{A}_X}$ disjoint from $\Delta$. Let $y \in Y \setminus f(S)$ be given. If $y \not\in \Delta'$, let $D_y \subset Y \setminus f(S)$ disjoint from $\Delta'$ be a disk containing $y$, fully contained in a trivializing neighbourhood of $f$ around $y$. Define $W_y := f^{-1}(D_y) \cong D_y \times F_y$, with projection map $p\colon W_y \to D_y$. Let $\eta_y'$ be an area form on $F_y$ compatible with the preimage orientation. Define $\eta_y = \lambda_y \rho_X^* p^* \eta_y'$, where $\lambda_y \in \R$ is chosen such that $\langle [F], c \rangle = \langle [F], \eta_y \rangle$. As $H_2(W_y;\R)$ is generated by $[F_y]$, it follows that $[\eta_y] = c|_{W_y} \in H^2(W_y;\R)$. But then $\eta_y$ tames $J$ on $\ker \varphi \cong \ker Tf$ for all $x \in W_y$, as the restriction of $(\rho_X^*)^{-1} \eta_y = \lambda_y p^* \eta_y'$ is an area form for that fiber.

If $y \in \Delta'$, the singular fiber $F_y$ either is indecomposable or consists of two irreducible components $F_y^\pm$ which satisfy $[F_y^+] \cdot [F_y^-] = 1$ and $[F_y^\pm]^2 = -1$, see \cite{GompfStipsicz99}. In the latter case, note that $0 < 1 = \langle c, [F] \rangle = \langle c, [F_y]\rangle = \langle c, [F_y^+]\rangle + \langle c, [F_y^-]\rangle$. If either term is nonpositive assume without loss of generality that $\langle c, [F_y^-]\rangle = r \leq 0$. Define $c' := c + (\frac12 - r) c_y^+$, where $c_y^+ \in H^2(X;\R)$ is a class dual to $[F_y^+]$. As $[F_y] \cdot [F_y^\pm] = 0$ we then have $\langle c', [F] \rangle = \langle c, [F] \rangle > 0$, and furthermore $\langle c', [F_y^+]\rangle = \langle c, [F_y^+] \rangle +(\frac12 - r) > 0$ and $\langle c', [F_y^-]\rangle = \frac12 > 0$. Moreover, as different fibers do not intersect, we have $c|_{W_{y'}} = c'|_{W_{y'}}$ for $y' \neq y$. After finitely many repetitions, at most once for each $y \in \Delta'$, one obtains a class, again denoted by $c$, pairing positively with every fiber component (see \cite[Exercise 10.2.19]{GompfStipsicz99}).

Return to $y \in \Delta'$ and let $\sigma$ be the standard symplectic form on $C_y$ given locally in real coordinates by $\sigma = dx_1 \wedge dy_1 + dx_2 \wedge dy_2$, where $z_i = x_i + i y_i$. As all fibers $F_y'$ in $C_y$ are holomorphic, $\rho_X^* \sigma|_{F_{y'} \cap C_y}$ tames $J$ for all $y' \in f(C_y)$, so that $\rho_X^* \sigma$ tames $\btwo{J}$ on $C_y$. Let $\sigma_y$ be an extension of $\sigma$ to $F_y$ as a positive area form with total area $\langle \sigma_y, [F_y] \rangle$ equal to $\langle c, [F_y] \rangle$. Let $p\colon W_y \to F_y$ be a retraction and let $f\colon C_y \to [0,1]$ be a smooth radial function so that $f \equiv 0$ in a neighbourhood of $x = f^{-1}(y) \in \Delta$ and $f \equiv 1$ in a neighbourhood of $\partial C_y$, which is smoothly extended to $W_y$ by being identically $1$ outside $C_y$. On the ball $C_y$, the form $\sigma$ is exact, say equal to $\sigma = d\alpha$ for $\alpha \in \Omega^1(C_y)$. Define a two-form $\eta_y'$ on $W_y$ by $\eta_y' := p^*(f\sigma_y) + d((1-f)\alpha)$, which is closed as $f\sigma_y$ is a closed area form on $F_y$. Near $x$ we have $f \equiv 0$ so that $\eta_y' = d\alpha = \sigma$. Set $\eta_y := \rho_X^* \eta_y'$. Then there $\eta_y = \rho_X^* \sigma$ tames $J$, hence in particular tames $J|_{\ker \varphi}$. Similarly, $\sigma_y$ is an area form on $F_y \setminus \{x\}$ for the orientation given by $J$ under the isomorphism by $\rho_X$. But then $\rho_X^* \sigma_y$ tames $J$ on $\ker \varphi_y \cong \ker Tf = TF_y$ on $F_y$, so that the same holds for $\eta_y$ as this condition is convex. By openness of the taming condition, shrinking $V_y$ and hence $W_y$ and possibly $C_y$ we can ensure that $\eta_y$ tames $J|_{\ker \varphi}$ on $W_y$. Finally, note that $[\eta_y'] = c|_{W_y} \in H^2(W_y;\R)$ by construction. For any point $y \in f(S)$, take the Lie algebroid two-form $\eta_y := \eta$ on the neigbourhood $W_y := S$ of $F_y$. We have now obtained the required neighbourhoods $W_y$ and forms $\eta_y$ for all $y \in Y$ to apply \autoref{prop:gluingfiberform} and obtain a Lie algebroid closed two-form again denoted by $\eta$ such that $\eta$ tames $J$ on $\ker \varphi_x$ for all $x \in X$. By \autoref{thm:lathurstontrick} we obtain an $\mc{A}_X$-symplectic structure on $X$.
\ep			
\chapter{Constructing log-symplectic structures}
\label{chap:constructingblogs}
\renewcommand{\b}[1][]{{^b}{#1}}
In this chapter we use the techniques developed in the previous chapter to obtain existence results for \blog{}s on total spaces of fibration-like maps. As a warm-up, we prove the $b$-analogue of Thurston's result for fibrations with two-dimensional fibers as \autoref{thm:bthurstonfibration}, giving the following (see \autoref{cor:fibrationblogs}).
\begin{thm2} Let $f\colon X^{2n} \to Y^{2n-2}$ be a fibration between compact connected manifolds. Assume that $Y$ admits a \blog with singular locus $Z_Y$, and that the generic fiber $F$ of $f$ is orientable and satisfies $[F] \neq 0 \in H_2(X;\R)$. Then $(X,Z_X)$ admits a \blog, where $Z_X = f^{-1}(Z_Y)$.
\end{thm2}
After introducing the notion of a $b$-hyperfibration (the $b$-analogue of a hyperpencil with empty base locus \cite{Gompf04}) in Section \ref{sec:bhyperfibration}, we show the following as \autoref{thm:bhyperfibration}.
\begin{thm2}\label{thm:introbhyperfibration} Let $f\colon (X,Z_X) \to (Y, Z_Y, \btwo{\omega}_Y)$ be a $b$-hyperfibration between compact connected $b$-oriented $b$-manifolds. Assume that there exists a finite collection $S$ of sections of $f$ intersecting all fiber components non-negatively and for each fiber component at least one section in $S$ intersecting positively. Then $(X,Z_X)$ admits a \blog.
\end{thm2}
Moreover, we recover the following result of Cavalcanti \cite{Cavalcanti17} linking \alf{}s to \blog{}s (see \autoref{thm:alflogsymp}) using our framework, as discussed in the introduction to this thesis.
\newtheorem*{thm:introblogalf}{Theorem \ref{thm:introblogalf}}
\begin{thm:introblogalf}[\cite{Cavalcanti17}] Let $f\colon X^4 \to \Sigma^2$ be an \alf between compact connected manifolds. Assume that the generic fiber $F$ is orientable and $[F] \neq 0 \in H_2(X;\R)$. Then $X$ admits a \blog.
\end{thm:introblogalf}

The contents of this chapter have appeared before in \cite{CavalcantiKlaasse16} and are joint with Gil Cavalcanti. However, we use somewhat different notation, and directly use the general results for arbitrary Lie algebroids as established in Chapter \ref{chap:constrasymp}.
\subsection*{Organization of the chapter}
In Section \ref{sec:blogsbtangent} we discuss the definition of \blog{}s and describe how they can be viewed as symplectic forms for the $b$-tangent bundle. We further discuss the existence of \blog{}s on surfaces. In Section \ref{sec:blogconstructions} we prove the $b$-analogue of Thurston's result, \autoref{thm:bthurstonfibration}, and \autoref{thm:alflogsymp} that four-dimensional \alf{}s give rise to \blog{}s. In Section \ref{sec:bhyperfibration} we define the notion of a $b$-hyperfibration and prove \autoref{thm:bhyperfibration}, that $b$-hyperfibrations lead to \blog{}s.
\section{\texorpdfstring{$b$}{b}-Fibrations and sections}
\label{sec:blogsbtangent}
In this section we discuss the $b$-geometry language used to study log-symplectic structures. We will not distinguish between log pairs and $b$-manifolds, which both are manifolds $X$ together with a hypersurface $Z \subseteq X$. Moreover, recall that a $b$-map $f\colon (X,Z_X) \to (Y,Z_Y)$ is either a transverse strong map of pairs, or a morphism of log divisors (see \autoref{rem:tsmopbmaps}). Finally, we recall \autoref{prop:logdivmorph}, which says that such maps induce Lie algebroid morphisms $(\varphi,f)\colon \mc{A}_{Z_X} \to \mc{A}_{Z_Y}$ between the respective log-tangent bundles.
\begin{defn} A \emph{$b$-fibration} is $b$-map $f\colon (X,Z_X) \to (Y,Z_Y)$ whose induced Lie algebroid morphism $\varphi$ is a Lie algebroid fibration between log-tangent bundles.
\end{defn}
\begin{rem}\label{rem:fibrintobfibr} It follows immediately from the definition that given a fibration $f\colon X \to Y$, one can turn it into a $b$-fibration by choosing a hypersurface $Z_Y \subset Y$ and considering the $b$-map $f\colon (X,Z_X) \to (Y,Z_Y)$, where $Z_X = f^{-1}(Z_Y)$.
\end{rem}
Given a $b$-manifold $(X,Z_X)$, we will call an orientation given to the log-tangent bundle a \emph{$b$-orientation}, and an almost-complex-structure for the log-tangent bundle a \emph{$b$-\acs}. Recall that a log-symplectic structure induces a $b$-orientation. We will only consider $b$-orientable $(X,Z_X)$, and \bacs{}s and log-symplectic forms inducing the same $b$-orientation.
\begin{rem} As remarked in Section \ref{sec:logtangentbundle}, a $b$-oriented $b$-manifold $(X,Z_X)$ gives an orientation to $TX$ away from $Z_X$. However, this orientation cannot come from an existing orientation on $X$ when $Z_X \neq \emptyset$. See also \autoref{prop:orlogseparating}.
\end{rem}
Given a map $f\colon X \to (Y,Z_Y)$, note that $Z_X$ is uniquely determined by $Z_Y$ and the requirement that $f$ is a $b$-map. Given a $b$-map $f\colon (X,Z_X) \to (Y,Z_Y)$, its level sets are either contained in $Z_X$ or are disjoint from it.
We now discuss the extent to which there is a difference between $Tf$-critical points and $\varphi$-critical points for the Lie algebroid morphism of \autoref{prop:logdivmorph}, given a $b$-map $f\colon (X,Z_X) \to (Y,Z_Y)$.
\begin{prop}\label{prop:kerbdf} Let $f\colon (X,Z_X) \to (Y,Z_Y)$ be a $b$-map. Then $\rho_{X,x}\colon \ker \varphi_x \to \ker Tf_x$ is an isomorphism for all $x \in X$.
\end{prop}
Consequently we can unambiguously speak of a critical point of $f$, without specifying whether we mean with respect to $\varphi$ or $Tf$. The essential ingredients are contained in the following lemma. Two vector spaces $V, V_1$ will be called a \emph{pair} if $V_1 \subset V$ is a subspace. A linear map $f\colon V \to W$ is a map between pairs $(V,V_1)$ and $(W,W_1)$ if $f(V_1) \subset W_1$.
\begin{lem}\label{lem:blinalg} Let $F\colon (V,V_1) \to (W,W_1)$ be a linear map between pairs such that under the projection maps ${\rm pr}_V\colon V \mapsto V / V_1$ and ${\rm pr}_W\colon W \mapsto W / W_1$, $F$ descends to an isomorphism $\overline{F}\colon V / V_1 \to W / W_1$. Assume that there are vector spaces $V'$, $W'$ and maps $\rho_V\colon V' \to V$, $\rho_W\colon W' \to W$, $F'\colon V' \to W'$ so that $F \circ \rho_V = \rho_W \circ F'$. Assume that $\im \rho_V = V_1$, $\im \rho_W = W_1$ and $F'\colon \ker \rho_V \to \ker \rho_W$ is an isomorphism. Then $\rho_V: \, \ker F' \to \ker F$ is an isomorphism.
\end{lem}
The situation is summarized by the following diagram, in which the rows are exact. The two vertical maps on the far left and right are assumed to be isomorphisms, while the conclusion of the lemma is that the top horizontal map is an isomorphism.
\begin{center}
	\begin{tikzpicture}
	\matrix (m) [matrix of math nodes, row sep=2.5em, column sep=2.5em,text height=1.5ex, text depth=0.25ex]
	{
		&                      & \ker F' & \ker F  &              &   \\
		0 & \ker \rho_V   & V'        & V         & V/V_1   & 0 \\
		0 & \ker \rho_W  & W'       & W        & W/W_1 & 0 \\};
	\path[-stealth]
	(m-1-3) edge node [above] {$\rho_V$} node [below] {$\cong$} (m-1-4)
	edge [right hook-latex] (m-2-3)
	(m-1-4) edge [right hook-latex] (m-2-4)
	(m-2-1) edge (m-2-2)
	(m-2-2) edge [right hook-latex] (m-2-3)
	(m-2-2) edge node [left] {$\cong$} node [right] {$F'$} (m-3-2)
	(m-2-3) edge node [right] {$F'$} (m-3-3)
	edge node [above] {$\rho_V$} (m-2-4)
	(m-2-4) edge node [right] {$F$} (m-3-4)
	edge node [above] {${\rm pr}_V$} (m-2-5)
	(m-2-5) edge (m-2-6)
	(m-2-5) edge node [right] {$\overline{F}$} node [left] {$\cong$} (m-3-5)
	(m-3-1) edge (m-3-2)
	(m-3-2) edge [right hook-latex] (m-3-3)
	(m-3-3) edge node [above] {$\rho_W$} (m-3-4)
	(m-3-4) edge node [above] {${\rm pr}_W$} (m-3-5)
	(m-3-5) edge (m-3-6);
	\end{tikzpicture}
\end{center}
\bp If $v' \in \ker F'$, then $F \rho_V (v') = \rho_W F' (v') = 0$, so that $\rho_V: \, \ker F' \to \ker F$. Given $v \in \ker F$ we have ${\rm pr}_W F(v) = {\rm pr}_W(0) = 0$, so by assumption ${\rm pr}_V(v) = 0$, hence $v \in V_1$. As $v \in V_1$, there exists a $v'_0 \in V'$ such that $\rho_V(v'_0) = v$, and $\rho_V^{-1}(v) = v'_0 + \ker \rho_V$. Consider a vector $v' = v'_0 + k$ for $k \in \ker \rho_V$, so that $\rho_V (v') = v$. Then $\rho_W F'(v') = F \rho_V(v') = F(v) = 0$, so that $F'(v') \in \ker \rho_W$. But then there exists a unique $k \in \ker \rho_V$ such that $F'(k) = - F'(v'_0)$, for which $F'(v') = 0$. Hence there exists a unique $v' \in \ker F' \cap \rho_V^{-1}(v)$. We conclude that $\rho_V\colon \ker F' \to \ker F$ is an isomorphism.
\ep
\bp[ of \autoref{prop:kerbdf}] Let $x \in X$ be given and denote $y = f(x)$. If $x \in X \setminus Z_X$ the statement follows as $\rho_X$ is an isomorphism away from $Z_X$. If $x \in Z_X$, we wish to apply \autoref{lem:blinalg} to the situation where $V = T_x X$, $V_1 = T_x Z_X$, $W = T_{y} Y$, $W_1 = T_{y} Z_Y$, $V' = T_x X(-\log Z_X)$, $W' = T_{y} Y(-\log Z_Y)$, $F = Tf_x$, $F' = \varphi_x$, $\rho_V = \rho_{X,x}$ and $\rho_W = \rho_{Y,y}$. As $f$ is a $b$-map, we have $f^{-1}(Z_Y) = Z_X$ so that $Tf_x(T_x {Z_X}) \subset T_{y} Z_Y$ making $F$ a map between pairs. Splitting $T_x X = N_x Z_X \oplus T_x Z_X$ pointwise we see that $T_x X / T_x Z_X \cong N_x Z_X$ and $Tf_x(T_x X) + T_y Z_Y = Tf_x(N_x Z_X) + T_y Z_Y$. Because $f$ is a transverse to $Z_Y$ we have $Tf_x(T_x X) + T_{y} Z_Y = T_{y} Y$. Similarly splitting $T_y Y = N_y Z_Y \oplus T_y Z_Y$ we conclude that $Tf_x(N_x Z_X) + T_y Z_Y = N_y Z_Y \oplus T_y Z_Y$ and we see that $Tf_x$ restricts to an isomorphism from $N_x Z_X \cong V / V_1$ to $N_y Z_Y \cong W / W_1$. By construction $Tf \circ \rho_X = \rho_Y \circ \varphi$, hence also pointwise at $x$ and $y$, and furthermore by definition of the log-tangent bundles we have $\im \rho_{X,x} = T_x Z_X$ and $\im \rho_{Y,y} = T_{y} Z_Y$. Note that $f^* \mathbb{L}_Y \cong \mathbb{L}_X$, because local defining functions for $Z_Y$ pull back to local defining functions for $Z_X$. Hence $\varphi_x\colon \mathbb{L}_{X,x} \to \mathbb{L}_{Y,y}$ is an isomorphism. Finally, $\ker \rho_{X,x} = \mathbb{L}_{X,x}$ and $\ker \rho_{Y,y} = \mathbb{L}_{Y,y}$, so that $\varphi_x\colon \ker \rho_{X,x} \to \ker \rho_{Y,y}$ is an isomorphism. By \autoref{lem:blinalg} we conclude that $\rho_{X,x}\colon \ker \varphi_x \to \ker Tf_x$ is an isomorphism.
\ep

The situation is summarized by the following diagram with exact rows.
\begin{center}
	\begin{tikzpicture}
	\matrix (m) [matrix of math nodes, row sep=2.5em, column sep=2.2em,text height=1.5ex, text depth=0.25ex]
	{	   &                               & \ker \varphi_x    & \ker Tf_x     &                             &   \\
		0 & \mathbb{L}_{X,x}   & T_x X(-\log Z_X)        & T_x X         & T_x X / T_x Z_X   & 0 \\
		0 & \mathbb{L}_{Y,y}   & T_y Y(-\log Z_Y)        & T_y Y          & T_y Y / T_y Z_Y   & 0 \\};
	\path[-stealth]
	(m-1-3) edge node [above] {$\rho_{X,x}$} node [below] {$\cong$} (m-1-4)
	edge [right hook-latex] (m-2-3)
	(m-1-4) edge [right hook-latex] (m-2-4)
	(m-2-1) edge (m-2-2)
	(m-2-2) edge [right hook-latex] (m-2-3)
	(m-2-2) edge node [left] {$\cong$} node [right] {$\varphi_x$} (m-3-2)
	(m-2-3) edge node [right] {$\varphi_x$} (m-3-3)
	edge node [above] {$\rho_{X,x}$} (m-2-4)
	(m-2-4) edge node [right] {$Tf_x$} (m-3-4)
	edge node [above] {${\rm pr}_x$} (m-2-5)
	(m-2-5) edge[anchor = south] (m-2-6)
	(m-2-5) edge node [right] {$\overline{Tf}_x$} node [left] {$\cong$} (m-3-5)
	(m-3-1) edge (m-3-2)
	(m-3-2) edge [right hook-latex] (m-3-3)
	(m-3-3) edge node [above] {$\rho_{Y,y}$} (m-3-4)
	(m-3-4) edge node [above] {${\rm pr}_y$} (m-3-5)
	(m-3-5) edge (m-3-6);
	\end{tikzpicture}
\end{center}
\begin{rem} The statement that $\overline{Tf}_x\colon N_x Z_X \to N_y Z_Y$ is an isomorphism for $y = f(x) \in Z_Y$ can be colloquially phrased as follows. As $f$ is transverse to $Z_Y$, the normal direction to $Z_Y$ at points in $Z_Y$ must be contained in the image of $TX$ under $Tf$. The fact that $f^{-1}(Z_Y) = Z_X$ then implies that it must in fact be obtained from the normal direction to $Z_X$. As both $N_x Z_X$ and $N_y Z_Y$ are one-dimensional subspaces and $Tf_x$ gives a surjection, it is an isomorphism.
\end{rem}
We next discuss sections of $b$-maps, in preparation of Section \ref{sec:bhyperfibration}.
\begin{prop}\label{prop:bmapsections} Let $f\colon (X,Z_X) \to (Y,Z_Y)$ be a $b$-map along with a section $s$. Then $s\colon (Y,Z_Y) \to (X,Z_X)$ is a $b$-map with induced Lie algebroid morphism $\psi$, and $\ker \varphi_x \oplus \psi_y(\mc{A}_{Z_Y,y}) = \mc{A}_{Z_X,x}$ for all $x \in X$, where $y = f(x)$.
\end{prop}
\bp By definition $f \circ s = {\rm id}_Y$, hence $s^{-1}(Z_X) = s^{-1}(f^{-1}(Z_Y)) = (f \circ s)^{-1}(Z_Y) = Z_Y$. Given finite-dimensional vector spaces $U, W, V', V$ with $V' \subset V$ so that $U \hookrightarrow V' \twoheadrightarrow W$ and $U \hookrightarrow V \twoheadrightarrow W$, we have $V' = V$ by counting dimensions. Let $y \in Z_Y$ and $x \in f^{-1}(y)$ be given. By definition $\ker Tf_x \hookrightarrow T_x X \twoheadrightarrow T_y Y$ using $Tf_x$. There is a surjection of $Ts_y(T_y Y)$ onto $T_y Y$. As $Ts_y(T_y Y) + T_x Z_X \subset T_x X$ and $\ker Tf_x \hookrightarrow T_x Z_X$ because $f^{-1}(Z_Y) = Z_X$, we see that $\ker Tf_x \hookrightarrow Ts_y(T_y Y) + T_x Z_X \twoheadrightarrow T_y Y$, so that $Ts_y(T_y Y) + T_x Z_X = T_x X$. We conclude that $s$ is a $b$-map. Similarly, consider $x \in X$ and denote $y = f(x)$. As $f$ is a $b$-map we have $\ker \varphi_x \hookrightarrow \mc{A}_{X,x} \twoheadrightarrow \mc{A}_{Y,y}$ using $\varphi_x$. Again there is a surjection of $\psi_y(\mc{A}_{Y,y})$ onto $\mc{A}_{Y,y}$. Note that $\ker \varphi_x \oplus \psi_y(\mc{A}_{Y,y}) \subset \mc{A}_{X,x}$ and furthermore $\ker \varphi_x \hookrightarrow \ker \varphi_x \oplus \psi_y (\mc{A}_{Y,y}) \twoheadrightarrow \mc{A}_{Y,y}$, so that $\ker \varphi_x \oplus \psi_y(\mc{A}_{Y,y}) = \mc{A}_{X,x}$.
\ep
Let $f\colon (X^4,Z_X) \to (\Sigma^2,Z_\Sigma)$ be a $b$-map for which $\ker \varphi$ is even-dimensional. For example, this situation arises when there is a log-symplectic form $\btwo{\omega}_\Sigma$ on $(\Sigma,Z_\Sigma)$ and a \bacs $\btwo{J}$ on $(X,Z_X)$ so that $\btwo{J}$ is $(\btwo{\omega}_\Sigma, f)$-tame. Because $f$ is transverse to $Z_\Sigma$, we conclude that $Z_\Sigma$ cannot contain critical values of $f$. Namely, if $x$ were a critical point of $f$, by a dimension count and \autoref{prop:kerbdf} we would have $\ker Tf_x = T_x X$, hence $Tf_x(T_x X) + T_{f(x)} Z_Y = T_{f(x)} Z_Y$, not $T_{f(x)} Y$. Said differently, such a map $f\colon X^4 \to \Sigma^2$ can only be turned into a $b$-map if $Z_\Sigma$ is disjoint from the set of critical values of $f$.

\section{Constructing \blog{}s}
\label{sec:blogconstructions}
In this section we apply the results of Chapter \ref{chap:constrasymp} to obtain existence results for \blog{}s. Let $X$ and $Y$ be compact connected manifolds, and assume that $Y$ is equipped with a log-Poisson structure. Using \autoref{prop:blogbsymp} we can view $Y$ as a $b$-manifold with a log-symplectic structure, obtaining a triple $(Y,Z_Y,\btwo{\omega}_Y)$.

Given a map $f\colon X \to Y$ such that $f$ is transverse to $Z_Y$, we can turn it into a $b$-map $f\colon (X,Z_X) \to (Y,Z_Y)$ by defining $Z_X := f^{-1}(Z_Y)$. We wish to equip $X$ with a \blog. We will use the following notation.
\bi
\item $F_y = f^{-1}(y)$ for $y \in Y$ is the level set, or \emph{fiber}, of $f$ over $y$;
\item $[F]$ is the homology class of a generic fiber.
\ei
Here a generic fiber is the inverse image of a regular value. This homology class will only be used when it is well defined and independent of the regular value.
\begin{rem} Given a $b$-map $f\colon (X,Z_X) \to (Y,Z_Y)$, the singular locus of $X$ is given by $Z_X = f^{-1}(Z_Y)$, so that it consists of fibers. This means the fibers of $f$, being natural candidates for log-symplectic submanifolds of $X$, never only hit $Z_X$.
\end{rem}
Using \autoref{thm:thurstonfibration} we obtain the $b$-version of Thurston's result for symplectic fiber bundles with two-dimensional fibers \cite{Thurston76}, adapting the proof by Gompf in \cite{Gompf01}.
\begin{thm}\label{thm:bthurstonfibration} Let $f\colon (X^{2n},Z_X) \to (Y^{2n-2}, Z_Y)$ be a $b$-fibration between compact connected $b$-manifolds. Assume that $(Y,Z_Y)$ is log-Poisson and that the generic fiber $F$ is orientable and $[F] \neq 0 \in H_2(X;\R)$. Then $(X,Z_X)$ is log-Poisson.
\end{thm}
This theorem has the following corollary, phrased without using the $b$-language. Given a log-Poisson structure $\pi$, we let $Z_\pi = (\wedge^n \pi)^{-1}(0)$ be its degeneracy locus.
\begin{cor}\label{cor:fibrationblogs} Let $f\colon X^{2n} \to Y^{2n-2}$ be a fibration between compact connected manifolds. Assume that $Y$ admits a log-Poisson structure $\pi$ and that the generic fiber $F$ is orientable and $[F] \neq 0 \in H_2(X;\R)$. Then $(X,Z_X)$ admits a log-Poisson structure for $Z_X = f^{-1}(Z_{\pi})$.
\end{cor}
\bp Let $Z_Y := Z_{\pi}$. Then by \autoref{prop:blogbsymp}, $(Y,Z_Y)$ is log-symplectic. By \autoref{rem:fibrintobfibr}, the map $f$ gives a $b$-fibration $f\colon (X,Z_X) \to (Y,Z_Y)$, where $Z_X = f^{-1}(Z_Y)$. Using \autoref{thm:bthurstonfibration} we conclude that $(X,Z_X)$ is log-Poisson.
\ep
For example, any homologically essential oriented surface bundle over a surface is log-Poisson (regardless of whether the base is orientable). The fact that any such base surface is log-Poisson is discussed in Section \ref{sec:logpoisson}.
\bp[ of \autoref{thm:bthurstonfibration}] We check the conditions of \autoref{thm:thurstonfibration}. If necessary, pass to a finite cover of $Y$ so that the fibers of $f$ are connected. We need only construct a closed $\mc{A}_{Z_X}$-two-form that is nondegenerate on $\ker \varphi$. This we will do locally using \autoref{prop:gluingfiberform}, by constructing the required neighbourhoods $W_y$ and forms $\eta_y$ for each fiber. Let $c \in H^2(X;\R)$ be through duality a class such that $\langle [F], c\rangle = 1$, using that $[F] \neq 0$. Given $y \in Y$, let $D_y \subset Y$ be an open disk containing $y$, fully contained in a trivializing neighbourhood of $f$ around $y$. Define $W_y := f^{-1}(D_y) \cong D_y \times F_y$. Using that $F_y$ is two-dimensional, choose an area form on $F_y$ inducing the preimage orientation of the fiber, and let $\overline{\eta}_y \in \Omega^2(W_y)$ be the pullback of this form via the projection $p\colon W_y \to F_y$. Because $\langle [F_y], \overline{\eta}_y \rangle = 1 = \langle [F_y], c\rangle$ and $H_2(W_y;\R)$ is generated by $[F_y]$,  it follows that $[\overline{\eta}_y] = c|_{W_y} \in H^2(W_y;\R)$. Now define $\eta_y := \rho_X^* \overline{\eta}_y$, and recall that by \autoref{prop:kerbdf} we know that $\rho_X$ gives a pointwise isomorphism between $\ker \varphi$ and $\ker Tf$. To check that $\eta_y$ tames $\btwo{J}$ on $\ker \varphi_x$ for $x \in W_y$, recall that there $\btwo{J}$ is defined via rotation. As $\eta_y$ is the pullback of the area form of a fiber, taming follows as its restriction to a fiber is an area form for that fiber. By \autoref{thm:thurstonfibration} we obtain a \blog on $(X,Z_X)$.
\ep
Next we revisit the result proven in \cite{Cavalcanti17} that any homologically essential \alf with orientable fibers gives rise to a \blog in dimension four, which can be chosen to be bona fide.
\begin{thm}[{\cite[Theorem 6.7]{Cavalcanti17}}] \label{thm:alflogsymp} Let $f\colon X^4 \to \Sigma^2$ be an \alf between compact connected manifolds. Assume that the generic fiber $F$ is orientable and $[F] \neq 0 \in H_2(X;\R)$. Then $X$ admits a \blog.
\end{thm}
\begin{rem} As remarked by Cavalcanti \cite{Cavalcanti17}, the assumption on homological essentialness cannot be dropped. Indeed, $S^4$ admits an \alf with orientable fibers (see \autoref{rem:alfs4}), however cannot admit a \blog by the cohomological obstruction of \autoref{thm:blogobstr}.
\end{rem}
This theorem should be viewed as a direct analogue of \autoref{thm:gompflfsymp}, that homologically essential four-dimensional Lefschetz fibrations provide symplectic structures. We prove it by first showing that an \alf gives rise to what we call a \emph{\blf}. Noting \autoref{prop:kerbdf} we use the following notation, given a map $f\colon X \to Y$ between manifolds (in what follows, $Y = \Sigma$).
\bi
\item $\Delta = {\rm Crit}(f) \subset X$ is the set of critical points of $f$;
\item $\Delta_y = \Delta \cap F_y$ for $y \in Y$ is the set of critical points of $f$ lying on the fiber $F_y$;
\item $\Delta' = {\rm Sing}(f) \subset Y$ is the set of singular values of $f$.
\ei
\begin{defn}\label{defn:blefschetzfibration} A \emph{\blf} is a $b$-map $f\colon (X^{2n},Z_X) \to (\Sigma^2, Z_\Sigma)$ between compact connected $b$-oriented $b$-manifolds so that for each critical point $x \in \Delta$ there exist complex coordinate charts compatible with orientations induced from the $b$-orientations centered at $x$ and $f(x)$ in which $f$ takes the form $f\colon \C^n \to \C$, $f(z_1, \dots, z_n) = z_1^2 + \dots + z_n^2$.
\end{defn}
\begin{rem}
	Given a \blf $f\colon (X,Z_X) \to (\Sigma,Z_\Sigma)$, the local model for $f$ around critical points $x \in \Delta$ implies $\ker Tf_x = T_x X$. Because $f$ is a $b$-map so is transverse to $Z_\Sigma$, we conclude that $Z_\Sigma$ and $\Delta'$ are disjoint. The $b$-orientation induces an orientation away from the singular locus so that it makes sense to demand compatibility of the charts. This discussion can be summarized by noting that a \blf is nothing more than a $b$-map whose induced Lie algebroid morphisms is a Lie algebroid Lefschetz fibration in the sense of \autoref{defn:lalf}.
\end{rem}

Alternatively, we could first define the notion of a $b$-\alf and then note that its critical values must be disjoint from the singular locus, so that we can further demand compatibility of the charts specifying the local model of $f$.
\begin{prop}\label{prop:alfblf} Let $f\colon X^{2n} \to \Sigma^2$ be an \alf between compact connected manifolds which is injective on critical points. Assume that the generic fiber $F$ is orientable and $[F] \neq 0 \in H_{2n-2}(X;\R)$. Then there exists a hypersurface $Z_\Sigma \subset \Sigma$ so that $f\colon (X^{2n},Z_X) \to (\Sigma^2,Z_\Sigma)$ is a \blf, where $Z_X = f^{-1}(Z_\Sigma)$.
\end{prop}
\bp We deal with orientations as in \cite[Theorem 6.7]{Cavalcanti17}. Fix an orientation for the generic fiber $F$. As $[F] \neq 0$, this forces $f\colon X \setminus \Delta \to \Sigma \setminus \Delta'$ to be an orientable fibration, which in turn orients all fibers, including the singular ones. We conclude that $X$ is orientable if and only if $\Sigma$ is. If they are, choose orientations and split $\Delta'$ into disjoint sets $\Delta'_+$ and $\Delta'_-$ according to the sign of the critical points. Then, pick a separating curve $\gamma \subset \Sigma$ disjoint from $\Delta'$ such that its interior contains all of $\Delta'_-$ and no points from $\Delta'_+$. If $\Sigma$ is not orientable, there instead exists a curve $\gamma \subset \Sigma$ so that $\Sigma \setminus \gamma$ is orientable, hence so is $X \setminus f^{-1}(\gamma)$. Choose orientations and then homotope $\gamma$ through negative critical values of $f$ so that all critical points are positive. Define $Z_\Sigma := \gamma$ and let $Z_X := f^{-1}(Z_\Sigma)$. Because $Z_\Sigma$ does not hit $\Delta'$, it is immediate that $f$ is a $b$-map from $(X,Z_X)$ to $(Y,Z_\Sigma)$. Moreover, the orientations we chose give the appropriate $b$-orientations. But then $f$ is a \blf.
\ep
\begin{rem}\label{rem:curvechoice} The curve $\gamma$ used in the previous proof is not unique. For example, in the orientable case we chose a separating curve, but we could just as well have chosen disjoint curves around each negative critical point separately. The effect of this is that the $b$-manifold structures that are used are not unique either.
\end{rem}
Given a four-dimensional \blf seen as a Lie algebroid Lefschetz fibration, we can use \autoref{thm:thurstonlfibration} to construct a log-symplectic structure on $X$.
\begin{thm}\label{thm:blfblog} Let $f\colon (X^4,Z_X) \to (\Sigma^2,Z_\Sigma)$ be a \blf between compact connected $b$-oriented $b$-manifolds which is injective on critical points. Assume that the generic fiber $F$ is orientable and $[F] \neq 0 \in H_2(X;\R)$. Then $(X,Z_X)$ admits a \blog.
\end{thm}
\bp We show the conditions of \autoref{thm:thurstonlfibration} are satisfied. Note that \autoref{prop:lfexactsequence} still holds for \alf{}s and hence also for \blf{}s. If necessary, use this exact sequence in homotopy to lift $f$ to a cover of $\Sigma$ so that $f$ has connected fibers. Fix an orientation for the generic fiber $F$, which orients $\ker Tf$ and $\ker \varphi$ at regular points, using \autoref{prop:kerbdf}. As $[F] \neq 0$, this forces $f\colon (X \setminus f^{-1}(\Delta'),Z_X \setminus f^{-1}(\Delta')) \to (Y\setminus \Delta',Z_Y\setminus \Delta')$ to be an orientable $b$-fibration. Let $c \in H^2(X;\R)$ be a class dual to $[F] \in H_2(X;\R)$, i.e.\ such that $\langle c, [F] \rangle = 1$. Let $\pi_{\mc{A}_{Z_\Sigma}} \in \Gamma(\mc{A}_{Z_\Sigma})$ be a transverse section specifying the $b$-orientation of $(\Sigma,Z_\Sigma)$. Then by \autoref{prop:standarddivideal}, $\rho_\Sigma(\pi_{\mc{A}_{Z_\Sigma}})$ is a log-Poisson structure on $(X,Z_\Sigma)$, and $\omega_{\mc{A}_{Z_\Sigma}} := \pi_{\mc{A}_{Z_\Sigma}}^{-1}$ is a log-symplectic form on $(\Sigma,Z_\Sigma)$.

We now need only construct the closed $\mc{A}_{Z_X}$-two-form as in \autoref{thm:thurstonlfibration}. However, this is very simple: take any two-form $\overline{\eta} \in \Omega^2_{\rm cl}(X)$ which is Poincar\'e dual to $c$, and define $\eta := \rho_X^* \overline{\eta} \in \Omega^2_{\rm cl}(\mc{A}_{Z_X})$. As the fibers of $f$ are connected, it is now immediate that $\eta$ satisfies the desired properties because $Z_X$ is disjoint from $\Delta$, noting \autoref{prop:kerbdf} to move between $\ker Tf$ and $\ker \varphi$.
\ep
\bp[ of \autoref{thm:alflogsymp}] This follows as an immediate corollary to \autoref{prop:alfblf} and \autoref{thm:blfblog}. We require that $f$ is injective on critical points so that critical values also obtain a sign, allowing us to group them appropriately. Here we note that we can always perturb $f$ so that it is injective on critical points. If one does not want to assume this one can proceed as in \cite[Lemma 3.3]{Gompf05}.
\ep
\begin{rem} The \blog from the proof of \autoref{thm:alflogsymp} has connected singular locus which may be empty. By \autoref{rem:curvechoice} we could have easily ensured that the singular locus has multiple components. Moreover, by \autoref{thm:addsinglocus} we can add components to the singular locus at will using a neighbourhood of a fiber. In particular we can always ensure the \blog is bona fide.
\end{rem}
It is an interesting question whether every \blog on a four-manifold can be obtained out of an \alf using the construction of \autoref{thm:alflogsymp}. This parallels the development by Gompf and Donaldson between Lefschetz pencils and symplectic structures on four-manifolds. The first thing to note is that in our construction the fibers are always contained in the eventual singular locus, or are disjoint from it. Moreover, all \blog{}s we construct are \emph{proper}, in that all connected components of their singular loci are compact and have a compact symplectic leaf. This implies compact log-symplectic four-manifolds $(X,Z_X,\pi)$ for which $\pi$ is not proper are not reached by our construction. Note however that one can deform $\pi$ to a proper \blog if all components of $Z_X$ are compact \cite[Theorem 3.6]{Cavalcanti17}. More serious is the fact that in our construction the singular locus fibers over a circle in the base with specified diffeomorphism type of the fiber. Hence our construction cannot create log-symplectic four-manifolds $X$ with a disconnected singular locus $Z_X$ with at least two components not fibering over $S^1$ with the same genus fiber. The following example shows this can indeed happen.
\begin{exa} Let $X^4 = \Sigma_g \times \Sigma_h$ be the product of compact surfaces of genera $g \geq 2$ and $h \geq 1$ carrying the product symplectic form, and consider the map $f\colon X \to \Sigma_h$ given by projection. Consider a copy of the torus $T^2 = S^1 \times S^1$ by picking an essential circle in both the base and the fiber. It is Lagrangian and homologically nontrivial, so that by a result of Gompf \cite[Lemma 1.6]{Gompf95}, by a slight perturbation of the symplectic structure in a neighbourhood of the torus we find a symplectic structure on $X$ for which $T^2$ is symplectic. Applying \autoref{thm:addsinglocus} we obtain a log-symplectic structure on $(X, T^3)$. Now use a fiber $\Sigma_g$ of $f$ disjoint from the torus and apply \autoref{thm:addsinglocus} once more to obtain a \blog on $(X, Z_X)$, where $Z_X = T^3 \sqcup \Sigma_g \times S^1$. 
	
	Note that $T^3$ cannot fiber over $S^1$ with fibers of genus other than one. Any fibration $p\colon T^3 \to S^1$ with fiber $F$ induces a long exact sequence in homotopy groups, a part of which reads $\pi_2(S^1) \to \pi_1(F) \to \pi_1(T^3) \to \pi_1(S^1)$, or more concretely $0 \to \pi_1(F) \hookrightarrow \Z^3 \to \Z$. This shows that $\pi_1(F)$ injects into the Abelian group $\Z^3$, hence must itself be Abelian and cannot have torsion. By counting its rank must be two, so that that the genus of $F$ must be one, and the fibers of $p$ are tori. Similarly, the product $Y := \Sigma_g \times S^1$ cannot fiber over $S^1$ with fiber $F$ being a torus. We have $b_1(Y) = 1 + b_1(\Sigma_g) = 2g +1 \geq 5$ as $g \geq 2$. However, if $F$ were a torus, $Y$ would be a mapping torus of $T^2$, hence $b_1(Y) \leq 1 + b_1(T^2) = 3$ which is a contradiction. By the discussion preceding this example we conclude that the log-symplectic structure on $(X,Z_X)$ cannot be obtained from the construction in \autoref{thm:alflogsymp}.
\end{exa}
\section{\texorpdfstring{$b$}{b}-Hyperfibrations}
\label{sec:bhyperfibration}
In this section we introduce a class of fibration-like maps with two-dimensional fibers that can be made to satisfy the conditions of \autoref{thm:lathurstontrick}. We call them \emph{$b$-hyperfibrations}, as they are the $b$-analogue of the notion of a hyperpencil with empty base locus, as introduced by Gompf in \cite{Gompf04}. After defining them we show that a $b$-hyperfibration satisfying a condition analogous to having homologically essential fibers gives rise to a \blog.

Let $E, F \to X$ be real vector bundles over a manifold $X$ and let $T\colon E \to F$ be a continuous bundle map. Call a point $x \in X$ \emph{regular} if $T_x\colon E_x \to F_x$ is surjective, and \emph{critical} otherwise. Let ${\rm Reg}(T)$ denote the spaces of regular points of $T$. Form the space $P \subset E$ by
\be
P = \ol{\bigcup\limits_{x \in {\rm Reg}(T)} \ker T_x},
\ee
and let $P_x = P \cap E_x$ for $x \in X$. Then $P_x = \ker T_x$ when $x$ is regular, and otherwise $P_x \subset \ker T_x$ consists of all limits of sequences of vectors at regular points which are annihilated by $T$.
\begin{defn}\label{defn:wrapped} A point $x \in X$ is called \emph{$T$-wrapped} if ${\rm span}_\R P_x$ has real codimension at most two in $\ker T_x$.
\end{defn}
Note that all points $x \in {\rm Reg}(T)$ are wrapped as then $P_x = \ker T_x$. In our applications $\ker T_x$ will be even-dimensional, so that the wrappedness condition is immediate unless ${\rm rank}(E) \geq 6$. Let $f\colon (X,Z_X) \to (Y,Z_Y)$ be a $b$-map between compact connected $b$-manifolds.  As before, let $\Delta$ and $\Delta_y$ be the set of critical points of the map $f\colon X \to Y$ and those on the fiber $F_y$ for $y \in Y$ respectively. We can then apply the above definition to the continuous or in fact smooth map $\varphi\colon \mc{A}_{Z_X} \to \mc{A}_{Z_Y}$, noting \autoref{prop:kerbdf}. Note that if $X$ is four-dimensional and $Y$ is a surface, every critical point of $f$ will be automatically $\varphi$-wrapped if $\ker \varphi$ is always even-dimensional. This is because $\ker \varphi$ is two-dimensional at regular points and the dimension of $\ker \varphi$ cannot exceed four at singular points. With this we can give the definition of a $b$-hyperfibration.
\begin{defn}\label{defn:bhyperfibration} A \emph{$b$-hyperfibration} is a $b$-map $f\colon (X^{2n},Z_X) \to (Y^{2n-2},Z_Y)$ between compact connected $b$-oriented $b$-manifolds so that there exists a \blog $\btwo{\omega}_Y$ on $(Y, Z_Y)$ and such that
	\bi
	\item[i)] each critical point $x \in \Delta$ is $\varphi$-wrapped;
	\item[ii)] for each critical point $x \in \Delta$, there exists a neighbourhood of $x$ and a \bacs $\btwo{J}_x$ on the $b$-manifold $(W_x, W_x \cap Z_X)$ such that $\btwo{J}_x$ is $(\btwo{\omega}_Y,f)$-tame;
	\item[iii)] each $y \in Y$ has only finitely many critical points $\Delta_y$ lying on its fiber $F_y$.
	\ei
\end{defn}
This should be compared with \cite[Definition 2.4]{Gompf04}. The definition does not require $X$ nor $Y$ to be orientable. When $Z_Y = \emptyset$ one almost recovers the definition by Gompf of a hyperpencil with empty base locus, over an arbitrary symplectic base.
\begin{rem} While a $b$-hyperfibration $f$ may have infinitely many critical points, note that regular points of $f$ are dense in $X$, arguing as in \cite[Theorem 2.11]{Gompf04}. If an open $W \subset X$ would consist entirely of critical points, choose a point $x_0 \in W$ which minimizes $\dim \ker \varphi_x$, and using \autoref{prop:kerbdf} note that $\ker \varphi \cong \ker Tf$ is a smooth distribution near $x_0$ as it can be realized as $\ker T(\pi \circ f)$ for a projection $\pi$. Then take a vector field in $\ker Tf$ and integrate it to obtain a curve of critical points all lying in a single fiber of $f$. This contradicts assumption iii) in the definition of a $b$-hyperfibration.
\end{rem}
The wrappedness of the critical points will be used to obtain a global $(\btwo{\omega}_Y,f)$-tame \bacs out of the locally existing ones. With this new notion in hand we can move on to the following result, which is the appropriate $b$-analogue of \cite[Theorem 2.11]{Gompf04}. Given a $b$-hyperfibration $f\colon (X,Z_X) \to (Y,Z_Y)$ and $y \in Y$ fixed, we refer to the closures of the connected components of $F_y \setminus \Delta_y$ as the \emph{components} of the fiber $F_y$.
\begin{thm}\label{thm:bhyperfibration} Let $f\colon (X,Z_X) \to (Y, Z_Y, \btwo{\omega}_Y)$ be a $b$-hyperfibration between compact connected $b$-oriented $b$-manifolds. Assume that there exists a finite collection $S$ of sections of $f$ intersecting all fiber components non-negatively and for each fiber component at least one section in $S$ intersecting positively. Then $(X,Z_X)$ admits a \blog.
\end{thm}
Note that the condition on the existence of such a collection $S$ of fibers of $f$ implies that each component of each fiber is homologically essential.
\begin{rem} In the above theorem, sections $s \in S$ are $b$-maps $s\colon (Y,Z_Y) \to (X,Z_X)$ by \autoref{prop:bmapsections} and furthermore that $\ker \varphi_x \oplus \psi_y(\mc{A}_{Y,y}) = \mc{A}_{X,x}$ for all $x \in X$, where $y = f(x)$. In the proof of \autoref{thm:bhyperfibration} we show a $b$-hyperfibration naturally gives rise to a \bacs which is $(\btwo{\omega}_Y,f)$-tame so that $\ker \varphi$ carries a $b$-orientation. By \autoref{prop:kerbdf} for smooth points $x \in X \setminus \Delta$ we then have an orientation for $T_x F_y = \ker Tf_x \cong \ker \varphi_x$. Both $\mc{A}_{X,x}$ and $\mc{A}_{Y,y}$ carry orientations, hence so does $\psi_y(\mc{A}_{Y,y})$. We can then define the positive intersection of $s \in S$ with $F_y$ by comparing the $b$-orientations on these tangent spaces in the usual way. Note that $s$ must intersect fibers in smooth points of $f$ as it is a section.
\end{rem}
The proof of \autoref{thm:bhyperfibration} will be modelled on Gompf's proof of \cite[Theorem 2.11]{Gompf04}. There will be an interplay between two types of singular behavior, namely that of the $b$-manifold structure and that of the fibration itself. The relation between these has been discussed before in \autoref{prop:kerbdf}. Note that in the case of a $b$-hyperfibration $f\colon (X,Z_X) \to (Y, Z_Y)$ with $Z_Y = \emptyset$, we cannot apply \cite[Theorem 2.11]{Gompf04} directly. Indeed, there is no base locus, but instead a set of sections $S$. This is akin to obtaining a \lf out of a Lefschetz pencil by blowing up the base locus.

For $V$ a real finite-dimensional vector space, let $\mc{B}_V \subset {\rm Aut}(V)$ be the open set of linear operators on $V$ with no real eigenvalues, and $\mc{J}_V \subset \mc{B}_V$ the set of complex structures on $V$ in either orientation. The following lemma is proven in \cite{Gompf04} as Corollary 4.2, loc.\ cit.
\begin{lem}\label{lem:acsretraction} Let $V, W$ be real finite-dimensional vector spaces. Then there exists a canonical real-analytic retraction $j\colon \mc{B}_V \to \mc{J}_V$, satisfying for all linear maps $T\colon V \to W$ such that $T A = B T$, that $T j(A) = j(B) T$ (whenever both sides are defined).
\end{lem}
Because the retraction in the previous lemma is canonical we can apply it pointwise to a continuously varying map, to again obtain a continuous map. Now let $E^{2n}, F^{2n-2} \to X$ be real oriented vector bundles over a compact manifold $X$. In what follows, a \emph{two-form} on a vector bundle is a continuously varying choice of skew-symmetric bilinear form on each fiber. The next proposition can be extracted from \cite[Lemma 3.2]{Gompf04}. We include a proof for completeness.
\begin{prop}\label{prop:acsconstruct} Let $T\colon E \to F$ be a continuous bundle map and $\omega_F$ a non-degenerate two-form on $F$. Assume that for all $x \in X$ there exists a neighbourhood $W_x$ of $x$ with an $(\omega_F, T)$-tame complex structure on $E|_{W_x}$. Assume that each critical point $x \in {\rm Crit}(T)$ is wrapped. Then there exists a continuous $(\omega_F, T)$-tame complex structure $J$ on $E$.
\end{prop}
\bp
Cover $X$ by open sets $W_\alpha$ equipped with complex structures $J_\alpha$ on $E|_{W_\alpha}$ which are $(\omega_F, T)$-tame. Let $\{\varphi_\alpha\}$ be a subordinate partition of unity, and define
\be
A := \sum_\alpha \varphi_\alpha J_\alpha\colon E \to E, \qquad\qquad B := \sum_\alpha \varphi_\alpha T_* J_\alpha\colon T(E) \to T(E),
\ee
so that $TA = BT$. Since $\ker T$ is $J_\alpha$-complex for all $\alpha$ because $J_\alpha$ is $(\omega_F, T)$-tame (see below \autoref{defn:taming}), $J_\alpha$ descends to a map $T_* J_\alpha\colon T(E) \to T(E)$, hence $B$ is well-defined. In order to apply \autoref{lem:acsretraction} we show that $A_x \in \mc{B}_{E_x}$ for all $x \in X$, i.e.\ that $A$ has no real eigenvalues. Let $\lambda$ be an eigenvalue of $A$ with eigenvector $v \in E$. Then $B T v = T A v = T \lambda v = \lambda T v$, so either $T v = 0$, or $T v$ is a $\lambda$-eigenvector for $B$. As each $T_* J_\alpha$ is $\omega_F$-tame, $B$ has no real eigenvalues on any fiber. Indeed, for $0 \neq w \in T(E)$, we have $\omega_F(w,w) = 0$. Hence if $B w = \lambda w$ for some $\lambda \in \R$, we have
\be
0 = \omega_F(w, \lambda w) = \omega_F(w, \sum_\alpha \varphi_\alpha T_* J_\alpha w) = \sum_\alpha \varphi_\alpha \omega_F(w, T_* J_\alpha w) > 0,
\ee
which is a contradiction. We conclude that any real eigenvector of $A$ must lie in $\ker T$. Let $x \in X$ be given. As $T$-regular points are always $T$-wrapped, and by hypothesis the same holds for all $T$-critical points, we know that $x$ is $T$-wrapped. Recall the subspace $P \subset \ker T \subset E$ used in \autoref{defn:wrapped}. We construct a decomposition ${\rm span}_\R P_x = \bigoplus_j \Pi_j$, with each $\Pi_j$ a real two-plane which is a $J_\alpha$-complex line for all $J_\alpha$ defined on $E_x$. Let $v \in P_x$ be given. By definition of $P$, there exists a sequence $(x_i)_{i \in \N}$ of $T$-regular points converging to $x$ and elements $v_i \in \ker T_{x_i}$ such that $v = \lim_{i\to\infty} v_i$. As the points $x_i$ are $T$-regular, the subspaces $\ker T_{x_i}$ are two-planes in $E_{x_i}$ oriented by the fiber-first convention. Pass to a subsequence so that the $\ker T_{x_i}$ converge to an oriented two-plane $\Pi \subset P_x$ containing $v$. Consider an open $W_\alpha$ containing $x$. Then there exists an $N_\alpha \in \N$ such that $x_i \in W_\alpha$ for all $i \geq N_\alpha$. But then for all $i \geq N_\alpha$, $\ker T_{x_i}$ is a $J_\alpha$-complex line, hence so is their limit $\Pi$. We conclude that $\Pi$ is a $J_\alpha$-complex line for each $J_\alpha$ defined at $x$. Proceed by induction to constructed multiple such real oriented two-planes $\Pi_j \subset P_x$ so that ${\rm span}_\R P_x = \bigoplus_j \Pi_j$, with each $\Pi_j$ being a $J_\alpha$-complex line for all $J_\alpha$ defined at $x$. Consider the quotient $Q_x := \ker T_x / {\rm span}_\R P_x$, which inherits an orientation from $\ker T_x$, which in turn is oriented as it is $J_\alpha$-complex for all $\alpha$ defined at $x$, all of which are $(\omega_F,T)$-tame. Then $Q_x$ inherits complex structures $\ol{J}_\alpha$ from each $J_\alpha$ defined at $x$, and these are all compatible with the orientation on $Q_x$. As $x$ is $T$-wrapped, $\dim_\C Q_x \leq 1$. But then there exists a fixed nondegenerate skew-symmetric bilinear form $\omega_x$ on $Q_x$ so that all $\ol{J}_\alpha$ are $\omega_x$-tame, as one can just pick an $\omega_x$ giving the orientation on $Q_x$. Consider the map $\ol{A}_x := \sum_\alpha \varphi_\alpha(x) \ol{J}_\alpha\colon Q_x \to Q_x$. Then $\ol{A}_x$ has no real eigenvalues on $Q_x$. As before, if $\ol{A}_x \ol{v} = \lambda \ol{v}$ for $0 \neq \ol{v} \in Q_x$ with $\lambda$ real, we would have
\be
0 = \omega_x(\ol{v},\lambda \ol{v}) = \omega_x(\ol{v}, \sum_\alpha \varphi_\alpha(x) \ol{J}_\alpha \ol{v}) = \sum_\alpha \varphi_\alpha(x) \omega_x(\ol{v}, \ol{J}_\alpha \ol{v}) > 0.
\ee
As $\ker T_x \cong {\rm span}_\R P_x \oplus Q_x$, we conclude that any real eigenvector of $A_x$ must lie in ${\rm span}_\R P_x = \bigoplus_j \Pi_j$. Construct a direct sum two-form $\wt{\omega}_x = \bigoplus_j \omega_j$ on ${\rm span}_\R P_x$ which tames each $J_\alpha$ at $x$. Then if $A_x v = \lambda v$ for $0 \neq v \in {\rm span}_\R P_x$ with $\lambda$ real and $v = \bigoplus_j v_j$ with respect to the direct sum decomposition of ${\rm span}_\R P_x$,
\be
0 = \wt{\omega}_x(v, \lambda v) = \sum_j \omega_j(v_j, A_x v_j) = \sum_{j,\alpha} \varphi_\alpha(x) \omega_j(v_j, J_\alpha(x) v_j) > 0.
\ee
We conclude that $A_x$ has no real eigenvalues, hence nor does $A$. By \autoref{lem:acsretraction} we obtain from $A$ a continuous complex structure $J = j(A)$ on $E$. As $T j(A) = j(B) T$ and convex combinations of $(\omega_F, T)$-tame endomorphisms are still tamed, the resulting \acs $J$ is $(\omega_F, T)$-tame.
\ep
\bp[ of \autoref{thm:bhyperfibration}] By \autoref{prop:blogbsymp} we need to show the existence of a log-symplectic form on $(X,Z_X)$. This is done by appealing to \autoref{thm:lathurstontrick} and \autoref{prop:gluingfiberform}, hence it suffices to construct local closed two-forms $\overline{\eta}_y$ around fibers so that the respective log-forms $\eta_y = \rho_X^* \overline{\eta}_y$ tame a global \bacs $\btwo{J}$ on $\ker \varphi$, and such that they are all cohomologous to the restriction of one global class $c \in H^2(X;\R)$. Construct around each point $x \in X$ a neighbourhood $W_x$ and a \bacs $\btwo{J}_x$ on $\mc{A}_{Z_X}|_{W_x}$ which is $(\btwo{\omega}_Y, f)$-tame. These exist by definition around critical points of $X$ and away from critical points these are constructed as in the proof of \autoref{thm:bthurstonfibration}, using that the fibers of $f$ are two-dimensional. Apply \autoref{prop:acsconstruct} to the situation where $E = \mc{A}_{Z_X}$, $F = f^*\mc{A}_{Z_Y}$, $T = \varphi$ and $\omega_F = f^* \btwo{\omega}_Y$, to obtain a global \bacs $\btwo{J}$ on $X$ which is $(\btwo{\omega}_Y, f)$-tame. Let $S$ be the finite set of sections in the hypothesis, and define $c$ to be the Poincar\'e dual of $\sum_{s \in S} [{\rm im}(s)] \in H_{2n-2}(X;\R)$. To apply Poincar\'e duality in the absence of an orientation on $X$, the images of the sections must be cooriented. However, this exactly means that the fibers of $f$ must be oriented. Regular fibers obtain an orientation through $\btwo{J}$ and \autoref{prop:kerbdf}, while singular fibers are oriented in their smooth locus.

Let $y \in Y$ be given. By the definition of a hyperfibration, $\Delta_y$ is finite. Let $C_y \subset X$ be the disjoint union of closed balls centered at each point in $\Delta_y$. Choose a closed two-form $\sigma \in \Omega^2(C_y)$ so that $\rho_X^* \sigma$ tames $\btwo{J}$ on $\ker \varphi|_{\Delta_y}$, noting that this is a condition at a finite set of points hence can easily be satisfied. Then, as taming $\btwo{J}$ on $\ker \varphi$ is an open condition by \autoref{prop:tamingopen}, after possibly shrinking the balls in $C_y$ we can assume that $\btwo{J}|_{\ker \varphi}$ is $\rho_X^* \sigma$-tame on the entirety of $C_y$. As by assumption $\btwo{J}$ is $(\omega_Y,f)$-tame, $F_y \setminus \Delta_y$ is a smooth noncompact $\btwo{J}$-holomorphic curve in $X \setminus \Delta_y$, whose complex $b$-orientation from $\btwo{J}$ agrees with its preimage $b$-orientation. By assumption there exists for each component of $F_y \setminus \Delta_y$ a section in $S$ intersecting that component positively. Choose $C_y$ so that $\partial C_y$ is transverse to $F_y$ and consider the intersecting circles $F_y \cap \partial C_y$. Connect each such circle to these points of intersection by a path in $F_y \setminus \Delta_y$. Let $C_y^0$ be the disjoint union of smaller concentric closed balls around $\Delta_y$ disjoint from these paths and with $\partial C_y^0$ transverse to $F_y$. Then each component $F_i$ of the compact surface $F_y \setminus {\rm int}(C_y^0)$ either lies fully inside ${\rm int}(C_y)$, or has a point of intersection with a section in $S$. In the latter case we will say that $F_i$ \emph{intersects $S$}. See the following figures for illustration.

\newcommand{\fibercurve}[2][black]{
	\draw[line width = #2, #1]
	plot[smooth, tension=.7] coordinates{($(a) + (-0.9,0)$) ($(a) + (-2.8,1)$) ($(a) + (-3.5,2.3)$) ($(a) + (-2.5,2.3)$) ($(a) + (-0.5,1)$) ($(a) + (-1.5,3.8)$) ($(a) + (-.5,4)$) ($(a) + (0.7,1.3)$) ($(a) + (2.3,1.2)$) ($(a) + (2.3,.7)$) (a)};
}
\begin{center}
	\begin{tikzpicture}[scale=1,x=1em,y=1em]
	\begin{scope}[xshift=-180pt]
	\coordinate (a) at (0,0);
	\fill (a) circle (0.25) node [label=below:{$\Delta_y$}]{};
	\draw[thin]
	(a) circle (2) node [label={[label distance = 0.9em]65:{$C_y^0$}}]{}
	(a) circle (3.5) node [label={[label distance = 2.5em]65:{$C_y$}}]{}
	($(a) + (-7,0)$) -- ($(a) + (-0.9,0)$) node [above, near start] {$F_y$}
	(a) -- ($(a) + (7,0)$);
	\fibercurve{0.4pt}
	\end{scope}
	\begin{scope}
	\coordinate (a) at (0,0);
	\fill (a) circle (0.25) node [label=below:{$\Delta_y$}]{};
	\draw[thin]
	(a) circle (2) node [label={[label distance = 0.9em]65:{$C_y^0$}}]{}
	(a) circle (3.5) node [label={[label distance = 2.5em]65:{$C_y$}}]{}
	($(a) + (-7,0)$) -- ($(a) + (-0.9,0)$) node [above, near start] {$F_i$}
	(a) -- ($(a) + (7,0)$);
	\fibercurve{0.4pt}
	\draw[very thick]
	($(a) + (-7,0)$) -- ($(a) + (-2,0)$)
	($(a) + (2,0)$) -- ($(a) + (7,0)$);
	\draw[thick]
	($(a) + (-3.1,2.075)$) -- ($(a) + (-2.8,2.825)$) 
	($(a) + (-1,3.85)$) -- ($(a) + (-1,4.65)$) 
	($(a) + (-4.5,-0.4)$) -- ($(a) + (-4.5,0.4)$) 
	($(a) + (+5,-0.4)$) -- ($(a) + (+5,0.4)$) node [label={[label distance = 0.5em]below:{$S$}}]{}; 
	\end{scope}
	\begin{scope}[even odd rule]
	\coordinate (a) at (0,0);
	\clip ($(a) + (-7,-3.5)$) rectangle ($(a) + (7,5)$)
	(a) circle (2);
	\fibercurve{1.4pt}
	\end{scope}
	\end{tikzpicture}
\end{center}

Let $W_y$ be the union of ${\rm int}(C_y^0)$ with a tubular neighbourhood rel boundary of $F_y \setminus {\rm int}(C_y^0)$ inside $X \setminus {\rm int}(C_y^0)$. Extend each $F_i$ to a closed oriented smooth surface $\wh{F}_i \subset W_y$ by arbitrarily attaching a surface inside $C_y^0$. Then the classes $[\wh{F}_i] \in H_2(W_y;\Z)$ form a basis for the homology of $W_y$. Indeed, contracting the whole neighbourhood $C_y^0$ to $\Delta_y$ we see that $F_y$ becomes homotopy equivalent to a wedge of the $F_i$. This $W_y$ is the desired neighbourhood of $F_y$ on which we construct the form $\overline{\eta}_y$ in the hypothesis of \autoref{prop:gluingfiberform}, using ideas similar to the proof of \autoref{thm:thurstonlfibration}. Since $F_y$ is $\btwo{J}$-holomorphic with $\btwo{J}|_{\ker \varphi}$ being $\rho_X^*\sigma$-tame on $C_y$, $\sigma|_{F_i \cap C_y}$ is a positive area form for each $i$. For each $F_i$ intersecting $S$, let $\sigma_i$ be an extension of $\sigma$ to $F_i$ as a positive area form with total area $\langle \sigma_i, [\wh{F}_i] \rangle$ equal to $\#(F_i \cap \im(S)) > 0$. Let $p$ denote the tubular neighbourhood projection onto $F_y \setminus {\rm int}(C_y^0)$ and let $f\colon C_y \to [0,1]$ be a smooth radial function defined on each ball around points in $\Delta_y$ so that $f \equiv 0$ in a neighbourhood of $C_y^0$ and $f \equiv 1$ in a neighbourhood of $\partial C_y$, which is smoothly extended to $W_y$ by being identically $1$ outside of $C_y$. See the following figure for illustration.

\begin{center}
	\begin{tikzpicture}[scale=1,x=1em,y=1em]
	\coordinate (a) at (0,0);
	
	\begin{scope}[even odd rule]
	\clip ($(a) + (-7.1,-3.5)$) rectangle ($(a) + (7.1,5)$)
	(a) circle (2);
	\draw[line width = 2.3mm]
	($(a) + (-7,0)$) -- ($(a) + (7,0)$) node [above, very near start] {$W_y$};
	\draw[line width = 1.5mm, white]
	($(a) + (-7,0)$) -- ($(a) + (7,0)$);
	\draw[thin]
	($(a) + (-7,0)$) -- ($(a) + (7,0)$);
	\end{scope}
	
	\begin{scope}[even odd rule]
	\coordinate (a) at (0,0);
	\clip ($(a) + (-7,-3.5)$) rectangle ($(a) + (7,5)$)
	(a) circle (2);
	
	\fibercurve{2.3mm}
	\fibercurve[white]{1.5mm}
	\fibercurve{0.4pt}
	\end{scope}
	
	\begin{scope}
	\fill (a) circle (0.25) node [label=below:{$\Delta_y$}]{};
	\draw[line width=0.4mm]
	(a) circle (2) node [label={[label distance = 0.85em]65:{$C_y^0$}}]{};
	\draw[thin]
	(a) circle (3.5) node [label={[label distance = 2.5em]65:{$C_y$}}]{}
	($(a) + (-7,0)$) -- ($(a) + (-0.9,0)$) 
	(a) -- ($(a) + (7,0)$);
	\fibercurve{0.4pt}
	\end{scope}
	\end{tikzpicture}
\end{center}

On the balls $C_y$, the form $\sigma$ is exact, say equal to $\sigma = d\alpha$ for $\alpha \in \Omega^1(C_y)$. Define a two-form $\overline{\eta}_y$ on $W_y$ by $\overline{\eta}_y := \sum_i p^* (f \sigma_i) + d((1-f) \alpha)$. In other words, $\sigma$ is extended by $0$ outside $C_y$, while the $\sigma_i$ are extended by $0$ inside $C_y^0$. For all functions $f \neq 0$ the form $f \sigma_i$ is a closed area form defined on the surface $F_i$, so that $\overline{\eta}_y$ is nonnegative and closed. On $W_y \cap C_y^0 = C_y^0$ we have $f \equiv 0$ so that $\overline{\eta}_y = d\alpha = \sigma$, hence there $\eta_y = \rho_X^* \overline{\eta}_y = \rho_X^* \sigma$ tames $\btwo{J}|_{\ker \varphi}$. Similarly, on $F_y \setminus \Delta_y$ the forms $\sigma_i$ are area forms for the orientation given by $\btwo{J}$. But then $\eta_y$ tames $\btwo{J}$ on $T_x F_y = \ker Tf_x \cong \ker \varphi_x$ for all $x \in F_y \setminus \Delta_y$ as this condition is convex. By openness of the taming condition, narrowing the tubular neighbourhood $p$ defining $W_y$ ensures that $\eta_y$ tames $\btwo{J}|_{\ker \varphi_x}$ for all $x \in W_y \setminus {\rm int}(C_y^0)$. But then $\eta_y$ tames $\btwo{J}$ on $\ker \varphi_x$ for all $x \in W_y$.

What remains is to show that $c|_{W_y} = [\overline{\eta}_y] \in H^2(W_y; \R)$. Recall that every component $F_i$ of $F_y \setminus {\rm int}(C_y^0)$ either intersects $S$ or lies in ${\rm int}(C_y)$. For those $F_i$ intersecting $S$ we have $\langle [\overline{\eta}_y], [\wh{F}_i] \rangle = \langle [\sigma_i], [\wh{F}_i] \rangle = \# (F_i \cap {\rm im}(S)) = \langle c, [\wh{F}_i] \rangle$. For $F_i$ disjoint from $S$ we know that $F_i \subset C_y$, but $\overline{\eta}_y$ is exact in $C_y$, so that we have $\langle [\overline{\eta}_y], [\wh{F}_i] \rangle = \langle 0, [\wh{F}_i] \rangle = 0 = \langle c, [\wh{F_i}] \rangle$. As $[\overline{\eta}_y]$ agrees with the class $c|_{W_y} \in H^2(W_y; \R)$ when evaluated on the classes [$\wh{F}_i$], which form a basis of $H_2(W_y;\Z)$, we see that $[\overline{\eta}_y] = c|_{W_y} \in H^2(W_y; \R)$ as desired. Applying \autoref{prop:gluingfiberform} and then \autoref{thm:lathurstontrick} we conclude that $(X,Z_X)$ admits a \blog.
\ep
\chapter{Constructing stable generalized complex structures}
\label{chap:constructingsgcs}
In this chapter we describe how to construct \sgcs{} out of \blog{}s using the techniques developed in Chapter \ref{chap:constrasymp}.

Recall from Section \ref{sec:scgs} that \sgcs{}s $\mc{J}$ can alternatively be viewed as symplectic forms in the elliptic tangent bundle defined using the complex log divisor $(K^*_\mc{J},s)$ (they are further of zero elliptic residue). This symplectic viewpoint allows one to apply symplectic techniques to the study of \sgcs{}s.

In this chapter we introduce a class of maps called boundary maps. These induce morphisms between the elliptic and log tangent bundles, and we use them to construct \sgcs{}s out of \blog{}s. A specific type of boundary map we call a \placeholder{} is shown to induce a Lie algebroid Lefschetz fibration between the elliptic and log-tangent bundle. This concept extends and formalizes the generalized Lefschetz fibrations of \cite{CavalcantiGualtieri09}. Using this correspondence we are able to prove the following result (\autoref{thm:lfibrationsgcs}) which was stated in the introduction to this thesis.
\newtheorem*{thm:sgcs}{Theorem \ref{thm:introsgcs}}
\begin{thm:sgcs} Let $f\colon (X^4,D) \to (\Sigma,Z)$ be a boundary Lefschetz fibration such that $(\Sigma,Z)$ carries a \blog. Assume that ${[F] \neq 0 \in H_2(X \setminus D;\R)}$ for the generic fiber $F$ and that $D$ is coorientable. Then $(X,D)$ admits a \sgcs.
\end{thm:sgcs}
A similar result is true in arbitrary dimension, using boundary fibrations instead (\autoref{thm:fibrationsgcs}). This result validates the definition of a boundary Lefschetz fibration as being the type of fibration-like map which should be linked with \sgcs{}s.

To produce concrete examples, one has to construct \placeholder{}s on explicit four-manifolds. This can be done using genus one Lefschetz fibrations over punctured surfaces with boundary monodromy given by powers of Dehn twists using a process we call completion (see \autoref{cor:lfcompletion}). In particular, we have the following corollary (\autoref{thm:g1lfsgcs}).
\begin{thm2} Let $f\colon X^4 \to D^2$ be a genus one Lefschetz fibration over the disk whose monodromy around the boundary is a power of a Dehn twist. Then all possible completions $\wt{f}\colon (\wt{X}, D) \to (D^2, \partial D^2)$ admit a \sgcs.
\end{thm2}
Together with Stefan Behrens and Gil Cavalcanti we have classified boundary Lefschetz fibrations over the disk, establishing the following.
\begin{thm}[\cite{BehrensCavalcantiKlaasse17}]\label{thm:classbdrylefschetz} Let $f\colon (X^4,D)\to (D^2,\partial D^2)$ be a relatively minimal boundary Lefschetz fibration. Then $X$ is diffeomorphic to one of the following:
	\begin{enumerate}
		\item $S^1 \times S^3$;
		\item $\#m (S^2 \times S^2)$, including $S^4$ for $m=0$;
		\item $\#m \C P^2 \# n \overline{\C P^2}$ with $m> n \geq 0$.
	\end{enumerate}
	In all cases the generic fiber is nontrivial in $H_2(X\setminus D;\R)$. In case \emph{(1)}, $D$ is coorientable, while in cases \emph{(2)} and \emph{(3)}, $D$ is coorientable if and only if $m$ is odd.
\end{thm}
This result can be combined with \autoref{thm:introsgcs} to equip the listed four-manifolds $X^4$ with coorientable $D$ with \sgcs{}s.

The contents of this chapter have appeared before in \cite{CavalcantiKlaasse17} and are joint with Gil Cavalcanti. We would like to thank Selman Akbulut, Stefan Behrens, Robert Gompf and Andr\'as Stipsicz for useful discussions during its preparation.
\subsection*{Organization of the chapter}
This chapter is built up as follows. In Section \ref{sec:boundarymaps} we define the normal Hessian of a map and introduce boundary maps. We define boundary fibrations and boundary Lefschetz fibrations and prove several normal form results for boundary maps (\autoref{prop:bdrymaplocalcoord}, \autoref{prop:fibratingpointwise}, and \autoref{prop:fibratingsemiglobal}). In Section \ref{sec:constrbdrlfs} we introduce a standard boundary fibration (\autoref{prop:standardbdryfibration}) and use this to obtain boundary Lefschetz fibrations out of genus one Lefschetz fibrations over a punctured surface using monodromy data (\autoref{cor:lfcompletion}) via a completion process. In Section \ref{sec:constrsgcs} we prove our main results: \autoref{thm:lfibrationsgcs}, that four-dimensional boundary Lefschetz fibrations give rise to \sgcs{}s; and \autoref{thm:fibrationsgcs}, that the same holds for boundary fibrations. Moreover, we show compatibility with existing fibrations over $T^2$ and $S^1$ can be achieved in the compact case (\autoref{cor:compatiblefibrations}). In Section \ref{sec:applications} we give examples of stable generalized complex manifolds constructed using our methods. In particular, we recover the examples $m \C P^2 \# n \overline{\C P}^2$ for $m$ odd from \cite{CavalcantiGualtieri09}.
\section{Boundary maps and \placeholder{}s}
\label{sec:boundarymaps}
In this section we introduce the notion of a boundary map, which is a map degenerating suitably on a submanifold. When this submanifold has codimension two these will supply us with morphisms from elliptic to log divisors, and hence with Lie algebroid morphisms from the respective elliptic to log-tangent bundles by \autoref{prop:elltologmorphism}. After this we define the notion of a \placeholder, which is a Lefschetz-type fibration that can be interpreted as a Lie algebroid \lf between these Lie algebroids.
\subsection{The normal Hessian}
Let $(X,D)$ be a \emph{pair}, i.e.\ $X$ is a manifold and $D \subseteq X$ is a submanifold.
\begin{defn} A \emph{map of pairs} $f\colon (X,D) \to (Y,Z)$ is a map $f\colon X \to Y$ such that $f(D) \subseteq Z$. A \emph{strong map of pairs} is a map of pairs $f\colon (X,D) \to (Y,Z)$ such that $f^{-1}(Z) = D$.
\end{defn}
Given a map of pairs $f\colon (X,D) \to (Y,Z)$, note that $Tf(TD) \subseteq TZ$, so that $Tf$ induces a map $\nu(Tf)\colon ND \to NZ$ between normal bundles. When $f$ is a strong map of pairs and $f$ is transverse to $Z$, the map $\nu(Tf$) is an isomorphism.

Consider a map of pairs $f\colon (X,D) \to (Y,Z)$ such that $\nu(Tf)\colon ND \to NZ$ is the zero map. Equivalently, one can assume that ${\rm im}\, Tf \subset TZ$. Let $z_1, \dots, z_\ell$ be local defining functions for $Z$ and consider their pullbacks $h_i := f^*(z_i)$ for $i = 1, \dots, \ell$. As $f(D) \subseteq Z$, the functions $h_i$ vanish on $D$. Moreover, because $\nu(Tf)$ is the zero map, the derivatives $dh_i$ vanish on $D$ as well. Consequently, we can consider their Hessians $H(h_i)\colon {\rm Sym}^2(TX) \to \R$, which descend to maps $H(h_i)\colon {\rm Sym}^2(ND) \to \R$. As the differentials $dz_i$ span $N^*Z$, these combine to give a map $H^\nu(f)\colon {\rm Sym}^2(ND) \to f^*(NZ)$, which one checks to be invariantly defined.
\begin{defn} Let $f\colon (X,D) \to (Y,Z)$ such that ${\rm im}\, Tf \subset TZ$. The \emph{normal Hessian} of $f$ along $D$ is the map $H^\nu(f): {\rm Sym}^2(ND) \to f^*(NZ)$ over $D$.
\end{defn}
When ${\rm codim}(Z) = 1$, the normal Hessian $H^\nu(f)$ can be viewed as the matrix of second partial derivatives of the coordinate function of $f$ normal to $Z$ in directions normal to $D$.
\subsection{Boundary maps}
Let $f\colon (X,D) \to (Y,Z)$ be a strong map of pairs with $Z$ a hypersurface and such that ${\rm codim}(D) \geq 2$. Then $f$ cannot be transverse to $Z$, as then $f^{-1}(Z) = D$ would be of codimension one. As $Z$ is of codimension one, $f$ being transverse to $Z$ is equivalent to $\nu(Tf)\colon ND \to NZ$ being nonzero, hence in this case $\nu(Tf)\colon ND \to NZ$ is the zero map. Equivalently we have ${\rm im}\, Tf \subset TZ$, so that the normal Hessian of $f$ along $D$ is well-defined.
\begin{defn} Let $f\colon (X,D) \to (Y,Z)$ be a strong map of pairs, $Z$ a hypersurface, and ${\rm codim}(D) \geq 2$. Then $f$ is a \emph{boundary map} if its normal Hessian $H^\nu(f)$ is definite along $D$.
\end{defn}
As $Z$ is a hypersurface, $NZ$ is one-dimensional. Because of this, $H^\nu(f)$ being definite makes sense, as locally it is a map $H^\nu(f)\colon {\rm Sym}^2(\R^d) \to \R$ where $d = {\rm codim}(D)$. The choice of the name will become clearer after establishing some properties (see \autoref{rem:bdrymapname}). We will also call $f$ a \emph{codimension-$k$ boundary map} if ${\rm codim}(D) = k$, and sometimes implicitly assume that ${\rm codim}(D) = 2$ when $\dim(X) = 4$. Indeed, the main reason for introducing the notion of a boundary map comes from \autoref{prop:bdrymapmorphdivisors} below, where ${\rm codim}(D) = 2$.
\begin{rem} Let $f\colon (X,D) \to (Y,Z)$ be a codimension-$k$ boundary map and $f''\colon (X',D') \to (X,D)$ and $f'\colon (Y,Z) \to (Y',Z')$ strong maps of pairs with $f''$ transverse to $D$ and $f'$ transverse to $Z'$. Then $f'' \circ f \circ f'\colon (X',D') \to (Y',Z')$ is a codimension-$k$ boundary map.
\end{rem}
\begin{rem} Assuming that ${\rm codim}(D) \geq 2$ is only done to ensure that $\nu(Tf)$ is the zero map, as is required for the definition of the normal Hessian. If ${\rm codim}(D) = 1$ yet this condition holds, it makes sense to talk about codimension-one boundary maps.
\end{rem}
Codimension-two boundary maps naturally give rise to morphisms of divisors between elliptic and log divisors respectively, and hence to Lie algebroid morphisms by \autoref{prop:elltologmorphism}.
\begin{prop}\label{prop:bdrymapmorphdivisors} Let $f\colon (X,D) \to (Y,Z)$ be a map of pairs with $Z$ a hypersurface and ${\rm codim}(D) = 2$. Then $f$ is a boundary map if and only if $I_D := f^* I_Z$ is an elliptic ideal and $f$ a morphism of divisors.
\end{prop}
In other words, a codimension-two boundary map uniquely specifies a compatible elliptic divisor structure on $D$. 
\bp Assume that $f$ is a boundary map and consider $I_D = f^* I_Z$. Let $z$ be a local defining function for $Z$ so that locally $I_Z = \langle z \rangle$, and hence $I_D = \langle f^*(z) \rangle$. As $f$ is a boundary map, $H^\nu(f)$ is definite, so that $f^*z$ specifies the germ of a definite Morse--Bott function around $D$. By the discussion above \autoref{prop:ellmorsebott} we see that $I_D$ is an elliptic ideal specifying an elliptic divisor structure on $D$ by \autoref{prop:locprincideal}. The map $f$ is a morphism of divisor by construction. Alternatively, one shows that for $(L,s)$ the log divisor determined by $Z$, the pair $(R,q) := (f^* L, f^* s)$ is an elliptic divisor with $D_q = D$. The converse is similar, using again that $f$ has definite normal Hessian if and only if $f^*(z)$ is locally Morse--Bott of even index around $D$, where $z$ is a local defining function for $Z$.
\ep
\begin{cor}\label{cor:bdrymaplamorphism} Let $f\colon (X,D) \to (Y,Z)$ be a codimension-two boundary map. Then $Tf$ induces a Lie algebroid morphism $(\varphi,f)\colon TX(-\log |D|) \to TY(- \log Z)$ for the divisor structures of \autoref{prop:bdrymapmorphdivisors}.
\end{cor}
The pointwise conclusion of the Morse--Bott lemma, \autoref{lem:morsebott}, provides a local form for boundary maps around points in $D$.
\begin{prop}\label{prop:bdrymaplocalcoord} Let $f\colon (X^n,D^{k}) \to (Y^m,Z^{m-1})$ be a boundary map and $x \in D$. Then around $x$ and $f(x) \in Z$ there exist coordinates $(x_1,\dots,x_n)$ and $(z,y_2,\dots,y_m)$ for which ${\{x_1 = \dots = x_{n-k} = 0\} = D}$ and $\{z = 0\} = Z$ such that for some $g\colon \R^n \to \R^{k-1}$ we have $f(x_1,\dots,x_n) = (x_1^2 + \dots + x_{n-k}^2, g(x_1,\dots,x_n))$.
\end{prop}
\bp Let $x \in D$ and $z$ be a local defining function for $Z$ around $f(x) \in Z$. As $f$ is a boundary map, the proof of \autoref{prop:bdrymapmorphdivisors} shows that $f^*(z)$ is a local Morse--Bott function of index zero around $x$, after possibly replacing $z$ by $-z$. By \autoref{lem:morsebott}, after trivializing $ND$, there exist coordinates $(x_1,\dots,x_n)$ of $X$ around $x$ with $\{x_1 = \dots = x_{n-k} = 0\} = D$ such that $f^*(z)(x_1,\dots,x_n) = x_1^2 + \dots + x_{n-k}^2$. Complete $z$ to a coordinate system $(z,y_2,\dots,y_m)$ of $Y$ around $f(x)$. Then $f(x_1,\dots,x_n) = (f^*(z)(x_1,\dots,x_n),g(x_1,\dots,x_n)) = (x_1^2 + \dots + x_{n-k}^2, g(x_1,\dots,x_n))$ for some $g\colon \R^n \to \R^{k-1}$.
\ep
Using either \autoref{prop:bdrymaplocalcoord} or the proof \autoref{prop:bdrymapmorphdivisors}, the Lie algebroid morphisms $(\varphi,f)\colon TX(-\log|D|) \to TY(-\log Z)$ obtained from boundary maps $f\colon (X,D) \to (Y,Z)$ using \autoref{cor:bdrymaplamorphism} have the following extra property: for any $k$-form $\alpha \in \Omega^k(\log Z)$, we have ${\rm Res}_\theta(\varphi^* \alpha) = 0$, because ${\rm Res}_\theta f^*(z) = 0$ for any local defining function $z$ for $Z$.
\subsection{Fibrating boundary maps}
We next introduce specific boundary maps by demanding submersiveness on $D$.
\begin{defn} A \emph{fibrating boundary map} is a boundary map $f\colon (X,D) \to (Y,Z)$ such that $f|_{D}\colon D \to Z$ is submersive.
\end{defn}
Note that it is not required that $D$ surjects onto $Z$. For fibrating boundary maps we can improve upon \autoref{prop:bdrymaplocalcoord}, obtaining a local form around points in $D$.
\begin{prop}\label{prop:fibratingpointwise} Let $f\colon (X^n,D^{k}) \to (Y^{m},Z^{m-1})$ be a fibrating boundary map. Then around points $x \in D$ and $f(x) \in Z$ there exist coordinates $(x_1,\dots,x_n)$ and $(z,y_2,\dots,y_m)$ with $\{x_1 = \dots = x_{n-k} = 0\} = D$ and $\{z = 0\} = Z$ such that $f(x_1,\dots,x_n) = (x_1^2 + \dots + x_{n-k}^2, x_{n-m+1}, \dots, x_{n})$.
\end{prop}
Hence we can simultaneously put both the components of $f$ normal and tangent to $D$ in standard form, obtaining a commuting diagram near $x$ and $f(x)$:
\begin{center}
	\begin{tikzpicture}
	\matrix (m) [matrix of math nodes, row sep=2.5em, column sep=2.5em,text height=1.5ex, text depth=0.25ex]
	{	 ND & NZ \\ D & Z\\};
	\path[-stealth]
	(m-1-1) edge node [above] {$f$} (m-1-2)
	(m-1-1) edge node [left] {${\rm pr}_D$} (m-2-1)
	(m-1-2) edge node [right] {${\rm pr}_Z$} (m-2-2)
	(m-2-1) edge node [above] {$f|_D$} (m-2-2);
	\end{tikzpicture}
\end{center}
Call a finite collection of functions an \emph{independent set} at $p$ if their differentials are everywhere linearly at $p$. By the implicit function theorem, an independent set can be completed to a coordinate system in a neighbouhood of $p$. Independence is preserved under pulling back along a submersion.
\bp Choose a tubular neighbourhood embedding $\Phi\colon NZ \to V \subset Y$ where ${\rm pr}_Z\colon NZ \to Z$ is the projection and let $z$ be a local defining function for $Z$ on an open subset $V' \subset V$ of $f(x)$. Let $U := f^{-1}(V') \subset X$. Choose coordinate functions $y_2,\dots,y_m\colon V' \to \R$ for $Z$. Then $\{y_2,\dots,y_m\}$ forms an independent set on $Z$, and because ${\rm pr}_Z$ is a submersion, the same is true for $\{z, {\rm pr}_Z^*(y_2), \dots, {\rm pr}_Z^*(y_m)\}$ on $Y$. By \autoref{prop:bdrymaplocalcoord}, after possibly shrinking $U$ and $V'$ and changing $z$ to $-z$, there exist tubular neighbourhood coordinates $(x_1,\dots,x_n)$ around $x$ such $f^*(z) = x_1^2 + \dots + x_{n-k}^2$. Consider $\{x_1,\dots,x_{n-k},f^* {\rm pr}_Z^*(y_j)\}$, which is an independent set on $X$, using submersiveness of $f|_D$. Complete this to a coordinate system $\{x_1,\dots,x_{n-m+1},f^* {\rm pr}_Z^*(y_j)\}$ on $X$, and relabel these as $(x_1,\dots,x_n)$. Using these coordinates on $X$ and the coordinates $(z,\pi_Z^*(y_j))$ on $Y$, the map $f$ is given by $f(x_1,\dots,x_n) = (x_1^2 + \dots + x_{n-k}^2, x_{n-m+1}, \dots, x_{n})$ as desired.
\ep
The normal form result of \autoref{prop:fibratingpointwise} immediately implies the following.
\begin{cor}\label{cor:fibratingsubm} Let $f\colon (X,D) \to (Y,Z)$ be a fibrating boundary map. Then $f$ is submersive in a punctured neighbourhood around $D$.
\end{cor}
Consequently, fibrating boundary maps have well-defined fibers of dimension $\dim(X) - \dim(Y)$ near $D$, and of dimension $\dim(X) - \dim(Y) - {\rm codim}(D) + 1$ on $D$. In particular, when $f$ is a fibrating boundary map, $D$ will be a fiber bundle over certain components of $Z$.
\begin{rem}\label{rem:fibrliealgdmorph} An alternative way of viewing the proof of \autoref{cor:fibratingsubm} when ${\rm codim}(D) = 2$ uses \autoref{prop:bdrymapmorphdivisors}. Namely, let $f\colon (X,D) \to (Y,Z)$ give rise to a Lie algebroid morphism $(\varphi,f)\colon TX(-\log |D|) \to TY(-\log Z)$. While $Tf\colon ND \to NZ$ is the zero map, $\varphi|_D$ is surjective because $f$ is fibrating. This is an open condition, so that $\varphi$ is surjective in a neighbourhood around $D$. On $X \setminus D$, the Lie algebroid $TX(-\log |D|)$ is isomorphic to $TX$, so $f$ is submersive there. Thus fibrating boundary maps give Lie algebroid morphisms which are fibrations near $D$.
\end{rem}
We proceed to obtain similar normal bundle commutativity around components of $D$ whose image is coorientable. To prove this, we use a result by Bursztyn--Lima--Meinrenken \cite{BursztynLimaMeinrenken16} on normal bundle embeddings, which we now describe. Let $M \subseteq X$ be a submanifold and $NM$ its normal bundle. Denote by $\mc{E}_M \in \mf{X}(NM)$ the associated Euler vector field and given $V \in \mf{X}(X)$ tangent to $M$, let $\nu(X) \in \mf{X}(NM)$ be its linear approximation.
\begin{defn} Let $M \subseteq X$ be a submanifold. A vector field $V \in \mf{X}(X)$ is \emph{Euler-like} along $M$ if it is complete, and satisfies $V|_M = 0$ and $\nu(V) = \mc{E}_M$.
\end{defn}
A \emph{strong tubular neighbourhood embedding} for a submanifold $M \subseteq X$ is an embedding $\Phi\colon NM \to X$ taking the zero section of $NM$ to $M$, and such that the linear approximation $\nu(\Phi)\colon (NM,M) \to (X,M)$ is the identity map. Euler-like vector fields give rise to strong tubular neighbourhood embeddings.
\begin{prop}[{\cite[Proposition 2.6]{BursztynLimaMeinrenken16}}]\label{prop:strongtubnhood} Let $M \subseteq X$ be a submanifold and $V \in \mf{X}(X)$ Euler-like along $M$. Then there exists a unique strong tubular neighbourhood embedding $\Phi\colon NM \to X$ such that $\Phi_*(\mc{E}_M) = V$.
\end{prop}
We use this result to construct compatible tubular neighbourhood embeddings for fibrating boundary maps around components of $D$ whose image is coorientable.
\begin{prop}\label{prop:fibratingsemiglobal} Let $f\colon (X,D) \to (Y,Z)$ be a fibrating boundary map and $D_j \subseteq D$ a connected component such that $f(D_j) =: Z_j \subset Z$ is coorientable. Then there exist a defining function $z$ for $Z_j$ and tubular neighbourhood embeddings $\Phi_{D_j}\colon \wt{U} \to U \subset X$ for $D_j$ and $\Phi_{Z_j}\colon\wt{V} \to V \subset Y$ for $Z_j$ such that $\Phi_{D_j}^* f^*(z) = Q_{f^*(z)} \in \Gamma(D_j; {\rm Sym}^2 N^* D_j)$ and the following diagram commutes:
	\begin{center}
		\begin{tikzpicture}
		\matrix (m) [matrix of math nodes, row sep=2.5em, column sep=2.5em,text height=1.5ex, text depth=0.25ex]
		{	 \wt{U} & U & V & \wt{V}\\ D_j & D_j & Z_j & Z_j\\};
		\path[-stealth]
		(m-1-1) edge node [above] {$\Phi_{D_j}$} (m-1-2)
		(m-1-1) edge node [left] {${\rm pr}_{D_j}$} (m-2-1)
		(m-2-1) edge node [above] {${\rm id}$} (m-2-2)
		(m-1-2) edge node [above] {$f$} (m-1-3)
		(m-2-2) edge node [above] {$f|_D$} (m-2-3)
		(m-1-4) edge node [above] {$\Phi_{Z_j}$} (m-1-3)
		(m-2-4) edge node [above] {${\rm id}$} (m-2-3)
		(m-1-4) edge node [right] {${\rm pr}_{Z_j}$} (m-2-4);
		
		\draw[right hook->]
		(m-2-2) edge (m-1-2)
		(m-2-3) edge (m-1-3);
		\end{tikzpicture}
	\end{center}
\end{prop}
When $Z$ is separating this result implies a normal form for $f$ around any point in $Z$ and its entire inverse image, as then a global defining function for $Z$ can be used.
\bp Let $z\colon V \to \R$ be a defining function for $Z_j$ and let $U' \subset f^{-1}(V)$ be the connected component containing $D_j$. Using \autoref{lem:morsebott} applied to $f^*(z)$, shrink $U'$ so that $U' = \Phi_{D_j}(\wt{U}')$ for some tubular neighbourhood embedding $\Phi_{D_j}\colon \wt{U}' \to U$ of $D_j$. Choose a tubular neighbourhood embedding $\Phi_{Z_j}\colon \wt{V} \to V$ for $Z_j$. For $x \in D_j$, use the proof of \autoref{prop:fibratingpointwise} (possibly changing $z$ to $-z$) to obtain an open $U_x \subset U'$ containing $x$ and coordinates $(x_1,\dots,x_n)$ so that $f^*(z) = x_1^2 + \dots + x_{n-k}^2$. Note that $U'$ is connected so that $f^*(z)$ has a fixed sign. Consider $U := \cup_{x \in D_j} U_x \subset U'$ and extract a finite subcover $\{U_\alpha\}_{\alpha \in I}$. On each set $U_\alpha$, define the vector field $\mc{E}_\alpha := x_1 \partial_{x_1} + \dots + x_{n-k} \partial_{x_{n-k}}$. It satisfies $\mc{L}_{\mc{E}_\alpha} f^*(z) = 2 f^*(z)$, and $f_* \Phi_{\alpha *} \mc{E}_\alpha = z \partial_z = \mc{E}_{Z_j}$. Let $\{\psi_\alpha\}_{\alpha \in I}$ be a partition of unity subordinate to $\{U_\alpha\}_{\alpha \in I}$ and define $\mc{E} := \sum_{\alpha \in I} \psi_{\alpha} \mc{E}_\alpha$ on $U'$. Then $\mc{E}|_{D_j} = \sum_{\alpha \in I} \psi_\alpha \mc{E}_\alpha|_{D_j} = \sum_{\alpha \in I} \psi_\alpha \cdot 0 = 0$, and $\nu(\mc{E}) = \sum_{\alpha \in I} \psi_\alpha \nu(\mc{E}_\alpha) = \sum_{\alpha \in I} \psi_{\alpha} \mc{E}_{D_j} = \mc{E}_{D_j}$. Ensure that $\mc{E}$ is complete by multiplying by a bump function and shrinking $U'$, so that $\mc{E} \in \mf{X}(U)$ is Euler-like along $D_j$. Moreover, $\mc{L}_{\mc{E}} f^*(z) = \sum_{\alpha \in I} \mc{L}_{\psi_\alpha \mc{E}_\alpha} f^*(z) = \sum_{\alpha \in I} \psi_\alpha \mc{L}_{\mc{E}_\alpha} f^*(z) = \sum_{\alpha \in I} \psi_\alpha 2 f^*(z) = 2 f^*(z)$. Finally, $f_* \mc{E} = \sum_{\alpha \in I} f_* (\psi_\alpha \mc{E}_\alpha) = \sum_{\alpha \in I} \psi_\alpha \mc{E}_Z = \mc{E}_Z$. Use \autoref{prop:strongtubnhood} and possibly shrink $U$ to obtain a tubular neighbourhood embedding $\Phi_{D_j}'\colon \wt{U} \to U$ such that $\Phi_* \mc{E}_{D_j} = \mc{E}$, $f_* \mc{E} = \mc{E}_Z$, and $\mc{L}_{\mc{E}} f^*(z) = 2 f^*(z)$. This embedding satisfies all desired properties. Moreover, $\Phi^{'*}_{D_j} f^*(z) = Q_{f^*(z)} \in \Gamma(D_j; {\rm Sym}^2(N^*D_j))$ as $g :=\Phi^{'*}_{D_j} f^*(z)$ is smooth and satisfies $\mc{L}_{\mc{E}} g = 2 g$.
\ep
We can now obtain topological information about the generic fibers of $f$ near $D$. The following result is the main reason for wanting compatible tubular neighbourhood embeddings.
\begin{prop}\label{prop:fiberbundle} Let $f\colon (X^n,D^{k}) \to (Y^m,Z^{m-1})$ be a fibrating boundary map with $Z$ separating. Denote the fiber of $f|_D\colon D \to Z$ by $F_D^{k-m+1}$. Then the fiber $F^{n-m}$ of $f$ near $D$ is an $S^{n-k-1}$-sphere bundle over $F_D$.
\end{prop}
\bp Apply \autoref{prop:fibratingsemiglobal} to each connected component of $D$ using a single global defining function $z$ for $Z$. Let $y \in \wt{V} \setminus Z$. As $f$ is submersive on $\wt{U} \setminus D$, consider $F_y = \Phi_{D}^{-1} \circ f^{-1} \circ \Phi_{Z}(y)$, which is of dimension $n-m$. Consider ${\rm pr}_Z(y) \in Z$ and its $(k-m+1)$-dimensional fiber $F_{D,{\rm pr}_Z(y)} = f|_D^{-1}({\rm pr}_Z(y))$. Because the tubular neighbourhood embeddings are compatible with $f$, we have $F_y = {\rm pr}_D^{-1}(F_{D,{\rm pr}_Z(y)})$, noting that ${\rm pr}_D$ is submersive. As a point in $\wt{V} \subset NZ$ is given by a point in $Z$ together with a distance, the fiber of ${\rm pr}_D\colon F_y \to F_{D,{\rm pr}_Z(y)}$ is given by a sphere, consisting of all points with fixed distance above the corresponding point in $D$.
\ep
\begin{cor}\label{cor:fibertori} Let $f\colon (X^4,D^2) \to (Y^2,Z^1)$ be a codimension-two fibrating boundary map. Then the generic fibers of $f$ near $D$ are unions of tori.
\end{cor}
\bp This follows from \autoref{prop:fiberbundle} using $n = 4$ and $k = m = 2$ by a dimension count. The generic fibers $F$ of $f$ near $D$ satisfy ${\rm dim}(F) = 4 - 2 = 2$, and are $S^1$-bundles over fibers of $f|_D$, which in turn are one-dimensional. The only such two-manifolds are unions of tori.
\ep
Recall that when $Z$ is separating it admits a global defining function $z \in C^\infty(Y)$. In this case, any boundary map $f$ to $(Y,Z)$ must map onto one component of $Y$ with respect to $Z$. More precisely, given such a $z$, consider $Y_+ = z^{-1}([0,\infty)) \subset Y$, a manifold with boundary given by $Z = z^{-1}(0)$. Then we have the following.
\begin{prop}\label{prop:globalbdrymap} Let $f\colon (X,D) \to (Y,Z)$ be a boundary map, and $Z$ separating. Then there exists a global defining function $z$ for $Z$ so that $f(X) \subset Y_+$, and $f$ defines a boundary map $f\colon (X,D) \to (Y', Z')$, where $Y' = Y_+ \cap f(X)$ and $Z' = Z \cap f(D)$.
\end{prop}
\bp Let $z$ be a defining function for $Z$. Then $f^*(z)$ is globally defined on $X$, with zero set $D = f^{-1}(Z)$. As $D$ has codimension at least two in $X$, its complement $X \setminus D$ is connected, so that $f^*(z)$ has fixed sign on $X \setminus D$. Consequently, by changing $z$ to $-z$ if necessary, the function $f^*(z)$ is non-negative. But then for all points $x \in X$, $z(f(x)) = f^*(z)(x) \geq 0$, so that $f(X) \subset Z_+$. Moreover, $f$ is a boundary map when restricting its codomain to its image.
\ep
\begin{rem}\label{rem:bdrymapname} The previous proposition explains the name `boundary map', as the defining condition specifies the behavior of $f$ normal to $Z$, the boundary of its restricted codomain.
\end{rem}
\subsection{Boundary (Lefschetz) fibrations}
We introduce further submersiveness assumptions.
\begin{defn} A \emph{boundary fibration} $f\colon (X,D) \to (Y,Z)$ is a fibrating boundary map such that $f|_{X \setminus D}\colon X \setminus D \to Y \setminus Z$ is a fibration.
\end{defn}
\begin{defn} A \emph{\placeholder} $f\colon (X^{2n},D) \to (\Sigma^{2},Z)$ is a fibrating boundary map such that $f|_{X \setminus D}\colon X \setminus D \to \Sigma \setminus Z$ is a Lefschetz fibration.
\end{defn}
The following is immediate from \autoref{cor:bdrymaplamorphism} together with \autoref{rem:fibrliealgdmorph}.
\begin{cor}\label{cor:bdryliealgdmorph} Let $f\colon (X^{2n},D) \to (Y,Z)$ be a codimension-two boundary (Lefschetz) fibration. Then $f$ gives rise to a Lie algebroid (Lefschetz) fibration, i.e.\ ${(\varphi,f)\colon TX(-\log|D|) \to TY(-\log Z)}$.
\end{cor}
The statement that $f$ induces a Lie algebroid morphism is to be interpreted as in \autoref{prop:elltologmorphism}, namely that there is a Lie algebroid morphism $(\varphi,f)$ such that $\varphi = Tf$ on sections. Moreover, the elliptic divisor structure on $D$ is the one induced from $f$ and $Z$.
By adapting the usual argument, we can ensure that boundary (Lefschetz) fibrations have connected fibers.
\begin{prop}\label{prop:blfconnfibers} Let $f\colon (X,D) \to (Y,Z)$ be a codimension-two boundary (Lefschetz) fibration, with $Z$ separating. Define $Y' := f(X)$ and $Z' := f(D)$. Then there exists a cover $g\colon \wt{Y}' \to Y'$ of $Y'$ and a boundary (Lefschetz) fibration $\wt{f}\colon (X,D) \to (\wt{Y}',Z')$ with connected generic fibers such that $f = g \circ \wt{f}$.
\end{prop}
\bp As $f$ is a global boundary map, by \autoref{prop:globalbdrymap} we obtain a boundary map $f\colon (X,D) \to (Y', Z')$ which immediately is also a boundary (Lefschetz) fibration. By definition, $f\colon X \setminus D \to Y' \setminus Z'$ is a (Lefschetz) fibration. Consequently, denoting its generic fiber by $F$, there is a sequence in homotopy $\pi_1(F) \to \pi_1(X \setminus D) \to \pi_1(Y' \setminus Z') \to \pi_0(F) \to 0$ \cite[Proposition 8.1.9]{GompfStipsicz99}. Applying Van Kampen's theorem for each connected component of $D$ shows that $\pi_1(X \setminus D)$ surjects onto $\pi_1(X)$ via the natural inclusion. Because $Z'$ is the boundary of $Y'$, we have $\pi_1(Y' \setminus Z') \cong \pi_1(Y')$. We obtain the following commutative diagram.
\begin{center}
	\begin{tikzpicture}
	\matrix (m) [matrix of math nodes, row sep=2.5em, column sep=2.5em,text height=1.5ex, text depth=0.25ex]
	{	 \pi_1(F) & \pi_1(X\setminus D) & \pi_1(Y' \setminus Z') & \pi_0(F) & 0\\ & \pi_1(X) & \pi_1(Y') & & \\};
	\path[-stealth]
	(m-1-1) edge node [left] {} (m-1-2)
	(m-1-2) edge node [left] {} (m-1-3)
	(m-1-3) edge node [left] {} (m-1-4)
	(m-1-4) edge node [left] {} (m-1-5)
	(m-2-2) edge node [above] {} (m-2-3)
	(m-1-3) edge node [left] {$\cong$} (m-2-3);
	\draw[>=stealth,->>]
	(m-1-2) edge node [left] {} (m-2-2);
	\end{tikzpicture}
\end{center}
The generic fiber $F$ of $f\colon X \setminus D \to Y' \setminus Z'$ is compact, hence $\pi_0(F)$ is finite. But then $f_*(\pi_1(X \setminus D)) \subset \pi_1(Y' \setminus Z')$ is a subgroup of finite index. Using the isomorphism $\pi_1(Y' \setminus Z') \cong \pi_1(Y')$, let $G$ denote the corresponding finite index subgroup inside $\pi_1(Y')$, and let $g\colon \wt{Y}' \to Y'$ be the associated cover. Then $f\colon X \to Y'$ lifts to $\wt{f}\colon X \to \wt{Y}'$ if and only if $f_*(\pi_1(X)) \subset G$, but this is an equality because $\pi_1(X \setminus D) \to \pi_1(X)$ surjects. For $\wt{f}$ we have that $\wt{f}\colon X \setminus D \to \wt{Y}' \setminus Z'$ induces a surjection $\wt{f}_*\colon \pi_1(X \setminus D) \to \pi_1(\wt{Y}' \setminus Z')$. But then $\pi_0(\wt{F})$ is trivial, so that the generic fibers $\wt{F}$ of $\wt{f}$ are connected.
\ep
By the above we can replace a boundary Lefschetz fibration on $(X,D)$ with $Z$ separating by one for which the generic fibers in $X \setminus D$ are connected. We next study the fibrating case, concluding that the fibers of $f|_D\colon D \to Z$ are connected if those of $f$ near $D$ are.
\begin{prop} Let $f\colon (X,D) \to (Y,Z)$ be a fibrating boundary map whose generic fibers near $D$ are connected, and $Z$ separating. Then the fibers of $f|_D\colon D \to Z$ are connected.
\end{prop}
\bp Using \autoref{prop:globalbdrymap}, replace $(Y,Z)$ by $(Y',Z')$ and let $V \subset Y'$ be an open neighbourhood of a point $y \in Z'$ so that $f|_{f^{-1}(V) \setminus D}$ is a fibration. Choose a curve $\gamma\colon [0,1] \to Y'$ such that $\gamma([0,1)) \in V \setminus Z'$, $f(1) = y$ and $\gamma$ is transverse to $Z'$. Then $M := f^{-1}(\gamma([0,1])) \subset f^{-1}(V)$ is a compact submanifold of $X$. Let $F$ be the generic fiber of $f$ near $D$ and denote $F_y = f^{-1}(y) = f|_D^{-1}(y)$. Let $F_{y,i}$ be the connected components of $F_y$, and let $U_i$ be disjoint opens around $F_{y,i}$ in $M$. Choose a sequence of points $(y_n)_{n \in \N} \subseteq \gamma([0,1])$ converging to $y$ such that $y_n \neq y$ for all $n \in \N$. For each $n$, consider the set $V_n := F_{y_n} \setminus \bigsqcup_i (U_i \cap F_{y_n})$. As each $U_i$ is open in $M$, the sets $U_i \cap F_{y_n}$ are disjoint opens in $F_{y_n}$.

Assume that $F_y$ is not connected. Then for $n$ large enough, there will be at least two indices $i$ for which $U_i \cap F_{y_n} \neq \emptyset$. As $F_{y_n}$ is connected, we conclude that $V_n \neq \emptyset$, as $F_n$ cannot be covered by disjoint opens. For each such $n$, let $x_n \in V_n$ be some element. By compactness of $M$, the sequence $\{x_n\}_{n \in \N}$ has a convergent subsequence, so that after relabeling we have $x_n \to x$ for some $x \in M$. By definition $x \in M \setminus (\bigcup_i U_i)$, so that $x \not\in U_i$ for all $i$, hence $x \not\in F_y$. However, by continuity we have $f(x_n) = y_n \to y = f(x)$, which is a contradiction. We conclude that $F_y$ is connected, so that the fibers of $f|_D\colon D \to Z'$ are connected.
\ep
\section{Constructing \placeholder{}s}
\label{sec:constrbdrlfs}
In this section we discuss how to obtain four-dimensional \placeholder{}s. In particular, we construct them out of genus one Lefschetz fibrations $f\colon X^4 \to \Sigma^2$ by surgery, replacing neighbourhoods of fibers of points $x \in \Sigma$ by a certain standard boundary map.

Given $n \in \Z$, let $p\colon L_n \to T^2$ be the complex line bundle over $T^2$ with first Chern class equal to $n \in H^2(T^2;\Z)$. Choose a Hermitian metric on $L_n$, which provides a fiberwise radial distance $r\colon {\rm Tot}(L_n) \to \R_{\geq 0}$. Further, let ${\rm pr}\colon T^2 \to S^1$ be the projection onto the first factor.
\begin{defn} The total space ${\rm Tot}(L_n)$ with the map $f\colon{\rm Tot}(L_n) \to S^1 \times \R_{\geq 0}$ given by $f(x) := ({\rm pr}(p(x)), r^2(x))$ for $x \in {\rm Tot}(L_n)$, is called the \emph{standard $n$-model}.
\end{defn}
The first thing to note is that $f$ as defined above is a boundary fibration.
\begin{prop}\label{prop:standardbdryfibration} The map $f\colon ({\rm Tot}(L_n), T^2) \to (S^1 \times \R_{\geq 0}, S^1 \times \{0\})$ is a boundary fibration, where $T^2$ is identified with the zero section in ${\rm Tot}(L_n)$.
\end{prop}
\bp Certainly $f$ is a strong map of pairs, and the codimension of $D = T^2$ inside ${\rm Tot}(L_n)$ is (at least) two. The normal Hessian of $H^\nu(f)$ of $f$ along $D$ is given by the constant map equal to $2$, which is clearly definite. We conclude that $f$ is a boundary map. Both the bundle projection $p$ as well as the projection ${\rm pr}$ are submersive, so that $f|_D$ is submersive, making $f$ a fibrating boundary map. It is immediate that $f$ is submersive when $r \neq 0$, i.e.\ on ${\rm Tot}(L_n) \setminus T^2$, so that $f$ is a boundary fibration.
\ep
Let $f\colon X \to Y$ be a smooth map that is a fibration in a neighbourhood of an embedded circle $\gamma \subset Y$. Recall that the \emph{monodromy} of $f$ around $\gamma$ is the element in the mapping class group of the fiber $F$ of $f^{-1}(\gamma) \to \gamma$ as a mapping torus.

Returning to the standard $n$-model, let $\varepsilon > 0$ be small and consider the circle $\gamma = S^1 \times \{\varepsilon\} \subset S^1 \times \R_{\geq 0}$. We compute the monodromy of $f\colon {\rm Tot}(L_n) \to S^1 \times \R_{\geq 0}$ around $\gamma$.
\begin{prop} \label{prop:monodromy}The monodromy of $f$ around $\gamma$ is the $-n$th power of a Dehn twist.
\end{prop}
\bp Let $M := f^{-1}(\gamma)$ be the $T^2$-bundle $f\colon M \to \gamma$. To compute the monodromy of $f$ around $\gamma$, note that using $p\colon M \to T^2$ we can view $M$ as the principal $\varepsilon$-sphere bundle of $L_n$. The hermitian metric gives rise to a connection $i \theta \in \Omega^1(M;\R)$. As $L_n \to T^2$ has Euler class equal to $n$, its curvature equals $\frac{d\theta}{2 \pi i} = n da \wedge db$, where $da$ and $db$ are generators of $H^1(T^2)$. Recall now that circle bundles are classified by their Euler class due to the Gysin sequence. Consider $M' := \R^3 / \Gamma$, with $\Gamma = \langle \alpha_1,\alpha_2,\alpha_3\rangle$ the integral lattice generated by the following group elements acting on $\R^3$: 
\be
	\alpha_1 = (x,y+1,z), \qquad \alpha_2 = (x,y,z+1), \qquad \alpha_3 = (x+1,y, z-n y).
\ee
The projection ${\rm pr}_{12}\colon \R^3 \to \R^2$ given by $(x,y,z) \mapsto (x,y)$ descends to a map $g\colon \R^3 / \Gamma \to \R^2 / \Z^2$. This is an $S^1$-bundle over $T^2$, with invariant one-forms $e_1 = dx$ and $e_2 = dy$, and connection one-form $e_3 = dz + n x dy$. As $d e_3 = n e_1 \wedge e_2$, we conclude that $\R^3 / \Gamma$ has Euler class equal to $n$, so that $\R^3 / \Gamma \cong M$, as circle bundles over $T^2$. Note that the map ${\rm pr}_1\colon \R^3 \to \R$ given by $(x,y,z) \mapsto x$ also descends and exhibits $M'$ as a $T^2$-bundle over $\R / \Z \cong S^1$. To obtain a $T^2$-bundle over $S^1$ out of $p\colon M \to T^2$, the choice of any projection $T^2 \to S^1$ gives isomorphic bundles. Hence we can assume that the projection ${\rm pr}_1$ makes the following diagram commute.
\begin{center}
	\begin{tikzpicture}
	\matrix (m) [matrix of math nodes, row sep=2.5em, column sep=2.5em,text height=1.5ex, text depth=0.25ex]
	{	M & T^2 & S^1 \\ M' & T^2 & S^1 \\};
	\path[-stealth]
	(m-1-1) edge node [above] {$p$} (m-1-2)
	(m-1-2) edge node [above] {${\rm pr}$} (m-1-3)
	(m-2-1) edge node [above] {${\rm pr}_{12}$} (m-2-2)
	(m-2-2) edge node [above] {${\rm pr}_{1}$} (m-2-3)
	(m-1-1) edge node [left] {$\cong$} (m-2-1)
	(m-1-2) edge node [left] {$\cong$} (m-2-2)
	(m-1-3) edge node [left] {$\cong$} (m-2-3);
	\end{tikzpicture}
\end{center}
We can compute the monodromy of $f$ around $\gamma$ by considering ${\rm pr}_1\colon M' \to S^1$. By our concrete description it is immediate that $M' = T^2 \times S^1 / (x,0) \sim (f(x),1)$, with
\be
	f \in {\rm Aut}(T^2) \text{ and } f_* = \begin{pmatrix} 1 & -n \\ 0 & 1 \end{pmatrix}\colon H_2(T^2) \to H_2(T^2).
\ee
We conclude that $M$ has monodromy equal to the $-n$th power of a Dehn twist.
\ep
We next describe the surgery process whereby we replace suitable neighbourhoods of points by the above standard models.
\begin{defn} A \emph{punctured surface} is an open surface $\Sigma$ obtained from a closed surface by removing a finite number of disks.
\end{defn}
A punctured surface naturally has a compact closure $\overline{\Sigma}$ by adding the circle boundaries of the removed disks, which has a natural isomorphism $j\colon \Sigma \stackrel{\cong}{\to} \overline{\Sigma} \setminus \partial \overline{\Sigma}$.
\begin{prop}\label{prop:completion} Let $f\colon X^4 \to \Sigma^2$ a proper map over a punctured surface such that $f$ is a genus one fibration in a neighbourhood of $\partial \Sigma$. Assume that the monodromy of $f$ around every boundary component is a power of a Dehn twist. Then there exists a compact elliptic pair $(\wt{X},|D|)$ with $i\colon X \stackrel{\cong}{\to} \wt{X} \setminus D$ and a fibrating boundary map $\wt{f}\colon (\wt{X},|D|) \to (\overline{\Sigma}, \partial \overline{\Sigma})$ such that the following diagram commutes:
	\begin{center}
		\begin{tikzpicture}
		\matrix (m) [matrix of math nodes, row sep=2.5em, column sep=2.5em,text height=1.5ex, text depth=0.25ex]
		{	X & (\wt{X},|D|) \\ \Sigma & (\overline{\Sigma}, \partial \overline{\Sigma}) \\};
		\path[-stealth]
		(m-1-1) edge node [left] {$f$} (m-2-1)
		(m-1-2) edge node [right] {$\wt{f}$} (m-2-2);
		
		\draw [right hook-latex]
		(m-1-1) edge node [above] {$i$} (m-1-2)
		(m-2-1) edge node [above] {$j$} (m-2-2);
		\end{tikzpicture}
	\end{center}
	In the above situation, we say $f\colon X \to \Sigma$ \emph{can be completed} to $\wt{f}\colon (\wt{X},|D|) \to (\overline{\Sigma}, \partial \overline{\Sigma})$. These completions are unique when the monodromy is not trivial.
\end{prop}
\bp At each boundary end $E_i$ of $\Sigma$, glue in the standard $n$-model ${\rm Tot}(L_{n_i})$, where $n_i \in \Z$ is such that ${\rm Mon}(\gamma_i) = \delta^{n_i}$, with $\gamma_i = \partial E_i$. Here $\delta$ denotes a Dehn twist. As the monodromies agree, the fibrations are isomorphic, hence can be glued together to a new fibration.

It only remains to argue that if the monodromy is nontrivial, the completion is unique. Since the completion of each end of $X$ is obtained by gluing  $L_n$ to the end of $X$, by identifying $M$, the  $\varepsilon$-sphere bundle of $L_n$, with the pre-image of a loop around a boundary component of $\Sigma$, we see that any  other completion can be obtained by precomposing the gluing map by a diffeomorphism of $M$. However, due to a result of Waldhausen \cite{Waldhausen67}  any such diffeomorphism  is isotopic to a fiber preserving diffeomorphism (covering a diffeomorphism of the base) and therefore extends to the $\varepsilon$-disc bundle of $L_n$. Hence the completion of the end is unique up to diffeomorphism.
\ep
\begin{cor}\label{cor:lfcompletion} Let $f\colon X^4 \to \Sigma^2$ be a genus one \lf over a punctured surface. Assume that the monodromy of $f$ around every boundary component is a power of a Dehn twist. Then $f$ can be completed to a codimension-two \placeholder $\wt{f}\colon (\wt{X},|D|) \to (\overline{\Sigma}, \partial \overline{\Sigma})$.
\end{cor}
\bp The Lefschetz singularities lie in the interior of $\Sigma$ so that $f$ is a genus one fibration near $\partial \Sigma$. After completing using \autoref{prop:completion}, $\wt{f}$ is a fibrating boundary map and a genus one Lefschetz fibration in the interior, so that it is a \placeholder.
\ep
Consequently, given a genus one \lf on $X$, removing several discs and their inverse image and then completing as in \autoref{cor:lfcompletion} gives a \placeholder on $\wt{X}$. Homological essentialness of generic fibers is preserved by the completion process.
\begin{prop}\label{prop:retainhomess} Let $f\colon X^4 \to \Sigma^2$ be a genus one \lf with boundary Lefschetz completion $\wt{f}\colon (\wt{X},|D|) \to (\overline{\Sigma}, \partial \overline{\Sigma})$. Then $f$ is homologically essential if and only if $\wt{f}$ is.
\end{prop}
A codimension-two boundary \lf on a four-manifold $X$ has singular fibers equal to the Euler characteristic of $X$, as is true for genus one Lefschetz fibrations.
\begin{prop}\label{prop:blfsingfibers} Let $f\colon (X^4, D) \to (\Sigma^2, Z)$ be a codimension-two boundary Lefschetz fibration. Then we have $\chi(X) = \mu$, with $\mu$ the number of singular fibers of $f|_{X\setminus D}$.
\end{prop}
\bp Given two open sets $U, V \subseteq X$ we have $\chi(U \cup V) = \chi(U) + \chi(V) - \chi(U \cap V)$. Let $U := X \setminus D$ and take $V \subseteq ND$ to be an open subset in a tubular neighbourhood of $D$. Then $V$ is homotopy equivalent to $D$, which is a union of tori, so that $\chi(V) = \chi(D) = 0$. Moreover, $U \cap V$ is deformation equivalent to $S^1 ND$, which is a principal circle bundle over $D$, hence $\chi(U \cap V) = 0$ as well. Recall that if a manifold $M$ admits a genus $g$ Lefschetz fibration, by \autoref{prop:lfeulerchar} its Euler characteristic satisfies $\chi(M) = 2 (2-2g) + \mu$, where $\mu$ is the number of singular fibers. Applying this to $U = X \setminus D$ with $g = 1$ we obtain that $\chi(X) = \chi(X \setminus D) = \mu$.
\ep
\section{Constructing \sgcs{}s}
\label{sec:constrsgcs}
In this section we combine the results from previous sections to construct \sgcs{}s using codimension-two boundary fibrations with two-dimensional fibers, and from codimension-two boundary \lf{}s in dimension four. For notational convenience, we call a boundary map $f\colon (X^{2n},D) \to (Y^{2n-2},Z)$ \emph{homologically essential} if the (generic) two-dimensional fibers $F$ (near $D$) specify nonzero homology classes ${[F] \neq 0 \in H_2(X\setminus D; \R)}$. We prove the following results stated in the introduction.
\begin{thm}\label{thm:lfibrationsgcs} Let $f\colon (X^4,D^2) \to (\Sigma^2,Z^1)$ be a homologically essential \placeholder between compact oriented manifolds with $D$ coorientable and $[Z] = 0 \in H_1(\Sigma;\Z_2)$. Then $(X,D)$ admits a \sgcs.
\end{thm}
\begin{thm}\label{thm:fibrationsgcs} Let $f\colon (X^{2n},D^{2n-2}) \to (Y^{2n-2},Z^{2n-3})$ be a homologically essential boundary fibration between compact oriented manifolds with $D$ coorientable such that $(Y,Z)$ admits a log-symplectic structure. Then $(X,D)$ admits a \sgcs.
\end{thm}
In the remainder of this section, any mention of an elliptic divisor structure on $(X,D)$, or the Lie algebroid $TX(-\log |D|)$, will refer to the structure induced by a codimension-two boundary map via \autoref{prop:bdrymapmorphdivisors}. Similarly, codimension-two boundary maps induce Lie algebroid morphisms between the appropriate elliptic and log-tangent bundles using \autoref{cor:bdryliealgdmorph} without further mention. Before proving the above two results, we first establish that we can create a suitable closed Lie algebroid two-form on $X$. 
\begin{prop}\label{prop:homessfibrbdrymap} Let $f\colon (X^{2n},D^{2n-2}) \to (Y^{2n-2},Z^{2n-3})$ be a homologically essential fibrating boundary map with $f(D)$ coorientable. Then there exists a closed two-form $\eta \in \Omega^2(\log |D|)$ with ${\rm Res}_q(\eta) = 0$ with $\eta|_{\ker \varphi}$ nondegenerate near $D$.
\end{prop}
\bp Note that $f$ is submersive with fibers $F$ around $D$ by \autoref{cor:fibratingsubm}. As $[F] \neq 0$, there exists $c \in H^2(X \setminus D)$ such that $\langle c, [F] \rangle = 1$. We construct local Lie algebroid two-forms around $D$ as in \autoref{prop:gluingfiberform}. By \autoref{thm:elllacohomology}, $H^2(X\setminus D)$ includes into $H_0^2(\log |D|)$. Let $\xi \in \Omega^2(\log |D|)$ be a closed two-form satisfying $[\xi] = c$ and ${\rm Res}_q(\xi) = 0$. By \autoref{prop:fibratingsemiglobal} there exist open neighbourhoods $U$ and $V$ around $D$ and $f(D)$ and a defining function for $f(D)$ on $V$ such that $f$ and $f|_D$ commute with the tubular neighbourhood projections ${\rm pr}_D$ and ${\rm pr}_Z$, and such that $f^*(z) = Q_{f^*(z)}$. Let $y \in f(D)$ and let $V' \subseteq V$ be a smaller open disk containing $y$ on which $NZ$ is trivial, and similarly so that $U' := f^{-1}(U)$. Then as in \autoref{prop:fibratingpointwise} we have $f(r e^{i\theta}, x, y) = (x, r^2)$. Moreover, in these coordinates we have $\ker \varphi = \langle \partial_\theta, \partial_y \rangle$. Let $\{V_i\}$ be a finite covering extracted from such open sets and let $U_i = f^{-1}(V_i)$, which together cover a neighbourhood of $D$. Set $U_0 = X \setminus D$. As $f\colon U_i \setminus D \to V_i \setminus Z$ is a $T^2$-bundle it is necessarily trivial, so that $H^2(U_i \setminus D) = H^2(T^2)$, with $H^1(T^2)$ generated by $\theta_i$ and $\gamma_i$ say.

Define $\eta_0 = \xi|_{X\setminus D}$ and $\eta_i \in \Omega^2(U_i;\log |D|)$ for $i \geq 1$ via $\eta_i = \lambda_i \theta_i \wedge \gamma_i$, where $\lambda_i \in \R$ is chosen such that $\int_F \eta_i = \int_F c$. Then $[\eta_i] = c|_{U_i} = [\xi|_{U_i}] \in H^2(U_i; \log |D|)$, so that there exist one-forms $\alpha_i$ such that $\eta_i = \xi + d \alpha_i$. As in \autoref{prop:gluingfiberform}, define a closed Lie algebroid two-form $\eta:= \xi + d (\sum_i (\psi_i\circ f) \alpha_{y_i})$ using a partition of unity $\{\psi_i\}$ subordinate to $\{V_i\}$. As ${\rm Res}_r(\eta_i) = 0$, we conclude immediately that ${\rm Res}_q(\eta_i) = 0$. Moreover, as each form $\alpha_{y_i}$ is smooth, this implies that ${\rm Res}_q(\eta) = 0$ as well. Finally, near $D$ the form $\eta$ is nondegenerate on $\ker \varphi$ as there it is given by the convex combination of forms $\eta = \sum_i (\psi \circ f) \eta_i$, where each $\eta_i$ is nondegenerate on $\ker \varphi = \langle \partial_\theta, \partial_y \rangle$ by construction.
\ep
We can now prove \autoref{thm:lfibrationsgcs} and \autoref{thm:fibrationsgcs} simultaneously, using \autoref{thm:thurstonfibration} or \autoref{thm:thurstonlfibration} respectively.
\bp[ of \autoref{thm:lfibrationsgcs} and \autoref{thm:fibrationsgcs}] Let $\omega_Y$ be a log-symplectic structure on $(Y,Z)$. As $Y$ is oriented, by \autoref{prop:blogbsymp} and \autoref{prop:orlogseparating} we know that $Z$ is separating, and in particular coorientable. Using \autoref{prop:globalbdrymap}, replace $(Y,Z)$ by $(Y',Z')$, and then using \autoref{prop:blfconnfibers}, lift to a cover $g\colon \wt{Y}' \to Y'$ so that $\wt{f}\colon (X,D) \to (\wt{Y}',Z')$ has connected fibers. It is immediate that $g^*\omega_Y$ defines a log-symplectic structure on $(\wt{Y}',Z)$, and moreover that $\wt{f}$ is a boundary Lefschetz fibration. The generic fibers of  $\wt{f}\colon X \to \wt{Y}'$ are either all homologically essential, or none are. Let $y \in Y' \setminus Z'$ be a regular point and let $y_1, \dots, y_n \in g^{-1}(y) \subset \wt{Y}'$ be its inverse images under the covering, with $F_{y_1}, \dots F_{y_n}$ their fibers under $\wt{f}$. Then for some choice of signs we have $[F_{y_1}] \pm \dots \pm [F_{y_n}] = [F]$. However, $[F] \neq 0$ by assumption, so there exists $i \in \{1,\dots,n\}$ such that $[F_{y_i}] \neq 0$. But then $[F_{y_i}] \neq 0$ for all $i \in \{1,\dots,n\}$, so $\wt{f}$ is homologically essential.

Apply \autoref{prop:homessfibrbdrymap} to the boundary (Lefschetz) fibration $\wt{f}$. This yields a global closed Lie algebroid two-form $\eta \in \Omega^2_{\rm cl}(X, \log |D|)$ such that $\eta|_{\ker \varphi}$ is nondegenerate near $D$. Moreover, its cohomology class $c = [\eta]$ pairs nontrivally with (generic) fibers in $X \setminus D$. Now apply either \autoref{thm:thurstonfibration} or \autoref{thm:thurstonlfibration}, possibly changing $\eta$ to a form $\eta'$ which agrees with $\eta$ in a neighbourhood of $D$. We obtain an elliptic symplectic structure $\omega_t = \varphi^* \omega_Y + t \eta'$ on $(X,|D|)$. As ${\rm Res}_q(\eta) = 0$ and $\eta'$ agrees with $\eta$ near $D$, we have ${\rm Res}_q(\eta') = 0$ as well. Using \autoref{prop:residuemaps} we see that ${\rm Res}_q(\omega_t) = {\rm Res}_q(\varphi^*\omega_Y) + t\, {\rm Res}_q(\eta') = 0$. But then the elliptic symplectic structure $\omega_t$ has zero elliptic residue. By \autoref{thm:sgcscorrespondence}, we conclude that $\omega_t$ for $t > 0$ small determines a \sgcs on $(X,|D|)$.
\ep
The \sgcs{}s on $(X,D)$ constructed in the proofs of \autoref{thm:lfibrationsgcs} and \autoref{thm:fibrationsgcs} arise as elliptic log-symplectic forms $\omega_X$ with zero elliptic residue through the correspondence of \autoref{thm:sgcscorrespondence}. As mentioned below \autoref{thm:sgcscorrespondence}, the three-form $H$ required for integrability of the \gcs is given by $H = {\rm Res}_r(\omega_X) \cup {\rm PD}_X[D]$. The two-form $\eta'$ introduced during the proof satisfies ${\rm Res}_r(\eta') = 0$. Using \autoref{prop:residuemaps}, together with the fact that $\omega_X = \varphi^* \omega_Y + t \eta'$ for some $t > 0$, we see that ${\rm Res}_r(\omega_X) = {\rm Res}_r(f^* \omega_Y) + t\, {\rm Res}_r(\eta') = f^* {\rm Res}_Z(\omega_Y)$.
\begin{rem}\label{rem:aboutresidues} Given a defining function $z$ for $Z$ we can locally write $\omega_Y = d \log z \wedge \alpha + \beta$, with $(\alpha,\beta) = ({\rm Res}_Z(\omega_Y), j_Z^* \omega_Y)$ the induced cosymplectic structure on $Z$. Then $f^*\omega_Y = d\log r \wedge f^*(\alpha) + f^*(\beta)$, and ${\rm Res}_r(f^*\omega_Y) = f^*({\rm Res}_Z(\omega_Y)) = f^*(\alpha)$.
\end{rem}
\subsection{Fibering over $T^2$ and $S^1$}
Recall from Section \ref{sec:scgs} that in the compact case, for a stable generalized pair $(X,D)$ and a log-Poisson pair $(Y,Z)$, the manifolds $D$ and $Z$ fiber over $T^2$ and $S^1$ respectively. After perturbing each structure, this can be achieved using one-forms on $D$ and $Z$ induced by the geometric structures. For a \sgcs $\mc{J}$ on X, these one-forms are $({\rm Res}_r,{\rm Res}_\theta)(\omega_X) \in \Omega^{1}(D) \times \Omega^1(D)$, where $\omega_X$ is the elliptic symplectic structure obtained from $\mc{J}$ using \autoref{thm:sgcscorrespondence} \cite{CavalcantiGualtieri15}. In the log-Poisson case, the one-form used instead is ${\rm Res}_Z(\omega_Y)$ (see \cite{GualtieriLi14, GuilleminMirandaPires14, OsornoTorres15, MarcutOsornoTorres14}), with $\omega_Y$ the log-symplectic structure induced by the log-Poisson structure using \autoref{prop:blogbsymp}.
\begin{rem} When $Z$ fibers over $S^1$ using $\alpha := {\rm Res}_Z(\omega_Y)$, the form $\alpha$ is Poincar\'e dual on $Z$ to the fiber $F_p$ of the induced fibration $p_Z\colon Z \to S^1$. In this case we obtain that $H = f^*({\rm PD}_Z[F_p]) \cup {\rm PD}_X[D]$.
\end{rem}
The \sgcs{}s we construct out of boundary (Lefschetz) fibrations $f\colon (X,D) \to (Y,Z)$  are such that $D$ fibers over $Z$. It is natural to ask whether this can be made compatible with the above fibrations $p_D\colon D \to T^2$ and $p_Z\colon Z \to S^1$. Note that $({\rm Res}_r,{\rm Res}_\theta)(\omega_X) = (f^*{\rm Res}_Z(\omega_Y), t\, {\rm Res}_\theta(\eta))$ as $\omega_X = \varphi^* \omega_Y + t \eta$, and ${\rm Res}_\theta(\varphi^* \omega_Y) = 0$ by the discussion below \autoref{prop:bdrymaplocalcoord}. This immediately shows we can make the fibrations compatible, as we can choose $\eta$ such that $t\, {\rm Res}_\theta(\eta)$  and $f^*{\rm Res}_Z(\omega_Y)$ determine a fibration, instead of a just foliation. As in Section \ref{sec:constrbdrlfs}, let $p\colon T^2 \to S^1$ be projection onto the first factor.
\begin{cor}\label{cor:compatiblefibrations} Under the assumptions of \autoref{thm:lfibrationsgcs} or \autoref{thm:fibrationsgcs}, given a \blog $\omega_Y$ on $(Y,Z)$, the elliptic two-form $\omega_X$ can be chosen such that $p \circ p_D = p_Z \circ f|_D$.
\end{cor}
In other words, we have a full commutative diagram around $D$ and $Z$ as follows:
\begin{center}
	\begin{tikzpicture}
	\matrix (m) [matrix of math nodes, row sep=2.5em, column sep=2.5em,text height=1.5ex, text depth=0.25ex]
	{	 ND & D & T^2 \\ NZ & Z & S^1\\};
	\path[-stealth]
	(m-1-1) edge node [above] {${\rm pr}_D$} (m-1-2)
	(m-1-1) edge node [left] {$f$} (m-2-1)
	(m-1-2) edge node [above] {$p_D$} (m-1-3)
	(m-1-2) edge node [left] {$f|_D$} (m-2-2)
	(m-2-1) edge node [above] {${\rm pr}_Z$} (m-2-2)
	(m-1-3) edge node [left] {$p$} (m-2-3)
	(m-2-2) edge node [above] {$p_Z$} (m-2-3);
	\end{tikzpicture}
\end{center}
\section{Examples and applications}
\label{sec:applications}
In this section we discuss some applications of the obtained results. We start with the following immediate consequence of the results in Sections \ref{sec:constrbdrlfs} and \ref{sec:constrsgcs}.
\begin{thm}\label{thm:g1lfsgcs} Let $f\colon X \to D^2$ be a genus one Lefschetz fibration over the disk whose monodromy around the boundary is a power of a Dehn twist. Then all possible completions $\wt{f}\colon (\wt{X}, D) \to (D^2, \partial D^2)$ admit a \sgcs.
\end{thm}
\bp Any Lefschetz fibration over the disk $D^2$ is homologically essential. By the monodromy assumption around the boundary, $f$ admits completions $\wt{f}$ using \autoref{cor:lfcompletion}. By \autoref{prop:retainhomess}, the map $\wt{f}$ is a homologically essential boundary Lefschetz fibration, so that $(\wt{X}, D)$ admits a \sgcs by \autoref{thm:lfibrationsgcs}, as the completion has $D$ coorientable.
\ep
We now turn to exhibiting \sgcs{}s on concrete four-manifolds using boundary (Lefschetz) fibrations.
\begin{exa} Let $I = [-1,1]$ with coordinate $t$ and view $S^2 \subseteq \R^3$ with north and south pole $p_N$ and $p_S$. Consider the standard height function $p\colon (S^2, p_N \cup p_S) \to (I,\partial I)$ given by $(x,y,z) \mapsto z$. The map $h\colon (I,\partial I) \to (\R, 0)$ given by $h(t) := 1 - t^2$ is a defining function for $\partial I$. Then $p^*(h)(x,y,z) = 1 - (1 - x^2 - y^2) = x^2 + y^2$, so that $p$ is a boundary map by the description of \autoref{prop:bdrymaplocalcoord}. In fact, $p$ is a boundary fibration. Let $\varphi\colon S^3 \to S^2$ be the Hopf fibration, making $p \circ \varphi\colon (S^3, S^1_N \cup S^1_S) \to (I,\partial I)$ a boundary fibration, where $S^1_N = \varphi^{-1}(p_N)$ and similarly for $S^1_S$. Finally, define the boundary fibration $f\colon (S^3 \times S^1, D) \to (I \times S^1, Z)$, where $Z = \partial I \times S^1$ and $D = (S^1_N \cup S^1_S) \times S^1 = f^{-1}(Z)$, by $f(x,\theta) = (p \circ \varphi)(x)$. Using \autoref{lem:blogsurface}, $(I \times S^1, Z)$ admits a \blog, while the fibers of $f$ are clearly homologically essential. By \autoref{thm:fibrationsgcs} applied to $f$, we conclude that $(S^1 \times S^3, D)$ admits a \sgcs\ whose type change locus has two connected components. This structure is integrable with respect to the zero three-form due to \autoref{rem:aboutresidues}.
\end{exa}
Using a slightly different map, we can ensure that the type change locus is connected.
\begin{exa} View $S^3 \subseteq \C^2$ with coordinates $(z_1,z_2)$ in the standard way and consider the map $\varphi\colon S^3 \to D^2$ given by $\varphi(z_1,z_2) = z_2$, viewing $D^2 \subseteq \C$. In the interior $D^2 \setminus \partial D^2$, the map $\varphi$ admits a section $s(z) = (\sqrt{1 - |z|^2},z)$. As $\varphi \circ s = {\rm id}$, we see that $\varphi$ is submersive on $S^3 \setminus D$, where $D$ is the circle $D = \varphi^{-1}(\partial D^2) = \{(0,z_2) \in \C^2: |z_2| = 1\}$. Moreover, $\varphi|_D$ is the identity from $D$ to $Z = \partial D^2$, which in particular implies that $\varphi|_D$ is submersive. Let $h(z) := 1 - |z|^2$ be a defining function for $\partial D^2$ on $D^2$. Then $\varphi^*(h)(z_1,z_2) = 1 - |z_2|^2$, which shows that $\varphi$ is a boundary map using \autoref{prop:bdrymaplocalcoord}. But then $\varphi$ is in fact a boundary fibration. Noting that $\varphi|_{S^3 \setminus D}$ admits a section, $\varphi\colon S^3 \setminus D \to D^2 \setminus D^2$ is a trivial $S^1$-bundle. Now define a map of pairs $f\colon (S^3\times S^1, D \times S^1 ) \to (D^2, \partial D^2)$ by $f(z_1,z_2,x) := \varphi(z_1,z_2)$. This is also a boundary fibration, which is homologically essential as $f|_{ (S^3 \setminus D) \times S^1}$ is the projection $S^1  \times (D^2 \setminus \partial D^2) \times S^1 \to (D^2, \partial D^2)$. Finally, $(D^2, \partial D^2)$ admits a \blog by \autoref{lem:blogsurface}. Using \autoref{thm:fibrationsgcs} we conclude that $(S^3 \times S^1,  D \times S^1)$ admits a \sgcs with connected type change locus and is integrable with respect to a nonzero degree three cohomology class.
\end{exa}
Before we consider actual boundary Lefschetz fibrations, it is convenient to reinterpret \autoref{prop:completion} with a Kirby calculus point of view. There, the completion of each end of $X$ was done by gluing in a copy of the total space of $L_n$, the disk bundle over the torus with Euler class $n$. Observe that ${\rm Tot}(L_n)$ is a four-manifold built out of one 0-handle, two 1-handles and one 2-handle. The process of capping off an end of $X$ with ${\rm Tot}(L_n)$ amounts to inverting the handle structure of ${\rm Tot}(L_n)$ and adding it to $X$. That is, to cap each end of $X$ we add one 2-handle, two 3-handles and one 4-handle. Further, due to inversion, the 2-handle of ${\rm Tot}(L_n)$ has its original core and co-core interchanged and once this 2-handle is added there is a unique way to complete with the 3- and 4-handles, up to diffeomorphism. Therefore we only have to describe how to attach the 2-handle, whose core is the circle fiber of the projection $M \to T^2$, where $M$ is the $\varepsilon$-sphere bundle of $L_n$. From the discussion in \autoref{prop:monodromy}, this fiber is a circle left fixed by the monodromy map, which, if the monodromy is nontrivial, is determined by the monodromy map and is essentially unique. Further, since the circle is the fiber of the projection $M \to T^2$, the framing, in double strand notation, is  given by the nearby fibers.
\begin{exa} Consider the genus one Lefschetz fibration over $D^2$ with one singular fiber with vanishing cycle $b$, where $b \in H^1(T^2)$ is a generator. Then the monodromy around $\partial D^2$ is $b$, which is clearly the power of a Dehn twist. The resulting completion $\wt{X}$ thus admits a \sgcs by \autoref{thm:g1lfsgcs}, and we are left with determining its diffeomorphism type. Draw a Kirby diagram for the trivial $T^2$-bundle over $D^2$ in the standard way (see \cite{GompfStipsicz99}) and attach a $-1$-framed two-handle to represent the vanishing cycle $b$. The completion process adds a single $0$-framed two-handle along the same cycle (see Figure \ref{fig:S1xS3}). Slide the $-1$-framed two-handle over this to split off a $-1$-framed unknot, representing a copy of $\overline{\C P}^2$. The remaining diagram collapses to a single one-handle after two further handle slides, which shows that $\wt{X} = S^1 \times S^3 \# \overline{\C P}^2$.
\begin{figure}[h!!]
	\begin{center}
		\includegraphics[height=3.5cm]{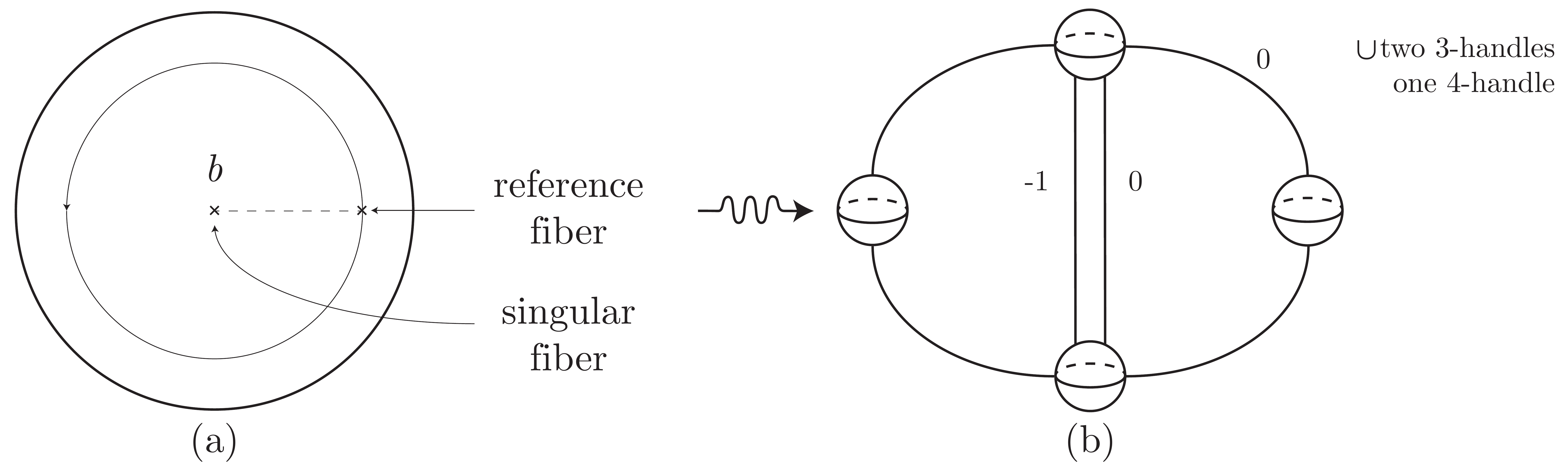} \caption{Figure (a) shows the base of a boundary Lefschetz fibration over the disk with a single nodal fiber, a choice of regular fiber for reference and loop around the boundary used to compute the monodromy and (b) shows the corresponding Kirby diagram.} \label{fig:S1xS3}
	\end{center}
\end{figure}
\end{exa}	
We can also recover the examples of \cite{CavalcantiGualtieri09}, showing that $m \C P^2 \# n \overline{\C P}^2$ admits a \sgcs if and only if $m$ is odd (i.e., if it admits an almost-complex structure).
\begin{exa} Let $m, n \in \N$ and assume that $m = 2k + 1$ is odd. Let $a, b \in H^1(T^2)$ be generators and consider $(D^2, \partial D^2)$ with clockwise assignment of singular fibers given by $a + (4k-1)b$, $a + 4 j b$ for $j$ from  $k-1$ up to $1-k$, and $a - (4k-1)b$, and $n$ copies of $b$. The monodromy around the boundary can be computed to be $(10k-1-n)b$. Therefore the associated genus one Lefschetz fibration admits a completion to $f\colon (\wt{X},D) \to (D^2, \partial D^2)$ (see Figure \ref{fig:An}). As is shown in \cite{CavalcantiGualtieri09}, we have $\wt{X} = m\C P^2\# n\overline{\C P}^2$, which hence has a \sgcs.
\begin{figure}[h!!]
	\begin{center}
		\includegraphics[height=4.5cm]{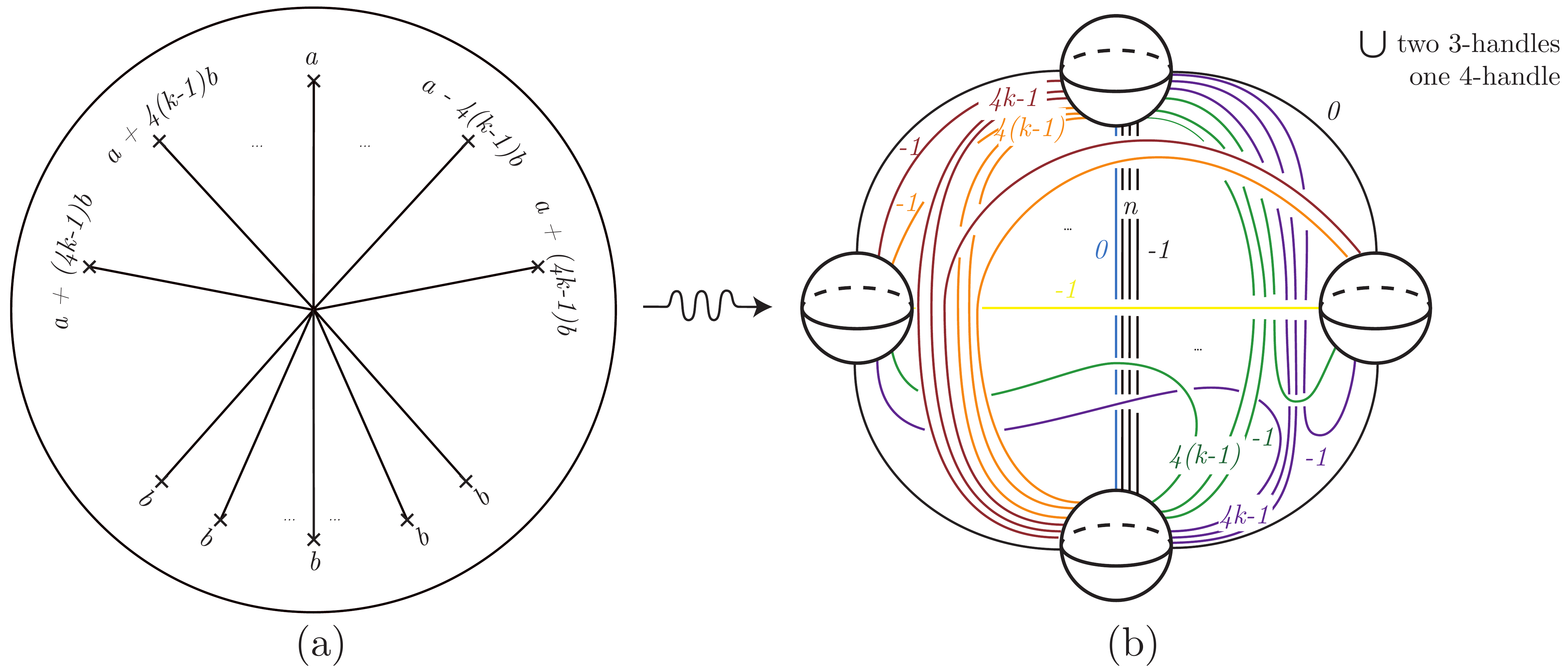} \caption{Figure (a) shows the base of a boundary Lefschetz fibration over the disc with $m+ n$ singular fibers and a choice of regular fiber for reference at the center of disk. Figure (b) shows the corresponding Kirby diagram.} \label{fig:An}
	\end{center}
\end{figure}
\end{exa}
Not every \sgcs comes from a boundary \lf, similarly to the case for symplectic structures and Lefschetz fibrations, or \blog{}s and achiral Lefschetz fibrations (see Section \ref{sec:blogconstructions}).
\begin{exa} The symplectic manifold $(\C P^2 , \omega_{\rm FS})$ carries a \sgcs with $D = \emptyset$. As $(\C P^2, \omega_{\rm FS})$ is not a symplectic fibration over any surface, this \sgcs can not be obtained through our construction.
\end{exa}
\begin{exa} Let $\Sigma_g$ be the genus $g$ surface and $g > 1$. Then $X = S^2 \times \Sigma_g$ has negative Euler characteristic so cannot admit a \placeholder by \autoref{prop:blfsingfibers}. Hence $X$ does not admit a \sgcs obtained through our methods. However, $X$ is symplectic hence carries a \sgcs with $D = \emptyset$.
\end{exa}
\chapter{Homotopical obstructions for \texorpdfstring{$\mc{A}$}{A}-symplectic structures}
\label{chap:homotopicalobstrs}
\fancypagestyle{empty}{%
	\fancyhf{}%
	\renewcommand\headrulewidth{0pt}%
	\fancyhead[RO,LE]{\thepage}%
}
In this chapter we discuss homotopical obstructions to the existence of $\mc{A}$-symplectic structures (see Chapter \ref{chap:asymplecticstructures}) on a given closed manifold $X$. We focus on $\mc{A}$-orientability and the existence of an $\mc{A}$-almost complex structure, i.e.\ a complex structure on $\mc{A}$. In particular, we introduce the discrepancy for four-dimensional log pairs and prove the following result (\autoref{cor:azsympstr}) mentioned in the introduction to this thesis.
\newtheorem*{thm:homtopobstr}{Theorem \ref{thm:introhomtopobstr}}
\begin{thm:homtopobstr} Let $X$ be a compact oriented four-manifold that admits a log-symplectic structure with singular locus $Z$. Then $b_2^+(X) + b_1(X) + f(X,Z)$ is odd.
\end{thm:homtopobstr}
Here $f(X,Z)$ is the discrepancy of the log pair $(X,Z)$ (see \autoref{defn:discrepancy}). This result can be used to obstruct the existence of \blog{}s on given four-manifolds, if one is willing to also prescribe the singular locus.
\subsection*{Organization of the chapter}
This chapter is built up as follows. In Section \ref{sec:stablebundleiso} we start by proving stable bundle isomorphisms for certain Lie algebroids $\mc{A}$ we have encountered before, notably the log-tangent bundle. In Section \ref{sec:compcharclasses} we then use these to compute the relevant characteristic classes of these Lie algebroids. We finish with Section \ref{sec:existanambu} and Section \ref{sec:existenceaacs} by using this information to obtain obstructions to the existence of $\mc{A}$-Nambu structures, and of $\mc{A}$-almost-complex structures, and discuss in some detail the case of the log-tangent bundle.
\section{Stable bundle isomorphisms}
\label{sec:stablebundleiso}
In this section we establish stable bundle isomorphisms for Lie algebroids obtained using a log pair $(X,Z)$, namely the log-tangent bundle $\mc{A}_Z$ of Section \ref{sec:logtangentbundle}, and more generally the $b^k$-tangent bundles $\mc{A}_Z^k$ for any $k$ (see Section \ref{sec:bktangentbundle}).
\begin{prop}\label{prop:logbundleiso} Let $(X,Z)$ be a log pair. Then $\mc{A}_Z \oplus L_Z \cong TX \oplus \underline{\R}$.
\end{prop}
\begin{rem} This was noted without proof in \cite{CannasDaSilva10, CannasDaSilvaGuilleminWoodward00} when $X$ is orientable, although they incorrectly state that in general one has $\mc{A}_Z \oplus \underline{\R} \cong TX \oplus L_Z$. However, when $NZ$ is trivial, so that $L_Z$ is trivial (see \autoref{prop:logdivnormalbundle}), this statement is equivalent to \autoref{prop:logbundleiso}. In particular, note that if $X$ is orientable and $(X,Z)$ admits a \blog, then $NZ$ is trivial by \autoref{prop:orlogseparating}.
\end{rem}
\bp The isomorphism locus of $\mc{A}_Z$ is $X \setminus Z$, and $L_Z$ is trivial over $X \setminus Z$. Hence, away from $Z$, $\mc{A}_Z$ is isomorphic to $TX$ and $L_Z$ is trivial, so that there we have an obvious isomorphism. Near $Z$ we define a bundle isomorphism as follows. Choose a tubular neighbourhood embedding of $Z$ and let $\{U_\alpha\}_\alpha$ be a trivializing cover of $L_Z$ (and hence $NZ$), with $z_\alpha$ the associated normal bundle coordinates.  Choose a metric on $NZ$ and consider the disk bundle in $NZ$ with radius $\pi/2$. Complete to a coordinate system $(z_\alpha,x_2, \dots, x_n)$ on $U_\alpha$. We define a map $\mc{A}_Z \oplus L_Z \to TX \oplus \underline{\R}$ on $U_\alpha$, or rather the disk bundle of radius $\pi/2$ inside $NZ \cap U_\alpha$. Generic sections of $\mc{A}_Z \oplus L_Z$ and $TX \oplus \underline{\R}$ are expressed as $(\lambda_1 z_\alpha \partial_{z_\alpha} + \sum_{2 \leq i \leq n} \lambda_{i} \partial_{x_i}, \lambda_{n+1} \cdot s_\alpha)$ and $(\mu_1 \partial_{z_\alpha} + \sum_{2\leq i \leq n} \mu_i \partial_{x_i}, \mu_{n+1} \cdot \mathbf{1})$, where $\lambda_i, \mu_i \in C^\infty(U_\alpha)$. The map is given by
\be
\begin{pmatrix} \mu_1 \\ \mu_{n+1} \\ \mu_2 \\ \dots \\ \mu_{n}\end{pmatrix} = \begin{pmatrix} z_\alpha \sin |z_\alpha| & \cos |z_\alpha| & 0 \\ - \cos |z_\alpha| & z_\alpha \sin |z_\alpha| & 0 \\ 0 & 0 & I_{n-1} \end{pmatrix} \begin{pmatrix} \lambda_1 \\ \lambda_{n+1} \\ \lambda_2 \\ \dots \\ \lambda_{n}\end{pmatrix},
\ee
where $I_{n-1}$ is the $(n-1) \times (n-1)$ identity matrix. When $|z_\alpha|$ is large, namely equal to $\pi/2$, the topleftmost part of the matrix is equal to $z_\alpha$ times the identity (and $\mathbf{1}$ is sent to $1 \cdot z_\alpha s_\alpha = s_0$), while it is equal to $\begin{pmatrix} 0 & 1 \\ -1 & 0 \end{pmatrix}$ when $|z_\alpha| = 0$, i.e.\ on $Z$. In particular this map can be extended continuously to $U_0$ as being $\pi/2$ times the identity (and also smoothly after introducing a radial bump function). Moreover, it is well-defined as on overlaps $U_\alpha \cap U_\beta$ we have $z_\alpha \partial_{z_\alpha} = z_\beta \partial_{z_\beta}$ so that $\partial_{z_\alpha} = \frac{z_\beta}{z_\alpha} \partial_{z_\beta} = g^{\beta}_\alpha \partial_{z_\beta}$, and also $s_\alpha = \frac{z_\beta}{z_\alpha} s_\beta = g^{\beta}_\alpha s_\beta$. The map has determinant equal to $z_\alpha^2 \sin^2 |z_\alpha| + \cos^2 |z_\alpha|$, which is positive for all values of $z_{\alpha}$, showing invertibility.
\ep
A similar thing can be done for the $b^k$-tangent bundles $\mc{A}_Z^k$.
\begin{prop}\label{prop:bkbundleiso} Let $(X,Z)$ be a log pair, $j_{k-1} \in J_Z^{k-1}$ a choice of $(k-1)$-jet, and $\mc{A}_Z^k$ the associated Lie algebroid. Then $\mc{A}_Z^k \oplus L_Z^k \cong TX \oplus \underline{\R}$.
\end{prop}
In particular, we claim that the statement is independent of $j_{k-1}$.
\bp Follow the same strategy as for \autoref{prop:logbundleiso}, except take $z_\alpha \in j_{k-1}$, replace $\lambda_1 z_\alpha \partial_{z_\alpha}$ by $\lambda_1 z_\alpha^k \partial_{z_\alpha}$, and define the map to be
\be
\begin{pmatrix} \mu_1 \\ \mu_{n+1} \\ \mu_2 \\ \dots \\ \mu_{n}\end{pmatrix} = \begin{pmatrix} z_\alpha^k \sin |z_\alpha| & \cos |z_\alpha| & 0 \\ - \cos |z_\alpha| & z_\alpha^k \sin |z_\alpha| & 0 \\ 0 & 0 & I_{n-1} \end{pmatrix} \begin{pmatrix} \lambda_1 \\ \lambda_{n+1} \\ \lambda_2 \\ \dots \\ \lambda_{n}\end{pmatrix}.
\ee
This map has determinant $z_\alpha^{2k} \sin^2 |z_\alpha| + \cos^2 |z_\alpha|$, again showing invertibility.
\ep
\begin{rem} Inspecting the proof of \autoref{prop:logbundleiso}, we see that more generally given a Lie algebroid $\mc{A} \to X$ and a $(\mc{B},Z)$-rescaling ${}^\mc{B} \mc{A}$ of $\mc{A}$ with ${\rm corank}(\mc{B}) = k$, we have ${}^\mc{B} \mc{A} \oplus L_Z \oplus \dots \oplus L_Z \cong \mc{A} \oplus \underline{\R} \oplus \dots \oplus \underline{\R}$, with $k$ copies of $L_Z$ and $\underline{\R}$. The case of \autoref{prop:logbundleiso} is when $\mc{A} = TX$ and $\mc{B} = TZ$, which has corank one. Compare this with \autoref{rem:larescaling}.
\end{rem}
Consider the case when $Z$ is separating, i.e.\ when $[Z]$ is trivial in homology with $\Z_2$-coefficients, as happens when $X$ is oriented and $Z$ is the vanishing locus of a \blog on $X$ (see \autoref{prop:orlogseparating}). Then from \autoref{prop:logbundleiso} and \autoref{prop:logdivstiefelwhitney} we see that $\mc{A}_Z$ and $TX$ are stably isomorphic. As we will see in the next section, having such an ``almost'' stable bundle isomorphism allows us to determine some of the relevant characteristic classes of $\mc{A}_Z$. However, even if $L_Z$ were trivial so that this is an actual stable bundle isomorphism, we cannot fully determine the isomorphism class of $\mc{A}_Z$ using \autoref{prop:logbundleiso}. More precisely, even if $X$ is orientable, the Euler class of $\mc{A}_Z$ cannot be recovered. This relates to the discussion on induced orientations in Section \ref{sec:logtangentbundle}, and the following proposition, which is remarked in \cite{CannasDaSilva10}.
\begin{prop} Let $E, F$ be two stably isomorphic orientable vector bundles of rank $2n$ over a $2n$-dimensional connected manifold. Then $E$ and $F$ are isomorphic if and only if their Euler classes agree.
\end{prop}
Recall that the Euler class of an oriented vector bundle $E \to X$ is given by $e(E) = {\rm PD}_{\Z}[{\rm im}(s)]$, where $s \in \Gamma(E)$ is a transverse section of $E$. The Euler characteristic of $E$ in case $X$ is compact and oriented is given by $\chi(E) = \langle e(E), [X] \rangle$, or by a signed count of the zeros of $s$. These signs are crucial, as they depend on the orientations involved. We will see in Section \ref{sec:existenceaacs} what this implies while comparing $\mc{A}_Z$ to $TX$.

\section{Computing characteristic classes}
\label{sec:compcharclasses}
In this section we compute characteristic classes of some of the Lie algebroids we have introduced. We will mainly be interested in $w_1$, $w_2$ and $p_1$, the first and second Stiefel--Whitney classes, and the first Pontryagin class. In preparation of this, recall the following properties of the Stiefel--Whitney classes, which can be found in e.g.\ \cite{MilnorStasheff74}. Note especially the behavior under tensor products.
\begin{prop}\label{prop:stiefelwhitney} Let $E, F \to X$ be real vector bundles of rank $m$ and $n$ respectively. Let $w\colon {\rm Vect}(X) \to H^\bullet(X;\Z_2)$ denote the full Stiefel--Whitney class. Then
	\bi
	\item[i)] $w(E \oplus F) = w(E) \cup w(F)$;
	\item[ii)] $w_1(E) = w_1(\det(E))$;
	\item[iii)] $\det(E \otimes F) \cong \det(E)^n \otimes \det(F)^m$;
	\item[iv)] $w_1(E \otimes F) = w_1(E) + w_1(F)$ if $m = n = 1$;
	\item[v)] $w_1(E \otimes F) = n w_1(E) + m w_1(F)$;
	\item[vi)] $w_2(E \otimes F) = w_2(E) + (m-1) w_1(F) \cup w_1(E) + \frac12 m(m-1) w_1(F)^2$ if $n = 1$;
	\item[vii)] $w_2(E \otimes F) = w_2(E) + w_1(F) \cup w_1(E)$ if $m = 4$ and $n = 1$.
	\ei
\end{prop}
\bp Part vii) immediately follows from vi). We prove v) and vi) as these relations are seldom used. For v), using properties ii)-iv) we have $w_1(E \otimes F) = w_1(\det(E\otimes F)) = w_1(\det(E)^n \otimes \det(F)^m) = n w_1(\det(E)) + m w_1(\det(F)) = n w_1(E) + m w_1(F)$. To show vi), we have to use a little more. Recall that in general we have
\be
	w(E \otimes F) = P_{m,n}(w_1(E),\dots,w_m(E),w_1(F),\dots,w_n(F)),
\ee
where $P_{m,n}$ are universal polynomials (depending only on $m$ and $n$). These polynomials can be characterized as follows. Let $\sigma_i, \tau_j$ denote the elementary symmetric polynomials in indeterminate variables $s_i$ and $t_j$. Then
\be
	P_{m,n}(\sigma_1,\dots,\sigma_m,\tau_1,\dots,\tau_n) = \prod_{i = 1}^m \prod_{j = 1}^n (1 + s_i + t_j). 
\ee
We now specify to the case where $n = 1$. Then we obtain
\begin{align*}
	P_{m,1}(\sigma_1,\dots,\sigma_m,\tau_1) &= \prod_{i = 1}^m (1 + s_i + t_1)\\
	&= \sum_{I \subseteq \{1,\dots,m\}} (1 + t_1)^{m - |I|} \cdot \prod_{i \in I} s_i\\
	&= \sum_{k = 0}^m (1 + \tau_1)^{m-k} \sigma_k.
\end{align*}
We now obtain the result by substituting the $w_1(F)$ for $\tau_1$ and $w_i(E)$ for the $\sigma_i$, as in degree one the above expression reads $n \tau_1 + \sigma_1$ (recovering v) in this case), while in degree two we obtain $\frac12 m(m-1) \tau_1^2 + (m-1) \tau_1 \sigma_1 + \sigma_2$.
\ep
In particular, the above proposition implies that $w_1(E \oplus F) = w_1(E) + w_1(F)$. We will also use the following properties of Pontryagin classes, see e.g.\ \cite{MilnorStasheff74}.
\begin{prop}\label{prop:pontryagin} Let $E, F \to X$ be vector bundles and denote by $p$ the full Pontryagin class. Then $2 p(E \oplus F) = 2 p(E) \cup p(F)$. Moreover, if $F$ has rank one, then $p(E \otimes F) = p(E)$.
\end{prop}
In other words, $p(E \oplus F)$ and $p(E) \cup p(F)$ agree up to two-torsion in $H^\bullet(X;\Z)$. Consider the log-tangent bundle $\mc{A}_Z = TX(-\log Z)$ associated to a log pair $(X,Z)$.
\begin{prop}\label{prop:logcharclass} Let $(X,Z)$ be a log pair. Then $w_1(\mc{A}_Z) = w_1(TX) + w_1(L_Z)$ and $w_2(\mc{A}_Z) = w_2(TX) + w_1(L_Z) \cup w_1(TX)$. If $X$ is orientable, then $w_2(\mc{A}_Z) = w_2(TX)$. If $X$ is moreover four-dimensional, then $p_1(\mc{A}_Z) = p_1(TX)$.
\end{prop}
\bp By \autoref{prop:logbundleiso} we have $\mc{A}_Z \oplus L_Z \cong TX \oplus \underline{\R}$, hence $w(\mc{A}_Z) = (1 + w_1(L_Z)) \cup w(TX)$ due to \autoref{prop:stiefelwhitney}.i). In degree one this gives $w_1(\mc{A}_Z) = w_1(TX) + w_1(L_Z)$ as desired. In degree two it follows that $w_2(\mc{A}_Z) = w_2(TX) + w_1(L_Z) \cup w_1(TX)$. Further, we see that $2 p_1(\mc{A}_Z) = 2 p_1(TX)$, as $p \equiv 1$ for line bundles. If $X$ is orientable then $w_1(TX) = 0$ so that $w_2(\mc{A}_Z) = w_2(TX)$ is immediate. Finally, if $X$ is also four-dimensional we know that $H^4(X;\Z)$ is equal to $\Z$, which in particular has no $2$-torsion, so that $p_1(\mc{A}_Z) = p_1(TX)$.
\ep
\begin{cor}\label{cor:orlogpairseparating} Let $(X,Z)$ be a log pair for which both $\mc{A}_Z$ and $TX$ are orientable. Then $Z$ is separating.
\end{cor}
\bp This is immediate from \autoref{prop:logbundleiso} as then ${\rm PD}_{\Z_2}[Z] = w_1(L_Z) = 0$, using \autoref{prop:logdivstiefelwhitney}.
\ep
Recall from \autoref{prop:zeroisotangent} that the zero tangent bundle satisfies $\mc{B}_Z \cong TX \otimes L_Z^*$. Using that $L_Z^* \cong L_Z$ and \autoref{prop:stiefelwhitney}.v) and vii) we obtain the following.
\begin{prop}\label{prop:zerocharclass} Let $(X,Z)$ be a log pair. Then $w_1(\mc{B}_Z) = w_1(TX) + n w_1(L_Z)$ and if $X$ is four-dimensional, $w_2(\mc{B}_Z) = w_2(TX) + w_1(L_Z) \cup w_1(TX)$.
\end{prop}
We have no direct use for $w_2(\mc{B}_Z)$ and $p_1(\mc{B}_Z)$ (when $X$ is four-dimensional), as by \autoref{prop:zerotangentnosymp} we know that $\mc{B}_Z$ does not admit Lie algebroid symplectic structures when $\dim X \geq 4$.

Similarly, for the scattering tangent bundle $\mc{C}_Z$ we have $\mc{C}_Z \cong \mc{A}_Z \otimes L_Z$ (see \autoref{prop:scatteringisolog}). This results in the following.
\begin{prop}\label{prop:scatcharclass} Let $(X,Z)$ be a log pair. Then $w_1(\mc{C}_Z) = w_1(TX) + (n+1) w_1(L_Z)$. If $X$ is four-dimensional, then $w_2(\mc{C}_Z) = w_2(TX)$. If $X$ is also orientable, then $p_1(\mc{C}_Z) = p_1(TX)$.
\end{prop}
\bp By \autoref{prop:logbundleiso} we have $\mc{A}_Z \oplus L_Z \cong TX \oplus \underline{\R}$. Upon tensoring with $L_Z$ we thus get $\mc{C}_Z \oplus L_Z^2 \cong \mc{B}_Z \oplus L_Z$. As $L_Z^2$ is canonically trivial, this gives $w(\mc{C}_Z) = (1 + w_1(L_Z)) \cup w(\mc{B}_Z)$ by \autoref{prop:stiefelwhitney}.i). In degree one this results in $w_1(\mc{C}_Z) = w_1(\mc{B}_Z) + w_1(L_Z) = w_1(TX) + (n+1)w_1(L_Z)$, using \autoref{prop:zerocharclass}. In degree two we see if $X$ is four-dimensional that $w_2(\mc{C}_Z) = w_2(\mc{B}_Z) + w_1(L_Z) \cup w_1(\mc{B}_Z) =(w_2(TX) + w_1(L_Z) \cup w_1(TX)) + w_1(L_Z) \cup (w_1(TX) + 4 w_1(L_Z)) = w_2(TX)$, regardless of whether $X$ is orientable. Assuming orientability, \autoref{prop:pontryagin} and \autoref{prop:logcharclass} gives the result regarding $p_1(\mc{C}_Z)$.
\ep
Given a log pair $(X,Z)$, we can form the $b^k$-tangent bundles $\mc{A}_Z^k$ using a choice of jet data $j_{k-1} \in J_Z^{k-1}$, see Section \ref{sec:bktangentbundle}. For $\mc{A}_Z^k$ we obtain the following using \autoref{prop:bkbundleiso}, in analogy with \autoref{prop:logcharclass} upon replacing $L_Z$ by $L_Z^k$.
\begin{prop}\label{prop:bkcharclass} Let $(X,Z)$ be a log pair and $j_{k-1} \in J_Z^{k-1}$ a $(k-1)$-jet at $Z$. Then $w_1(\mc{A}_Z^k) = w_1(TX) + k w_1(L_Z)$ and $w_2(\mc{A}_Z^k) = w_2(TX) + k w_1(L_Z) \cup w_1(TX)$. If $X$ is orientable, then $w_2(\mc{A}_Z^k) = w_2(TX)$. If $X$ is moreover four-dimensional, then $p_1(\mc{A}_Z^k) = p_1(TX)$.
\end{prop}

\section{Existence of \texorpdfstring{$\mc{A}$}{A}-Nambu structures}
\label{sec:existanambu}
In this section we discuss the existence of Nambu structures on the Lie algebroids $\mc{A}_Z^k$, $\mc{B}_Z$ and $\mc{C}_Z$ associated to log pairs $(X,Z)$. This also settles when these Lie algebroids admit symplectic structures in dimension two (see \autoref{rem:asympranktwo}), and gives an obstruction to their existence in arbitrary dimensions. In particular, we give simple topological proofs of results in \cite{MirandaPlanas16}, only using a local description to obtain the isomorphism of \autoref{prop:bkbundleiso}. Recall that for a log pair $(X,Z)$, we denote by $\mc{A}_Z^k$ any $b^k$-tangent bundle constructed out of a choice of jet data, see Section \ref{sec:bktangentbundle}.
\begin{thm}\label{thm:akznambustructure} Let $(X,Z)$ be a compact log pair. Then $X$ admits an $\mc{A}_Z^k$-Nambu structure if and only if $w_1(TX) + k {\rm PD}_{\Z_2}[Z] = 0$.
\end{thm}
\bp This follows immediately from \autoref{prop:anambustructure} and \autoref{prop:bkcharclass}, using \autoref{prop:logdivstiefelwhitney} to determine $w_1(L_Z)$.
\ep
\begin{cor}[\cite{MirandaPlanas16}] Let $(X,Z)$ be a compact log pair. If $X$ admits an $\mc{A}_Z^k$-Nambu structure, then $X$ is orientable if and only if $k$ is even or $[Z] = 0$.
\end{cor}
\bp By \autoref{thm:akznambustructure} we have that $w_1(TX) + k {\rm PD}_{\Z_2}[Z] = 0$. As orientability of $X$ is equivalent to the condition $w_1(TX) = 0$, the result follows.
\ep
We obtain similar results for the Lie algebroids $\mc{B}_Z$ and $\mc{C}_Z$.
\begin{thm} Let $(X,Z)$ be an $n$-dimensional compact log pair. Then $X$ admits a $\mc{B}_Z$-Nambu structure if and only if $w_1(TX) + n {\rm PD}_{\Z_2}[Z] = 0$, and a $\mc{C}_Z$-Nambu structure if and only if $w_1(TX) + (n+1) {\rm PD}_{\Z_2}[Z] = 0$.
\end{thm}
\begin{cor} Let $(X,Z)$ be an $n$-dimensional compact log pair. If $X$ admits a $\mc{B}_Z$- respectively $\mc{C}_Z$-Nambu structure, then $X$ is orientable if and only if either $[Z] = 0 \in H_{n-1}(X;\Z_2)$, or $n$ is even (for $\mc{B}_Z)$, or $n$ is odd (for $\mc{C}_Z$).
\end{cor}
\section{Existence of \texorpdfstring{$\mc{A}$}{A}-almost-complex structures}
\label{sec:existenceaacs}
In this section we specify to the case when $X$ is four-dimensional. We recall a classical result by Wu \cite{Wu52}, see also \cite{HirzebruchHopf58,LePotier}, characterizing when an oriented rank four vector bundle $E \to X$ admits a complex structure, usually only stated for $E = TX$.
\begin{thm}[Wu, \cite{Wu52}]\label{thm:wuacs} Let $X$ be a compact oriented four-manifold and $E \to X$ an oriented Euclidean vector bundle of rank four. Then $E$ admits a complex structure if and only if there exists a class $c \in H^2(X;\Z)$ such that $c \pmod 2 \equiv w_2(E) \in H^2(X;\Z_2)$ and $c^2 = p_1(E) + 2 e(E)$.
\end{thm}
If the class $c$ in the above theorem exists, the complex structure on $E$ will satisfy $c_1(E) = c$.
This result is relevant as follows: if $X$ admits an $\mc{A}$-symplectic structure for a certain Lie algebroid $\mc{A} \to X$, it also admits an $\mc{A}$-almost complex structure, i.e.\ $\mc{A}$ can be equipped with a complex structure (see \autoref{prop:asympaacs}). Consequently, for Lie algebroids with dense isomorphism locus (which hence satisfy ${\rm rank}(\mc{A}) = \dim X = 4$), we can obstruct the existence of $\mc{A}$-symplectic structures using \autoref{thm:wuacs}.

We now specify to the Lie algebroids associated to log pairs $(X,Z)$. As it is necessary that both $TX$ and $\mc{A}_Z$ are orientable, by \autoref{prop:logbundleiso} and \autoref{prop:logdivstiefelwhitney} we see that $Z$ must be nullhomologous, i.e.\ $[Z] = 0$. Denote the space of all such hypersurfaces by ${\rm Hyp}_0(X)$.
\begin{defn}\label{defn:discrepancy} Given $Z \in {\rm Hyp}_0(X)$ and orientations on both $\mc{A}_Z$ and $TX$, the \emph{discrepancy} $f(X,Z)$ of $Z$ is defined as the difference $2 f(X,Z) = e(\mc{A}_Z) - e(TX)$.
\end{defn}
\begin{rem} As the Euler class of an oriented vector bundle reduces mod 2 to its top Stiefel--Whitney class, we see that the difference in Euler classes indeed is even. Namely, Stiefel--Whitney classes are stable invariants, so that $e(\mc{A}_Z) = e(TX) \pmod 2$ in this case by \autoref{prop:logbundleiso}.
\end{rem}
There is also a more geometric description of the discrepancy. If $\mc{A}_Z$ is oriented and $X$ is orientable, any choice of orientation for $TX$ cannot agree with the orientation on the isomorphism locus $X \setminus Z$ induced by $\mc{A}_Z$ (see Section \ref{sec:logtangentbundle}). After a choice of orientation on $TX$, we can write $X \setminus Z = X_+ \sqcup X_-$, where $X_\pm$ denote the subsets where the orientations from $TX$ and $\mc{A}_Z$ do or do not agree. Then the discrepancy is given by $f(X,Z) = - \langle e(TX), [X_-] \rangle = - \chi(X_-)$.

We now state the obstruction to the existence of an $\mc{A}_Z$-almost-complex structure.
\begin{thm}\label{thm:azalmostcplx} Let $(X,Z)$ be a compact oriented four-dimensional log pair with $Z \in {\rm Hyp}_0(X)$, such that $X$ has an $\mc{A}_Z$-almost-complex structure. Then $c_1^2(\mc{A}_Z) = 3 \sigma(X) + 2 \chi(X) + 4 f(X,Z)$, and $b_2^+(X) + b_1(X) + f(X,Z)$ is odd.
\end{thm}
\bp This is very similar to when $Z = \emptyset$, see \cite[Theorem 1.4.13]{GompfStipsicz99}. As $\mc{A}_Z$ has a complex structure, by \autoref{thm:wuacs} we have $c_1^2(\mc{A}_Z) = p_1(\mc{A}_Z) + 2 e(\mc{A}_Z) = p_1(\mc{A}_Z) + 2 e(TX) + 4 f(X,Z)$. Using \autoref{prop:logcharclass} we get $c_1(\mc{A}_Z) \mod 2 \equiv w_2(\mc{A}_Z) = w_2(TX)$, so that $c_1(\mc{A}_Z)$ is characteristic. By van der Blij's lemma (\cite[Lemma 1.2.20]{GompfStipsicz99}) we obtain that $c_1^2(\mc{A}_Z) \equiv \sigma(X) \pmod 8$. Again using \autoref{prop:logcharclass} we have $p_1(\mc{A}_Z) = p_1(TX)$, the latter integrating to $3 \sigma(X)$ by the Hirzebruch signature theorem. We conclude that $3\sigma(X) + 2\chi(X) + 4 f(X,Z) \equiv \sigma(X) \pmod 8$, hence $\sigma(X) + \chi(X) + 2 f(X,Z) \equiv 0 \pmod 4$, so that $b_1(X) + b_2^+(X) + f(X,Z) \equiv 1 \pmod 2$. 
\ep
Out of this result we can immediately draw the following consequence.
\begin{cor}\label{cor:azsympstr} Let $(X,Z)$ be a compact oriented four-dimensional log pair with an $\mc{A}_Z$-symplectic structure. Then $b_2^+(X) + b_1(X) + f(X,Z)$ is odd.
\end{cor}
We can determine the parity of $f(X,Z)$ for example in the following situation.
\begin{cor}\label{cor:aznosympconn} Let $(X,Z)$ be a compact oriented four-dimensional log pair which is $\mc{A}_Z$-almost-complex, such that $X$ is not almost-complex. Then $f(X,Z)$ is odd, and $(X_- \cup (X_+ \# \C P^2), Z)$ does not admit an $\mc{A}_Z$-symplectic structure.
\end{cor}
\bp If $X$ is not almost-complex, then $b_2^+(X) + b_1(X) \equiv 0 \pmod 2$, while as $(X,Z)$ is $\mc{A}_Z$-almost-complex we obtain from \autoref{thm:azalmostcplx} that $b_2^+(X) + b_1(X) + f(X,Z) \equiv 1 \pmod 2$. We conclude that $f(X,Z) \equiv 1 \pmod 2$. If we perform a connected sum with $\C P^2$ in the subset $X_+$ to form the manifold $X' = X_- \cup (X_+ \# \C P^2)$, we see that $b_2^+(X') = b_2^+(X) + 1$ while $f(X',Z) = f(X,Z)$. Hence then $b_2^+(X') + b_1(X') + f(X',Z) \equiv 0 \pmod 2$, so that by \autoref{cor:azsympstr} we see that $(X',Z)$ does not admit an $\mc{A}_Z$-symplectic structure.
\ep
We finish by giving a simple example of how to apply the above results.
\begin{exa} Let $X = 2 \C P^2 \# \overline{\C P}^2$. Then $(X,Z)$ admits a \blog with $Z = S^1 \times S^2$, see \cite{Cavalcanti17}, and $b_2^+(X) = 2$ and $b_1(X) = 0$. Hence $X$ is not almost-complex while $(X,Z)$ admits an $\mc{A}_Z$-almost-complex structure, so that $f(X,Z)$ must be odd. If we consider the manifold $X' = 3 \C P^2 \# \overline{\C P}^2$ obtained by taking connected sum in $X_+$, then by \autoref{cor:aznosympconn} we conclude that $(X',Z)$ does not admit a \blog.
\end{exa}
\begin{rem} Note that \autoref{cor:azsympstr} is an obstruction for the log pair $(X,Z)$, and not just the manifold $X$ itself. This is unlike the cohomological obstructions due to M{\u{a}}rcu{\c{t}}--Osorno Torres \cite{MarcutOsornoTorres14two} and Cavalcanti \cite{Cavalcanti17} mentioned in Section \ref{sec:logsympstr}. As we make direct use of the Lie algebroid $\mc{A}_Z$ associated to $(X,Z)$, such a dependence on $Z$ seems inevitable for obstructions of this type.
\end{rem}
A similar thing can be done in dimension four for the bundles $\mc{A}_Z^k$, where the parity of $k$ determines whether the signs in determining the Euler class of $\mc{A}_Z^k$ agree with $TX$ or not. Indeed, the signs agree if and only if $k$ is even, so that for $k$ odd we obtain the same obstruction for $\mc{A}_Z^k$ as for $\mc{A}_Z$. On the other hand, the bundles $\mc{A}_Z^k$ for $k$ even admit an almost-complex structure if and only $TX$ does.			
\chapter{Splitting theorems for \texorpdfstring{$\mc{A}$}{A}-Lie algebroids}
\label{chap:splittingtheorems}
In this chapter we present in part ongoing work with Melina Lanius \cite{KlaasseLanius17} on extending results by Burszyn--Lima--Meinrenken \cite{BursztynLimaMeinrenken16} to their $\mc{A}$-counterparts.

More precisely, we establish splitting results for $\mc{A}$-anchored vector bundles and $\mc{A}$-Lie algebroids, in full analogy with \cite{BursztynLimaMeinrenken16}. Our original motivation was to obtain splitting results for $\mc{A}$-Dirac and $\mc{A}$-Poisson structures, which forms the basis for future work. Our strategy of proof will heavily follow that of Bursztyn--Lima--Meinrenken in \cite{BursztynLimaMeinrenken16}. In the interest of brevity we have chosen to assume the reader is familiar with that work. We have tried to stay as close as possible to the notation of \cite{BursztynLimaMeinrenken16}, except with replacing the notation for a manifold $M$ by $X$, and a vector field $X$ by $V$, as is consistent with the rest of this thesis. Unlike the notation we used earlier, in this chapter $N$ will denote a submanifold of $X$. We rely on the functoriality of the constructions in \cite{BursztynLimaMeinrenken16} and show that the ingredients needed to obtain their splitting results can be chosen compatibly in our context. The following is the main result of this chapter, as mentioned in the introduction. It is phrased in imprecise terms (see \autoref{thm:aanchored} and \autoref{thm:aliealgebroid} for precise versions), and is the $\mc{A}$-analogue of \cite[Theorems 3.13 and 4.1]{BursztynLimaMeinrenken16}, to which it reduces when $\mc{A} = TX$.
\newtheorem*{thm:introsplitting}{Theorem \ref{thm:introsplitting}}
\begin{thm:introsplitting} Let $\mc{A} \to X$ be a Lie algebroid and $E_\mc{A}$ either an $\mc{A}$-involutive $\mc{A}$-anchored vector bundle, or an $\mc{A}$-Lie algebroid. Moreover, let $i\colon N \hookrightarrow X$ be an $\mc{A}$-transversal for $E_\mc{A}$. Then $E_\mc{A}$ is isomorphic to its $\mc{A}$-linear approximation in a neighbourhood of $N$.
\end{thm:introsplitting}
\subsection*{Organization of the chapter}
This chapter is built up as follows. In Section \ref{sec:splittingbasic} we recall basic notions from \cite{BursztynLimaMeinrenken16} and define the $\mc{A}$-analogue of involutivity. In Section \ref{sec:atransversals} we introduce the notion of an $\mc{A}$-transversal, along which an $\mc{A}$-anchored vector bundle can be linearized. In Section \ref{sec:aautolifts} we discuss bundle morphisms which preserve the structure that is present, and moreover when vector fields can be lifted to such automorphisms. We finish with Section \ref{sec:splittingtheorems}, which proves the splitting theorems for $\mc{A}$-anchored vector bundles and $\mc{A}$-Lie algebroids mentioned above.
\section{Basic notions}
\label{sec:splittingbasic}
In this section we recall several basic notions that we will need. Let $(X,N)$ be a \emph{pair}, i.e.\ $X$ is a manifold and $N \subseteq X$ is a submanifold. Denote by $\nu_N = \nu(X,N)$ the normal bundle to $N$, and let $\mf{X}(X) = \Gamma(TX)$. Let us recall \autoref{defn:pairsandmaps}.
\begin{defn} A \emph{map of pairs} $f\colon (X',N') \to (X,N)$ is a map $f\colon X' \to X$ such that $f(N') \subseteq N$. It is a \emph{strong map of pairs} if $f^{-1}(N) = N'$.
\end{defn}
Given a map of pairs $f\colon (X',N') \to (X,N)$, note that $Tf(TN') \subseteq TN$, so that $Tf$ induces a map $\nu(f)\colon \nu_{N'} \to \nu_N$ by applying the normal bundle functor.
\begin{prop}\label{prop:strongmapiso} Let $f\colon (X',N') \to (X,N)$ be a strong map of pairs such that $f$ is transverse to $N$. Then $\nu(f)\colon \nu(X',N') \to \nu(X,N)$ is a fiberwise isomorphism.
\end{prop}
Let $V \in \mf{X}(X)$ be a vector field tangent to $N$. Then $V\colon (X,N) \to (TX,TN)$ and $\nu(V)\colon \nu(X,N) \to \nu(TX,TN) \cong T\nu(X,N)$, which we interpret as a vector field $\nu(V)$ on $\nu(X,N)$.
\begin{defn} The vector field $\nu(V) \in \mf{X}(\nu_N)$ is the \emph{linear approximation} of $V$.
\end{defn}
The \emph{tangent lift} $V_T \in \mf{X}(TX)$ of a vector field $V \in \mf{X}(X)$ is obtained by applying the tangent functor $T$ to the map $V\colon X \to TX$, as $V_T = J \circ T V$ where $J$ is the canonical involution on $T(TX)$. If $V$ is tangent to $N$, then $V_T$ is tangent to $TN$, and $\nu(V_T) = \nu(V)_T$.
\subsection{Euler-like vector fields}
Let $N \subseteq X$ be a submanifold and denote by $\mc{E} \in \mf{X}(\nu_N)$ the associated Euler vector field. Note that if $V \in \mf{X}(X)$ vanishes on $N$, it is also tangent to $N$.
\begin{defn}[{\cite[Definition 2.6]{BursztynLimaMeinrenken16}}] A vector field $V \in \mf{X}(X)$ is \emph{Euler-like along $N$} if it is complete, satisfies $V|_N = 0$, and its linear approximation at $N$ satisfies $\nu(V) = \mc{E} \in \mf{X}(\nu_N)$.
\end{defn}
Note that if $V$ is Euler-like along $N$, then $\nu(V_T) = \nu(V)_T = \mc{E}_T \in \mf{X}(TX)$, which is the Euler vector field for $TN$ (see \cite[Example 2.10]{BursztynLimaMeinrenken16}). Hence then $V_T$ is Euler-like along $TN$.
\begin{defn} A \emph{strong tubular neighbourhood embedding} of $N$ into $X$ is an embedding $\Phi\colon (\nu_N,N) \to (X,N)$ with linear approximation $\nu(\Phi)\colon \nu(\nu_N,N) \cong \nu_N \to \nu_N$ the identity.
\end{defn}
Euler-like vector fields give rise to strong tubular neighbourhood embeddings.
\begin{prop}[{\cite[Proposition 2.7]{BursztynLimaMeinrenken16}}] Let $V \in \mf{X}(X)$ be Euler-like along $N$. Then there exists a unique strong tubular neighbourhood embedding $\Phi\colon \nu_N \to X$ such that $\Phi_*(\mc{E}_N) = V$.
\end{prop}
\subsection{$\mc{A}$-involutivity}
\begin{defn} An \emph{anchored vector bundle} $(E,\rho_E)$ is \emph{involutive} if $\rho_E(\Gamma(E))$ is a Lie subalgebra of $\Gamma(TX)$.
\end{defn}
If $E$ is equipped with a map $[\cdot,\cdot]\colon \Gamma(E) \times \Gamma(E) \to \Gamma(E)$ such that $\rho_E([\sigma,\tau]) = [\rho_E(\sigma), \rho_E(\tau)]$ for all $\sigma,\tau \in \Gamma(E)$, then $\rho_E(\Gamma(E))$ is a subalgebra so that $(E,\rho_E)$ is involutive. As a consequence, Lie algebroids are involutive anchored vector bundles. The corresponding $\mc{A}$-analogue for a Lie algebroid $\mc{A} \to X$ naturally presents itself.
\begin{defn} An \emph{$\mc{A}$-anchored vector bundle} $(E_\mc{A}, \varphi_\mc{A})$ is \emph{$\mc{A}$-involutive} if the image $\varphi_{E_\mc{A}}(\Gamma(E_\mc{A}))$ is a Lie subalgebra of $\Gamma(\mc{A})$.
\end{defn}
\begin{prop}\label{prop:ainvolinvol} Let $E_\mc{A}$ be $\mc{A}$-involutive. Then $E_\mc{A}$ is involutive.
\end{prop} 
\bp Note that $\rho_{E_{\mc{A}}}(\Gamma(E_{\mc{A}})) = \rho_{\mc{A}}(\varphi_{E_\mc{A}}(\Gamma(E_\mc{A}))$. Hence, given $\sigma, \tau \in \Gamma(E_\mc{A})$ we have that $[\rho_{E_\mc{A}}(\sigma),\rho_{E_\mc{A}}(\tau)] = \rho_{\mc{A}}([\varphi_{E_\mc{A}}(\sigma),\varphi_{E_\mc{A}}(\tau)]_\mc{A}) = \rho_\mc{A}(\varphi_{E_\mc{A}}(u)) = \rho_{E_\mc{A}}(u)$, for $u \in \Gamma(E_\mc{A})$ such that $[\varphi_{E_\mc{A}}(\sigma),\varphi_{E_\mc{A}}(\tau)]_\mc{A} = \varphi_{E_\mc{A}}(u)$, which exists by $\mc{A}$-involutivity.
\ep
Essentially, the above proposition follows because $\mc{A}$ is always involutive. Note that $\mc{A}$-Lie algebroids are always $\mc{A}$-involutive.
\begin{rem} It is not true in general that if $E_\mc{A}$ is involutive, it is also $\mc{A}$-involutive. While for $\sigma,\tau \in \Gamma(E_\mc{A})$ we have $[\varphi_{E_\mc{A}}(\sigma),\varphi_{E_\mc{A}}(\tau)]_\mc{A} = v$ for some $v \in \Gamma(\mc{A})$, it is not implied that $v = \varphi_{E_\mc{A}}(u)$ for some $u \in \Gamma(E_\mc{A})$. We hence get $[\rho_{E_\mc{A}}(\sigma),\rho_{E_\mc{A}}(\tau)] = \rho_{\mc{A}}([\varphi_{E_\mc{A}}(\sigma),\varphi_{E_\mc{A}}(\tau)]_\mc{A}) = \rho_\mc{A}(v)$, which is not sufficient for involutivity of $E_\mc{A}$.
\end{rem}
\section{\texorpdfstring{$\mc{A}$}{A}-transversals}
\label{sec:atransversals}
Let $\mc{A} \to X$ be a Lie algebroid and $E_{\mc{A}}$ an $\mc{A}$-anchored vector bundle. Note that $E_{\mc{A}}$ is also an anchored vector bundle. Recall that for Lie algebroids $\mc{A}$ we have the notion of a map being transverse to $\mc{A}$, i.e.\ transverse to its anchor $\rho_\mc{A}$. The same definition works for anchored vector bundles $E$. With this we can define the notion of a transversal.
\begin{defn} A submanifold $i\colon N \hookrightarrow X$ is a \emph{transversal} for $E$ if $i$ is transverse to $E$.
\end{defn}
Given a transversal $N$ for $\mc{A}$, we can form the pullback Lie algebroid $i^! \mc{A} \to N$. In this case, the pullback admits a simpler description, as a Lie subalgebroid of $\mc{A}$ supported on $N$.
\begin{prop}\label{prop:transversal} Let $i\colon N \to X$ be a transversal for $\mc{A}$. Then $i^! \mc{A} = \rho_\mc{A}^{-1}(TN) \subseteq \mc{A}|_N$.
\end{prop}
\bp By definition $i^! \mc{A}$ consists of all pairs $(v,V) \in f^*\mc{A} \times TN$ for which $\rho_\mc{A}(v) = Ti(V)$. As $Ti$ is injective, the condition on $v$ is precisely that $\rho_\mc{A}(v) \in TN$.
\ep
We next extend the notion of transversality to that of $\mc{A}$-transversality.
\begin{defn} Let $(\varphi,f)\colon (N,\mc{B}) \to (X,\mc{A})$ be a Lie algebroid morphism. Then $\varphi$ is \emph{$\mc{A}$-transverse} to $E_{\mc{A}}$ if $\varphi$ is transverse to $\varphi_{E_\mc{A}}$.
\end{defn}
Note that $f$ is not necessarily transverse to $\rho_{E_\mc{A}}$ as well. To define the appropriate notion of transversal for $\mc{A}$-anchored vector bundles, let $N$ be a transversal for $\mc{A}$ and consider the induced Lie subalgebroid $\mc{B} = i^! \mc{A}$ of $\mc{A}$ supported on $N$, with $(\iota,i)\colon i^! \mc{A} \to \mc{A}$ the associated Lie algebroid morphism.
\begin{defn} Let $i \colon N \to X$ be a transversal for $\mc{A}$. Then $N$ is an \emph{$\mc{A}$-transversal} for $E_\mc{A}$ if $\iota$ is $\mc{A}$-transverse to $E_\mc{A}$.
\end{defn}
Note that it follows that every $\mc{A}$-transversal for $E_{\mc{A}}$ is also a transversal for $E_{\mc{A}}$.
\begin{prop}\label{prop:atranstrans} Let $(\iota,i)\colon (N, i^! \mc{A}) \to (X,\mc{A})$ be an $\mc{A}$-transversal for $E_\mc{A}$. Then $N$ is a transversal for $E_\mc{A}$.
\end{prop}
\bp As $N$ is a transversal for $\mc{A}$, we have over $f(N) \subseteq X$ that $Ti(TN) + \rho_\mc{A}(\mc{A}) = TX$. Moreover, $\iota$ is transverse to $\varphi_{E_\mc{A}}$, so that $\iota(i^! \mc{A}) + \varphi_{E_\mc{A}}(E_\mc{A}) = \mc{A}$. Applying the anchor $\rho_\mc{A}$, we see that $\rho_\mc{A}(\iota(i^! \mc{A})) + \rho_\mc{A}(\varphi_{E_\mc{A}}(E_\mc{A})) = \rho_\mc{A}(\mc{A})$, so that by \autoref{prop:transversal} we have $\rho_\mc{A}(\rho_\mc{A}^{-1}(TN)) + \rho_{E_\mc{A}}(E_\mc{A}) = \rho_\mc{A}(\mc{A})$. But $\rho_\mc{A}(\rho_\mc{A}^{-1}(TN)) \subseteq TN$, so that $Ti(TN) + \rho_{E_\mc{A}}(E_\mc{A}) = Ti(TN) + \rho_\mc{A}(\rho_\mc{A}^{-1}(TN)) + \rho_{E_\mc{A}}(E_\mc{A}) = Ti(TN) + \rho_\mc{A}(\mc{A}) = TX$ as desired.
\ep
Let $\mc{A} \to X$ be a Lie algebroid and $(\varphi,f)\colon (N,\mc{B}) \to (X,\mc{A})$ be $\mc{A}$-transverse to an $\mc{A}$-anchored vector bundle $E_{\mc{A}}$. We can then define the $\mc{A}$-pullback of $E_\mc{A}$ along $(\varphi,f)$, which will be a $\mc{B}$-anchored vector bundle on $N$.
\begin{defn} The \emph{$\mc{A}$-pullback} of $E_{\mc{A}}$ along $(\varphi,f)$ is the $\mc{B}$-anchored vector bundle $\varphi^! E_\mc{A} \to N$ defined as the following fiber product, with $\mc{B}$-anchor the projection to $\mc{B}$,
	\begin{equation*}
		\varphi^! E_\mc{A} = \{(V,v)\in f^* E_{\mc{A}} \times \mc{B} \, | \, \varphi_{E_\mc{A}}(V) = \varphi(v)\}.
	\end{equation*}
\end{defn}
If $E_\mc{A}$ is in addition an $\mc{A}$-Lie algebroid, the pullback $\varphi^! E_\mc{A}$ inherits a unique $\mc{B}$-Lie algebroid structure, with the diagonal map $\varphi^! E_\mc{A} \to f^* E_\mc{A} \times \mc{B}$ the inclusion as a Lie subalgebroid.
\begin{rem} In the special case where $N = X \times Q$ and $\mc{B} = \mc{A} \oplus \mc{C}$ with $(\varphi,f)\colon (N,\mc{B}) \to (X,\mc{A})$ the obvious projections, we have $\varphi^! E_\mc{A} = E_\mc{A} \times \mc{C}$.
\end{rem}
For the next lemma, see \cite{HigginsMackenzie90} regarding direct products of Lie algebroids.
\begin{lem}\label{lem:directprodinvol} Let $E_\mc{A} \to X$ be an $\mc{A}$-involutive $\mc{A}$-anchored vector bundle and $F_\mc{B} \to N$ a $\mc{B}$-involutive $\mc{B}$-anchored vector bundle, and let $\mc{C} := \mc{A} \oplus \mc{B} \to X \times N$. Then $E_\mc{A} \oplus F_\mc{B}$ is a $\mc{C}$-involutive $\mc{C}$-anchored vector bundle.
\end{lem}
\bp As $\rho_\mc{C} = (\rho_\mc{A}, \rho_\mc{B})$, the bundle $E_\mc{A} \oplus F_\mc{B}$ is $\mc{C}$-anchored with anchor $(\rho_{E_\mc{A}}, \rho_{F_\mc{B}})$ and $\mc{C}$-anchor $\varphi_\mathcal{C}=(\varphi_{E_\mathcal{A}}, \varphi_{F_\mc{B}})$, i.e. $E_\mc{A} \oplus F_\mc{B} \stackrel{\varphi_\mc{C}}{\to} \mc{A} \oplus \mc{B} \stackrel{\rho_\mc{C}}{\to} T(X \times N) \cong TX \oplus TN$. The bracket $[\cdot,\cdot]_\mathcal{C}$ on the Lie algebroid $\mathcal{C}$ is generated as follows. Let $\sigma, \tau \in \Gamma(\mc{C})$ be decomposable, i.e.\ $\sigma = g v \oplus i v'$ and $\tau = h w \oplus j w'$ for $g,h,i,j \in \mathcal{C}^\infty(X\times N )$, $v, w \in\Gamma(\mathcal{A})$ and $v', w' \in \Gamma(\mathcal{B})$. Then we have:
\begin{align*}
[\sigma,\tau]_\mc{C} &= \left(gh[v,w]_\mc{A} + (\mc{L}_{\rho_\mc{C}(\sigma)} h) \cdot w - (\mc{L}_{\rho_\mc{C}(\tau)} g) \cdot v,\right.\\
\qquad&\qquad\qquad \left.ij[v',w']_\mc{B} + (\mc{L}_{\rho_\mc{C}(\sigma)} j) \cdot w' - (\mc{L}_{\rho_\mc{C}(\tau)} i) \cdot v'\right).
\end{align*}
Here by e.g.\ $\mc{L}_{\rho_\mc{C}(\sigma)} h$ we mean that $\rho_\mc{C}(\sigma) = g \rho_\mc{A}(v) \oplus i \rho_\mc{B}(v')$ acts on $h \in C^\infty(X \times N)$ by pulling back the vector fields $\rho_\mc{A}(v)$ and $\rho_\mc{B}(v')$ to $X \times N$ along the projections to $X$ and $N$ respectively. Note that for a section  $u \in \Gamma(E_\mc{A} \oplus F_\mc{B})$ of the form $u = (e,f)$ we have $\varphi_\mc{C}(u) = (\varphi_{E_\mathcal{A}}(e),\varphi_{F_\mathcal{B}}(f))$. As these sections appropriately generate $\Gamma(E_\mc{A} \oplus F_\mc{B})$, to show $\mathcal{C}$-involutivity of $E_\mathcal{A}\oplus F_\mathcal{B}$, it suffices to check the case when we have $u, v \in \Gamma(E_\mc{A} \oplus F_\mc{B})$ of this form, with $[\varphi_\mathcal{C}(u),\varphi_{\mathcal{C}}(v)]_\mc{C} = [(g\varphi_{E_\mathcal{A}}(e),i\varphi_{F_\mathcal{B}}(f)), (h\varphi_{E_\mathcal{A}}(e'),j\varphi_{F_\mathcal{B}}(f'))]_\mc{C}$. As $E_\mc{A}$ is $\mc{A}$-involutive we have that $[\varphi_{E_\mathcal{A}}(e),\varphi_{E_\mathcal{A}}(e')]_\mc{A} \in \varphi_{E_\mathcal{A}}( \Gamma(E_\mathcal{A}))$, and similarly because $F_\mathcal{B}$ is $\mathcal{B}$-involutive also $[\varphi_{F_\mathcal{B}}(f),\varphi_{F_\mathcal{B}}(f')]_\mc{B} \in \varphi_{F_\mathcal{B}}(\Gamma(F_\mathcal{B}))$. By the description of $[\cdot,\cdot]_\mc{C}$ above we see that $[\varphi_\mathcal{C}(u),\varphi_{\mathcal{C}}(v)]\in \varphi_\mathcal{C}(\Gamma(E_\mathcal{A}\oplus F_\mathcal{B}))$, so that $E_\mc{A} \oplus F_\mc{B}$ is $\mc{C}$-involutive as desired.
\ep
Note that $\mc{A}$-involutivity is preserved under $\mc{A}$-pullback along $\mc{A}$-transverse maps.
\begin{prop}\label{prop:ainvolatransverse} Let $(\varphi,f)\colon (N,\mc{B}) \to (X,\mc{A})$ be an $\mc{A}$-transverse map to an $\mc{A}$-anchored vector bundle $E_\mc{A} \to X$. If $E_\mc{A}$ is $\mc{A}$-involutive, then $\varphi^! E_\mc{A}$ is $\mc{B}$-involutive.
\end{prop}
\bp This is similar to \cite[Proposition 3.12]{BursztynLimaMeinrenken16}. Let ${\rm Gr}(f) = \{(x, f(x)) \in M \times N \, | \, x \in X\} \subseteq X \times N$ denote the graph of $f$. We begin by identifying $\varphi^!E_{\mathcal{A}}$ with a subbundle of $E_\mathcal{A}\times \mathcal{B} \to X \times N$ over ${\rm Gr}(f)$. Recall that $\varphi^!E_{\mathcal{A}} = \{(V,v)\in E_{\mathcal{A}}\times \mathcal{B} \, | \, \varphi_{E_\mathcal{A}}(V)=\varphi(v)\}$. Let $\mc{C} := \mc{A} \times \mc{B} \to X \times N$ be the direct product of Lie algebroids. The bundle $E_\mc{A} \times \mc{B} \to X \times N$ is $\mc{C}$-anchored with $\mc{C}$-anchor $\Phi = (\varphi_{E_\mc{A}}, {\rm id})$. Let $F$ be the subbundle of $(E_{\mathcal{A}}\times \mathcal{B})|_{{\rm Gr}(f)}$ consisting of all elements whose image under $\Phi$ lies in $\mc{F} := \{(Y,v)\in\mathcal{A}\times \mathcal{B} \, | \, Y = \varphi(v)\} \subseteq \mc{A} \times \mc{B}$. This condition ensures that $F$ consists of all elements $(V,v)$ such that $\varphi_{E_\mathcal{A}}(V)=\varphi(v)$, which is precisely the condition defining $\varphi^!E_\mathcal{A}$, so that $\varphi^! E_\mc{A}$ and $F$ can be identified.

Next, let $\sigma_1, \sigma_2 \in \Gamma(F)$ be given and extend them to sections $\tau_1,\tau_2$ of $E_{\mathcal{A}}\times \mathcal{B}$, using that ${\rm Gr}(f)$ is a closed submanifold of $X\times N$. Note that $[\Phi(\tau_1), \Phi(\tau_2)]_\mc{C}$ restricts over ${\rm Gr}(f)$ to lie in $\mc{F}$: the bracket $[\cdot,\cdot]_\mc{C}$ is given for $\tau_1 = g v \oplus i v'$ and $\tau_2 = h w \oplus j w'$ as in \autoref{lem:directprodinvol} by
\begin{align*}
	[\tau_1,\tau_2]_\mc{C} &= \left(gh[v,w]_\mc{A} + (\mc{L}_{\rho_\mc{C}(\tau_1)} h) \cdot w - (\mc{L}_{\rho_\mc{C}(\tau_2)} g) \cdot v,\right.\\
	\qquad&\qquad\qquad \left.ij[v',w']_\mc{B} + (\mc{L}_{\rho_\mc{C}(\tau_1)} j) \cdot w' - (\mc{L}_{\rho_\mc{C}(\tau_2)} i) \cdot v'\right).
\end{align*}
Thus it suffices to check the case when $[\Phi(\tau_1),\Phi(\tau_2)]_\mc{C} = [(g v,iv'), (h w,j w')]_\mc{C}$. In this case, since $E_\mathcal{A}$ is $\mathcal{A}$-involutive, we obtain $[v, w]_\mc{A} \in \varphi_{E_\mathcal{A}}(\Gamma(E_\mathcal{A}))$ and consequently $[\Phi(\tau_1),\Phi(\tau_2)]_\mc{C}|_{{\rm Gr}(f)}$ is a section of $\mc{F}$ because $\varphi$ is a Lie algebroid morphism. Note now that by \autoref{lem:directprodinvol}, the bundle $E_\mc{A} \times \mc{B}$ is $\mc{C}$-involutive. Hence $[\Phi(\tau_1),\Phi(\tau_2)]_\mc{C}$ lifts to a section $\tau\in\Gamma(E_{\mathcal{A}}\times \mathcal{B})$ such that $\Phi(\tau) = [\Phi(\tau_1),\Phi(\tau_2)]_\mc{C}$, and we just showed that $\Phi(\tau)|_{{\rm Gr}(f)} \in \Gamma(\mc{F})$. We conclude that the restriction $\sigma := \tau|_{{\rm Gr}(f)}$ is a section of $F$ satisfying $[\Phi(\sigma_1),\Phi(\sigma_2)] = \Phi(\sigma)$. Hence $\Phi(\Gamma(F)) \subseteq F$ is a Lie subalgebra and $\varphi^! E_\mc{A}$ is $\mc{B}$-involutive.
\ep
We mainly use the $\mc{A}$-pullback construction for $\mc{A}$-transversals $i\colon N \to X$ for $E_\mc{A}$, obtaining the $i^!\mc{A}$-anchored vector bundle $\iota^! E_\mc{A}$. As in \autoref{prop:transversal}, in this case the $\mc{A}$-pullback admits a simpler description. Note that by \autoref{prop:atranstrans}, $N$ is also a transversal for $E_\mc{A}$.
\begin{prop}\label{prop:iiotaea} Let $N$ be an $\mc{A}$-transversal for $E_\mc{A}$. Then $\iota^! E_\mc{A} = i^! E_\mc{A} = \rho_{E_\mc{A}}^{-1}(TN)$.
\end{prop}
\bp The Lie algebroid $\mc{B} \to N$ that is used is given by $i^! \mc{A} = \rho_\mc{A}^{-1}(TN)$. Consequently, as $\iota^! E_\mc{A}$ consists of all pairs $(V,v) \in f^* E_\mc{A} \times i^! \mc{A}$ for which $\varphi_{E_\mc{A}}(V) = \iota(v)$, the condition on $V$ is that $\varphi_{E_\mc{A}}(V) \in \rho_\mc{A}^{-1}(TN)$. But this is precisely that $V \in \varphi_{E_\mc{A}}^{-1}(\rho_\mc{A}^{-1}(TN)) = \rho_{E_\mc{A}}^{-1}(TN)$.
\ep
In view of the above proposition, we will always write $i^! E_{\mc{A}}$ instead of $\iota^! E_{\mc{A}}$ from now on. In conclusion, $\mc{A}$-transversals $i\colon N \to X$ for $E_\mc{A}$ are both transversals for $\mc{A}$ and for $E_\mc{A}$, and allow us to obtain an $i^! \mc{A}$-anchored vector bundle $i^! E_\mc{A}$ over $N$.
We now define the bundle which will provide the approximation of an $\mc{A}$-anchored vector bundle $E_\mc{A}$ along an $\mc{A}$-transversal $N$.
\begin{defn} Let $i\colon N \to X$ be an $\mc{A}$-transversal for an $\mc{A}$-anchored vector bundle $E_\mc{A}$. The \emph{$\mc{A}$-linear approximation} of $E_\mc{A}$ along $N$ is the $p^! i^! \mc{A}$-anchored vector bundle $p^! i^! E_\mc{A}$.
\end{defn}
If $E_\mc{A}$ is moreover an $\mc{A}$-Lie algebroid, the $\mc{A}$-linear approximation will be a $p^! i^! \mc{A}$-Lie algebroid. If $E_\mc{A}$ is an $\mc{A}$-Dirac structure relative to $\eta_\mc{A}$, then $p^! i^! E_\mc{A}$ is an $p^! i^! \mc{A}$-Dirac structure relative to $p^* i^* \eta_\mc{A}$ (note our abuse of notation in using e.g.\ $i$ instead of $\iota$, c.f.\ \autoref{prop:iiotaea}). We thus have a diagram of fiber products.
\begin{center}
	\begin{tikzpicture}
	\matrix (m) [matrix of math nodes, row sep=2.5em, column sep=2.5em,text height=1.5ex, text depth=0.25ex]
	{	p^! i^! E_\mc{A} & i^! E_\mc{A} & E_\mc{A} \\ p^! i^! \mc{A} & i^! \mc{A} & \mc{A} \\ T\nu(X,N) & TN & TX
		\\};
	\path[-stealth]
	(m-1-1) edge (m-1-2)
	(m-1-2) edge (m-1-3)
	(m-1-1) edge (m-2-1)
	(m-1-2) edge (m-2-2)
	(m-1-3) edge (m-2-3)
	(m-2-1) edge (m-2-2)
	(m-2-2) edge (m-2-3)
	(m-2-1) edge (m-3-1)
	(m-2-2) edge (m-3-2)
	(m-2-3) edge (m-3-3)
	(m-3-1) edge (m-3-2)
	(m-3-2) edge (m-3-3)
	;
	\end{tikzpicture}
\end{center}
\begin{prop} Let $i\colon N \to X$ be an $\mc{A}$-transversal for an $\mc{A}$-involutive $\mc{A}$-anchored vector bundle $E_\mc{A}$. Then $p^! i^! E_\mc{A}$ is $p^! i^! \mc{A}$-involutive.
\end{prop}
\bp Follows from \autoref{prop:ainvolatransverse}, as $i$ is $\mc{A}$-transverse and similarly for $p$.
\ep
Note that $\nu(\mc{A}, i^! \mc{A})$ is an anchored vector bundle over the manifold $\nu(X,N)$, with anchor $\nu(\rho_\mc{A})\colon \nu(\mc{A},i^! \mc{A}) \to \nu(TX,TN) \cong T \nu(X,N)$.
\begin{prop} Let $i\colon N \to X$ be a transversal for $\mc{A}$. Then the anchor map $\nu(\rho_\mc{A})\colon \nu(\mc{A},i^!\mc{A}) \to \nu(TX,TN)$ is a fiberwise isomorphism.
\end{prop}
\bp Note that $\rho_\mc{A}\colon (\mc{A},i^! \mc{A}) \to (TX,TN)$ is a strong map of pairs by \autoref{prop:transversal}. Because $N$ is a transversal for $\mc{A}$, we have $Ti(TN) + \rho_\mc{A}(\mc{A}) = TX$. But then $\rho_\mc{A}$ is transverse to $TN$, so that the result follows from \autoref{prop:strongmapiso}.
\ep
A similar result holds for $\mc{A}$-transversals to $\mc{A}$-anchored vector bundles.
\begin{prop} Let $i\colon N \to X$ be an $\mc{A}$-transversal for an $\mc{A}$-anchored vector bundle $E_\mc{A}$. Then $\nu(\varphi_{E_\mc{A}})\colon \nu(E_\mc{A},i^! E_\mc{A}) \to \nu(\mc{A},i^!\mc{A})$ is a fiberwise isomorphism.
\end{prop}
The situation for these bundles is thus summarized by the following diagram.
\begin{center}
	\begin{tikzpicture}
	\matrix (m) [matrix of math nodes, row sep=2.5em, column sep=2.5em,text height=1.5ex, text depth=0.25ex]
	{	\nu(E_\mc{A}, i^! E_\mc{A}) & i^! E_\mc{A} \\ \nu(\mc{A}, i^! \mc{A}) & i^! \mc{A} \\ \nu(TX, TN) & TN
		\\};
	\path[-stealth]
	(m-1-1) edge node[above] {$p_{i^! E_\mc{A}}$} (m-1-2)

	(m-1-1) edge node[left] {$\nu(\varphi_{E_\mc{A}})$} node[right] {$\cong$} (m-2-1)
	(m-1-2) edge node[left] {$\varphi_{E_\mc{A}}$} (m-2-2)

	(m-2-1) edge node[above] {$p_{i^! \mc{A}}$} (m-2-2)

	(m-2-1) edge node[left] {$\nu(\rho_\mc{A})$} node[right] {$\cong$} (m-3-1)
	(m-2-2) edge node[left] {$\rho_\mc{A}$} (m-3-2)

	(m-3-1) edge node[above] {$p_{TN}$} (m-3-2)
	;
	\draw [->] (m-1-2) [bend left] edge node [right] {$\rho_{E_\mc{A}}$} (m-3-2)
	(m-1-1) [bend right=80] edge node [left] {$\nu(\rho_{E_\mc{A}})$} node[right] {$\cong$} (m-3-1)
	;
	\end{tikzpicture}
\end{center}
As in \cite[Lemma 3.8]{BursztynLimaMeinrenken16}, the $\mc{A}$-linear approximation gets its name from the following.
\begin{prop} There is a canonical isomorphism of $p^! i^! \mc{A}$- and $\nu(\mc{A},i^! \mc{A})$-anchored vector bundles $p^! i^! E_\mc{A} \cong \nu(E_\mc{A}, i^! E_\mc{A})$.
\end{prop}
\bp By \cite[Lemma 3.8]{BursztynLimaMeinrenken16}, there is a canonical isomorphism \emph{as anchored vector bundle} given by ${(p_{i^! E_\mc{A}}, \nu(\rho_{E_\mc{A}}))\colon \nu(E_\mc{A},i^! E_\mc{A}) \stackrel{\cong}{\to} p^! i^! E_\mc{A}}$ and similarly we have $(p_{i^! \mc{A}}, \nu(\rho_{\mc{A}}))\colon \nu(\mc{A}, i^! \mc{A}) \stackrel{\cong}{\to} p^! i^! \mc{A}$. As applying $\nu$ is functorial, we get that $p^! i^! E_\mc{A} \cong \nu(E_\mc{A}, i^! E_\mc{A})$ as $p^! i^! \mc{A}$- and $\nu(\mc{A},i^! \mc{A})$-anchored vector bundles as well.
\ep

\section{\texorpdfstring{$\mc{A}$}{A}-automorphisms and lifts}
\label{sec:aautolifts}
In this section we discuss the required automorphisms for $\mc{A}$-bundles and show that $\mc{A}$-involutive $\mc{A}$-anchored vector bundles admit an appropriate $\mc{A}$-lift of their anchor.

Let ${\rm pr}_E\colon E \to X$ be a vector bundle and denote by $\mf{aut}(E) \subseteq \mf{X}(E)$ the Lie algebra of infinitesimal bundle automorphisms. Any infinitesimal automorphism $\wt{V} \in \mf{aut}(E)$ restricts to a vector field $V \in \mf{X}(X)$ such that $\wt{V} \sim_{{\rm pr}_E} V$. We say that $\wt{V}$ \emph{lifts} $V$. Given $\wt{V} \in \mf{aut}(E)$ there is an associated operator $D\colon \Gamma(E) \to \Gamma(E)$ such that $D(f \sigma) = f D(\sigma) + (\mc{L}_V f) \cdot \sigma$ for all $f \in C^\infty(X)$, $\sigma \in \Gamma(E)$. Conversely, any linear operator $D$ with this property corresponds to a unique $\wt{V} \in \mf{aut}(E)$ lifting $V$.

Recall from \cite{BursztynLimaMeinrenken16} that for an anchored vector bundle $E \to X$, its Lie algebra of infinitesimal automorphism $\mf{aut}_{\rm AV}(E) \subseteq \mf{aut}(E)$ consists of all $\wt{V} \in \mf{aut}(E)$ such that $\wt{V} \sim_{\rho_E} V_T$. Equivalently, the operator $D\colon \Gamma(E) \to \Gamma(E)$ corresponding to $\wt{V}$ satisfies $\rho_E(D(\tau)) = [V, \rho_E(\tau)]$ for all $ \tau \in \Gamma(E)$.
\begin{defn} An element $V \in \mf{X}(X)$ is \emph{liftable} to $\mf{aut}_{\rm AV}(E)$ if there exists an element $\wt{V} \in \mf{aut}_{\rm AV}(E)$ called the \emph{lift} such that $\wt{V} \sim_{{\rm pr}_E} V$. We say $\rho_E$ \emph{admits a lift} if there exists a map $\wt{\rho_E}\colon \Gamma(E) \to \mf{aut}_{\rm AV}(E)$ called its \emph{lift} such that ${\rm pr}_E \circ \wt{\rho}_E = \rho_E$.
\end{defn}
By \cite[Proposition 3.17]{BursztynLimaMeinrenken16}, the existence of a lift for $\rho_E$ is equivalent to involutivity of $E$. If $E$ further carries a compatible bracket $[\cdot,\cdot]_E$ making it into a Lie algebroid, the associated Lie algebra of infinitesimal automorphisms $\mf{aut}_{\rm LA}(E) \subseteq \mf{aut}_{\rm AV}(E)$ consists of all $\wt{V} \in \mf{aut}_{\rm AV}(E)$ such that the corresponding operator $D\colon \Gamma(E) \to \Gamma(E)$ is a derivation of the Lie bracket, i.e.\ $D([\sigma,\tau]_E) = [D\sigma,\tau]_E + [\sigma, D \tau]_E$ for all $\sigma,\tau \in \Gamma(E)$. A Lie algebroid $E$ has a \emph{canonical lift} $\overline{\rho}_E\colon \Gamma(E) \to \mf{aut}_{\rm LA}(E)$ of $\rho_E$, given by the operators $D_\sigma = [\sigma, \cdot]_E$ for $\sigma \in \Gamma(E)$.
\begin{rem} If $E = TX$, the canonical lift of $\rho_E = {\rm id}_{TX}$ is the tangent lift.
\end{rem}
\subsection{$\mc{A}$-anchored vector bundles}
Let $\mc{A} \to X$ be a Lie algebroid and $E_\mc{A}$ an $\mc{A}$-anchored vector bundle. Let $V \in \mf{X}(X)$ be given, and let $\sigma_{\mc{A}} \in \Gamma(\mc{A})$ be such that $\rho_\mc{A}(\sigma_\mc{A}) = V$. Then $V$ has a canonical lift $V_\mc{A} = \overline{\rho}_{\mc{A}}(\sigma_{\mc{A}}) \in \mf{aut}_{\rm LA}(\mc{A})$ using the canonical lift of $\mc{A}$. Denote by $\mf{aut}_{\mc{A}-{\rm AV}}(E_{\mc{A}})$ the Lie algebra of infinitesimal $\mc{A}$-anchored vector bundle automorphisms of $E_{\mc{A}}$, consisting of all $\wt{V} \in \mf{aut}_{\rm AV}(E_\mc{A})$ such that $\wt{V} \sim_{\varphi_{E_\mc{A}}} V_{\mc{A}} \in \mf{aut}_{\rm LA}(\mc{A})$, where ${\rm pr}_{E_\mc{A}}(\wt{V}) = V$. In terms of operators, for $\sigma \in \Gamma(E_\mc{A})$, $D_\sigma$ should satisfy $\varphi_{E_{\mc{A}}}(D_\sigma(\tau)) = D_{\sigma_{\mc{A}}}(\varphi_{E_{\mc{A}}}(\tau)) = [\sigma_\mc{A}, \varphi_{E_\mc{A}}(\tau)]_{\mc{A}}$ for all $\tau \in \Gamma(E_{\mc{A}})$, where $\sigma_\mc{A} = \varphi_{E_\mc{A}}(\sigma) \in \Gamma(\mc{A})$.

We now define the concept of liftability for $\mc{A}$-anchored vector bundles.
\begin{defn} An element $V \in \mf{X}(M)$ is \emph{$\mc{A}$-liftable} to $\mf{aut}_{\mc{A}-{\rm AV}}(E_\mc{A})$ if there exists an element $\wt{V} \in \mf{aut}_{\mc{A}-{\rm AV}}(E_\mc{A})$ called the \emph{$\mc{A}$-lift} such that $\wt{V} \sim_{{\rm pr}_{E_\mc{A}}} V$. We say $\rho_{E_\mc{A}}$ \emph{admits an $\mc{A}$-lift} if there exists an \emph{$\mc{A}$-lift} $\wt{\rho}_{E_\mc{A}}\colon \Gamma(E_\mc{A}) \to \mf{aut}_{\mc{A}-{\rm AV}}(E_\mc{A})$ such that ${\rm pr}_{E_\mc{A}} \circ \wt{\rho}_{E_\mc{A}} = \rho_{E_\mc{A}}$.
\end{defn}
Hence, an $\mc{A}$-lift $\wt{\rho}_{E_\mc{A}}$ is a lift of $\rho_{E_\mc{A}}$ that is compatible with the canonical lift $\overline{\rho}_\mc{A}$ of $\rho_\mc{A}$. By \autoref{prop:ainvolinvol} and \cite[Proposition 3.17]{BursztynLimaMeinrenken16}, an $\mc{A}$-involutive $\mc{A}$-anchored vector bundle admits a lift of $\rho_{E_\mc{A}}$. We can characterize when an $\mc{A}$-lift exists.
\begin{prop}\label{prop:aanchoredalift} Let $E_\mc{A}$ be an $\mc{A}$-anchored vector bundle. Then there exists an $\mc{A}$-lift $\wt{\rho}_{E_\mc{A}}$ of $\rho_{E_\mc{A}}$ if and only if $E_\mc{A}$ is $\mc{A}$-involutive.
\end{prop}
Our proof of the above proposition is based on \cite[Proposition 3.17]{BursztynLimaMeinrenken16} and uses a specific type of connection. Recall that a \emph{$\rho_E$-connection} for an anchored vector bundle $E \to X$ is a bilinear map $\nabla\colon \Gamma(E) \times \Gamma(E) \to \Gamma(E)$, $(\sigma,\tau) \mapsto \nabla_\sigma(\tau)$, which is $C^\infty(X)$-linear in $\sigma$ and satisfies $\nabla_\sigma(f \tau) = f \nabla_\sigma(\tau) + (\rho_E(\sigma) f) \tau$, for all $\sigma, \tau \in \Gamma(E)$, $f \in C^\infty(X)$. Any connection $\nabla'$ on $E$ defines a $\rho_E$-connection $\nabla$ via $\nabla_\sigma := \nabla'_{\rho_E(\sigma)}$ for $\sigma \in \Gamma(E)$, so that $\rho_E$-connections always exist. Any $\rho_{E}$-connection $\nabla$ has a \emph{torsion tensor} $T_\nabla \in \Gamma(\wedge^2 E^* \otimes TX)$ given by $T_\nabla(\sigma,\tau) = \rho_E(\nabla_\sigma \tau) - \rho_E(\nabla_\tau \sigma) - [\rho_E(\sigma), \rho_E(\tau)]$ for $\sigma,\tau \in \Gamma(E)$. We say $\nabla$ is \emph{torsion-free} if $T_\nabla \equiv 0$. Given an $\mc{A}$-anchored vector bundle $E_\mc{A}$, a $\rho_{E_\mc{A}}$-connection further has an \emph{$\mc{A}$-torsion tensor} $T_\nabla^\mc{A} \in \Gamma(\wedge^2 E_\mc{A}^* \otimes \mc{A})$ given by $T_\nabla^\mc{A}(\sigma,\tau) = \varphi_{E_\mc{A}}(\nabla_\sigma \tau) - \varphi_{E_\mc{A}}(\nabla_\tau \sigma) - [\varphi_{E_\mc{A}}(\sigma), \varphi_{E_\mc{A}}(\tau)]_\mc{A}$ for $\sigma,\tau \in \Gamma(E_\mc{A})$. We then say $\nabla$ is \emph{$\mc{A}$-torsion-free} if $T_\nabla^\mc{A} \equiv 0$. It is immediate that $\rho_\mc{A} \circ T_\nabla^\mc{A} = T_\nabla$, as $\rho_{E_\mc{A}} = \rho_\mc{A} \circ \varphi_{E_\mc{A}}$. This implies that an $\mc{A}$-torsion-free $\rho_{E_\mc{A}}$-connection is also torsion-free.
\bp[ of \autoref{prop:aanchoredalift}] Let $\nabla$ be a $\rho_{E_\mc{A}}$-connection on $E_\mc{A}$ and assume that $E_\mc{A}$ is $\mc{A}$-involutive. Similarly to \cite[Proposition 3.17]{BursztynLimaMeinrenken16}, we can then lift its $\mc{A}$-torsion tensor $T_\nabla^\mc{A}$ to a tensor $S_\nabla \in \Gamma(\wedge^2 E_\mc{A}^* \otimes E_\mc{A})$, so that $\varphi_{E_\mc{A}} \circ S_\nabla = T_\nabla^\mc{A}$. This will also be a lift of the torsion tensor $T_\nabla$, as $\rho_{E_\mc{A}} \circ S_\nabla = \rho_\mc{A} \circ (\varphi_{E_\mc{A}} \circ S_\nabla) = \rho_\mc{A} \circ T_\nabla^\mc{A} = T_\nabla$.  Define a new $\rho_{E_\mc{A}}$-connection $\overline{\nabla}$ by $\overline{\nabla}_\sigma \tau = \nabla_\sigma \tau - \frac12 S_\nabla(\sigma,\tau)$ for $\sigma,\tau \in \Gamma(E_\mc{A})$. Then $\overline{\nabla}$ is torsion-free, $T_{\overline{\nabla}} \equiv 0$, but crucially it is also $\mc{A}$-torsion-free, $T_{\overline{\nabla}}^\mc{A} \equiv 0$. The operators $D_\sigma\colon \Gamma(E_\mc{A}) \to \Gamma(E_\mc{A})$ for $\sigma \in \Gamma(E_\mc{A})$ given by $D_\sigma \tau = \overline{\nabla}_\sigma \tau - \overline{\nabla}_\tau \sigma$ for $\tau \in \Gamma(E_\mc{A})$ defines a lift $\wt{\rho}_{E_\mc{A}}$ of $\rho_{E_\mc{A}}$ as in \cite[Proposition 3.17]{BursztynLimaMeinrenken16}. However, as the $\mc{A}$-torsion tensor of $\overline{\nabla}$ vanishes, the operators $D_\sigma$ satisfy $\varphi_{E_\mc{A}}(D_\sigma \tau) = [\varphi_{E_\mc{A}}(\sigma), \varphi_{E_\mc{A}}(\tau)]_\mc{A}$, so that $\wt{\rho}_{E_\mc{A}}$ is an $\mc{A}$-lift for $\rho_{E_\mc{A}}$. Finally, the operators $D_\sigma$ satisfy the additional property that $D_{f\sigma} \tau = f D_\sigma \tau - (\rho_{E_\mc{A}}(\tau) f) \sigma$ for all $f \in C^\infty(X)$. Conversely, if $\wt{\rho}_{E_\mc{A}}$ is an $\mc{A}$-lift of $\rho_{E_\mc{A}}$, then the corresponding operators $D_\sigma$ satisfy $\varphi_{E_\mc{A}}(D_\sigma \tau) = [\varphi_{E_\mc{A}}(\sigma), \varphi_{E_\mc{A}}(\tau)]_\mc{A}$ for all $\sigma,\tau \in \Gamma(E_\mc{A})$, showing $\mc{A}$-involutivity of $E_\mc{A}$.
\ep
\subsection{$\mc{A}$-Lie algebroids} Let $E_\mc{A}$ be an $\mc{A}$-Lie algebroid. Then the appropriate infinitesimal automorphisms of $E_\mc{A}$ form the Lie algebra $\mf{aut}_{\mc{A}-{\rm LA}}(E_\mc{A}) = \mf{aut}_{\mc{A}-{\rm AV}}(E_\mc{A}) \cap \mf{aut}_{\rm LA}(E_{\mc{A}})$ of all $\mc{A}$-anchored vector bundle automorphisms which are also Lie algebroid morphisms. In other words, the corresponding operator $D_\sigma$ for $\sigma \in \Gamma(E_\mc{A})$ should be a derivation of $[\cdot,\cdot]_{E_\mc{A}}$. Note that $E_\mc{A}$ is a Lie algebroid so carries a canonical lift $\overline{\rho}_{E_\mc{A}}\colon \Gamma(E_\mc{A}) \to \mf{aut}_{\rm LA}(E_\mc{A})$ of $\rho_{E_\mc{A}}$ given by the operators $D_\sigma = [\sigma, \cdot]_{E_\mc{A}}$ for $\sigma \in \Gamma(E_\mc{A})$.
\begin{prop} The canonical lift $\overline{\rho}_{E_\mc{A}}$ is an $\mc{A}$-lift of $\rho_{E_\mc{A}}$.
\end{prop}
\bp Let $\sigma, \tau \in \Gamma(E_\mc{A})$ be given. Because $\varphi_{E_\mc{A}}$ is a Lie algebroid morphism, we have $\varphi_{E_\mc{A}}(D_\sigma(\tau)) = \varphi_{E_\mc{A}}([\sigma, \tau]_{E_\mc{A}}) = [\varphi_{E_\mc{A}}(\sigma), \varphi_{E_\mc{A}}(\tau)]_\mc{A} = [\sigma_{\mc{A}}, \varphi_{E_\mc{A}}(\tau)]_{\mc{A}} = D_{\sigma_\mc{A}}(\varphi_{E_\mc{A}}(\tau))$. We conclude that $\overline{\rho}_{E_\mc{A}}$ is an $\mc{A}$-lift of $\rho_{E_\mc{A}}$.
\ep
Consequently, the canonical lift maps $\overline{\rho}_{E_\mc{A}}\colon \Gamma(E_\mc{A}) \to \mf{aut}_{\mc{A}-{\rm LA}}(E_\mc{A})$. We always use the canonical lift to construct infinitesimal $\mc{A}$-Lie algebroid automorphisms.
\section{Splitting theorems}
\label{sec:splittingtheorems}
In this section we prove the $\mc{A}$-analogues of splitting theorems obtained in \cite{BursztynLimaMeinrenken16}. Let $i\colon N \to X$ be an $\mc{A}$-transversal for an $\mc{A}$-anchored vector bundle $E_\mc{A}$, giving an $\mc{A}$-linear approximation $p^! i^! E_\mc{A}$ of $E_\mc{A}$ along $N$. The following is the $\mc{A}$-analogue of \cite[Theorem 3.13]{BursztynLimaMeinrenken16}.
\begin{thm}\label{thm:aanchored} Let $\wt{V} \in \mf{aut}_{\mc{A}-{\rm AV}}(E_\mc{A})$ be vanishing along $i^! E_\mc{A}$ such that $V = \rho_{E_\mc{A}}(\wt{V})$ is Euler-like along $N$. Then there exists a unique isomorphism of $p^! i^! \mc{A}$- and $\mc{A}$-anchored vector bundles $\wt{\psi}\colon p^! i^! E_\mc{A} \to E_\mc{A}|_U$ with base map $\psi\colon \nu_N \to U \subseteq X$ a strong tubular neighbourhood embedding.
\end{thm}
\bp It follows that $V_\mc{A} = \varphi_{E_\mc{A}}(\wt{V}) \in \mf{aut}_{{\rm LA}}(\mc{A})$ vanishes along $i^! \mc{A}$ because the map $\varphi_{E_\mc{A}}\colon (E_\mc{A}, i^! E_\mc{A}) \to (\mc{A}, i^! \mc{A})$ is a strong map of pairs, and moreover we have $V = \rho_\mc{A}(V_\mc{A})$. As a consequence, by \cite[Theorems 3.13 and 4.1]{BursztynLimaMeinrenken16}  we obtain isomorphisms $\wt{\psi}\colon p^! i^! E_\mc{A} \to E_\mc{A}|_U$ and $\psi_\mc{A}\colon p^! i^! \mc{A} \to \mc{A}|_U$ as anchored vector bundles and Lie algebroids respectively. The functorial properties of the results loc.\ cit.\ applied to the map of pairs ${\rm id}_X \colon (X,N) \to (X,N)$ covering the morphism of anchored vector bundles $\varphi_{E_\mc{A}}\colon E_\mc{A} \to \mc{A}$ satisfying $\wt{V} \sim_{\varphi_{E_\mc{A}}} V_\mc{A}$ shows that $\varphi_{E_\mc{A}} \circ \wt{\psi} = \psi_\mc{A}$, so that $\wt{\psi}$ is an isomorphism of $p^! i^! \mc{A}$- and $\mc{A}$-anchored vector bundles.
\ep
If $E_\mc{A}$ is moreover an $\mc{A}$-Lie algebroid, its $\mc{A}$-anchor $\varphi_{E_\mc{A}}$ is a Lie algebroid morphism. The proof of \autoref{thm:aanchored} then immediately gives the following $\mc{A}$-analogue of \cite[Theorem 4.1]{BursztynLimaMeinrenken16}, establishing a splitting result for $\mc{A}$-Lie algebroids.
\begin{thm}\label{thm:aliealgebroid} Let $\wt{V} \in \mf{aut}_{\mc{A}-{\rm LA}}(E_\mc{A})$ be vanishing along $i^! E_\mc{A}$ such that $V = \rho_{E_\mc{A}}(\wt{V})$ is Euler-like along $N$. Then there exists a unique isomorphism of $p^! i^! \mc{A}$- and $\mc{A}$-Lie algebroids $\wt{\psi}\colon p^! i^! E_\mc{A} \to E_\mc{A}|_U$ with base map $\psi\colon \nu_N \to U \subseteq X$ a strong tubular neighbourhood embedding.
\end{thm}
The infinitesimal automorphisms $\wt{V}$ in the statements of Theorems \ref{thm:aanchored} and \ref{thm:aliealgebroid} exist due to the following, which are the $\mc{A}$-analogues of \cite[Theorem 3.13.a)]{BursztynLimaMeinrenken16}.
\begin{lem}\label{lem:aanchored} Let $E_\mc{A}$ be an $\mc{A}$-involutive $\mc{A}$-anchored vector bundle and $i\colon N \to X$ an $\mc{A}$-transversal for $E_\mc{A}$. Then there exists $\wt{V} \in \mf{aut}_{\mc{A}-{\rm AV}}(E_\mc{A})$ vanishing along $i^! E_\mc{A}$ such that $V = \rho_{E_\mc{A}}(\wt{V})$ is Euler-like along $N$.
\end{lem}
\bp As $N$ is a transversal for $E_\mc{A}$ by \autoref{prop:atranstrans}, using \cite[Lemma 3.9]{BursztynLimaMeinrenken16} there exists a section $\varepsilon \in \Gamma(E_\mc{A})$ vanishing at $i^! E_{\mc{A}}$ such that $V = \rho_{E_\mc{A}}(\varepsilon)$ is Euler-like along $N$. As $E_\mc{A}$ is $\mc{A}$-involutive, by \autoref{prop:aanchoredalift} there exists an $\mc{A}$-lift $\wt{\rho}_{E_\mc{A}}\colon \Gamma(E_\mc{A}) \to \mf{aut}_{\mc{A}-{\rm AV}}(E_\mc{A})$ of $\rho_{E_\mc{A}}$. Using the proof of \cite[Theorem 3.13.a)]{BursztynLimaMeinrenken16}, the vector field $\wt{V} := \wt{\rho}_{E_\mc{A}}(\varepsilon)$ on $E_\mc{A}$ satisfies the desired properties.
\ep
A similar statement for $\mc{A}$-Lie algebroids follows from the proof of \autoref{lem:aanchored}, using the canonical $\mc{A}$-lift $\overline{\rho}_{E_\mc{A}}\colon \Gamma(E_\mc{A}) \to \mf{aut}_{\mc{A}-{\rm LA}}(E_\mc{A})$ of $\rho_{E_\mc{A}}$.
\begin{lem}\label{lem:aliealgebroid} Let $E_\mc{A}$ be an $\mc{A}$-Lie algebroid and $i\colon N \to X$ an $\mc{A}$-transversal for $E_\mc{A}$. Then there exists $\wt{V} \in \mf{aut}_{\mc{A}-{\rm LA}}(E_\mc{A})$ vanishing along $i^! E_\mc{A}$ such that $V = \rho_{E_\mc{A}}(\wt{V})$ is Euler-like along $N$.
\end{lem}		
\backmatter
\bibliographystyle{hyperamsplain-nodash}
\cleardoublepage
\phantomsection
\addcontentsline{toc}{chapter}{Bibliography}
\bibliography{thesis}
\cleardoublepage
\chapter*{Summary\\[1ex]\hspace{1em}\large\emph{an overview for non-specialists}}
\addcontentsline{toc}{chapter}{Summary}
\markboth{Summary}{}
\label{chap:summary}

The title of this thesis is ``Geometric structures and Lie algebroids''. Here we provide an overview of its contents aimed at non-specialists, and in particular explain the words that are used in the title.

\subsection*{Differential geometry}

In this thesis we study specific types of mathematical structures on spaces that are called \emph{manifolds}. A manifold is a space that locally looks like Euclidean space, $\R^n$ for some nonnegative integer $n$, and the field studying them is called \emph{differential geometry}. The number $n$ in the local picture is called the \emph{dimension} of the manifold. One-dimensional examples of manifolds include the line and the circle, with the sphere and the donut giving examples in dimension two.

While it is often helpful to visualize a manifold using our three-dimensional viewpoint, this is no longer possible when considering manifolds which are, say, four-dimensional. The abstract concept of a manifold is ubiquitous in modern mathematics and physics. For example, many physical theories have a geometric flavour and admit a natural description using the language of differential geometry, without the manifolds that are used necessarily directly referring to physical space.

\subsection*{Geometric structures}

While manifolds are interesting in their own right, one often equips them with further structure before proceeding. Indeed, manifolds should be thought of as an abstract canvas on which the subject of study takes place. These extra geometric structures are usually put not on the manifold $X$ itself, but rather on its \emph{tangent bundle} $TX$, which is an auxiliary manifold keeping track of velocities of curves on $X$.

One example of a geometric structure is called a \emph{symplectic structure}. These originally arose from physics as they provide the language to discuss Hamiltonian dynamics in classical mechanics. A symplectic structure $\omega$ on a given manifold $X$ allows for the measurement of ``volume'' of all of its even-dimensional submanifolds (this concept of volume is not exactly what one perhaps expects, but it does agree with intuition in dimension two). A symplectic structure further provides a so-called Poisson bracket, which gives rise to the Hamiltonian dynamics on $X$ mentioned before. This is in fact a very special type of Poisson bracket, as the symplectic structure is maximally nondegenerate. The notion of a \emph{Poisson structure} $\pi$ naturally generalizes a symplectic structure, focusing on the Poisson bracket it gives rise to. As will be important later on, a Poisson structure can be viewed as a certain kind of map
\be
	\pi^\sharp\colon T^*X \to TX,
\ee
going from the \emph{cotangent bundle} $T^*X$ to the tangent bundle $TX$. With the tangent bundle describing velocities, its dual, the cotangent bundle, keeps track of momentum. For Poisson structures $\pi$ that are symplectic, the map $\pi^\sharp$ is an \emph{isomorphism}, giving a nice correspondence between the two bundles and allowing us to identify them. In particular, it can then be inverted (so that we write $\pi = \omega^{-1}$ as a shorthand), giving rise a map in the opposite direction, as follows.
\be
	\omega^\flat\colon TX \to T^*X.
\ee
However, general Poisson structures allow for this map to become more degenerate, and in extreme cases even zero. In this thesis we are mainly interested in Poisson structures which are close to being symplectic, as measured by the map $\pi^\sharp$ above being well-behaved. Indeed, we consider Poisson structures where this map is almost-everywhere an isomorphism, only failing to be on a small portion of the manifold.

A more complicated geometric structure that is also of interest to us is that of a \emph{generalized complex structure}. While we cannot here define them precisely, they should be thought of as having both a symplectic and a \emph{complex} flavour.\footnote{Roughly speaking, a complex structure describes the concept of ``clockwise rotation by $90$ degrees'', although there need not be a way to measure degrees.} As with Poisson structures, in this thesis we are particularly interested in those examples that are close to being symplectic. Underlying any \gcs $\mc{J}$ is a Poisson structure, $\pi_\mc{J}$, and it is this Poisson structure that measures the distance $\mc{J}$ has from being symplectic in the same manner as before.

\subsection*{Lie algebroids}

While the geometric structures we put on our manifold $X$ naturally live on its tangent bundle $TX$, there is a slight inconvenience. Namely, a general Poisson structure $\pi$, even if it is close to being symplectic, does not give an isomorphism between $T^*X$ and $TX$. Because of this, $\pi$ cannot be inverted smoothly to $\omega = \pi^{-1}$. If this would be possible, we could use powerful techniques from symplectic geometry to study $\pi$, but alas. In general, its inverse $\omega$ is singular, and does not exist on the entire manifold. 

This problem can be resolved in favorable cases by replacing the tangent bundle $TX$ by another bundle, called a \emph{Lie algebroid} $\mc{A}$. This bundle comes equipped with a map $\rho_\mc{A}\colon \mc{A} \to TX$ relating it to $TX$. The Lie algebroids we use are those for which this map is again almost-everywhere an isomorphism, failing to be so exactly where the map $\pi^\sharp$ is not an isomorphism either. In this sense the degenerate behavior of $\pi$ is absorbed in that of $\mc{A}$ using the process of \emph{lifting}. Indeed, as we show it is sometimes possible to take our not-symplectic Poisson structure $\pi$, and lift it to its analogue on $\mc{A}$, called an \emph{$\mc{A}$-Poisson structure $\pi_\mc{A}$}. Diagrammatically this results in the following.
\begin{center}
	\begin{tikzpicture}
	\matrix (m) [matrix of math nodes, row sep=2.5em, column sep=2.5em,text height=1.5ex, text depth=0.25ex]
	{	\mc{A}^* & \mc{A} \\ T^*X & TX \\};
	\path[-stealth]
	(m-1-1) edge node [above] {$\pi_\mc{A}^\sharp$} (m-1-2)
	(m-2-1) edge node [left] {$\rho_\mc{A}^*$} (m-1-1)
	(m-2-1) edge node [above] {$\pi^\sharp$} (m-2-2)
	(m-1-2) edge node [right] {$\rho_\mc{A}$} (m-2-2);
	\end{tikzpicture}
\end{center}
This diagram should be familiar from the cover. Once we have lifted appropriately, the $\mc{A}$-Poisson structure $\pi_\mc{A}$ now \emph{is} symplectic, and we can use the aforementioned symplectic techniques to study $\pi_\mc{A}$, and hence our original Poisson structure $\pi$. We will not attempt to explain the results obtained in this thesis using this process here. However, after reading up to this point, the reader is invited to revisit the \hyperref[chap:intro]{Introduction}, especially the first two pages.

\selectlanguage{dutch}
{\let\cleardoublepage\relax \chapter*{Samenvatting\\[1ex]\hspace{1em}\large\emph{een overzicht voor niet-experts}}}
\addcontentsline{toc}{chapter}{Samenvatting}
\markboth{Samenvatting}{}
\label{chap:samenvatting}

De titel van dit proefschrift luidt ``Meetkundige structuren en Lie algebro\"iden''. Hier geven we een overzicht van de inhoud voor niet-experts. In het bijzonder zullen we de woorden uit de titel toelichten.\\

\subsection*{Differentiaalmeetkunde}

In dit proefschrift bestuderen we specifieke typen wiskundige structuren op ruimten die \emph{vari\"eteiten} genoemd worden. Een vari\"eteit is een ruimte die er lokaal uitziet als de Euclidische ruimte, $\R^n$ voor een zeker geheel getal $n$. Het gebied binnen de wiskunde dat hen bestudeert wordt \emph{differentiaalmeetkunde} genoemd. Het getal $n$ in het lokale plaatje heet de \emph{dimensie} van de vari\"eteit. Eendimensionale voorbeelden van vari\"eteiten zijn de lijn en de cirkel, waarbij de bol en de donut voorbeelden geven in twee dimensies.

Hoewel het vaak handig is om vari\"eteiten te visualiseren vanuit ons driedimensionaal oogpunt, is dit niet meer mogelijk bij vari\"eteiten die, zeg, vierdimensionaal zijn. Het abstracte concept van een vari\"eteit is alomtegenwoordig in de moderne wis- en natuurkunde. Zo hebben bijvoorbeeld vele natuurkundige theorie\"en een meetkundig aspect en kunnen ze natuurlijk in de taal van de differentiaalmeetkunde beschreven worden. Hierbij hoeven de gebruikte vari\"eteiten niet noodzakelijkerwijs direct te verwijzen naar de fysieke ruimte.

\subsection*{Meetkundige structuren}

Hoewel vari\"eteiten op zichzelf interessant zijn, beschouwt men vaak vari\"eteiten met extra structuur. In deze zin kunnen vari\"eteiten gezien worden als een abstract canvas waarop het onderzoeksonderwerp zich afspeelt. Deze extra meetkundige structuren worden meestal niet direct op de vari\"eteit $X$ zelf aangelegd, maar in plaats daarvan op zijn \emph{raakbundel} $TX$, een hulp-vari\"eteit die snelheden van krommen op $X$ bijhoudt.

Een voorbeeld van zo'n meetkundige structuur wordt een \emph{symplectische structuur} genoemd. Deze komen oorspronkelijk uit de natuurkunde, waar ze de taal verschaffen om te spreken over Hamiltoniaanse dynamica in klassieke mechanica. Een symplectische structuur $\omega$ op een gegeven vari\"eteit $X$ leidt tot het kunnen meten van het ``volume'' van al zijn evendimensionale deelvari\"eteiten (dit concept van volume is niet precies wat men misschien verwacht, maar komt wel overeen met intu\"itie in twee dimensies). Een symplectische structuur geeft verder een zogeheten Poisson haakje, welke leidt tot de eerdergenoemde Hamiltoniaanse dynamica op $X$. Dit is een erg speciaal type Poisson haakje, omdat de symplectische structuur maximaal niet-gedegeneerd is. Het concept van een \emph{Poisson structuur} is een natuurlijke generalisatie van een symplectische structuur, waarbij men de nadruk legt op het bijbehorende Poisson haakje. Zoals verderop van belang zal zijn kan een Poisson structuur gezien worden als een bepaald soort afbeelding
\be
\pi^\sharp\colon T^*X \to TX,
\ee
van de \emph{coraakbundel} $T^*X$ naar de raakbundel $TX$. Waar de raakbundel snelheden beschrijft, houdt zijn duale, de coraakbundel, het concept van momentum bij. Voor Poisson structuren $\pi$ die symplectisch zijn is de afbeelding $\pi^\sharp$ een \emph{isomorfisme}, en geeft een fijne correspondentie tussen de twee bundels, zodat ze met elkaar ge\"identificeerd kunnen worden. In het bijzonder kan deze afbeelding dan ge\"inverteerd worden (en schrijven we ookwel $\pi = \omega^{-1}$), zodat er ook een afbeelding in de andere richting bestaat, als volgt.
\be
\omega^\flat\colon TX \to T^*X.
\ee
Echter, algemene Poisson structuren laten toe dat deze afbeelding meer gedegenereerd wordt, en in extreme gevallen zelfs nul. In dit proefschrift zijn we met name ge\"interesseerd in Poisson structuren die vrijwel symplectisch zijn, zoals gemeten wordt door het gedrag van de afbeelding $\pi^\sharp$. In het bijzonder beschouwen we Poisson structuren waar deze afbeelding vrijwel overal een isomorfisme is, en alleen hierin faalt op een klein deel van de vari\"eteit.

Een meer ingewikkelde meetkundige structuur welke ons ook interesseert is dat van een \emph{gegeneraliseerde complexe structuur}. Hoewel we deze hier niet precies kunnen defini\"eren, kan men over ze denken alsof ze zowel een symplectische als een \emph{complexe} kant hebben.\footnote{Grof gezegd beschrijft een complexe structuur het concept van ``90 graden rotatie met de klok mee'', al hoeft er geen manier te zijn om graden te kunnen meten.} Zoals voor Poisson structuren zijn we in dit proefschrift in het bijzonder ge\"interesseerd in die voorbeelden welke vrijwel symplectisch zijn. Aan elke gegeneraliseerde complexe structuur $\mc{J}$ ligt een Poisson structure $\pi_{\mc{J}}$ ten grondslag, en het is deze Poisson structuur die net als voorheen meet in hoeverre $\mc{J}$ niet-symplectisch is.

\subsection*{Lie algebro\"iden}

Hoewel de meetkundige structuren welke we op onze vari\"eteit $X$ aanleggen van nature leven op zijn raakbundel $TX$, is er een kleine tegenslag. Een algemene Poisson structuur $\pi$, zelfs wanneer hij dichtbij een symplectische is, geeft namelijk geen isomorfisme tussen $T^* X$ en $TX$. Hierdoor kan $\pi$ niet op een gladde manier ge\"inverteerd worden tot $\omega = \pi^{-1}$. Als dit mogelijk zou zijn zouden we krachtige technieken uit de symplectische meetkunde kunnen gebruiken om $\pi$ te bestuderen, maar helaas. In het algemeen is zijn inverse $\omega$ singulier, en bestaat niet op de gehele vari\"eteit.

Dit probleem kan in gunstige gevallen opgelost worden door de raakbundel $TX$ te vervangen door een andere bundel, genaamd een \emph{Lie algebro\"ide} $\mc{A}$. Deze bundel is uitgerust met een afbeelding $\rho_\mc{A}\colon \mc{A} \to TX$ die hem relateert aan $TX$. Het soort Lie algebro\"iden welke we gebruiken zijn zodanig dat deze afbeelding weer vrijwel overal een isomorfisme is, en hierin faalt precies wanneer de afbeelding $\pi^\sharp$ dat doet. In deze zin wordt het gedegenereerde gedrag van $\pi$ geabsorbeerd in dat van $\mc{A}$ met behulp van het proces van \emph{liften}. Zoals we laten zien is het soms mogelijk om onze niet-symplectische Poisson structure $\pi$ te liften naar zijn analogon op $\mc{A}$, ofwel een \emph{$\mc{A}$-Poisson structuur} $\pi_\mc{A}$. In een diagram resulteert dit in het volgende.
\begin{center}
	\begin{tikzpicture}
	\matrix (m) [matrix of math nodes, row sep=2.5em, column sep=2.5em,text height=1.5ex, text depth=0.25ex]
	{	\mc{A}^* & \mc{A} \\ T^*X & TX \\};
	\path[-stealth]
	(m-1-1) edge node [above] {$\pi_\mc{A}^\sharp$} (m-1-2)
	(m-2-1) edge node [left] {$\rho_\mc{A}^*$} (m-1-1)
	(m-2-1) edge node [above] {$\pi^\sharp$} (m-2-2)
	(m-1-2) edge node [right] {$\rho_\mc{A}$} (m-2-2);
	\end{tikzpicture}
\end{center}
Dit diagram zou bekend moeten zijn van de omslag van dit proefschrift. Nadat we op de juiste manier hebben gelift zal de $\mc{A}$-Poisson structuur $\pi_\mc{A}$ nu \emph{wel} symplectisch zijn. Hierdoor kunnen we dan de eerdergenoemde symplectische technieken gebruiken om $\pi_\mc{A}$ te bestuderen, en daarmee ook onze originele Poisson structuur $\pi$. We zullen niet proberen om hier de resultaten uit dit proefschrift behaald via dit proces te beschrijven. Echter, na tot dit punt te hebben gelezen, wordt de lezer uitgenodigd om de \hyperref[chap:intro]{Introductie} nogmaals te bezoeken, in het bijzonder de eerste twee pagina's.

\selectlanguage{english}
\chapter*{Acknowledgements}
\addcontentsline{toc}{chapter}{Acknowledgements}
\markboth{Acknowledgements}{}

As my advisor and many before him have said, completing a Ph.D.\ is never a solo effort. I would like to briefly thank some who have helped me along the way.

First I would like to thank Gil for his willingness to accept me as his Ph.D.\ student after our interaction during my Master thesis project. Gil, thank you for your guidance and effort invested in me, our many meetings and talks, and your open door policy. Marius, thanks for being my promotor and our mathematical discussions, but also for the barbecues and fostering such a pleasant atmosphere within the research group.

To the members of the reading committee, Erik van den Ban, Christian Blohmann, Rui Loja Fernandes, Marco Gualtieri, and Eckhard Meinrenken: thank you for the time and effort invested in reading this thesis, and for your comments and suggestions.

I would like to thank especially my roommates Arjen and Joey for their mathematical and LaTeX help, but most of all their friendship and camaraderie. Special thanks also go to our next-door neighbours, Davide and Mike.

Then I would thank the group of people that have been part of the Poisson geometry group and Friday Fish seminars during my stay, including (in alphabetical order) Boris, Daniele, David, Dima, Florian S., Florian Z., Francesco, Ionu\c{t}, Jo\~{a}o, Joost, Lauran, Maria, Matias, Ori, Panos, Pedro, Roy, and Stefan. Thank you for the many seminars we've had together and for the sense of being part of a mathematical family. I would also like to thank the other Ph.D.\ candidates I've met either at Utrecht or through the many GQT schools: Anshui, Arie, Bram, Dana, Danilo, Du\v{s}an, Jules, KaYin, Marcelo, Peter, Shan, Sebastian, Valentijn, and many others.

I further thank the other staff at the Mathematical Institute with which I have had a chance to interact, particularly Carel Faber, Martijn Kool, Andr\'e Henriques, Ieke Moerdijk, and Fabian Ziltener. I am also grateful to C\'ecile, Jean, Ria, and Wilke, but also Carin, Sylvia, and Barbara, for taking care of things and ensuring everything runs smoothly.

Then, I would like to thank my family and friends for their support and distraction. Finally, I mention Kirsten: this is just the beginning of our journey together.

{\let\cleardoublepage\relax \chapter*{Curriculum Vit\ae}}
\addcontentsline{toc}{chapter}{Curriculum Vit\ae}

Ralph Leonard Klaasse was born on April 1, 1990 in Nijmegen, the Netherlands.

\vspace{\baselineskip}

\noindent In June 2007 he graduated high school with the profiles N\&T and N\&G from Dominicus College in Nijmegen. In September 2007 he began his studies in Applied Mathematics at Eindhoven University of Technology. After one year he further began studies in Applied Physics, and obtained bachelor degrees with cum laude distinction in both disciplines in October 2011 and September 2012 respectively. While at Eindhoven University of Technology, he participated in the Honors Horizon program for excellence. His bachelor thesis was titled ``Protein-mediated instability in curved membranes'', and was supervised by prof.~dr.~M.~A.~Peletier and dr.~C.~Storm.

\vspace{\baselineskip}

\noindent He then moved to Utrecht University and started the master program Mathematical Sciences in September 2011. In August 2013 he obtained a master degree with cum laude distinction. His master thesis was titled ``Seiberg--Witten theory for symplectic manifolds'', and was supervised by dr.~G.~R.~Cavalcanti.

\vspace{\baselineskip}

\noindent In September 2013 he started his Ph.D.\ candidacy at Utrecht University under supervision of dr.~G.~R.~Cavalcanti. While at Utrecht University, he was a member of the Geometric Structures group led by prof.~dr.~M.~N.~Crainic. During February and March 2017, he visited prof.~dr.~R.~L.~Fernandes at the University of Illinois at Urbana--Champaign. He was awarded the KWG prize for Ph.D.\ candidates at the Dutch mathematical congress in April 2017. He went on to defend his Ph.D.\ thesis, titled ``Geometric structures and Lie algebroids'', on August 29, 2017.

\vspace{\baselineskip}

\noindent Starting from September 2017, he will be a Postdoctoral Fellow at the Universit\'e Libre de Bruxelles, in the research group of dr.~J.~Fine.

\cleardoublepage
\thispagestyle{empty}\pagenumbering{gobble}\null
\includepdf[pages={1}]{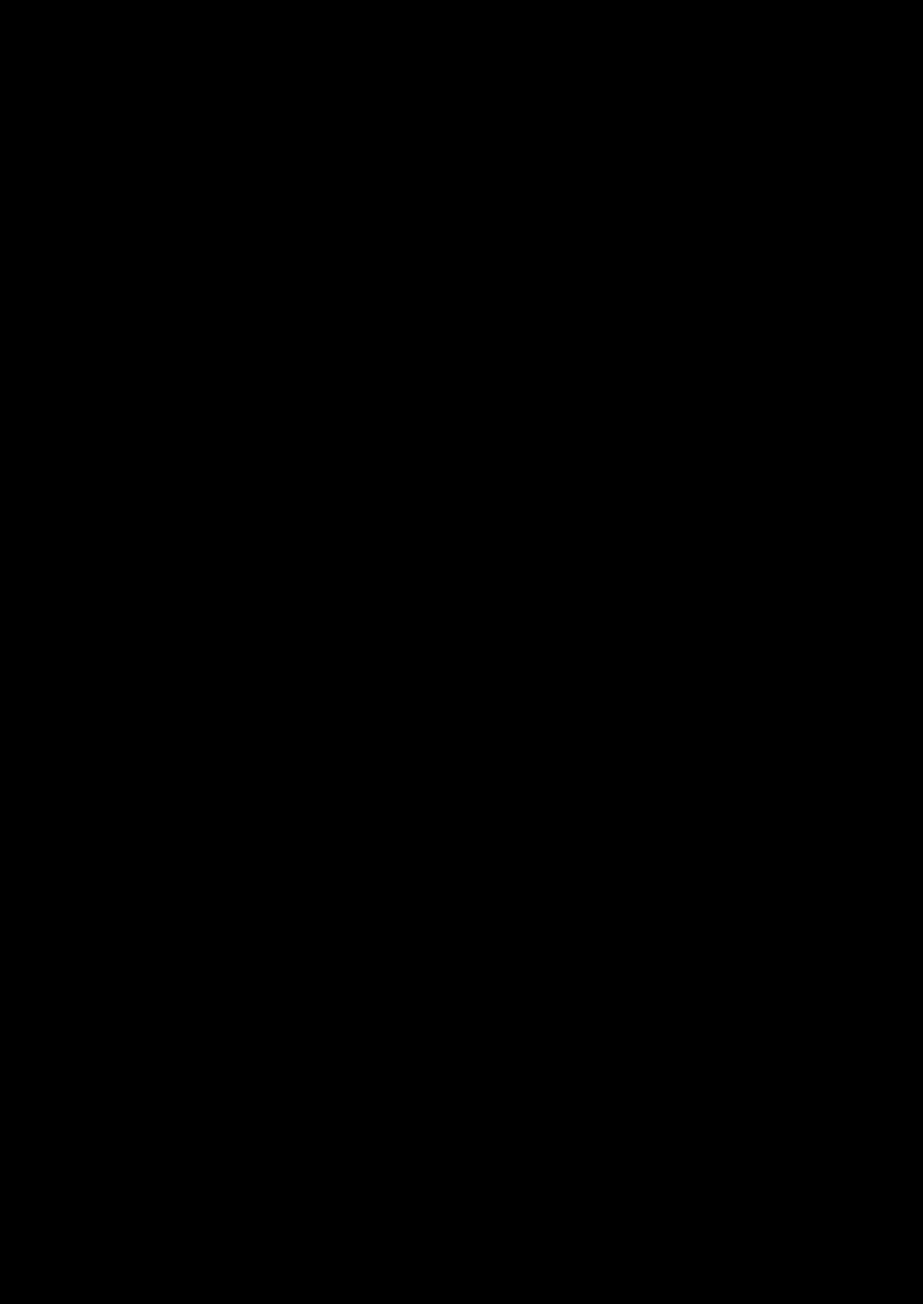}
\end{document}